\documentclass{article}
\usepackage{latexsym}
\usepackage{amsmath,amssymb}
\usepackage{txfonts}
\usepackage[dvips]{graphicx, color}

\usepackage{color}

\begin{document}
\title{\bf Stability and Arithmetic}
\author{\bf Lin WENG}
\date{\it \footnotesize{Dedicated to the memory of my father Jiahua WENG\\
21 September 1933 -- 12 April 2009}}
\maketitle

\begin{abstract}
Stability plays a central role in arithmetic.
In this article, we explain some basic ideas and 
present certain constructions for our studies. 
There are two aspects: namely,	
general Class Field Theories for Riemann surfaces 
using semi-stable parabolic bundles
and for $p$-adic number fields using what we 
call semi-stable filtered $(\varphi,N;\omega)$-modules;
and non-abelian zeta functions for  function fields 
over finite fields using semi-stable bundles 
and for number fields using semi-stable lattices. 
\end{abstract}


\vskip 1cm
\noindent
\centerline{C\scriptsize{ONTENTS}}

\vskip 0.5cm
\noindent
{Introduction}\hfill 7
\vskip 0.20cm
\noindent
{Part A. Guidances from Geometry}
\vskip 0.20cm
\noindent
{Chapter I. Micro Reciprocity Law in Geometry}\hfill 13
\vskip 0.20cm
\noindent
\S1\ \ Narasimhan-Seshadri Correspondence\hfill 13

1.1\ \ Uniformization\hfill 13

1.2\ \ Narasimhan-Seshadri Correspondence\hfill 13

\noindent
\S2\ \ Micro Reciprocity Law\hfill 14

2.1\ \ Weil's Program\hfill 14

2.2\ \ Micro Reciprocity Law\hfill 15
\vskip 0.20cm
\noindent
{Chapter II. CFTs in Geometry}\hfill 17
\vskip 0.20cm
\noindent
\S3\ \ Arithmetic CFT: Class Field Theory\hfill 17

\noindent
\S4\ \ Geometric CFT: Conformal Field Theory\hfill 17
\eject
\vskip 0.40cm
\noindent
{Part B. High Rank Zeta Functions and Stability}
\vskip 0.20cm
\noindent
{Chapter III. High Rank Zeta Functions}\hfill 19
\vskip 0.20cm
\noindent
\S5\ \ Function Fields\hfill 19

5.1\ \ Definition and Basic Properties\hfill 19

5.2\ \ Global High Rank Zetas via Euler Products\hfill 20

\noindent
\S6\ \ Number Fields\hfill 21

6.1\ \ Stability of $\mathcal O_F$-Lattices\hfill 21

6.2\ \ Geo-Arithmetical Cohomology\hfill 22

6.3\ \ High Rank Zetas\hfill 23
\vskip 0.20cm
\noindent
{Chapter IV. Geometric Characterization of Stability}\hfill 25
\vskip 0.20cm
\noindent
\S7\ \ Upper Half Space Model\hfill 25

7.1\ \ Upper Half Plane\hfill 25

7.2\ \ Upper Half Space\hfill 25

7.3\ \ Rank Two $\mathcal O_F$-Lattices: Upper Half Space Model\hfill 26

\noindent
\S8\ \ Cusps\hfill 26

8.1\ \ Definition\hfill 26

8.2\ \ Cusp and Ideal Class Correspondence\hfill 27

8.3\ \ Stablizer Groups of Cusps\hfill 27

8.4\ \ Fundamental Domain for $\Gamma_\infty$ on
 ${\mathcal H}^{r_1}\times{\mathbb H}^{r_2}$\hfill 28

\noindent
\S9\ \ Fundamental Domain\hfill 29

9.1\ \ Siegel Type Distance\hfill 29

9.2\ \ Fundamental Domain\hfill 30

\noindent
\S10\ \  Stability in Rank Two\hfill 31

10.1\ \ Stability and Distances to Cusps\hfill 31

10.2\ \ Moduli Space of Rank 1\ \ Semi-Stable $\mathcal O_F$-Lattices\hfill 31
\vskip 0.20cm
\noindent
{Chapter V. Algebraic Characterization of Stability}\hfill 33
\vskip 0.20cm
\noindent
\S11\ \ Canonical Filtrations\hfill 33

11.1\ \ Canonical Filtrations\hfill 33

11.2\ \ Examples of Lattices\hfill 34

\noindent
\S12\ \ Algebraic Characterization\hfill 34

12.1\ \ A GIT Principle\hfill 34

12.2\ \ Micro-Global Relation for Geo-Ari Truncation\hfill 34
\newpage
\vskip 0.20cm
\noindent
{Chapter VI. Analytic Characterization of Stability}\hfill 36
\vskip 0.20cm
\noindent
\S13\ \ Arthur's Analytic Truncation\hfill 36

13.1\ \ Parabolic Subgroups\hfill 36

13.2\ \ Logarithmic Map\hfill 36

13.3\ \ Roots, Coroots, Weights and Coweights\hfill 37

13.4\ \ Positive Cone and Positive Chamber\hfill 37

13.5\ \ Partial Truncation and First Estimations\hfill 38

\noindent
\S14\ \ Reduction Theory\hfill 38

14.1\ \ Langlands' Combinatorial Lemma\hfill 38

14.2\ \ Langlands-Arthur's Partition: Reduction Theory\hfill 39

\noindent
\S15\ \ Arthur's Analytic Truncation\hfill 40

15.1\ \ Definition\hfill 40

15.2\ \ Basic Properties\hfill 41

15.3\ \ Truncation $\Lambda^T{\bold 1}$\hfill 41

\noindent
\S16\ \ Analytic Characterization of Stability\hfill 42

16.1\ \ A Micro Bridge\hfill 42

16.2\ \ Analytic Truncations and Stability\hfill 43
\vskip 0.20cm
\noindent
{Chapter VII. Non-Abelian $L$-Functions}\hfill 44
\vskip 0.20cm
\noindent
\S17\ \ High Rank Zetas and Eisenstein Series\hfill 44

17.1\ \ Epstein Zeta Functions and High Rank Zetas\hfill 44

17.2\ \ Rankin-Selberg Method: An Example with $SL_2$\hfill 44

\noindent
\S18\ \ Non-Abelian $L$-Functions\hfill 46

18.1\ \ Automorphic Forms and Eisenstein Series\hfill 46

18.2\ \ Non-Abelian L-Functions\hfill 50

\noindent
\S19\ \ Basic Properties of Non-Abelian $L$-Functions\hfill 51

19.1\ \ Meromorphic Extension and Functional Equations\hfill 51

19.2\ \ Holomorphicity and Singularities\hfill 52
\vskip 0.20cm
\noindent
{Chapter VIII. Symmetries and the Riemann Hypothesis}\hfill 53
\vskip 0.20cm
\noindent
\S20\ \  Abelian Parts of High Rank Zetas\hfill 53

20.1\ \ Analytic Studies of High Rank Zetas\hfill 53

20.2\ \ Advanced Rankin-Selberg and Zagier Methods\hfill 54

20.3\ \ Discovery of Maximal Parabolics: SL, Sp and $G_2$\hfill 55

\noindent
\S21\ \ Abelian Zetas for $(G,P)$\hfill 56

21.1\ \ Definition\hfill 56

21.2\ \ Conjectural FE and the RH\hfill 57

\noindent
\S22\ \ Abelian Parts of High Rank Zetas\hfill 57
\eject
\vskip 0.50cm
\noindent
{Part C. General CFT and Stability}
\vskip 0.20cm
\noindent
{Chapter IX. $l$-adic Representations for $p$-adic Fields}\hfill 59
\vskip 0.20cm
\noindent
\S23\ \ Finite Monodromy and Nilpotency\hfill 59

23.1\ \ Absolute Galois Group and Its pro-$l$ Structures\hfill 59

23.2\ \ Finite Monodromy\hfill 60

23.3\ \ Unipotency\hfill 60
\vskip 0.20cm
\noindent
{Chapter X. Primary Theory of $p$-adic Representations}\hfill 62
\vskip 0.20cm
\noindent
\S24\ \ Preliminary Structures of Absolute Galois Groups\hfill 62

24.1\ \ Galois Theory: A $p$-adic Consideration\hfill 62

24.2\ \ Arithmetic Structure: Cyclotomic Character\hfill 62

24.3\ \ Geometric Structure: Fields of Norms\hfill 62

\noindent
\S25\ \ Galois Representations: Characteristic $p$-theory\hfill 64

25.1\ \ $\mathbb F_p$-Representations\hfill 64

25.2\ \ Etale $\varphi$-modules\hfill 65

25.3\ \ Characteristic $p$ Representation and Etale $\varphi$-Module\hfill 65

\noindent
\S26\ \ Lifting to Characteristic Zero\hfill 66

26.1\ \ Witt Vectors and Teichm\"uller Lift\hfill 66

26.2\ \ $p$-adic	Representations of Fields of Characteristic 0\hfill 67

26.3\ \ $p$-adic Representations and Etale $(\varphi,\Gamma)$-Modules\hfill 68
\vskip 0.20cm
\noindent
{Chapter XI. $p$-adic Hodge and Properties of Periods}\hfill 69
\vskip 0.20cm
\noindent
\S27\ \ Hodge Theory over $\mathbb C$\hfill 69

\noindent
\S28\ \ Admissible Galois Representations\hfill 70

\noindent
\S29\ \ Basic Properties of Various Periods\hfill 70

29.1\ \ Hodge-Tate Periods\hfill 70

29.2\ \ de Rham Periods\hfill 71

29.3\ \ Crystalline Periods\hfill 71

29.4\ \ Semi-Stable Periods\hfill 72

\noindent
\S30\ \  Hodge-Tate, de Rham, Semi-Stable and Crystalline Reps\hfill 73

30.1\ \ Definition\hfill 73

30.2\ \ Basic Structures of $\mathbb D_{\bullet}(V)$\hfill 73

30.3\ \ Relations among Various $p$-adic Representations\hfill 74

30.4\ \ Examples\hfill 75

\noindent
\S31\ \ $p$-adic Hodge Theory\hfill 76
\eject
\vskip 0.20cm
\noindent
{Chapter XII. Fontaine's Rings of Periods}\hfill 78
\vskip 0.20cm
\noindent
\S32\ \ The Ring of de Rham Periods $\mathbb B_{\mathrm{dR}}$\hfill 78

\noindent
\S33\ \ The Ring of Crystalline Periods $\mathbb B_{\mathrm{crys}}$\hfill 80

\noindent
\S34\ \ The Ring of Semi-Stable Periods $\mathbb B_{\mathrm{st}}$\hfill 81

\vskip 0.20cm
\noindent
{Chapter XIII. Micro Reciprocity Laws and General CFT}\hfill 82
\vskip 0.20cm
\noindent
\S35\ \ Filtered $(\varphi,N)$-Modules and Semi-Stable Representations\hfill 82

35.1\ \ Definition\hfill 82

35.2\ \ Weak Admissibility and Semi-Stablility\hfill 84

\noindent
\S36\ \ Monodromy Theorem for $p$-adic Galois Representations\hfill 85

\noindent
\S37\ \ Semi-Stability of Filtered $(\varphi,N;\omega)$-Modules\hfill 85

37.1\ \ Weak Admissibility = Stability and of Slope Zero\hfill 85

37.2\ \  Ramifications\hfill 86

37.3\ \  $\omega$-Structures\hfill 87

37.4\ \ Semi-Stability of Filtered $(\varphi,N;\omega)$-Modules\hfill 88

\noindent
\S38\ \ General CFT for $p$-adic Number Fields\hfill 88

38.1\ \ Conjectural Micro Reciprocity Law\hfill 88

38.2\ \ General CFT for $p$-adic Number Fields\hfill 89
\vskip 0.20cm
\noindent
{Chapter XIV. GIT Stability, Moduli and Invariants}\hfill 90
\vskip 0.20cm
\noindent
\S39\ \ Moduli Spaces\hfill 90

\noindent
\S40\ \  Polarizations and Galois Cohomology\hfill 91

\noindent
\S41\ \ Iwasawa Cohomology and Dual Exp Map\hfill 92

41.1\ \ Galois Cohomology\hfill 92

41.2\ \ $(\varphi,\Gamma)$-Modules and Galois Cohomology\hfill 93

41.3\ \ Iwasawa Cohomology $H_{\mathrm{Iw}}^i(K,V)$\hfill 93

41.4\ \ Two Descriptions of $H_{\mathrm{Iw}}^i(K,V)$\hfill 93

41.5\ \ Dual Exponential Maps\hfill 95
\vskip 0.20cm
\noindent
{Chapter XV. Two Approaches to Conjectural MRL}\hfill 96
\vskip 0.20cm
\noindent
\S42\ \ Algebraic and Geometric Methods\hfill 96

\noindent
\S43\ \ MRL with Limited Ramifications\hfill 96

43.1\ \ Logarithmic Map\hfill 97

43.2\ \ Basic Structures of $\mathbb B_{\mathrm{crys}}^{\varphi=1}$\hfill 97

43.3\ \ Rank One Structures\hfill 98

\noindent
\S44\ \ Filtration of Invariant Lattices\hfill 99

\eject
\vskip 0.5cm
\noindent
\S45\ \ Sen-Tate Theory and Its Generalizations\hfill 100

45.1\ \ Sen's Method\hfill 100

45.2\ \ Sen's Theory for $\mathbb B_{\mathrm{dR}}$\hfill 101

45.3\ \ Overconvergency\hfill 101

\noindent
\S46\ \ $p$-adic Monodromy Theorem\hfill 102

\noindent
\S47\ \ Infinitesimal, Local and Global\hfill 104

47.1\ \ From Arithmetic to Geometry\hfill 104

47.2\ \ From Infinitesimal to Global\hfill 104

\noindent
\S48\ \ Convergent $F$-isocrystals and Rigid Stable $F$-Bundles\hfill 105

48.1\ \ Rigid Analytic Space\hfill 105

48.2\ \ Convergent $F$-Isocrystals\hfill 105

48.3\ \ Integrable and Convergent Connections\hfill 106

48.4\ \ Frobenius Structure\hfill 106

48.5\ \ Unit-Root $F$-Isocrystals\hfill 107

48.6\ \ Stability of Rigid $F$-Bundles\hfill 107

\noindent
\S49\ \ Overconvergent $F$-Isocrystals, Log Geometry and Stability\hfill 108

49.1\ \ Overconvergent Isocrystals\hfill 108

49.2\ \ $p$-adic Representations with Finite Local Monodromy\hfill 109

49.3\ \ Logarithmic Rigid Analytic Geometry\hfill 109

\vskip 0.20cm
\noindent{References}\hfill 112
\vfill
\eject
\vskip 1cm
\centerline{\large\bf{Introduction}}
\vskip 1cm

\noindent
In the past a decade or so, importance of stability, which originally 
appeared and has played key roles in algebraic geometry, was gradually 
recognized by many people working in arithmetic. As typical examples, we now 
have 

\noindent
(i) Existence theorem and reciprocity law of a non-abelian class field theory for function fields over complex numbers, based on Seshadri's work of semi-stable parabolic bundles over Riemann surfaces;	

\noindent
(ii) High rank zeta functions for global fields,
defined as natural integrations over moduli spaces of semi-stable bundles/lattices; and 

\noindent
(iii) Characterization of the so-called semi-stable representations for absolute Galois groups of $p$-adic number fields, in terms of weakly admissible filtered $(\varphi,N)$-modules,  or better,	 semi-stable filtered $(\varphi,N)$-modules of slope zero.

Along with this line,  as an integrated part of our Program on Geometric Arithmetic, in this article, we explain some basic ideas and present certain constructions  using stability to study two non-abelian aspects of arithmetic, one at	 a micro level and	
the other on large scale.
\vskip 0.30cm
\noindent
{\bf I) Micro Level} 
\vskip 0.30cm
\noindent
We, at this micro level, want to give a characterization for each 
individual Galois representation. For this, first we, according to the 
associated base field and coefficients, classify Galois representations 
into four types, namely, $v$-adic/adelic representations for local/global (number) fields. As such, then our aim becomes to expose some totally independent structures from which the original Galois representations
can be reconstructed.

In general, arbitrary Galois  representations are too complicated to have 
clearer structures, certain natural restrictions should be imposed. 
In this direction, we, as a natural continuation of existing theories of Galois representations, 
choose to make the following rather standard 
restrictions: 

\begin{table}[htbp]
\label{tab:galois}
\begin{center}
\begin{tabular}{|c|c|c|}  \hline
Fields$\backslash$Coefs & $v$-Adic & Adelic\\ \hline
Local & Fin Monodromy \&  Nilpotency & Compatible System\\ \hline
Global & + Finite Ramification & + Admissible System\\ \hline
\end{tabular}
\end{center}
\end{table}
To be more precise, they are:
\vskip 0.20cm
\noindent
{\bf (i)  $v$-adic Galois Representations} for
\vskip 0.20cm
\noindent
{\bf (i.a)} {\it  Local Field $K_w$}:  Here Galois
representations $\rho_{w,v}: G_{K_w}\to GL_n(F_v)$ involved are
for the absolute Galois group $G_{K_w}$	 of a local $w$-adic number 
field $K_w$ with coefficients in a fixed $v$-adic number field $F_v$.
Motivated by 

\noindent
($\alpha$) Grothendieck's Monodromy Theorem for $v$-adic Galois representations of $w$-adic number fields, where $v\nparallel w$, i.e., $v$ and $w$ are with different residual characteristics; and

\noindent
($\beta$) Fontaine$\|$Berger's Monodromy Theorem for $v$-adic Galois representations of $w$-adic number fields, with $v\|w$, i.e., $v$ and $w$ are with the same residual characteristics,
\vskip 0.15cm
\noindent 
we assume that

\noindent
({\bf pST})  {\it $\rho_{w,v}$ is potentially semi-stable.}

Clearly, when $v\nparallel w$, this is equivalent to the following

\noindent
({\bf pFM\&U}) {\it $\rho_{w,v}$ is potentially of finite monodromy and unipotent.}

\noindent
In other words, we assume that there exists a finite Galois extension $L_{w'}/K_w$ such that for the induced Galois representation $\rho_{w',v}:G_{L_{w'}}\to GL_n(F_v)$, the image of the associated ramification group $I_{L_{w'}}$ is both finite and nilpotent.
\vskip 0.20cm
\noindent
{\bf (i.b)} {\it Global Field $K$}: Here Galois
representations $\rho_{K,v}: G_{K}\to GL_n(F_v)$ involved are
for the absolute Galois group $G_{K}$  of a global number 
field $K$ with coefficients in a fixed $v$-adic number field $F_v$.
Motivated by etale cohomology theory of algebraic varieties,
we assume that 

\noindent
({\bf pST}) {\it For all local completions $K_w$, the associated local 
$v$-adic representations $\rho_{w,v}: G_{K_w}\to GL_n(F_v)$ 
satisfies condition {\bf pST} of (i.a);} and

\noindent
({\bf Unr}) {\it For almost all $w$, the associated $v$-adic representations
$\rho_{w,v}: G_{K_w}\to GL_n(F_v)$ are unramified.}
\vskip 0.20cm
\noindent
{\bf (ii)  Adelic Galois Representations} for
\vskip 0.20cm
\noindent
{\bf (ii.a)} {\it Local Field $K_w$}: Here Galois
representations $\rho_{w,\mathbb A_F}: G_{K_w}\to GL_n(\mathbb A_F)$ involved are
for the absolute Galois group $G_{K_w}$	 of a $w$-adic number 
field $K_w$ with coefficients in the adelic space $\mathbb A_F$ 
associated to a number field $F$. Continuity of $\rho_{w,\mathbb A_F}$ 
proves to be too loose. Stronger algebraic condition should be imposed. 
Motivated by Grothendieck's  etale cohomology theory of 
algebraic varieties, and 
Deligne's solution to the Weil conjecture when $v\nparallel w$, 
together with  Katz-Messing's modification when $v\|w$,
we assume that 

\noindent
({\bf Unr}) {\it For almost all $v$ (in coefficients), the associated $v$-adic representation $\rho_{w,v}:G_{K_w}\to GL_n(F_v)$ are unramified;} and 

\noindent
({\bf Inv}) {\it For all $v$, i.e., for $v$ satisfying either $v\|w$ or 
$v\nparallel w$, 
the associated characteristic polynomials of the Frobenius induced from
 $\rho_{w,v}$ are the same, particularly, independent of $v$.} 

\noindent 
We call such a representation a {\it thick one}, as the invariants do not depend on the coefficients chosen.
\vskip 0.20cm
\noindent
{\it Remark.} The compatibility conditions stated here are standard. (See e.g. [Se2], [Hi],  [Tay].)
However, from our point of view, the {\bf Inv} condition appears 
to be to practical
-- yes, it is very convenient and extremely useful to impose 
the independence for the associated characteristic polynomials 
of Frobenius; on the other hand, 
this independence should not be the cause 
but rather an ultimate	goal. In other words, 
it would be much better if the Inv condition can be replaced by 
other principles, e.g., certain compatibility from class field theory. (See e.g., [Kh1,2,3].) 
We leave the details to the reader.

\vskip 0.2cm
\noindent
{\bf (ii.b)} {\it Global Field $K$}: 
Here Galois
representations $\rho_{K,\mathbb A_F}: G_{K}\to 
GL_n(\mathbb A_F)$ involved are
for the absolute Galois group $G_{K}$  of a global number 
field $K$ with coefficients in the adelic space $\mathbb A_F$ 
associated to a number field $F$. 
As above, only continuity of $\rho_{w,\mathbb A_F}$ 
appears to be too weak to get a good theory. 
Much stronger algebraic conditions should be imposed.
Certainly, there are two different directions to be considered, 
namely, the horizontal one consisting of places $w$ of $K$, 
and the vertical one consisting of places $v$ of coefficients field $F$. From ii.a),  we assume that

\noindent
({\bf Comp}) {\it For every fixed  place $w$ of $K$, the induced representation
$\rho_{w,\mathbb A_F}:G_{K_w}\to\mathrm{GL}_n({\mathbb A_F})$ forms a compatible system.} 

\noindent
As such, the corresponding theory is a thick one. Hence, 
by {\bf Inv}, we are 
able to select good representatives for $\rho_{w,\mathbb A_F}$, e.g., the induced
$\rho_{w,v}:G_{K_w}\to\mathrm{GL}_n({F_v})$ where $v\|w$.
In this language, we then further assume that
the admissible conditions for the other direction $v$ can be read 
from these selected $\rho_{w,v}$, $v\|w$. More precisely,
we assume that

\noindent
({\bf dR}) {\it All  $\rho_{w,v}$, $v\|w$, are of de Rham type;}

\noindent
({\bf Crys}) {\it For almost all $w$ and $v$,  $\rho_{v,w}$ are crystalline.}

For this reason, we may form what we call the {\it anleric ring}
$$\mathbb B_{\mathbb A}:=\prod\,'\Big(\mathbb B_{\mathrm{dR}},\mathbb B_{\mathrm{crys}}^+\Big),$$ where $\mathbb B_{\mathrm{dR}}$ denotes the ring of de Rham periods, and 
$\mathbb B_{\mathrm{crys}}^+$ the ring of crystalline periods, and 
$\prod'$ means the restricted product. As such, the final global condition
we assume is the following:

\noindent
({\bf Adm}) {\it $\{\rho_{w,v}\}_{v\|w}$ are $\mathbb B_{\mathbb A}$-admissible.}

\noindent
Even this admissibility is not clearly stated due to \lq the lack of space',
which will be discussed in details elsewhere,
one may sense  it say via determinant formalism from 
abelian CFT, (see e.g., the reformulation by Serre for rank one case ([Se2]) and
the conjecture of Fontaine-Mazur on 
geometric representations ([FM], see also [Tay]).
For the obvious reason, we will call such a representation a 
{\it thin} one.
\vskip 0.30cm
With the restrictions on Galois representations stated, 
let us next turn our attention to their characterizations. 
Here by a characterization, 
we mean a certain totally independent but intrinsic structure from which
the original Galois representation can be reconstructed.	
 There are two different approaches,  analytic one and	algebraic one. 
\vskip 0.20cm
\noindent
$\bullet$ {\bf Analytic One} Here the objects seeking are supposed to be equipped with analytic structures such as  connections and residues (at least for $v$-adic representations).  Good examples are the related works of Weil on flat bundles, of Seshadri on logarithmic
unitary flat  bundles, and of Dwork on $p$-adic differential equations;

\noindent
$\bullet$ {\bf Algebraic One} Here the structures involved are supposed to be purely alebraic. Good examples are Mumford's semi-stable bundles, Seshadri's parabolic bundles, Fontaine's various rings of periods, and semi-stable filtered $(\varphi,N;\omega)$-modules. 
We will leave the details to the main text. Instead, let me point out that 
for  local theories, when $l\not=p$, 
we should equally have	$l$-adic analogues
$\mathbb B_{\mathrm{total}}, \,\mathbb B_{\mathrm{pFM\&N}},\,
\mathbb B_{\mathrm{ur}}$
of Fontaine's $p$-adic ring of de Rham, semi-stable, crystalline periods, 
namely, $\mathbb B_{\mathrm {dR}},  \mathbb B_{\mathrm {st}},\, \mathbb B_{\mathrm {crys}}$, respectively.  Practically, this is possible due to the following reasons.

\noindent
$\bullet$  {\it Hodge-Tate Filtration}: Since every $l$-adic representation, 
$l\not=p$, is geometric. Hence, it can be realized in terms of etale 
cohomology over which by the comparison theorem there is a natural 
Hodge-Tate filtration structure;

\noindent
$\bullet$ {\it Monodromy Operator}: This is a direct consequence of Grothendieck's Monodromy Theorem for $l$-adic Galois Representations;

\noindent
$\bullet$ {\it Frobenius}, or equivalently, {\bf Dieudonne Filtration}: 
This should be put into the context that Weil's conjecture works in both 
$l$-adic and $p$-adic settings mentioned above;

\noindent
$\bullet$ {\it Ramifications}, or equivalently, {\bf $\omega$-structures}: 
This may be read from the so-called theory of breaks and conductors for $l$-adic Galois representations. For details see e.g., the main text and Chapter 1 of [Ka2].

\noindent
To uniform the notation,
denote the corresponding rings of periods 
in both $l$-adic theory and $p$-adic theory 
by  $\mathbb B_{\mathrm{dR}}, \,\mathbb B_{\mathrm{st}},\,\mathbb B_{\mathrm{ur}}^+$. Accordingly, for adelic representations of local fields,  we then can formulate a huge {\it anleric ring} $\mathbb B_{\mathbb A}:=\prod\,'\Big(\mathbb B_{\mathrm{dR}},\mathbb B_{\mathrm{ur}}^+\Big),$ {\it of	 adelic periods},  namely, the restricted product of $\mathbb B_{\mathrm{dR}}$ with respect to
$\mathbb B_{\mathrm{ur}}^+$. In this language, the algebraic condition
for thin adelic Galois representations of global fields along with 
the vertical direction may also be stated as:
\vskip 0.20cm
\noindent
({\bf Adm}) {\it It is $\mathbb B_{\mathbb A}$-admissible.}

\vskip 0.30cm
\noindent
{\bf II) Large Scale} 
\vskip 0.20cm
\noindent
A characterization of each individual Galois representation in terms of
pure algebraic structures
 may be called	a Micro Reciprocity Law, MRL for short,	 as it exposes
an intrinsic connection between Galois representations and certain 
algebraic aspects of the base fields. Assuming such a MRL, 
we then are in a position to  understand the
mathematics involved in a global way.
There are also two different approaches, at least when the coefficients are local. Namely, the categorical 
theoretic one,
based on the fact that Galois representations selected automatically
form a Tannakian category,
 and the moduli theoretic one, based on the fact that the associated 
 algebraic structures admit GIT stability interpretations.
 (In the case when the coefficients are global adelic spaces, the 
 existing standard
Tannakian category theory and  GIT should be extended.
Indeed, as pointed out by Hida, it is already an interesting problem to see whether our restricted adelic Galois representations form a Tannakian category: After all, the forgetful functor now is not
to the category of finite vector spaces over local fields but to that of
 adelic spaces.)
\vskip 0.30cm
\noindent
$\bullet$ {\bf Tannakian Categories} The main aim here is to offer a general Class Field Theory, CFT for short, for the associated base field. Roughly speaking, this goes as follows, at least when the coefficients are local fields.
With the Micro Reciprocity Law,
we then can get a clone Tannakian category, consisting of certain intrinsically defined pure algebraic objects associated to the base fields, for the Tannakian category
consisting of selected	Galois representations. As a direct consequence of the finite
monodromy and nilpotence, using the so-called finitely generated sub-Tannakian categories and automorphism groups of the associated restrictions of the fiber functors,
 one then can establish an existence theorem and a global reciprocity law for all finite (non-abelian) extensions of the base fields
 so as to obtain a general CFT for them.
As one may expect here,	 much refined results can be obtained. Indeed,
 via a certain truncation process, not only the associated Galois groups 
 but the whole system of high ramification groups can be reproduced. 
For details, see  Part C. 
\vskip 0.30cm
\noindent
$\bullet$ {\bf Moduli Spaces} From the MRL, Galois representations
selected can be characterized by intrinsically defined algebraic structures 
associated to based fields. These algebraic structures are further expected to be able to put together to form well-controlled moduli spaces. Accordingly, we have certain geometric objects to work with.	The importance of such geometric spaces can hardly be overestimated since,
 with such spaces,
we can introduce intrinsic (non-abelian) invariants for the base fields. Good examples are high rank zeta functions and their associated abelian parts. For details, see Part B.
\vskip 0.30cm
To achieve this, we clearly need to have a good control of objects 
selected. As usual, this is quite delicate:
If the selection is too restrictive, then there might not be enough
information involved; on the other hand, it should not be too loose, as 
otherwise,  it is too complicated to see structures in a neat manner, even 
we know many things are definitely there. (The reader can sense this
from our current studies of the Langlands Program.)
It is for the purpose of overcoming such difficulties 
that we introduce the following
\vskip 0.30cm
\noindent
{\bf Key} to the Success: {\bf Stability} 
\vskip 0.20cm
\noindent
This is supposed to be a  condition which helps us to make {\it good 
selections} and hence to get nice portions among all possibilities.	
Particularly, for the algebraic objects selected, we then expect
to establish a general MRL  (using them) so that the Tannakian 
category formalism can be applied and a general CFT can be established;
and to construct moduli spaces (for them) so that intrinsic invariants 
can be introduced naturally. This condition is {\it Stability}.
In accordance with what said above, 
as a general principle of selection,  the condition of stability then should be

\noindent
(a) algebraic, (b) intrinsic, and (c) rigid, 

\noindent
so that, with it, we can

\noindent
(i)  have a nice characterization of Galois representations in terms of semi-stable algebraic structures;

\noindent
(ii) form a Tannakian category for these semi-stable objects; and hence

\noindent
(iii) construct natural moduli spaces.

\noindent
Good examples are for (parabolic) bundles, filtered $(\varphi,N;\omega)$-modules, etc.
For details, please see Parts A, B, and C in the main text.
\vskip 0.30cm
This paper consists of three parts. In Part A, we indicate how a general non-abelian CFT for Riemann surfaces can be established using Tannakian category theory based on Seshadri's work on semi-stable parabolic bundles. This serves as a general guidance for our discussions in later parts. In Part B, we, motivated by yet another CFT, the conformal field theory, for Riemann surfaces, discussed in Part A, make an intensive
study on non-abelian invariants, namely, the high rank
zetas for global fields defined using stability. Along with the course, we give a geometric characterization for rank two semi-stable lattices using generalized Siegel type distances between moduli points and cusps, an analytic characterization of stability using Arthur's truncation, and a definition of general non-abelain $L$-functions using Langlands' theory of  Eisenstein series and spectral decompositions. 
In addition, we also briefly recall abelian zetas associated to $(G,P)$, with $G$ reductive 
groups and $P$ their maximal parabolic groups, which may be viewed as abelian parts of our non-abelian zetas. These abelian parts, naturally related with constant terms of 
Eisenstein series are expected to help us to understand
the hidden role played by symmetry in the Riemann Hypothesis.
Finally, in Part C, we outline a program aiming at establishing a general CFT for $p$-adic number fields. 
Key points are the notion of semi-stable filtered $(\varphi,N;\omega)$-modules of slope zero and a conjectural 
Micro Reciprocity Law claiming that there is a natural one-to-one and onto correspondence between de Rham 
representations and semi-stable modules of slope zero. 
Key ingredients of Fontaine's theory of $p$-adic Galois representations are recalled as well.
\vskip 12.0cm
\noindent
{\bf Acknowledgements}: I would like to thank Deninger and
Hida for their keen interests and huge supports during the long periods of
preparations of this paper: our visits to M\"unster in Sept-Oct 2004,6,8	
were very crucial 
to the studies of zetas explained in Part B; and a series personal notes 
on General CFT written at UCLA in March-April, 2007,8,9 is essential to
Part C. Special thanks also due to anonymous referees for their careful 
readings and detailed suggestions.

\noindent
This work is partially supported by JSPS.
\eject
\centerline{\Large{\bf Part A. Guidances from Geometry}}
\vskip 1.0cm
\centerline{\Large\bf Chapter I. Micro Reciprocity Law in Geometry}

\section{Narasimhan-Seshadri Correspondence}
\subsection{Uniformization}
Let $M$ be a compact Riemann surface of genus $g$ and $M^o\hookrightarrow M$ a punctured Riemann surface with $M\backslash M^o:=\{P_1,P_2,\dots,P_N\}$. Assume that $2g-2+N>0$ so that by uniformization theorem there exists a Fuchsian group of first type $\Gamma\subset \mathrm{PSL}(2,\mathbb R)$ and the associated universal covering map $$(\pi^o:\frak H\to \Gamma\backslash \frak H\simeq M^o)	
\hookrightarrow (\pi:\frak H^+\to \Gamma\backslash \frak H^+\simeq M)$$ where $\frak H$ denotes the usual upper half plane and $\frak H^+$ denotes the extended upper half plane, namely, $\frak H$ together with cusps associated to $(M^o,M)$, or better, to $\Gamma$.

\subsection{(Narasimhan-)Seshadri Correspondence}
Let $\rho:\pi_1(M^o,*)\to\mathrm{GL}(V)$ be a {\it unitary representation} of the fundamental group $\pi_1(M^o,*)(\simeq\Gamma)$ of $M^o$. For simplicity, assume that it is irreducible. Then we know that $\rho$ satisfies the {\it finite monodromy} property at all $P_i$'s.  This then implies that there exists a finite Galois covering $$\pi':M'\to M$$ of compact Riemann surfaces
ramified possibly at $P_i$'s such that $\rho$ naturally induces a unitary representation $$\rho':\pi_1(M',*)\to\mathrm{GL}(V)$$ of the fundamental group of the {\it compact} Riemann surface $M'$ on $V$. As such, by the uniformization theorem, we obtain a {\it unitary flat bundle} over $M'$ equipped with a natural action of the Galois group $\mathrm{Gal}(\pi')$, namely, the four-tuple 
$$\Big(M',E_{\rho'}:=\big(\pi_1(M',*),\rho'\big)\big\backslash \big(\frak H^{(+)}\times V\big),\nabla_{\rho'};
\mathrm{Gal}(\pi')\Big).$$
One checks that the $\mathrm{Gal}(\pi')$-invariants of the
direct image of the differentials of $M'$ with coefficients in $E_{\rho'}$
coincides with the logarithmic differentials on $(M,Z)$ with coefficients in $E_\rho$,
namely, $$\Big(\pi_*'\big(E_{\rho'}\otimes\Omega_{M'}^1\big)\Big)^{\mathrm{Gal}(\pi')}
=E_\rho\otimes\Omega_M^1(\log Z)$$ where $Z=P_1+P_2+\dots+P_N$
denotes the reduced branch divisor on $M$. Consequently, we then obtain a {\it logarithmic unitary flat bundle} $\big(E_\rho,\nabla_\rho(\log Z)\big)$ on the compact Riemann surface $M$. Thus by using $\mathrm{Res}_{P_i}\nabla_\rho(\log Z)$, that is, by taking {\it residues} of  logarithmic unitary connection $\nabla_\rho(\log Z)$ at $P_i$'s, we then obtain Seshadri's {\it parabolic structures} on the fibers of $E_\rho$, which is nothing but the quotient bundle $\big(\pi_1(M,*),\rho\big)\big\backslash \big(\frak H^{(+)}\times V\big)$, at punctures $P_i$'s. As such, 
an important discovery of Seshadri is that the parabolic bundle obtained then is stable of degree zero. More strikingly, the converse is correct as well. Namely, any stable parabolic bundle of degree zero can be constructed in this manner.

\section{Micro Reciprocity Law}
\subsection{Weil's Program}

This result of Seshadri, obtained with the help of Metha ([MS]), is in fact 
motivated by an earlier fundamental work of Narasimhan-Seshadri 
([NS]), which claims that
there is a natural one-to-one and onto correspondence between irreducible
unitary representations of fundamental group $\pi_1(M,*)$ of compact Riemnn surface $M$ and stable bundles of degree zero on $M$.
In this sense, Seshadri's result on parabolic bundles above is a generalization
of Narasimhan-Seshadri's work from compact Riemann surfaces to punctured Riemann surfaces, in which vector bundles are replaced by parabolic bundles.

In (algebraic) geometry, Narasimhan-Seshadri's work then leads to a natural moduli space for irreducible unitary representations for fundamental groups of compact Riemann surfaces via Mumford's Geometric Invariant Theory, GIT for short. Indeed, by Narasimhan-Seshadri's result, 
it suffices to consider that for stable bundles of degree zero. While being stable and 
of degree zero for vector bundle are  conditions in terms of intersection theory,
it can be shown that this condition is equivalent to a certain GIT-stability. As such, via GIT quotient technique of Mumford ([M]), we can naturally realize the moduli space of stable bundles of degree zero on a compact Riemann surfaces as a  quasi-projective variety.
Moreover, following GIT, a natural compactification can be made by adding the so-called semi-stable points, which in terms of bundles means (Seshadri classes of) semi-stable vector bundles of degree zero. As Seshadri class corresponding naturally to equivalence class of unitary representations of  fundamental group of the compact Riemann surface in question (modulo unipotency, or better after taking semi-simplification), this then gives also an algebraic construction for moduli spaces of these 
 representations of fundamental groups.

However, moduli spaces of semi-stable bundles of degree zero over  compact Riemann surfaces in general are singular. It was once a central problem to resolve these singularities in a natural manner. In terms of what was happened, there were in fact two different approaches, one of which due to Seshadri. It is this work of Seshadri that leads to the notion of parabolic bundles.

Before the notion of parabolic bundles, Seshadri also studied the so-called $\pi$-{\it bundles} ([S2]), a notion introduced by Grothendieck ([G]). In particular, 
Seshadri's main discovery may be stated as that there is a natural one-to-one and onto correspondence between the so-called $\pi$-bundles and bundles with parabolic structures (say, when $\pi$ is a finite ramified covering).
For more details, see e.g., Biswas related work on orbifold bundles and parabolic bundles ([Bis]).
\vskip 0.20cm
Despite their huge successes in (algebaric) geometry, these fundamental works on stability have not made any serious 
impact in arithmetic (see however Nori's basic work ([Nor]) on fundamental groups via 
Tannakian category, 
even in which stability plays no role): until the time around the beginning of 90's of last century, 
the above works had been largely ignored by mathematicians working in arithmetic.
This is in fact very much unfortunate and shows us how interesting mathematics is exposed 
as human being's activities. By contrast, as we now know, not just as a result 
these works play a central 
role in establishing a general non-abelian class field theory for Riemann surfaces, 
or the same, for function fields over complex numbers, 
but, all these works  are  generalizations of 
Weil's pioneer	work  claiming that the assignment
$\rho\leftrightarrow E_\rho$ (resp.  $\rho\leftrightarrow (E_\rho,\nabla_\rho)$) gives a canonical one-to-one and onto correspondence between irreducible representations of fundamental groups
of compact Riemann surfaces and indecomposible degree zero bundles	
(resp. and indecomposible flat bundles) on the associated Riemann 
surfaces. And in history, it was 
\vskip 0.20cm
\noindent
(i) aiming at establishing a general CFT for Riemann surfaces that 
motivated Weil to prove such a result in his master piece on generalization of 
abelian functions ([We1]); 
And 

\noindent
(ii) clearly with arithmetic applications in mind that Grothendieck gave a Bourbaki seminar explaining Weil's work in which
the notion of $\pi$-bundles was introduced ([G]). 
\vskip 0.20cm
This unfortunate situation has been graduatelly changed. Say, 
at the end of 90's, There was a short note [W1]. This note is a rediscovery of
Weil's program, starting with a crucial observation that the above correspondences of Weil, 
Narasimhan-Seshadri and Seshadri can be viewed as a kind of reciprocity law; 
after all, 

\noindent
(a) the correspondences are relating fundamental groups (reading as analogue of Galois groups) 
with certain intrinsic algebraic structures (reading as
non-abelian analogues and generalizations of ideal classes); and 

\noindent
(b) by using parabolic structures, ramification information can be taken care of completely. 

\noindent
Along with such a line, naturally, these works on stability then further leads to the part of 
our Program aiming at establishing a general CFT for various fields (using stability) [W1].

\subsection{Micro Reciprocity Law}

Seshadri's fundamental works may be summarized as the follows.
\vskip 0.20cm
\noindent
{\bf Theorem.} {\it Let $(M^0,M)$ be a punctured Riemann surface. Then we have}
\vskip 0.20cm
\noindent
(i) {\bf Micro Reciprocity Law} ((Weil, Mumford, Narasimhan-Seshadri,) Seshadri)
\vskip 0.20cm

{\it There exists a natural one-to-one and onto correspondence} 
\vskip 0.20cm
\centerline{$\Big\{\, irreducible\ unitary\ representations\ of\ \pi_1(M^o,*)\,\Big\}$}

\centerline{$\Updownarrow$}

\centerline{$\Big\{\,stable\ parabolic\ bundles\ of\ degree\ zero\ on\ (M^o,M)\,\Big\};$}
\eject
\vskip 0.20cm
\noindent
(ii) {\bf Ramifications versus Parabolic Structures} ((Grothendieck), Seshadri) 

{\it There exists a natural one-to-one and onto correspondence

\centerline{$\Big\{\,vector\ bundles\ W/M'\ with\ compatible\ action\ Gal(M'/M)\,\Big\}$}

\centerline{$\Updownarrow$}

\centerline{$\Big\{\,parabolic\ bundles\  E_*/(M^o,M)\ with\ compatible\ parabolic\ weights\,\Big\}$}

\noindent
such that 

\noindent
(i) the correspondence induces a natural one on sub-objects of $W$ and of $E_*(W)$;
and 

\noindent
(ii) the degrees satisfy the relation 
$$\qquad\mathrm{deg}(W)=\mathrm{deg}(M'/M)\cdot\mathrm{par.deg}\big(E_*(W)\big).$$}

\vfill
\eject
\centerline{\Large\bf Chapter II. CFTs in Geometry}
\vskip 0.80cm
\section{Arithmetic CFT: Class Field Theory}
Building on the above detailed micro study of individual representation of fundamental groups 
of Riemann surfaces and hence individual semi-stable parabolic bundle, we can 
 study them from a more global point of view.
There are two approaches, one using category theory and the other using moduli theory.

As a starting point of the category approach,  let us first consider the 
category consisting of
semi-stable parabolic bundles of (parabolic) degree zero over $(M^o,M)$.
Note that, as building blocks of general semi-stable objects, stable ones are very rigid. 
That is to say, there is no non-trivial morphisms between two stable objects, 
a fact corresponding to Schur's Lemma in representation
theory (for irreducible representations). Consequently, 
we conclude that
the just formed category admits	 much finer structures: 
It is clearly abelian, has a tensor product
structure and admits a natural functor $\mathbb F$ to the category of finite dimensional vector spaces (the fibers of base bundles to a fixed point of $M^o$). Thus, from the rigid properties mentioned above, which guarantees the faithfulness  of the functor just mentioned, we see that the category is in fact Tannakian. Denote it by $\Big(\mathbb{PV}_{M^o,M}^{\mathrm{ss};0};\mathbb F\Big)$. Then, 
from Tannakian category theory, 
we obtain the following main theorem of CFT for Riemann surfaces, or 
the same, for function fields over complex numbers;
\vskip 0.20cm
\noindent
{\bf Main Theorem of Arithmetic CFT} ([W1])

\noindent
$\bullet$ ({\bf Existence}) {\it There exists a canonical one-to-one and onto correspondence

\centerline{$\Big\{Finitely\,Generated\,SubTannakian\,Cats\,
\big(\Sigma,\mathbb F|_{\Sigma}\big)\,of\,\Big(\mathbb{PV}_{M^o,M}^{\mathrm{ss};0};
\mathbb F\Big)\Big\}$}

\centerline{$\Big\Updownarrow\ \Pi$}

\centerline{$\Big\{\,Finite\ Galois\ Coverings\ M'\to (M^o,M)\,\Big\}$}

\noindent
  which induces naturally an isomorphism}
	
\noindent
$\bullet$ ({\bf Reciprocity Law})
$$\mathrm{Aut}^{\otimes}\Big(\Sigma,\mathbb F|_{\Sigma}\Big)\simeq
\mathrm{Gal}\Big(\Pi(\Sigma,\mathbb F|_{\Sigma})\Big).$$

\section{Geometric CFT: Conformal Field Theory}

Here we give some most important aspects of the second global approach, namely the one using 
moduli spaces. As a starting point, for a fixed compact Riemann surface $M$, denote by 
$\mathcal M_M(r,0)$  the moduli spaces of rank $r$ semi-stable bundles of 
 degree zero on $M$. (Recall that then we squeeze semi-stable bundles into their associated 
 Seshadri classes, defined using graded pieces of the associated Jordan-H\"older filtrations.)
Over such moduli spaces, we can construct many global invariants. Analytically we may expect 
that a still ill-defined Feymann integral would give us something interesting. We will not 
pursue this line further, instead, let us start with an algebraic construction.

Since each moduli point corresponds to a semi-stable vector bundle, it makes sense to 
talk about the associated cohomology groups. As such, then we may form the so-called 
Grothendieck-Mumford determinant line of cohomologies, i.e., the alternative tensor products 
of determinants of cohomologies. Consequently, if we move our moduli points over all moduli spaces, 
we can glue the above determinant lines to form the so-called Grothendieck-Mumford determinant line 
bundles $\lambda_M$ on $\mathcal M_M(r,0)$.
Note that the Picard group of $\mathcal M_M(r,0)$ is isomorphic to $\mathbb Z$, we see that a suitable multiple of $\lambda_M$ is indeed very ample. (For all this, we in fact need to restrict ourselves only to the stable part. Let us assume it was the case now while leaving the details on how to fix it
to the literatures, or better to the reader.) It then makes sense to talk about the 
$\mathbb C$-vector space
$H^0\big(\mathcal M_M(r,0),\lambda_M^{\otimes n}\big)$ 
(for $n$ sufficiently away from 0).
This is a finite dimensional vector space naturally associated to $M$, whose
dimension is given by the so-called Verlinde formula.

The most interesting and certain a very deep point is somehow we expect that
the space itself $H^0\big(\mathcal M_M(r,0),\lambda_M^{\otimes n}\big)$, also called 
{\it conformal blocks}, 
does not really very much related with the complex structure on the compact Riemann surface $M$ 
used. (For more details on Conformal Field Theory, intitated by Belavin-Polykov-Zamodolochikov, see
e.g., [US].)
 More precisely,
let us now move $M$ in $\mathcal M_g\hookrightarrow \overline{\mathcal M_g}$, the moduli space of compact Riemann surfaces of genus $g=g(M)$
 and its stable compactification of  Deligne-Mumford ([DM]). Denote by $\Delta_{\mathrm{bdy}}$
  the boundary of $\mathcal M_g$, which is a normal crossing divisor by Deligne-Mumford theory.
  Then the conformal blocks form a natural vector bundle 
$\Pi_*\Big(\lambda_M^{\otimes n}\Big)\Big|_{\mathcal M_g}$
on $(\mathcal M_g\hookrightarrow) \overline{\mathcal M_g}$, with which, we may state the following: 
\vskip 0.20cm
\noindent
{\bf Main Theorem in Geometric CFT:} (Tsuchiya-Ueno-Yamada, see also [Hi])
{\it There exists a projectively flat logarithmic connection on the bundle 
$\Pi_*\Big(\lambda_M^{\otimes n}\Big)\Big|_{\mathcal M_g}$ over
$(\mathcal M_g,\Delta_{\mathrm{bdy}})$.}

\vfill
\eject
\centerline{\Large{\bf Part B. High Rank Zeta Functions	 and Stability}}
\vskip 0.80cm
\centerline{\Large\bf Chapter III. High Rank Zeta Functions}
\section{Function Fields}
\subsection{Definition and Basic Properties}
Let $C$ be a regular, geometrically irreducible projective curve of genus $g$ defined 
over $\mathbb F_q$, the finite field with $q$ elements and $\mathcal M_{C,r}$ the moduli
space of semi-stable bundles of rank $r$ over $C$. These spaces are projective varieties.
So following Weil, we may try to attach them with the standard Artin-Weil zeta functions. However, there is another more intrinsic way. Namely, instead of simply viewing these moduli spaces 
as algebraic varieties, 
we here want to fully use the moduli aspects of these spaces
by viewing rational points of these varieties as rational bundles: This is possible at least for 
the stable part by a work of Harder-Narasimhan on Brauer groups ([HN]). 
Accordingly, for each rational moduli point, we can have a very natural weighted count. 
All this then leads to	the following
\vskip 0.20cm
\noindent
{\bf Definition.} (Weng) {\it The rank $r$ zeta function for $C/\mathbb F_q$ is defined by
$$\zeta_{C,\mathbb F_q;r}(s):=\sum_{V\in[V]\in\mathcal M_{C,r}}
\frac{q^{h^0(C,V)}-1}{\#\mathrm{Aut}(V)}\cdot\Big(q^{-s}\Big)^{\mathrm{deg}(V)},\qquad\mathrm{Re}(s)>1.$$
Here as usual, $[V]$ denotes the Seshadri class of (a rational) semi-stable bundle $V$, and $\mathrm{Aut}(V)$ denotes the automorphism group of $V$.}
\vskip 0.20cm
By semi-stable condition, the summation above is only taken over the part of moduli 
space whose points have non-negative degrees. Thus by the duality, Riemann-Roch and a Clifford 
type lemma for semi-stable bundles, we then can expose the following basic properties for our zeta functions of curves.
\vskip 0.20cm
\noindent
{\bf Zeta Facts}(Weng) (0) {\it	 $\zeta_{C,1,{\bf F}_q}(s)$ 
is nothing but the classical Artin zeta function $\zeta_C(s)$ for curve $C$.}

\noindent
(1) {\it  $\zeta_{C,r,{\bf F}_q}(s)$ 
is well-defined for ${\rm
Re}(s)>1$, and admits a meromorphic continuation to the whole complex 
$s$-plane;}

\noindent
(2) ({\bf Rationality}) {\it Set $t:=q^{-s}$ and introduce the non-abelian $Z$-function of $C$ by 
$$\zeta_{C,r,{\bf F}_q}(s)=:Z_{C,r,{\bf F}_q}(t):=
\sum_{V\in [V]\in {\mathcal M}_{C,r}(d),d\geq
0}{{q^{h^0(C,V)}-1}\over {\#{\rm Aut}(V)}}\cdot t^{d(V)}, \ \ |t|<1.$$ Then 
 there exists a polynomial $P_{C,r,{\bf F}_q}(s)\in {\bf Q}[t]$ such that 
$$Z_{C,r,{\bf F}_q}(t)={{P_{C,r,{\bf F}_q}(t)}\over {(1-t^r)(1-q^rt^r)}};$$}

\noindent
(3) ({\bf Functional Equation}) {\it Set the rank $r$
non-abelian $\xi$-function}\linebreak
$\xi_{C,r,{\bf F}_q}(s)$ {\it by 
$$\xi_{C,r,{\bf F}_q}(s):=\zeta_{C,r,{\bf F}_q}(s)\cdot (q^{s})^{r(g-1)}.$$	
Then
$$\xi_{C,r,{\bf F}_q}(s)=\xi_{C,r,{\bf F}_q}(1-s).$$}

\noindent
{\it Remarks.} (1) ({\bf Count in Different Ways}) The above weighted count is designed for 
all rational semi-stable bundles, motivated by Harder-Narasimhan's interpretation on Siegel's 
work about Tamagawa numbers ([HN]). As such, even the moduli space is used, 
it does not really play a key role as all elements in a Seshadri class are counted. 
For this reason,  modifications for the definition of high rank zetas 
can be given, say, count only one within a fixed Seshadri class, or count only what are called 
strongly semi-stable bundles, etc...

\noindent
(2) ({\bf Stratifications and Cohomological Interpretations}) Deninger once asked whether there 
was a cohomological interpretation
for our zeta functions. There is a high possibility for it: We expect that our	
earlier works on refined Brill-Noether loci would play a key role here, since 
refined Brill-Noether loci induce  natural 
stratifications on moduli spaces. Thus, following Grothendieck's work on
cohomological interpretation of Weil's zeta functions,
 what we have to do next is to expose a certain weighted fixed point formula.

\subsection{Global High Rank Zetas via Euler Products}

Let ${\mathcal C}$ be a regular, reduced, irreducible projective curve of genus
$g$ defined over a number field $F$. Let $S_{\rm bad}$ be the collection 
of all infinite places and
these finite places of $F$ at which ${\mathcal C}$  does not have good 
reductions. 
As usual, a place $v$ of $F$ is called good if $v\not\in S_{\rm bad}$.
For any good place $v$ of $F$,	 the $v$-reduction of 
${\mathcal C}$, denoted as
${\mathcal C}_v$, gives a regular, reduced, irreducible projective curve 
defined 
over the residue field
$F(v)$ of $F$ at $v$. Denote the cardinal number of $F(v)$ by $q_v$.	
Then,  we obtain
the associated rank $r$ non-abelian zeta function
$\zeta_{{\mathcal C}_v,r,{\bf F}_{q_v}}(s)$. Moreover, from the rationality of 
$\zeta_{{\mathcal
C}_v,r,{\bf F}_{q_v}}(s)$, there exists a degree $2rg$ polynomial 
$P_{{\mathcal C}_v,r,{\bf F}_{q_v}}(t)\in {\bf Q}[t]$
such that
$$Z_{{\mathcal C}_v,r,{\bf F}_{q_v}}(t)={{P_{{\mathcal C}_v,r,{\bf F}_{q_v}}(t)}\over {(1-t^r)(1-q^rt^r)}}.$$
Clearly, $P_{{\mathcal C}_v,r,{\bf F}_{q_v}}(0)\not=0.$ Set
$$\tilde P_{{\mathcal C}_v,r,F(v)}(t):={{P_{{\mathcal C}_v,r,F(v)}(t)}\over
 {P_{{\mathcal C}_v,r,F(v)}(0)}}.$$
 
\noindent
{\bf Definition.} (Weng) {\it The rank $r$ non-abelian 
zeta function $\zeta_{{\mathcal C},r,F}(s)$ of ${\mathcal C}$ over $F$ is defined as the 
following Euler product}
$$\zeta_{{\mathcal C},r,F}(s)
:=\prod_{v:{\rm good}}{1\over{
\tilde P_{{\mathcal C}_v,r,{\bf F}_{q_v}}(q_v^{-s})}},\hskip 2.0cm 
{\rm Re}(s)\gg 0.$$

Clearly, when $r=1$, $\zeta_{{\mathcal C},r,F}(s)$ coincides with the 
classical Hasse-Weil zeta function
for $C$ over $F$.
\vskip 0.20cm
\noindent
{\bf Conjecture.} {\it For a regular, reduced, geometrically irreducible projective 
curve ${\mathcal C}$ of
genus $g$ defined over a number field $F$,  its associated rank 
$r$ global non-abelian
zeta function
$\zeta_{{\mathcal C},r,F}(s)$  admits a meromorphic continuation to the whole 
complex $s$-plane.}
\vskip 0.20cm
\noindent
Recall that even  when $r=1$, i.e., for the classical Hasse-Weil zeta 
functions, this statement, as a part of a series of high profile conjectures
is still open. On the other hand, we have the following
\vskip 0.30cm
\noindent
{\bf Proposition.} ([W4]) {\it	When ${\rm Re}(s)> 1+g+(r^2-r)(g-1)$, 
$\zeta_{{\mathcal C},r,F}(s)$ converges.}
\vskip 0.20cm
Like in the theory for abelian zeta functions, we want to use our
non-abelian zeta functions  to study  non-abelian aspect of
arithmetic of curves. For this purpose, completed zetas, or better, local factors for \lq bad' 
places, should be introduced:	

\noindent
(i) For $\Gamma$-factors, motivated by the local rationality, we take these associated to
 $\zeta_F(rs)\cdot\zeta_F\big(r(s-1)\big)$, where $\zeta_F(s)$ denotes
the standard Dedekind zeta function for $F$; and 

\noindent
(ii) for finite bad factors, first choose a semi-stable model 
for ${\mathcal C}$ so as to get a semi-stable reduction for curves at bad places.
Then, either
\noindent 
(a) use Seshadri's moduli spaces of semi-stable parabolic bundles as suggested in [W4]; or

\noindent
(b) use moduli space of semi-stable bundles over nodal curvces, as pointed out by Seshadri.

For the time being, even we know that each produces local factors for 
singular fibers, usually polynomials with degree lower than $2rg$, but
we do not know which one is right. To test them,
we propose the following 
functional equation.
\vskip 0.30cm
\noindent
{\bf Working Hypothesis.} {\it The completed zeta function 
$\xi_{{\mathcal C},r,F}(s)$ of ${\mathcal C}/F$ admits a unique
meromorphic continuation to the whole complex $s$-plane and satisfies 
the functional equation 
$$\xi_{{\mathcal C},r,F}\big(s\big)=\varepsilon\cdot 
\xi_{{\mathcal C},r,F}\Big(1+{1\over r}-s\Big)$$ with $|\varepsilon|=1$.}

\section{Number Fields}

\subsection{Stability of $\mathcal O_F$-Lattices}

Let $F$ be a number field with $\mathcal O_F$ the ring of integer and $\Delta_F$ the 
discriminant. By definition, an $\mathcal O_F$-{\it lattice} $\Lambda$ of rank $r$ 
consists of a pair $(P,{\bf \rho})$, where $P$ is a rank $r$ projective $\mathcal O_F$-module 
and ${\bf \rho}$ is a metric on the space $\big(\mathbb R^{r_1}\times\mathbb C^{r_2}\big)^r=
 \big(\mathbb R^r\big)^{r_1}\times\big(\mathbb C^r\big)^{r_2}$, where $r_1$ (resp. $r_2$) 
 denotes the number of real embeddings (resp. complex embeddings) of $F$. 
 Recall that,  being projective, there exists a fractional idea $\frak a$ of $F$ such that $P\simeq\mathcal O_F^{r-1}\oplus\frak a$. Particularly, the natural inclusion $O_F^{r-1}\oplus\frak a\hookrightarrow F^r$ induces a natural embedding of $P$ into $
\big(\mathbb R^{r_1}\times\mathbb C^{r_2}\big)^r$ via the compositions
$$P\simeq O_F^{r-1}\oplus\frak a\hookrightarrow F^r\hookrightarrow
\Big(\mathbb R^{r_1}\times\mathbb C^{r_2}\Big)^r\simeq \big(\mathbb R^r\big)^{r_1}\times
\big(\mathbb C^r\big)^{r_2}.$$ As such, then the image of $P$  naturally offers us a lattice 
$\Lambda$ in the metrized space $\Big(\big(\mathbb R^r\big)^{r_1}\times\big(\mathbb C^r\big)^{r_2},
{\bold \rho}\Big)$.

By definition, an $\mathcal O_F$-lattice is called {\it semi-stable}
if for all sub-$\mathcal O_F$-lattice $\Lambda_1$ of $\Lambda$, we have
$$\mathrm{Vol}(\Lambda_1)^{\mathrm{rank}(\Lambda)}\geq \mathrm{Vol}(\Lambda)^{\mathrm{rank}(\Lambda_1)},$$ where
the volume $\mathrm{Vol}(\Lambda)$ of $\Lambda$ is usually called the covolume of $\Lambda$, namely,
$$\mathrm{Vol}(\Lambda):=\mathrm{Vol}\Big(\Big(\big(\mathbb R^r\big)^{r_1}\times\big(\mathbb C^r\big)^{r_2},{\bold \rho}\Big)\Big/\Lambda\Big).$$

Denote by $\mathcal M_{F,r}$ the moduli space of semi-stable $\mathcal O_F$ lattices of rank $r$,
i.e., the space of isomorphism classes of semi-stable $\mathcal O_F$ lattices of rank $r$. 
Then there is a natural topological structure on $\mathcal M_{F,r}$. In fact there is a much 
finer structure on it;
Denote by $\mathcal M_{F,r}[T]$ the volume $T$ part of $\mathcal M_{F,r}$, i.e., 
the part consisting of isomorphisms classes of rank $r$ semi-stable $\mathcal O_F$-lattices 
of volume $T$, then 
\vskip 0.20cm
\noindent
(i) there is a natural decomposition 
$$\mathcal M_{F,r}=\bigcup_{T\in\mathbb R_{>0}}\mathcal M_{F,r}[T];$$ Moreover, 
\vskip 0.20cm
\noindent
(ii) for each fixed $T$, $\mathcal M_{F,r}[T]$ is compact; and 
\vskip 0.20cm
\noindent
(iii) there are natural measures $d\mu$ on $\mathcal M_{F,r}$ such that 
$$d\mu=d\mu\Big|_{\mathcal M_{F,r}[|\Delta_F|^{\frac{r}{2}}]}\times \frac{dT}{T}.$$
(The compactness of $\mathcal M_{F,r}[T]$ is the main reason why we use the stability 
condition in the study of non-abelian zetas in [W5].)

\subsection{Geo-Arithmetical Cohomology}

Let $\Lambda$ be an $\mathcal O_F$-lattice. Then define its {\it geo-arthmetical
cohomology groups} by 
$$H^0(F,\Lambda):=\Lambda,\qquad\mathrm{and}\qquad H^1(F,\Lambda):=\Big(\mathbb R^{r_1}\times\mathbb C^{r_2}\Big)^r\Big/\Lambda.$$ As such, unlike in algebraic geometry and/or in arithmetic geometry, cohomological groups $H^i$ are not vector spaces, but locally compact topological groups. 

With this simple but genuine definition, then the basic properties such as the duality and the Riemann-Roch theorem can be realized as follows;
\vskip 0.20cm
\noindent
{\bf Pontrjagin Duality} (Weng) {\it There is a natural topological isomorphism
$$H^1(F,\Lambda)\simeq\widehat{H^0(F,\omega_F\otimes\Lambda^\vee)}
$$ where $\omega_F:=(\frak d_F,\rho_{\mathrm{st}})$ denotes the differential lattice of $F$, 
namely, the (rank one) projective module given by the standard differential module $\frak d_F$ of $\mathcal O_F$, and the metric given by the standard metric $\rho_{\mathrm{st}}$ on $\mathbb R^{r_1}\otimes\mathbb C^{r_2}$}.
\vskip 0.20cm
Moreover, since $H^{i=0,1}(F,\Lambda)$ are locally compact topological groups, we can apply 
Fourier analysis to introduce quantitive invariants
for them ([F]), say, for $h^0$, or better for $exp(h^0)$, counting each element 
${\bold x}\in H^0(F,\Lambda)$, (which is nothing but the lattice $\Lambda$ itself,) 
with weight of the Gaussian distribution	
$$e^{-\pi\sum_{\sigma:\mathbb R}\|{\bold x}\|_{\rho_\sigma}
-2\pi\sum_{\tau:\mathbb C}\|{\bold x}\|_{\rho_\tau}}.$$
(As such, this definition then coincides with the one previously introduced by 
van der Geer and Schoof, for which an arithmetic analogue of effectivity is used ([GS]).)
\vskip 0.20cm
\noindent
{\bf Geo-Arithmetical Riemann-Roch Theorem.} (Weng) 
{\it For an $\mathcal O_F$-lattice $\Lambda$,} 
$$h^0(F,\Lambda)-h^1(F,\Lambda)=\mathrm{deg}(V)-\frac{\mathrm{rank}(\Lambda)}{2}\cdot \log\big|\Delta_F\Big|.$$

Our Riemann-Roch is a direct consequence of the Fourier inversion formula, reflecting 
the topological Pointrjagin duality above, and the standard Poission summation formula. 
So it has its roots in Tate's Thesis ([Ta1]), even our result is not really there.

In the above RR, $\mathrm{deg}(V)$ denotes what we call Arakelov degree of $V$.
In fact, in Arakelov geometry, there is the following
\vskip 0.20cm
\noindent
{\bf Arakelov Riemann-Roch Theorem.} (See e.g. [L1,2,3]) 
$$-\log\Big(\mathrm{Vol}(\Lambda)\Big)=\mathrm{deg}(V)-\frac{\mathrm{rank}(\Lambda)}{2}\cdot\log\big|\Delta_F\Big|.$$

From this, it is simple to see that the above definition of ours for semi-satble
$\mathcal O_F$ lattices is equivalent to the following definition in [St1]:
an $\mathcal O_F$-lattice is semi-stable if
for all sub-$\mathcal O_F$-lattice $\Lambda_1$ of $\Lambda$, we have
$$\frac{\mathrm{deg}(\Lambda_1)}{\mathrm{rank}(\Lambda_1)}\leq\frac{\mathrm{rank}(\Lambda)}{\mathrm{rank}(\Lambda)},$$ an arithmetic-geometric analogue of the slope stability condition of Mumford for vector bundles  over compact Riemann surfaces:
A vector bundle $V$ over a compact Riemann surface $M$ is semi-stable
if for all subbundles $V_1$,
$$\frac{\mathrm{deg}(V_1)}{\mathrm{rank}(V_1)}\leq\frac{\mathrm{deg}(V)}{\mathrm{rank}(V)}.$$
\subsection{High Rank Zetas}
With the above preperation, we are ready to state the following
\vskip 0.20cm
\noindent
{\bf Definition.} (Weng) {\it The  rank $r$ zeta function of $F$ is defined by}
$$\xi_{F,r}(s):=\Big(\Big|\Delta_F\Big|^s\Big)^{\frac{r}{2}}
\cdot\int_{\mathcal M_{F,r}}\Big(e^{h^0(F,\Lambda)}-1\Big)
\cdot\Big(e^{-s}\Big)^{\mathrm{deg}(\Lambda)} d\mu(\Lambda),\ \mathrm{Re}(s)>1.$$

Tautologically, from the duality and the geo-arithmetical Riemann-Roch, we obtain 
the following standard properties for the high rank zeta functions (see however [We2]):
\vskip 0.20cm
\noindent
{\bf Zeta Facts.}  (Weng) (0) (Iwasawa) {\it $\xi_{F,1}(s)\buildrel.\over{=}\xi_{F}(s)$, the completed Dedekind zeta for $F$;}

\noindent
(1) ({\bf Meromorphic Extension}) {\it Non-abelian zeta function 
$$\xi_{F,r}(s):=\Big(|\Delta_F|^{{r\over
2}}\Big)^s\int_{\Lambda\in {\mathcal M}_{F,r}}
\Big(e^{h^0(F,\Lambda)}-1\Big)\big(e^{-s}\big)^{\mathrm{deg}(\Lambda)}
\cdot d\mu$$ converges
absolutely and uniformly when ${\rm Re}(s)\geq 1+\delta$ for any
$\delta>0$. Moreover, $\xi_{F,r}(s)$ admits a unique meromorphic continuation 
to the whole complex $s$-plane;}

\noindent
(2) ({\bf Functional Equation}) {\it The extended $\xi_{F,r}(s)$ satisfies the
functional equation  $$\xi_{F,r}(s)=\xi_{F,r}(1-s);$$}

\noindent
(3) ({\bf Singularities}) {\it The extended $\xi_{F,r}(s)$ has two singularities, 
all simple poles, at $s=0\,1$, with
$$\mathrm{Res}_{s=0}\,\xi_{F,r}(s)\,=\,-\mathrm{Res}_{s=0}\,\xi_{F,r}(s)
\,=\,\mathrm{Vol}\Big({\mathcal	 M}_{F,r}[|\Delta_F|^{r\over 2}]\Big).$$}

\vfill
\eject
\centerline{\Large\bf Chapter IV. Geometric Characterization of Stability}
\vskip 0.80cm
\noindent
Here we give an example on how to characterize stability 
in geometric terms. More precisely, 
in this chapter, we will offer a characterization of
semi-stable rank two $\mathcal O_F$-lattices in terms of a Siegel type distance to cusps.
We will present the materials in a classical way in which many fundamental results 
of algebraic number theory will be used. The main results are listed in \S 8 and \$ 9.

\section{Upper Half Space Model}
\subsection{Upper Half Plane}

As usual, denote by $${\mathcal H}:=\{z=x+iy\in \mathbb C: x\in\mathbb R, y\in\mathbb R_+^*\},$$
the upper half plane.
The group $SL(2,\mathbb R)$ naturally acts on $\mathcal H$ via: 
$$M\,z:=\frac{az+b}{cz+d},\qquad\forall 
M=\left(\begin{matrix} a&b\\ c&d\end{matrix}\right)\in 
SL(2,\mathbb R),\ \ z\in {\mathcal H}.$$
The stablizer of $i=(0,1)\in {\mathcal H}$  is equal to 
$SO(2):=\{A\in O(2):\det A=1\}$. Since this action 
 is transitive,	 we can identify
the quotient $SL(2,\mathbb R)/SO(2)$ with $\mathcal H$	by the 
quotient map induced from $SL(2,\mathbb R)\to {\mathcal H},\ 
g\mapsto g\cdot i.$

${\mathcal H}$ admits the real line $\mathbb R$ as its boundary. 
Consequently, to compactify it, we add on it the real projective line 
 $\mathbb P^1(\mathbb R)$ with $\infty=\left[\begin{matrix} 1\\ 
0\end{matrix}\right]$. Naturally, the above action of $SL(2,\mathbb R)$ 
also extends to $\mathbb P^1(\mathbb R)$ via 
 $$\left(\begin{matrix} a&b\\ c&d\end{matrix}\right)\left[\begin{matrix} x\\ 
y\end{matrix}\right]=
 \left[\begin{matrix} ax+by\\ cx+dy\end{matrix}\right].$$ 
 
\subsection{Upper Half Space}

Similarly,  3-dimensional hyperbolic space is defined to be 
$$\begin{aligned}{\mathbb H}:=&\mathbb C\times]0,\infty[\,=\,\Big\{(z,r):z=x+iy\in \mathbb C, 
r\in\mathbb R_+^*\Big\}\\
=&\Big\{(x,y,r):x,y \in\mathbb R, r\in\mathbb R_+^*\Big\}.\end{aligned}$$
We will think of ${\mathbb H}$ as a subset of 
Hamilton's quaternions with $1,\,i,\,j,\,k$  the standard 
$\mathbb R$-basis. Write points $P$ in 
 ${\mathbb H}$ as $$P=(z,r)=(x,y,r)=z+rj\qquad\mathrm{where}\  \
z=x+iy,\  j=(0,0,1).$$

The natural action of $SL(2,\mathbb C)$ on ${\mathbb H}$ and on its boundary 
$\mathbb P^1(\mathbb C)$ may be described as follows:
We represent an element of $\mathbb P^1(\mathbb C)$ by 
$\left[\begin{matrix} x\\ y\end{matrix}\right]$ where
 $x,y\in\mathbb C$ with $(x,y)\not=(0,0)$. Then the action of the matrix 
  $M=\left(\begin{matrix} a&b\\ c&d\end{matrix}\right)\in SL(2,\mathbb C)$ 
on $\mathbb P^1(\mathbb C)$ is defined to be
$$\left[\begin{matrix} x\\ y\end{matrix}\right]\mapsto \left(\begin{matrix} a&b\\ c&d\end{matrix}\right)\left[\begin{matrix} x\\ 
y\end{matrix}\right]:=
\left[\begin{matrix} ax+by\\ cx+dy\end{matrix}\right].$$ Moreover, if we 
represent points $P\in {\mathbb H}$ as 
quaternions whose fourth component equals zero, then the action of 
$M$ on ${\mathbb H}$ is defined to be 
$$P\mapsto M\,P:=(aP+b)(cP+d)^{-1},$$
where the inverse on the right is taken in the skew field of quaternions. 

Furthermore, with this action, the stablizer of $j=(0,0,1)\in 
{\mathbb H}$ in $SL(2,\mathbb C)$  is equal to 
$SU(2):=\{A\in U(2):\det A=1\}$. Since the 
action of $SL(2,\mathbb C)$ on ${\mathbb H}$ is transitive, 
we obtain also a natural identification ${\mathbb H}\simeq
SL(2,\mathbb C)/SU(2)$ via the quotient map induced from 
$SL(2,\mathbb C)\to {\mathbb H},\ \ g\mapsto g\cdot j.$ 

\subsection{Rank Two $\mathcal O_F$-Lattices: Upper Half Space Model}
Identify
${\mathcal H}$ with $SL(2,\mathbb R)/SO(2)$ and ${\mathbb H}$ with 
$SL(2,\mathbb C)/SU(2)$. Denote by $\mathcal M_{F,2;\frak a}$ the 
moduli space of semi-stable lattices of rank two whose associated projective models
are isomorphic to $\mathcal O_F\oplus\frak a$ for a certain ideal $\frak a$, 
and denote its volume $T$ part by $\mathcal M_{F,2;\frak a}[T]$. Make the identification 
$$\mathcal M_{F,2;\frak a}\big[N(\frak a)\cdot \Delta_F\big]\simeq 
\bigg(SL(\mathcal O_F\oplus\frak a)\Big\backslash 
\Big({\mathcal H}^{r_1}\times{\mathbb H}^{r_2}\Big)\bigg)_{\mathrm {ss}},$$ 
where as usual ss means the subset consisting of points 
corresponding to rank two semi-stable $\mathcal O_F$-lattices
in the quotient space \newline
$SL(\mathcal O_F\oplus\frak a)\Big\backslash 
\bigg(\Big(SL(2,\mathbb R)/SO(2)\Big)^{r_1}\times \Big(SL(2,\mathbb C)/SU(2)
\Big)^{r_2}\bigg).$

Hence clearly, if the metric on $\mathcal O_F\oplus\frak a$ is given by
$g=(g_\sigma)_{\sigma\in S_\infty}$ with $g_\sigma\in SL(2,F_\sigma)$, 
then the corresponding points on the 
right hand side is $g(\mathrm{ImJ})$ with 
$\mathrm {ImJ}:=(i^{(r_1)},j^{(r_2)})$, i.e., the point given by 
$(g_\sigma\bold\tau_\sigma)_{\sigma\in S_\infty}$ where 
$\tau_\sigma=i_\sigma:=(0,1)$ if $\sigma$ is real and 
 $\tau_\sigma=j_\sigma:=(0,0,1)$ if $\sigma$ is complex.	

\section{Cusps}
\subsection{Definition}
The working site now is shifted to the space 
$SL(\mathcal O_F\oplus\frak a)\Big\backslash\Big({\mathcal H}^{r_1}
\times{\mathbb H}^{r_2}\Big)$.
Here the action of $SL(2,\mathcal O_F\oplus \frak a)$ is 
via the action of $SL(2,F)$ on ${\mathcal H}^{r_1}
\times{\mathbb H}^{r_2}$. More precisely, $F^2$ admits 
natural embeddings $F^2\hookrightarrow 
\Big(\mathbb R^{r_1}\times\mathbb C^{r_2}\Big)^2\simeq 
\big(\mathbb R^2\big)^{r_1}\times\big(\mathbb C^2\big)^{r_2}$ so that
 $\mathcal O_F\oplus \frak a$ naturally embeds into 
$\big(\mathbb R^2\big)^{r_1}\times\big(\mathbb C^2\big)^{r_2}$ 
as a rank two $\mathcal O_F$-lattice. As such, 
$SL(\mathcal O_F\oplus\frak a)$ acts on the image of 
  $\mathcal O_F\oplus \frak a$ in $\big(\mathbb R^2\big)^{r_1}
\times\big(\mathbb C^2\big)^{r_2}$ as automorphisms. Our task 
here is to understand the cusps of this action of 
$SL(\mathcal O_F\oplus\frak a)$ on 
${\mathcal H}^{r_1}\times{\mathbb H}^{r_2}$. For this, we go as follows.

First, the space ${\mathcal H}^{r_1}\times{\mathbb H}^{r_2}$ 
admits a natural boundary 
$\mathbb R^{r_1}\times\mathbb C^{r_2}$, in which the field $F$ 
is imbedded via Archmidean places of $F$: 
$F\hookrightarrow \mathbb R^{r_1}
\times\mathbb C^{r_2}$. Consequently, $\mathbb P^1(F)
\hookrightarrow \mathbb P^1(\mathbb R)^{r_1}\times
\mathbb P^1(\mathbb C)^{r_2}$ with $\left[\begin{matrix} 1\\ 
0\end{matrix}\right]:=\infty\mapsto (\infty^{(r_1)},\infty^{(r_2)})$. 
As usual, via fractional linear transformations, $SL(2,\mathbb R)$ 
acts on	 $\mathbb P^1(\mathbb R)$, and $SL(2,\mathbb C)$ acts on 
 $\mathbb P^1(\mathbb C)$, hence so does $SL(2,F)$ on 
$$\mathbb P^1(F)\hookrightarrow \mathbb P^1(\mathbb R)^{r_1}\times 
\mathbb P^1(\mathbb C)^{r_2}.$$ Being a 
discrete subgroup of $SL(2,\mathbb R)^{r_1}\times SL(2,\mathbb C)^{r_2}$,
for the action of $SL(\mathcal O_K\oplus \frak a)$ 
 on $\mathbb P^1(F)$, we call the corresponding 
orbits (of $SL(\mathcal O_F\oplus\frak a)$ on $\mathbb P^1(F)$) 
the {\it cusps}. Very often we also call their associated representatives cusps.

\subsection{Cusp and Ideal Class Correspondence}

With this, we have the following fundamental result	
rooted back to Maa$\beta$.
\vskip 0.30cm
\noindent
{\bf Cusp and Ideal Class Correspondence.} (Maa$\beta$)
{\it There is a natural bijection $\Pi$ between the ideal class group 
$CL(F)$ of $F$ and the cusps $\mathcal C_\Gamma$ of 
$\Gamma=SL(\mathcal O_F\oplus\frak a)$ 
acting on ${\mathcal H}^{r_1}\times
{\mathbb H}^{r_2}$ given by
$$\mathcal C_\Gamma\to CL(F),\qquad \left[\begin{matrix} \alpha\\ 
\beta\end{matrix}\right]\mapsto 
\Big[\mathcal O_F\,\alpha+\frak a\,\beta\Big].$$}

Easily, one checks that the inverse map $\Pi^{-1}$ is given as follows: 
For a fractional ideal $\frak b$, by Chinese Reminder Theorem,
choose $\alpha_\frak b,
\beta_\frak b\in F$ such that $\mathcal O_F\cdot\alpha_\frak b+
\frak a\cdot\beta_\frak b=\frak b$; Define $\Pi^{-1}([\frak b])$ 
simply by the class of the point 
$\left[\begin{matrix} \alpha_\frak b\\ 
\beta_\frak b\end{matrix}\right]$ in 
$SL(2,\mathcal O_F\oplus a)\Big\backslash \mathbb P^1(F)$. Recall also that there 
always exists $M_{ \left[\begin{matrix} \alpha\\ \beta\end{matrix}\right]}
:= \left(\begin{matrix} 
\alpha&\alpha^*\\ \beta&\beta^*\end{matrix}\right)\in SL(2,F)$ such that 
$M_{\left[\begin{matrix} \alpha\\ 
\beta\end{matrix}\right]}\cdot\infty=
\left[\begin{matrix} \alpha\\ \beta\end{matrix}\right].$

\subsection{Stablizer Groups of Cusps}

Recall that under the Cusp-Ideal Class Correspondence, there are exactly $h$
inequivalence cusps $\eta_i,\, i=1,2,\ldots,h$, where $h:=\#\mathrm{CL}(F)$. Moreover, if we write the cusp 
$\eta:=\eta_i=\left[\begin{matrix} \alpha_i\\ \beta_i\end{matrix}\right]$ for 
suitable $\alpha_i,\,\beta_i\in F$, then the associated ideal class 
is exactly the one for the fractional ideal 
 $\mathcal O_F\alpha_i+\frak a\beta_i=:\frak b_i$. 
Denote the stablizer group of $\eta$ in	 $SL(\mathcal O_F\oplus\frak a)$ by 
$\Gamma_{\eta}.$ 
\vskip 0.20cm
\noindent
{\bf Lemma}. ([W-2,5])	{\it  The associated 
\lq lattice' for the cusp $\eta$ is given by  $\frak a\frak b^{-2}$. Namely,
$$A^{-1}\Gamma_\eta A=\Bigg\{\left(\begin{matrix} u&z\\ 
0&u^{-1}\end{matrix}\right)
:u\in U_F, z\in\frak a\frak b^{-2}\Bigg\},$$
where $U_F$ denotes the group of units of $F$.}
\vskip 0.30cm
Set $\Gamma_\eta':=\bigg\{A\left(\begin{matrix} 1&z\\ 
0&1\end{matrix}\right)A^{-1}:z\in\frak a\frak b^{-2}\bigg\},$ 
Then $$\Gamma_\eta=\Gamma_\eta'\times \bigg\{A\left(\begin{matrix} u&0\\ 
0&u^{-1}\end{matrix}\right)A^{-1}:u\in U_F\bigg\}.$$
Note that also componentwisely, $\left(\begin{matrix} u&0\\ 
0&u^{-1}\end{matrix}\right)z=\frac{uz} {u^{-1}}=u^2z$. So, in 
practice, what we really get is the following decomposition
$$\Gamma_\eta=\Gamma_\eta'\times U_F^2$$ with 
$$U_F^2\simeq \bigg\{A\cdot\left(\begin{matrix} u&0\\ 0&u^{-1}\end{matrix}
\right)\cdot A^{-1}:u\in U_F\bigg\}\,\simeq\, 
\bigg\{A\left(\begin{matrix} 1&0\\ 
0&u^2\end{matrix}\right)A^{-1}:u\in U_F\bigg\}.$$

\subsection{Fundamental Domain for $\Gamma_\infty$ on
 ${\mathcal H}^{r_1}\times{\mathbb H}^{r_2}$}

We are now ready to construct a fundamental domain for 
the action of $\Gamma_\eta\subset SL(\mathcal O_F\oplus\frak a)$ on 
${\mathcal H}^{r_1}\times{\mathbb H}^{r_2}$. 
 This is based on a construction of 
a fundamental domain for the action of $\Gamma_\infty$ on
 ${\mathcal H}^{r_1}\times{\mathbb H}^{r_2}$. More precisely, 
with an element 
$A=\left(\begin{matrix}\alpha&\alpha^*\\ \beta&\beta^*\end{matrix}\right)
\in SL(2,F)$ (always exists!), we have

\noindent
i) $A\cdot\infty=\left[\begin{matrix}\alpha\\ \beta\end{matrix}\right];$ and

\noindent
ii) The isotropy group of $\eta$ in $A^{-1}SL(\mathcal O_F\oplus\frak a)A$ is 
generated by translations $\bold\tau\mapsto\bold\tau+z$ 
with $z\in \frak a\frak b^{-2}$ and by dilations 
$\bold\tau\mapsto u\bold\tau$ where $u$ runs through the group
 $U_F^2$. 

\noindent
(Here, we use $A,\,\alpha,\,\beta,\,\frak b$ as
  running symbols for $A_i,\, \alpha_i,\,\beta_i,\,\frak b_i:=
\mathcal O_F\alpha_i+\frak a\beta_i$.)

Consider then the map $$\begin{matrix}
\mathrm {ImJ}:\ {\mathcal H}^{r_1}\times
{\mathbb H}^{r_2}&\to& \mathbb R_{>0}^{r_1+r_2},\\
(z_1,\cdots,z_{r_1};P_1,\cdots,P_{r_2})&\mapsto&
 (\Im(z_1),\cdots,\Im(z_{r_1});J(P_1),\cdots,J(P_{r_2})),\end{matrix}$$ 
where if $z=x+iy\in{\mathcal H}$, resp. $P=z+rj\in{\mathbb H}$,
 we set $\Im(z)=y$, resp. $J(P)=r$. It induces a map 
$$\Big(A^{-1}\cdot\Gamma_\eta\cdot A\Big)\Big\backslash \Big({\mathcal H}^{r_1}
\times{\mathbb H}^{r_2}\Big)\to U_F^2\Big\backslash 
\mathbb R_{>0}^{r_1+r_2},$$ which exhibits $\Big(A^{-1}\cdot
\Gamma_\eta\cdot A\Big)\Big\backslash \Big({\mathcal H}^{r_1}
\times{\mathbb H}^{r_2}\Big)$ as a torus bundle over $U_F^2\Big\backslash 
\mathbb R_{>0}^{r_1+r_2}$ with fiber the $n=r_1+2r_2$ 
dimensional torus $\Big(\mathbb R^{r_1}\times\mathbb C^{r_2}\Big)\Big/
\frak a\frak b^{-2}$. Having factored out the action of the 
translations, we only have to construct a fundamental domain 
for the action of $U_F^2$ on $\mathbb R_{>0}^{r_1+r_2}$. For this, we 
look first at the action of $U_F^2$ on the norm-one hypersurface
$\bold S:=\Big\{y\in \mathbb R_{>0}^{r_1+r_2}:N(y)=:\prod_iy_i=1\Big\}$. By taking 
logarithms, it is transformed bijectively into a trace-zero hyperplane
 which is isomorphic to the space $\mathbb R^{r_1+r_2-1}$
$$\begin{aligned}\bold S&\buildrel\log\over\to	\mathbb R^{r_1+r_2-1}:=
\Big\{(a_1,\cdots a_{r_1+r_2})\in \mathbb R^{r_1+r_2}:\sum a_i=0\Big\},\\
y&\mapsto\qquad \Big(\log y_1,\cdots,\log y_{r_1+r_2}\Big),\end{aligned}$$ 
where the action of $U_F^2$ on $\bold S$ is carried out 
over an action on $\mathbb R^{r_1+r_2-1}$ by the translations
 $a_i\mapsto a_i+\log\varepsilon^{(i)}$. By Dirichlet's 
Unit Theorem ([L1], [Ne]), the logarithm transforms $U_F^2$ into a 
lattice in $\mathbb R^{r_1+r_2-1}$. Accordingly, the exponential map transforms 
a fundamental domain, e.g., a fundamental parallelopiped, for
 this action back into a fundamental domain $\bold S_{U_F^2}$
 for the action of $U_F^2$ on $\bold S$. The cone over 
$\bold S_{U_F^2}\ $, that is,  $\ \mathbb R_{>0}\cdot
\bold S_{U_F^2}\subset \mathbb R_{>0}^{r_1+r_2}$, is then a 
fundamental domain for the action of $U_F^2$ on
$\mathbb R_{>0}^{r_1+r_2}$. Denote by $\mathcal T$ a fundamental 
domain for the action of the translations by elements of
 $\frak a\frak b^{-2}$ on $\mathbb R^{r_1}\times\mathbb C^{r_2}$,
and set $$\mathrm{ReZ}\,\Big(z_1,\cdots, z_{r_1};
P_1,\cdots,P_{r_2}\Big):=\Big(\Re(z_1),\cdots,\Re(z_{r_1});Z(P_1),
\cdots,Z(P_{r_2})\Big)$$ with
$\Re(z):=x$, resp. $Z(P):=z$ if $z=x+iy\in{\mathcal H}$, resp. 
$P=z+rj\in{\mathbb H}$, 
then what we have just said proves the following
\vskip 0.30cm
\noindent
{\bf Proposition.} ([W-2,5]) {\it A fundamental domain for the action of 
$A^{-1}\Gamma_\eta A$ on ${\mathcal H}^{r_1}\times
{\mathbb H}^{r_2}$ is given by}
$$\bold E:=\Big\{\bold\tau\in {\mathcal H}^{r_1}\times
{\mathbb H}^{r_2}:\mathrm{ReZ}\,(\bold\tau)\in {\mathcal T},\  
\mathrm{ImJ}\,(\bold\tau)\in  \mathbb R_{>0}\cdot
\bold S_{U_F^2}\Big\}.$$

For later use, we also set $\mathcal F_\eta:=A_\eta^{-1}\cdot \bold E$.

\section{Fundamental Domain}
\subsection{Siegel Type Distance}

Guided by Siegel's discussion on totally real fields [Sie] and the 
discussion above,	
we are now ready to construct  fundamental domains for 
 $SL(\mathcal O_F\oplus\frak a)
\Big\backslash \Big(\mathcal H^{r_1}\times\mathbb H^{r_2}\Big).$ 

As the first step, we generalize Siegel's \lq distance to cusps'. 
For this, recall that for a cusp $\eta=\left[\begin{matrix}\alpha\\ 
\beta\end{matrix}\right]\in \mathbb P^1(F)$, by the Cusp-Ideal Class
Correspondence, we obtain a natural  ideal class associated to
the fractional ideal	
$\frak b:=\mathcal O_F\cdot\alpha+\frak a\cdot\beta$. Moreover, by assuming 
that $\alpha,\beta$ are all contained in $\mathcal O_F$, 
as we may, we know that the corresponding stablizer group $\Gamma_\eta$ 
is given by
$$A^{-1}\cdot \Gamma_\eta\cdot A=\bigg\{\gamma=\begin{pmatrix} u&z\\ 0&u^{-1}
\end{pmatrix}\in\Gamma: u\in U_F, z\in\frak a\frak b^{-2}\bigg\},$$
where $A\in SL(2,F)$ satisfying $A\infty=\eta$ which may be further chosen in 
the form
$A=\begin{pmatrix} \alpha&\alpha^*\\ \beta&\beta^*\end{pmatrix}\in SL(2,F)$
so that $\mathcal O_F\beta^*+\frak a^{-1}\alpha^*=\frak b^{-1}$.

Now for $\bold\tau=(z_1,\ldots,z_{r_1};P_1,\cdots,P_{r_2})
\in {\mathcal H}^{r_1}\times{\mathbb H}^{r_2}$, set
$$N(\tau):=N\Big(\mathrm{ImJ}(\tau)\Big)=\prod_{i=1}^{r_1}
\Im (z_i)\cdot\prod_{j=1}^{r_2} J(P_j)^2=\Big(y_1\cdot\ldots\cdot y_{r_1}\Big)
\cdot\Big(v_1\cdot\ldots\cdot v_{r_2}\Big)^2.$$
Then  for all 
$\gamma=\begin{pmatrix} a&b\\ c&d
\end{pmatrix}\in SL(2,F)$,  $$N\Big(\mathrm{ImJ}(\gamma\cdot\tau)\Big)
=\frac{N(\mathrm{ImJ}(\bold\tau))}{\|N(c\bold \tau+d)\|^2}.\eqno(*)$$
(Note that here only the second row of $\gamma$ appears.)
Moreover, following [W-2,5], define the {\it reciprocal distance 
$\mu(\eta,\bold\tau)$  from the point
$\bold\tau\in {\mathcal H}^{r_1}\times{\mathbb H}^{r_2}$  to the cusp 
$\eta=\left[\begin{matrix}\alpha\\ \beta\end{matrix}\right]$  in 
$\mathbb P^1(F)$} by 
$$\begin{aligned}\mu(\eta,\bold\tau):=&
N\Big(\frak a^{-1}\cdot (\mathcal O_F\alpha+\frak a\beta)^2\Big)\\
&\times\frac{\Im(z_1)\cdots \Im(z_{r_1})\cdot J(P_1)^2\cdots 
J(P_{r_2})^2}{\prod_{i=1}^{r_1}|(-\beta^{(i)}z_i+
\alpha^{(i)})|^2\prod_{j=1}^{r_2}
\|(-\beta^{(j)}P_j+\alpha^{(j)})\|^2}\\
=&\frac{1}{N(\frak a\frak b^{-2})}\cdot\frac{N(\mathrm{ImJ}(\bold\tau))}
{\|N(-\beta\bold \tau+\alpha)\|^2}.\end{aligned}$$
\vskip 0.20cm
\noindent
{\bf {Lamma 1.}} ([W-2,5]) (i) {\it $\mu$ is well-defined;}

\noindent
(ii) {\it $\mu$ is invariant under the action of 
$SL(\mathcal O_F\oplus\frak a)$. That is to say,
$$\mu(\gamma\eta,\gamma\bold\tau)=
\mu(\eta,\bold\tau),\qquad\forall \gamma\in SL(\mathcal O_F\oplus\frak a).$$}

\noindent
(iii) {\it There exists a positive constant $C$ 
depending only on $F$ and $\frak a$ 
such that if $\mu(\eta,\bold\tau)>C$ and 
$\mu(\eta',\bold\tau)>C$ for $\bold\tau\in	
{\mathcal H}^{r_1}\times{\mathbb H}^{r_2}$ and 
$\eta,\,\eta'\in\mathbb P^1(F)$, then $\eta=\eta'$.}

\noindent
(iv) {\it There exists a positive real number $T:=T(F)$ 
depending only on $F$ such that for 
$\bold\tau\in {\mathcal H}^{r_1}\times {\mathbb H}^{r_2}$, 
there exists a cusp $\eta$ such that $\mu(\eta,\bold\tau)>T$.}
\vskip 0.30cm
Now for the cusp $\eta=\left[\begin{matrix}\alpha\\ 
\beta\end{matrix}\right]\in\mathbb P^1(F)$, define the 
\lq sphere of influence' of $\eta$ by 
$$F_\eta:=\Big\{\bold\tau\in {\mathcal H}^{r_1}\times
{\mathbb H}^{r_2}:\mu(\eta,\bold\tau)\geq 
\mu(\eta',\bold\tau),\forall\eta'\in \mathbb P^1(F)\Big\}.$$

\noindent
{\bf{Lemma 2}.} ([W-2,5]) {\it The action of $SL(\mathcal O_F\oplus\frak a)$ 
in the interior $F_\eta^0$ of $F_\eta$ reduces to that 
of the isotropy group $\Gamma_\eta$ of $\eta$, i.e., if 
$\bold\tau$ and $\gamma\bold\tau$ both belong to $F_\eta^0$, 
then $\gamma\bold\tau=\bold\tau.$}
\vskip 0.30cm
Consequently, we arrive at the following way to 
decompose the orbit space $SL(\mathcal O_F\oplus\frak a)
\Big\backslash \Big({\mathcal H}^{r_1}\times{\mathbb H}^{r_2}\Big)$ into
$h$ pieces glued in some way along pants of their boundary. 
\vskip 0.20cm
\noindent
{\bf Proposition.} ([W-2,5]) {\it Let 
$i_\eta:\Gamma_\eta\Big\backslash F_\eta\hookrightarrow 
SL(\mathcal O_F\oplus\frak a)\Big\backslash\Big({\mathcal H}^{r_1}
\times{\mathbb H}^{r_2}\Big)$ denote the natural map. Then 
$$SL(\mathcal O_F\oplus\frak a)\Big\backslash 
\Big({\mathcal H}^{r_1}\times{\mathbb H}^{r_2}\Big)=\bigcup_\eta 
i_\eta\Big(\Gamma_\eta\Big\backslash F_\eta\Big),$$
where the union is taken over a set of $h$ cusps 
representing the ideal classes of $F$. Each piece corresponds
to an ideal class of $F$.}
\vskip 0.30cm
Note that the action of $\Gamma_\eta$ on ${\mathcal H}^{r_1}
\times{\mathbb H}^{r_2}$ is free. Consequently, all fixed 
points of $SL(\mathcal O_F\oplus\frak a)$ on ${\mathcal H}^{r_1}
\times{\mathbb H}^{r_2}$ lie on the boundaries of $F_\eta$.
\subsection{Fundamental Domains}
We can give a more precise description of the fundamental domain, 
based on our understanding of that for stablizer groups of
cusps. To state it, denote by $\eta_1,\,\ldots,\,\eta_h$ 
inequivalent cusps for the action of $SL(\mathcal O_F\oplus\frak a)$ on 
${\mathcal H}^{r_1}\times{\mathbb H}^{r_2}$. Choose $A_{\eta_i}\in SL(2,F)$
such that $A_{\eta_i}\infty=\eta_i,\ i=1,2,\ldots,h$. Write $\bold S$ for
the norm-one hypersurface
$\bold S:=\Big\{y\in \mathbb R_{>0}^{r_1+r_2}:N(y)=1\Big\}$, and 
$\bold S_{U_F^2}$ for
 the action  of $U_F^2$ on $\bold S$. Denote by $\mathcal T$ a fundamental 
domain for the action of the translations by elements of
 $\frak a\frak b^{-2}$ on $\mathbb R^{r_1}\times\mathbb C^{r_2}$, 
and $$\bold E:=\Big\{\bold\tau\in {\mathcal H}^{r_1}\times
{\mathbb H}^{r_2}:
\mathrm{ReZ}\,(\bold\tau)\in {\mathcal T},\  
\mathrm{ImJ}\,(\bold\tau)\in  \mathbb R_{>0}\cdot
\bold S_{U_F^2}\Big\}$$ 
a fundamental domain for the action of 
$A_\eta^{-1}\Gamma_\eta A_\eta$ on ${\mathcal H}^{r_1}\times
{\mathbb H}^{r_2}$.  The intersections of $\bold E$ 
with $i_\eta(F_\eta)$ are connected. Consequently,
we have the following
\vskip 0.20cm
\noindent
{\bf Proposition.$'$} (Siegel, Weng) (1) {\it $D_\eta:=A_\eta^{-1} {\bold E}\cap F_\eta$ 
is a fundamental domain for the action of $\Gamma_\eta$ on $F_\eta$;}

\noindent
(2) {\it There exist $\alpha_1,\,\cdots,\,\alpha_h\in 
SL(\mathcal O_F\oplus\frak a)$ such that 
$\cup_{i=1}^h\alpha(D_{\eta_i})$ is connected and hence a
 fundamental domain for $SL(\mathcal O_F\oplus\frak a)$.}
\vskip 0.30cm
That is to say, a fundamental domain may be given as 
$S_Y\cup \mathcal F_1(Y_1)\cup\cdots\cup \mathcal F_h(Y_h)$ 
with $S_Y$ bounded, $\mathcal F_i(Y_i)=A_i\cdot \widetilde  {\mathcal F}_i(Y_i)$
and $$\widetilde{\mathcal F}_i(Y_i):=\Big\{\bold\tau\in {\mathcal H}^{r_1}
\times{\mathbb H}^{r_2}:\mathrm{ReZ}(\bold\tau)\in\Sigma,\,
\mathrm{ImJ}(\bold\tau)\in\mathbb R_{>T}\cdot \bold S_{U_F^2}\Big\}.$$
Moreover, all $\mathcal F_i(Y_i)$'s are	 
disjoint from each other when $Y_i$ are sufficiently large.
\section{Stability in Rank Two}
\subsection{Stability and Distances to Cusps}

Define now the {\it distance of $\tau$ to the cusp $\eta$} by
$$d(\eta,\tau):=\frac{1}{\mu(\eta,\tau)}\geq 1.$$ 
Then, with the use of a crucial result of Tsukasa Hayashi [Ha], 
 we are ready to state the following  fundamental result, which
exposes a beautiful intrinsic relation between stability and the 
distance to cusps.
\vskip 0.30cm
\noindent 
{\bf {Theorem.}} (Weng)	 {\it The lattice $\Lambda$ is semi-stable 
if and only if the distances  of corresponding point 
$\tau_\Lambda\in \mathcal H^{r_1}\times{\mathbb H}^{r_2}$ 
to all cusps are all bigger or equal to 1.}

\subsection{Moduli Space of Rank Two Semi-Stable $\mathcal O_F$-Lattices}

For a rank two 
$\mathcal O_F$-lattice $\Lambda$, denote by $\tau_\Lambda\in
 \mathcal H^{r_1}\times\mathbb H^{r_2}$ the corresponding module point. 
Then, by the 
previous subsection,  $\Lambda$ is semi-stable if and only if 
for all cusps $\eta$,	
$d(\eta, \tau_\Lambda):=\frac{1}{\mu(\eta, \tau_\Lambda)}$ 
are bigger than or equal to 1.
This then leads to the consideration of the following
truncation of the fundamental domain $\mathcal D$ of $SL(\mathcal O_F\oplus
\frak a)\Big\backslash \Big(\mathcal H^{r_1}
\times\mathbb H^{r_2}\Big)$:
For $T\geq 1$, denote by $$\mathcal D_T:=\Big\{\tau\in\mathcal D:
d(\eta,\tau)\geq T^{-1},\ \forall\mathrm{cusp}\ \eta\Big\}.$$
 
The space $\mathcal D_T$ may be precisely described in terms of $\mathcal D$ 
and certain neighborhood of cusps. To explain this, we first recall the 
following

\noindent
{\bf{Lemma}.} ([W-2,5]) {\it For a cusp $\eta$, denote by $$X_\eta(T):=\Big\{
\tau\in \mathcal H^{r_1}\times\mathbb H^{r_2}: d(\eta,\tau)<T^{-1}\Big\}.$$
Then for $T\geq 1$, $$X_{\eta_1}(T)\cap X_{\eta_2}(T)\not=\emptyset\qquad
\Leftrightarrow\qquad  \eta_1=\eta_2.$$}

With this, we are ready to state the following
\vskip 0.20cm
\noindent
{\bf Theorem.} (Weng) {\it There is a natural identification between}

\noindent
(a) {\it  moduli space of rank two semi-stable 
$\mathcal O_F$-lattices of volume $N(\frak a)\cdot|\Delta_F|$
with underlying projective module 
$\mathcal O_F\oplus\frak a$;} and 

\noindent
(b) {\it truncated compact domain 
$\mathcal D_1$ consisting of points in the fundamental domain $\mathcal D$
whose distances to all cusps are bigger than 1.}
\vskip 0.20cm
In other words, the truncated compact domain  $\mathcal D_1$
is obtained from the fundamental domain $\mathcal D$
of $SL(\mathcal O_F\oplus
\frak a)\Big\backslash \Big(\mathcal H^{r_1}
\times\mathbb H^{r_2}\Big)$ by delecting the disjoint open neighborhoods
$\cup\cup_{i=1}^h\mathcal F_i(1)$ associated to inequivalent cusps 
$\eta_1,\eta_2,\ldots,\eta_h$,
where $\mathcal F_i(T)$ denotes the neighborhood of $\eta_i$ consisting
of $\tau\in\mathcal D$ whose distance to $\eta_i$ is strictly 
less than $T^{-1}$.
\vfill
\eject
\centerline
{\Large\bf Chapter V. Algebraic Characterization of Stability}

\section{Canonical Filtrations}
\subsection{Canonical Filtrations}
Following Lafforgue [Laf], we call an abelian category ${\mathcal A}$ 
together with two additive
morphisms $$\mathrm {rk}:{\mathcal A}\to {\mathbb N},
\qquad \mathrm{deg}:{\mathcal A}\to {\mathbb R}$$ a {\it category with slope 
structure}. In particular, for non-zero $A\in {\mathcal A}$, 
\vskip 0.30cm
\noindent
({\bf 1}) define the {\it slope} of $A$ by 
$\mu(A):={{\mathrm{deg}(A)}\over {\mathrm{rk}A}};$
\vskip 0.30cm 
\noindent
({\bf 2}) If $0=A_0\subset A_1\subset\cdots\subset A_l=A$ is a filtration 
of $A$ in 
${\mathcal A}$ with $\mathrm {rk}(A_0)<\mathrm {rk}(A_1)<\cdots<
\mathrm {rk}(A_l)$, define the {\it associated polygon} 
to be the function $[0,\mathrm {rk}A]\to {\mathbb R}$ such that

\noindent
(i) its values at 0 and $\mathrm {rk}(A)$ are 0;

\noindent
(ii)  it is affine on the intervals 
$[\mathrm {rk}(A_{i-1}),\mathrm {rk}(A_i)]$ 
with slope $\mu(A_i/A_{i-1})-\mu(A)$;
\vskip 0.30cm
\noindent
({\bf 3}) If $\frak a$ is a collection of subobjects of $A$ in ${\mathcal A}$, 
then 
$\frak a$ is said to be {\it nice} if

\noindent
(i) $\frak a$ is stable under intersection and finite summation;

\noindent
(ii) $\frak a$ is Noetherian, i.e., every increasing chain 
of elements
in $\frak a$ has a maximal element in $\frak a$;

\noindent
(iii) if $A_1\in\frak a$ then $A_1\not=0$ if and only if 
$\mathrm {rk}(A_1)\not=0$; and

\noindent
(iv) for $A_1,A_2\in\frak a$ with $\mathrm {rk}(A_1)=\mathrm {rk}(A_2)$. 
Then
$A_1\subset A_2$ is proper implies that $\mathrm {deg}(A_1)<
\mathrm {deg}(A_2);$
\vskip 0.30cm
\noindent
({\bf 4}) For any nice $\frak a$, set 
$$\begin{aligned}\mu^+(A):=&\mathrm {sup}\,\Big\{\mu(A_1):A_1\in\frak a,
\mathrm {rk}(A_1)\geq 1\Big\},\\ 
\mu^-(A):=&\mathrm {inf}\,\Big\{\mu(A/A_1):A_1\in\frak a,
\mathrm {rk}(A_1)< \mathrm {rk}(A) \Big\}.\end{aligned}$$ 
Then we say $(A,\frak a)$ is {\it semi-stable} if $\mu^+(A)=\mu(A)=\mu^-(A)$. 
Moreover if $\mathrm {rk}(A)=0$, 
set also $\mu^+(A)=-\infty$ and $\mu^-(A)=+\infty$.
\vskip 0.30cm
\noindent
{\bf{Proposition 1}.} ([Laf]) {\it Let ${\mathcal A}$ be a 
category with slope structure, $A$ an object in ${\mathcal A}$ and 
$\frak a$ a nice family of subobjects of $A$ in ${\mathcal A}$. Then}

\noindent
(1) ({\bf Canonical Filtration}) {\it $A$ admits a unique filtration $0=
\overline A_0\subset \overline A_1\subset\cdots\subset\overline A_l=A$
with elements in $\frak a$  such that}

\noindent
(i) {\it $\overline A_i,0\leq i\leq k$ are maximal in $\frak a$;}

\noindent 
(ii) {\it $\overline A_i/\overline A_{i-1}$ are semi-stable; and}

\noindent 
(iii) $\mu(\overline A_1/\overline A_{0})>
\mu(\overline A_2/\overline A_{1}>\cdots>
\mu(\overline A_k/\overline A_{k-1})$;

\noindent
(2) ({\bf Boundness}) {\it All polygons of filtrations of $A$ with elements 
in $\frak a$ are bounded from above by $\overline p$, where 
 $\overline p:=\overline p^A$ is the associated polygon for the canonical 
filtration in (1);}

\noindent
(3) {\it For any $A_1\in \frak a, \mathrm {rk}(A_1)\geq 1$ implies
$\mu(A_1)\leq \mu(A)+\frac{\overline p(\mathrm {rk}(A_1))}
{\mathrm {rk}(A_1)};$}

\noindent
(4) {\it The polygon $\overline p$ is convex with maximal slope 
$\mu^+(A)-\mu(A)$ 
and minimal slope $\mu^-(A)-\mu(A)$;}

\noindent
(5) {\it If $(A',\frak a')$ is another pair, and $u:A\to A'$ is a homomorphism 
such that $\mathrm {Ker}(u)\in\frak a$ and $\mathrm{Im}(u)\in\frak a'$. Then
$\mu^-(A)\geq\mu^+(A')$ implies that $u=0$.}
\vskip 0.30cm
\noindent
This results from a Harder-Narasimhan type filtration 
consideration. A detailed proof may be found at pp. 87-88 in [Laf].
(There are some interesting approaches related to the topics here
in literatures. For examples, [An2], [Ch].)

\subsection{Examples of Lattices}

As an example, we have the following	

\noindent
{\bf {Proposition 2}.} ([W2,3]) {\it  Let $F$ be a number field. Then}

\noindent
(1) {\it the abelian category of hermitian vector sheaves on 
$\mathrm {Spec}\,{\mathcal O}_F$ together with the natural rank and the 
Arakelov degree is a category with slopes;}

\noindent
(2) {\it For any hermitian vector sheaf $(E,\rho)$, $\frak a$ consisting 
of pairs $(E_1,\rho_1)$ with $E_1$ sub vector sheaves of $E$ and 
$\rho_1$ the restrictions of $\rho$, forms a nice family.}
\vskip 0.30cm
\noindent
Indeed, (1) is obvious, while (2) is a direct consequence of 
the following standard facts: 

\noindent
(i) For a fixed $(E,\rho)$,
$\Big\{\mathrm {deg}(E_1,\rho_1):(E_1,\rho_1)\in \frak a\Big\}$ is discrete 
subset of ${\mathbb R}$; and 

\noindent
(ii) for any two sublattices $\Lambda_1,\, \Lambda_2$	
of $\Lambda$,
 $$\mathrm{Vol}\Big(\Lambda_1/(\Lambda_1\cap\Lambda_2)\Big)\geq
\mathrm{Vol}\Big((\Lambda_1+\Lambda_2)/\Lambda_2\Big).$$

Consequently, there exists canonical filtrations  of Harder-Narasimhan
 type for 
hermitian vector sheaves over $\mathrm {Spec}\,{\mathcal O}_F$.
Recall that hermitian vector sheaves over $\mathrm {Spec}\,{\mathcal O}_F$
are ${\mathcal O}_F$-lattices in $({\mathbb R}^{r_1}\times 
{\mathbb C}^{r_2})^{r=\mathrm{rk}(E)}$ in the language of Arakelov theory:
Say, corresponding $\mathcal O_F$-lattices are 
induced from their $H^0$ via the natural embedding
$F^r\hookrightarrow ({\mathbb R}^{r_1}\times {\mathbb C}^{r_2})^r$ where $r_1$ 
(resp. $r_2$) denotes the real (resp. complex) embeddings of $F$.

\section{Algebraic Characterization}

\subsection{A GIT Principle}
In  Geometric Invariant Theory ([M], [Kem], [RR]), 
a fundamental principle, the Micro-Global Principle, 
claims that if a point is not GIT stable then there exists
a parabolic subgroup which destroys the corresponding stability.

In the setting of $\mathcal O_F$-lattices, 
even we do not have a proper definition of GIT stability for lattices, 
in terms of intersection stability, an analogue of the Micro-Global Principle 
does hold.

\subsection{Micro-Global Relation for Geo-Ari Truncations}

Let $\Lambda=\Lambda^g$ be a rank $r$ lattice associated to 
$g\in\mathrm{GL}_r(\mathbb A)$ and  $P$ a parabolic 
subgroup.  Denote the sublattices 
filtration associated to  $P$ by $$0=\Lambda_0\subset
\Lambda_1\subset\Lambda_2\subset\cdots\subset\Lambda_{|P|}=\Lambda.$$ 
Assume that $P$ corresponds to the partition $I=(d_1,d_2,\cdots,d_{n=:|P|})$. 
Consequently, we have $$\mathrm{rk}(\Lambda_i)=r_i:=d_1+d_2+\cdots+d_i,\qquad 
\mathrm{for}\ i=1,2,\cdots,|P|.$$ 
Let $p,q:[0,r]\to {\mathbb R}$ be two polygons 
such that $p(0)=q(0)=p(r)=q(r)=0$. Then following Lafforgue, 
we say $q$ is {\it bigger than} $p$ {\it with respect to} $P$ and denote it by 
$q>_Pp$, if $q(r_i)-p(r_i)>0$ for all $i=1,\cdots,|P|-1.$ 
Introduce also	the characteristic function $\bold 1(\overline p^*\leq p)$ by
$$\bold 1(\overline p^g\leq p)
=\begin{cases} 1,&\text{if $\overline p^g\leq p$;}\\
0,&\text{otherwise}.\end{cases}$$ Recall that for a parabolic subgroup 
$P$, $p_P^g$ denotes the polygon induced by $P$ for (the lattice corresponding 
to) the element $g\in G(\mathbb A)$.
\vskip 0.30cm
\noindent
{\bf {Fundamental Relation.}} (Lafforgue, Weng) 
{\it Let $p:[0,r]\to 
{\mathbb R}$ be a fixed convex polygon such that $p(0)=p(r)=0$. Then we have
$$\bold 1(\overline p^g\leq p)=\sum_{P:\, \text{stand\, 
para}}(-1)^{|P|-1}\sum_{\delta\in
P(F)\backslash G(F)}\bold 1(p_P^{\delta g}>_Pp)
\qquad\forall g\in G(\mathbb A).$$}

\vfill
\eject
\centerline{\Large\bf Chapter VI. Analytic Characterization of Stability}
\section{Arthur's Analytic Truncation}
\subsection{Parabolic Subgroups}

Let $F$ be a number field with $\mathbb A=\mathbb A_F$ the ring of adeles. 
Let $G$ be a connected reductive group defined over $F$. Recall that a 
subgroup $P$ of $G$ is called {\it parabolic} if $G/P$ is a complete algebraic 
variety. Fix a minimal $F$-parabolic subgroup $P_0$ of $G$ with its unipotent 
radical $N_0=N_{P_0}$ and fix a $F$-Levi subgroup $M_0=M_{P_0}$ of $P_0$ so as 
to have a Levi decomposition $P_0=M_0N_0$. An $F$-parabolic subgroup $P$ is 
called {\it standard} if it contains $P_0$.  For such a parabolic subgroup $P$, 
there exists a unique  Levi subgroup $M=M_P$ containing $M_0$ which we 
call the 
{\it standard Levi subgroup} of $P$. Let $N=N_P$ be the unipotent radical.
 Let us agree to use the term parabolic subgroups and Levi subgroups to 
denote standard $F$-parabolic subgroups and standard Levi subgroups 
repectively, unless otherwise is stated. 

Let $P$ be a parabolic subgroup of $G$. Write $T_P$ for the maximal split 
torus in the center of $M_P$ and $T_P'$ for the maximal quotient split 
torus of $M_P$. Set $\tilde{\frak a}_P:=X_*(T_P)\otimes \mathbb R$ and 
denote its real dimension
by $d(P)$, where $X_*(T)$ is the lattice of 1-parameter subgroups in the torus 
$T$. Then it is known that $\tilde{\frak a}_P=X_*(T_P')\otimes \mathbb R$ as 
well. The two descriptions of $\tilde {\frak a}_P$ 
show that if $Q\subset P$ is a 
parabolic subgroup, then there is a canonical injection 
$\tilde {\frak a}_{P}\hookrightarrow \tilde{\frak a}_{Q}$ and a 
natural surjection	
$\tilde {\frak a}_{Q}\twoheadrightarrow \tilde{\frak a}_{P}$.	
We thus obtain a canonical decomposition 
$\tilde {\frak a}_{Q}=\tilde{\frak a}_Q^{P}\oplus \tilde {\frak a}_{P}$ 
for a certain subspace $\tilde{\frak a}_Q^{P}$ of $\tilde {\frak a}_{Q}$. 
In particular, $\tilde{\frak a}_G$ is a summand of 
$\tilde {\frak a}=\tilde{\frak a}_P$ for all $P$. Set 
$\frak a_{P}:= \tilde{\frak a}_{P}/\tilde {\frak a}_{G}$ and 
$\frak a_{Q}^P:= \tilde{\frak a}_Q^{P}/ \tilde{\frak a}_{G}.$ Then we have 
$$\frak a_{Q}=\frak a_Q^{P}\oplus\frak a_{P}$$
and $\frak a_P$ is canonically identified as a subspace of 
$\frak a_Q$. Set $\frak a_0:=\frak a_{P_0}$ and 
$\frak a_0^P=\frak a_{P_0}^P$ then we also have 
$\frak a_0=\frak a_0^P\oplus \frak a_P$ for all $P$.

\subsection{Logarithmic Map}

For a real vector space $V$, write $V^*$ its dual space over $\mathbb R$. Then
dually we have the spaces $\frak a_0^*, \frak a_P^*, 
\Big(\frak a_0^P\Big)^*$ and hence the decompositions 
$$\frak a_0^*=\Big(\frak a_0^Q\Big)^*\oplus
\Big(\frak a_Q^P\Big)^*\oplus\frak a_P^*.$$ In particular, 
$\frak a_P^*=X(M_P)\otimes \mathbb R$ with $X(M_P):=\mathrm {Hom}_F\Big(M_P,GL(1)\Big)$ i.e., collection of characters on 
$M_P$. It is known that
$\frak a_P^*=X(A_P)\otimes \mathbb R$ where $A_P$ denotes the split 
component of the center of $M_P$. Clearly, if $Q\subset P$, then 
$M_Q\subset M_P$ while $A_P\subset A_Q$. Thus via restriction,
the above two expressions of $\frak a_P^*$ also naturally induce
an injection $\frak a_P^*\hookrightarrow \frak a_Q^*$ and a sujection 
$\frak a_Q^*\twoheadrightarrow
\frak a_P^*$, compactible with the decomposition $\frak a_Q^*=
\Big(\frak a_Q^P\Big)^*\oplus\frak a_P^*.$

Every $\chi=\sum s_i\chi_i$ in $\frak a_{P,\mathbb C}^*:=
\frak a_P^*\otimes\mathbb C$ determines a morphism 
$P(\mathbb A)\to \mathbb C^*$ by $p\mapsto p^\chi:=\prod|\chi_i(p)|^{s_i}.$ 
Consequently, we have a natural logarithmic map $H_P:P(\mathbb A)\to\frak a_P$ 
defined by $$\langle H_P(p),\chi\rangle=p^\chi,\qquad 
\forall \chi\in\frak a_P^*.$$ The 
kernel of $H_P$ is denoted by $P(\mathbb A)^1$ and we set 
$M_P(\mathbb A)^1:=P(\mathbb A)^1\cap M_P(\mathbb A)$.
 
Let also $A_+$ be the set of $a\in A_P(\mathbb A)$ such that 

\noindent
(1) $a_v=1$ for all finite places $v$ of $F$; and 

\noindent
(2) $\chi(a_\sigma)$ 
is a positive number independent of infinite places $\sigma$ of $F$ 
for all $\chi\in X(M_P)$. 

\noindent
Then $M(\mathbb A)=A_+\cdot M(\mathbb A)^1$. 

\subsection{Roots, Coroots, Weights and Coweights}

We now introduce standard bases for above spaces and their duals. Let 
$\Delta_0$ and $\widehat\Delta_0$ be the subsets of simple roots and simple 
weights in $\frak a_0^*$ respectively. (Recall that elements of 
$\widehat\Delta_0$ are non-negative linear combinations of elements in 
$\Delta_0$.) Write $\Delta_0^\vee$ (resp. $\widehat\Delta^\vee_0$) for 
the basis of $\frak a_0$ dual to $\widehat\Delta_0$ (resp. $\Delta_0$). 
Being the dual of the collection of simple weights (resp. of simple roots), 
$\Delta_0^\vee$ (resp. $\widehat\Delta_0^\vee$) 
is the set of coroots (resp. coweights).

For every $P$, let $\Delta_P\subset\frak a_0^*$ be the set of non-trivial 
{\it restrictions} of elements of $\Delta_0$ to $\frak a_P$. Denote the dual 
basis of $\Delta_P$ by $\widehat\Delta_P^\vee$.
For each $\alpha\in\Delta_P$, let $\alpha^\vee$ be the projection of 
$\beta^\vee$ to $\frak a_P$, where $\beta$ is the root in $\Delta_0$ 
whose restriction to $\frak a_P$ is $\alpha$. Set 
$\Delta_P^\vee:=\Big\{\alpha^\vee:\alpha\in\Delta_P\Big\}$, and define the 
dual basis of $\Delta_P^\vee$ by $\widehat\Delta_P$. 

More generally, if $Q\subset P$, write $\Delta_Q^P$ to denote the {\it subset} 
$\alpha\in\Delta_Q$ appearing in the action of $T_Q$ in the unipotent radical 
of $Q\cap M_P$. (Indeed, $M_P\cap Q$ is a parabolic subgroup of $M_P$ with 
nilpotent radical $N_Q^P:=N_Q\cap M_P$. Thus $\Delta_Q^P$ is simply the set 
of roots of the parabolic subgroup $(M_P\cap Q,A_Q)$. And one checks that the 
map $P\mapsto \Delta_Q^P$ gives a natural bijection between parabolic 
subgroups $P$ containin $Q$ and subsets of $\Delta_Q$.)
Then $\frak a_P$ is the subspace of $\frak a_Q$ annihilated by 
$\Delta_Q^P$. Denote by $(\widehat\Delta^\vee)_Q^P$ the dual of $\Delta_Q^P$. 
Let $(\Delta_Q^P)^\vee:=\Big\{\alpha^\vee:\alpha\in \Delta_Q^P\Big\}$ 
and denote by $\widehat\Delta_Q^P$ the dual of $(\Delta_Q^P)^\vee$.
\subsection{Positive Cone and Positive Chamber}

Let $Q\subset P$ be two parabolic subgroups of $G$.
We extend the linear functionals in $\Delta_Q^P$ and $\widehat\Delta_Q^P$ to
 elements of the dual space $\frak a_0^*$ by means of the canonical 
projection from $\frak a_0$ to $\frak a_Q^P$ given by the decomposition 
$\frak a_0=\frak a_0^Q\oplus\frak a_Q^P\oplus\frak a_P$. Let $\tau_Q^P$ 
be the characteristic function of the {\it positive chamber} 
$$\begin{aligned}\Big\{H\in\frak a_0:&\langle\alpha,H\rangle>0\ \mathrm{for\ all}\ 
\alpha\in\Delta_Q^P\Big\}\\
=&\frak a_0^Q\oplus
\Big\{H\in\frak a_Q^P:\langle\alpha,H\rangle>0\ \mathrm{for\ all}\ 
\alpha\in\Delta_Q^P\Big\}
\oplus\frak a_P\end{aligned}$$ and let $\widehat\tau_Q^P$ be the characteristic 
function of the {\it positive cone} 
$$\begin{aligned}\Big\{H\in\frak a_0:&\langle\varpi,H\rangle>0\ 
\mathrm{for\ all}\ \varpi\in\widehat\Delta_Q^P\Big\}\\
=&\frak a_0^Q\oplus
\Big\{H\in\frak a_Q^P:\langle\varpi,H\rangle>0\ 
\mathrm{for\ all}\ \varpi\in\widehat\Delta_Q^P\Big\}
\oplus\frak a_P.\end{aligned}$$ Note that elements in 
$\widehat\Delta_Q^P$ are non-negative linear combinations of elements in 
$\Delta_Q^P$, we have $$
\widehat\tau_Q^P\geq \tau_Q^P.$$ 

\subsection{Partial Truncation and First Estimations}

Denote $\tau_P^G$ and $\widehat\tau_P^G$ simply by $\tau_P$ and 
$\widehat\tau_P$.
\vskip 0.30cm
\noindent
{\bf Basic Estimation.} (Arthur)  {\it Suppose that we are given a parabolic 
subgroup $P$, and a Euclidean norm $\|\cdot\|$ on $\frak a_P$. Then 
there are constants $c$ and $N$ such that for all $x\in G(\mathbb A)^1$ 
and $X\in\frak a_P$,
$$\sum_{\delta\in P(F)\backslash G(F)}\widehat\tau_P\Big(H(\delta x)-X\Big)
\leq c\Big(\|x\|e^{\|X\|}\Big)^N.$$} Moreover, {\it the sum is finite}.

\vskip 0.30cm
As a direct consequence,  we have the following

\noindent
{\bf Corollary.} ([Ar2,3]) {\it Suppose that $T\in\frak a_0$ and $N\geq 0$. Then
 there exist constants $c'$ and $N'$ such that for any function $\phi$ 
on $P(F)\backslash G(\mathbb A)^1$, and $x,y\in G(\mathbb A)^1$,
$$\sum_{\delta\in P(F)\backslash G(F)}\Big|\phi(\delta x)\Big|
\cdot\widehat\tau_P
\Big(H(\delta x)-H(y)-X\Big)$$ is bounded by} $$c'\|x\|^{N'}\cdot
\|y\|^{N'}\cdot
\sup_{u\in G(\mathbb A)^1}\Big(|\phi(u)|\cdot \|u\|^{-N}\Big).$$

\section{Reduction Theory}
\subsection{Langlands' Combinatorial Lemma}

If $P_1\subset P_2$, following Arthur [Ar2], set
$$\sigma_1^2(H):=\sigma_{P_1}^{P_2}:=\sum_{P_3:P_2\supset P_2}
(-1)^{\mathrm{dim}(A_3/A_2)}\tau_1^3(H)\cdot\hat\tau_3(H),$$ for 
$H\in\frak a_0$. Then we have 
\vskip 0.30cm
\noindent
{\bf Lemma 1.} ([Ar2]) If $P_1\subset P_2$, $\sigma_1^2$ {\it is a characteristic 
function of the subset of $H\in\frak a_1$} such that

\noindent
(i) $\alpha(H)>0$ for all $\alpha\in\Delta_1^2$;

\noindent
(ii) $\sigma(H)\leq 0$ for all $\sigma\in\Delta_1\backslash \Delta_1^2$; and

\noindent
(iii) $\varpi(H)>0$ for all $\varpi\in\hat\Delta_2$.

As a spacial case, with $P_1=P_2$, we get the following important consequence: 
\vskip 0.30cm
\noindent
{\bf Langlands' Combinatorial Lemma.} {\it If $Q\subset P$ are parabolic 
subgroups, then for all $H\in\frak a_0$,}
$$\begin{aligned}\sum_{R:Q\subset R\subset P}
(-1)^{\mathrm{dim} (A_R/A_P)}\tau_Q^R(H)\hat\tau_R^P(H)=&\delta_{QP};\\
\sum_{R:Q\subset R\subset P}(-1)^{\mathrm{dim} (A_Q/A_R)}\widehat\tau_Q^R(H)
\tau_R^P(H)=&\delta_{QP}.\end{aligned}$$

Suppose now that $Q\subset P$ are parabolic subgroups. Fix a vector
$\Lambda\in\frak a_0^*$. Let $$\varepsilon_Q^P(\Lambda):=
(-1)^{\#\{\alpha\in\Delta_Q^P:\Lambda(\alpha^\vee)\leq 0\}},$$ 
and let $$\phi_Q^P(\Lambda,H),\qquad H\in\frak a_0,$$ be the characteristic 
function of the set  $$\Big\{
H\in\frak a_0: 
\begin{matrix}
\varpi(H)>0, &\mathrm{if}\ \Lambda(\alpha^\vee)\leq 0\\ 
\varpi(H)\leq 0,&\mathrm{if}\ \Lambda(\alpha^\vee)>0
\end{matrix},
\forall\alpha\in\Delta_Q^P
\Big\}.$$

\noindent
{\bf Lemma 2.}([Ar2,3]) With the same notation as above,
$$\sum_{R:Q\subset R\subset P}\varepsilon_Q^R(\Lambda)\cdot
\phi_Q^R(\Lambda,H)\cdot\tau_R^P(H)=\begin{cases} 0, &\text{if 
$\Lambda(\alpha^\vee)\leq 0,\ 
\exists \alpha\in\Delta_Q^P$}\\
1,&\text{otherwise}\end{cases}.$$

\subsection{Langlands-Arthur's Partition: Reduction Theory}

Our aim here is to derive Langlands-Arthur's partition of 
$G(F)\backslash G(\mathbb A)$ into disjoint subsets, one 
for each (standard) parabolic subgroup.

To start with, suppose that $\omega$ is a compact subset of 
$N_0(\mathbb A)M_0(\mathbb A)^1$ and that $T_0\in-\frak a_0^+$. 
For any parabolic subgroup $P_1$, introduce the associated {\it Siegel set}
$\frak s^{P_1}(T_0,\omega)$ as the collection of
$$pak,\qquad p\in\omega,\ a\in A_0(\mathbb R)^0,\ k\in K,$$
where $\alpha\Big(H_0(a)-T_0\Big)$ is positive for each 
$\alpha\in\Delta_0^1$. Then from classical reduction theory,
we conclude that {\it for sufficiently big $\omega$ and sufficiently small
$T_0$, $G(\mathbb A)=P_1(F)\cdot \frak s^{P_1}(T_0,\omega).$}
\vskip 0.30cm
Suppose now that $P_1$ is given. Let $\frak s^{P_1}(T_0,T,\omega)$ be the set 
of $x$ in $\frak s^{P_1}(T_0,\omega)$ such that 
$\varpi\Big(H_0(x)-T\Big)\leq 0$ for each $\varpi\in\hat\Delta_0^1$. 
Let $F^{P_1}(x,T):=F^1(x,T)$ be the characteristic function of the 
set of $x\in G(\mathbb A)$ such that $\delta x$ belongs to 
$\frak s^{P_1}(T_0,T,\omega)$ for some $\delta\in P_1(F)$.

As such, $F^1(x,T)$ is left $A_1(\mathbb R)^0N_1(\mathbb A)M_1(F)$-invariant, 
and can be regarded as the characteristic function of the 
projection of $\frak s^{P_1}(T_0,T,\omega)$ onto 
$A_1(\mathbb R)^0N_1(\mathbb A)M_1(F)\backslash G(\mathbb A)$, 
a compact subset of the quotient space
$A_1(\mathbb R)^0N_1(\mathbb A)M_1(F)\backslash 
G(\mathbb A)$.
\vskip 0.20cm
For example,  $F(x,T):=F^G(x,T)$ admits the following more direct description
which will play a key role in our study of Arthur's periods:

If $P_1\subset P_2$ are (standard) parabolic subgroups, we write 
$A_1^\infty:=A_{P_1}^\infty$ for $A_{P_1}(\mathbb A)^0$, the identity 
component of $A_{P_1}(\mathbb R)$, and $$A_{1,2}^\infty:=
A_{P_1,P_2}^\infty:=A_{P_1}\cap M_{P_2}(\mathbb A)^1.$$ Then the logarthmic map
$H_{P_1}$ maps $A_{1,2}^\infty$ isomorphically onto $\frak a_1^2$, 
the orthogonal complement of $\frak a_2$ in $\frak a_1$. If $T_0$ and 
$T$ are points in $\frak a_0$, set $A_{1,2}^\infty(T_0,T)$ to be the set
$$\Big\{a\in A_{1,2}^\infty:\alpha\Big(H_1(a)-T\Big)>0,\ 
\alpha\in\Delta_1^2;\
\varpi\Big(H_1(a)-T\Big)<0,\,\varpi\in\hat\Delta_1^2\Big\},$$ where
$\Delta_1^2:=\Delta_{P_1\cap M_2}$ and $\hat\Delta_1^2
:=\hat\Delta_{P_1\cap M_2}.$ In particular, for $T_0$ such that $-T_0$ 
is suitably regular,
{\it $F(x,T)$ is the characteristic function of the compact subset of 
$G(F)\backslash G(\mathbb A)^1$ obtained by projecting 
$$N_0(\mathbb A)\cdot M_0(\mathbb A)^1\cdot A_{P_0,G}^\infty(T_0,T)\cdot K$$ 
onto $G(F)\backslash G(\mathbb A)^1.$}

All in all, we arrive at the following

\noindent
{\bf {Arthur's Partition.}} (Arthur) {\it Fix $P$ and let $T$ be any suitably 
point in $T_0+\frak a_0^+$. Then} $$\sum_{P_1:P_0\subset P_1\subset P}
\sum_{\delta\in P_1(F)\backslash G(F)}F^1(\delta x)\cdot\tau_1^P
\Big(H_0(\delta x)-T\Big)=1\qquad\forall x\in G(\mathbb A).$$

\section{Arthur's Analytic Truncation}
\subsection{Definition}
Following Arthur, we make the following
\vskip 0.20cm
\noindent 
{\bf Definition.} (Arthur) {\it Fix a suitably regular point $T\in\frak a_0^+$. 
If $\phi$ is a continuous 
function on $G(F)\backslash G(\mathbb A)^1$, define Arthur's analytic 
trunction $\Big(\Lambda^T\phi\Big)(x)$ to be the function
$$\Big(\Lambda^T\phi\Big)(x):=
\sum_P(-1)^{\mathrm{dim} (A/Z)}\sum_{\delta\in P(F)\backslash G(F)}
\phi_P(\delta x)\cdot\hat\tau_P\Big(H(\delta x)-T\Big),$$ where
$$\phi_P(x):=\int_{N(F)\backslash N(\mathbb A)}\phi(nx)\,dn$$ denotes the 
constant term of $\phi$ along $P$, and the sum is over all (standard) 
parabolic subgroups.}
\vskip 0.20cm
The main purpose for introducing analytic truncation is to give a natural way
to construct integrable functions: even from the example of $GL_2$, we 
know that automorphic forms are generally not integrable over the total 
fundamental domain $G(F)\backslash G(\mathbb A)^1$ mainly due to the fact 
that in the Fourier expansions of such functions,  constant 
terms are only of moderate growth (hence not integrable). Thus in order 
to naturally obtain  integrable functions, we should truncate the original 
function along the cuspidal regions by removing constant terms. 
Simply put, Arthur's analytic truncation is a 
well-designed divice in which constant terms are tackled in such a way 
that  different levels of parabolic subgroups are suitably counted at 
the corresponding cuspidal region so that the whole truncation will not 
be overdone while there will be no parabolic subgroups left untackled.

Note that all parabolic subgroups of $G$ can be obtained from standard 
parabolic subgroups by taking conjugations with elements from 
$P(F)\backslash G(F)$. So we have:
\vskip 0.30cm
\noindent
(a) $\displaystyle{\Big(\Lambda^T\phi\Big)(x)=\sum_P(-1)^{\mathrm{dim} (A/Z)}
\phi_P(x)
\cdot\hat\tau_P\Big(H(x)-T\Big),}$ {\it where the sum is over all, both
standard and non-standard, parabolic subgroups};
\vskip 0.30cm
\noindent
(b) {\it If $\phi$ is a cusp form, then $\Lambda^T\phi=\phi$};

This is because by definition, all constant terms along proper $P: P\not=G$ 
are zero.
Moreover, as a direct consequence of the Basic Estimation for partial 
truncation, we have
\vskip 0.30cm
\noindent
(c) {\it If $\phi$ is of moderate growth} in the sense that there exist some 
constants $C,N$ such that $\Big|\phi(x)\Big|\leq c\|x\|^N$ for all 
$x\in G(\mathbb A)$, {\it then so is $\Lambda^T\phi$}.

\subsection{Basic Properties}
 
Recall that an element $T\in\frak a_0^+$ is called 
{\it sufficiently regular}, if for any 
$\alpha\in\Delta_0$, $\alpha(T)\gg 0$.
Fundamental properties of Arthur's analytic truncation 
may be summarized as follows:
\vskip 0.20cm
\noindent
{\bf Proposition.} (Arthur) {\it For sufficiently regular $T$ in $\frak a_0$,} 

\noindent
(1) {\it Let $\phi:G(F)\backslash 
G(\mathbb A)\to\mathbb C$ be a locally $L^1$ function. Then
$$\Lambda^T\Lambda^T\phi(g)=\Lambda^T\phi(g)$$ for almost all $g$. 
If $\phi$ is also locally bounded,  then the above is true for all $g$};

\noindent
(2) {\it Let $\phi_1,\,\phi_2$ be two locally $L^1$ functions on 
$G(F)\backslash G(\mathbb A)$. Suppose that 
$\phi_1$ is of moderate growth and $\phi_2$ is	rapidly decreasing. Then
$$\int_{Z_{G(\mathbb A)}G(F)\backslash G(\mathbb A)}\overline{\Lambda^T
\phi_1(g)}\cdot \phi_2(g)\,dg
=\int_{Z_{G(\mathbb A)}G(F)\backslash G(\mathbb A)}\overline{\phi_1(g)}\cdot
\Lambda^T\phi_2(g)\,dg;$$}

\noindent
(3) {\it Let $K_f$ be an open compact subgroup of $G(\mathbb A_f)$, and 
$r, r'$ are 
two positive real numbers. Then there exists a finite subset 
$\Big\{X_i:i=1,2,\cdots,N\Big\}\subset\frak u$, the universal enveloping 
algebra of the Lie algebra associated to $G(\mathbb A_\infty)$, such that the 
following is satisfied: Let $\phi$ be a smooth function on 
$G(F)\backslash G(\mathbb A)$, right invariant under $K_f$ and let 
$a\in A_{G(\mathbb A)},\ g\in G(\mathbb A)^1\cap S$. Then 
$$\Big|\Lambda^T\phi(ag)\Big|\leq \|g\|^{-r}\sum_{i=1}^N
\sup\Big\{|\delta(X_i)\phi(ag')|\,\|g'\|^{-r'}:g'\in G(\mathbb A)^1\Big\},$$
where $S$ is a	Siegel domain with respect to $G(F)\backslash G(\mathbb A).$}

\subsection{Truncation $\Lambda^T{\bold 1}$}

To go further, let us give a much more detailed study of Authur's analytic 
truncation for the constant function ${\bold 1}$. Fix a sufficiently regular $T\in\frak a_0$. 
Introduce the truncated subset
$\Sigma(T):=\Big(Z_{G({\mathbb A})}G(F)\backslash G({\mathbb A})\Big)_T$ of 
the space $G(F)\backslash G(\mathbb A)^1$ by
$$\Sigma(T):=\Big(Z_{G({\mathbb A})}G(F)\backslash G({\mathbb A})\Big)_T:=
\Big\{g\in Z_{G({\mathbb A})}G(F)\backslash G({\mathbb A}):
\Lambda^T{\bf 1}(g)=1\Big\}.$$ We claim that
$\Sigma(T)$ or the same 
$\Big(Z_{G({\mathbb A})}G(F)\backslash G({\mathbb A})\Big)_T$, is compact. 
In fact, 
much stronger result is correct. Namely, we have the following
\vskip 0.20cm
\noindent
{\bf Lemma.}  (Arthur) {\it For sufficiently regular $T\in\frak a_0^+$, 
$\Lambda^T{\bold 1}(x)=F(x,T).$ That is to say,
$\Lambda^1{\bold 1}$ is the characteristic function of the compact subset
$\Sigma(T)$ of $G(F)\backslash G(\mathbb A)^1$ obtained by projecting 
$N_0(\mathbb A)\cdot M_0(\mathbb A)^1\cdot A_{P_0,G}^\infty(T_0,T)\cdot K$
onto $G(F)\backslash G(\mathbb A)^1.$} 

\section{Analytic Characterization of Stability}
\subsection{A Micro Bridge}

For simplicity, we in this subsection work only with 
the field of rationals $\mathbb Q$ and use mixed
languages of adeles and lattices. Also,	 
without loss of generality, we assume that $\mathbb Z$-lattices are of 
volume one. Accordingly, set $G=SL_r$. 

For a rank $r$ lattice $\Lambda$ of volume one,	 denote the sublattices 
filtration associated to  a parabolic subgroup $P$ by $$0=\Lambda_0\subset
\Lambda_1\subset\Lambda_2\subset\cdots\subset\Lambda_{|P|}=\Lambda.$$ 
Assume that $P$ corresponds to the partition $I=(d_1,d_2,\cdots,d_{|P|})$.
A polygon $p:[0,r]\to\mathbb R$ is called {\it normalized} 
if $p(0)=p(r)=0$. For a (normalized) polygon $p:[0,r]\to\mathbb R$,
define the associated (real) character $T=T(p)$ of 
$M_0$ by the condition that $$\alpha_i(T)=\Big[p(i)-p(i-1)\Big]-
\Big[p(i+1)-p(i)\Big]$$ 
for all $i=1,2,\cdots, r-1.$ Then one checks that $T(p)$ coincides with
$$\Big(p(1), p(2)-p(1),
\cdots, p(i)-p(i-1),\cdots,p(r-1)-p(r-2),-p(r-1)\Big).$$

Now take $g=g(\Lambda)\in G(\mathbb A)$. Denote its 
lattice by $\Lambda^g$, and its induced filtration from $P$ by
$$0=\Lambda_0^{g,P}\subset \Lambda_1^{g,P}\subset\cdots\subset 
\Lambda_{|P|}^{g,P}=\Lambda^g.$$ Consequently, the polygon 
$p_P^g=p_P^{\Lambda^g}:[0,r]\to\mathbb R$ is characterized by

\noindent
(1) $p_P^g(0)=p_P^g(r)=0$;

\noindent
(2) $p_P^g$ is affine on $[r_i,r_{i+1}]$, $i=1,2,\cdots,|P|-1$; and

\noindent
(3) $p_P^g(r_i)=\mathrm{deg}\big(\Lambda_i^{g,P}\big)-r_i\cdot
\frac{\mathrm{deg}\big(\Lambda^{g}\big)}{r},\ i=1,2,\cdots,|P|-1$.

\noindent
Note that the volume of $\Lambda$ is assumed to be one, therefore 
(3) is equivalent to 

\noindent
(3)$'$ $p_P^g(r_i)=\mathrm{deg}\big(\Lambda_i^{g,P}\big),\ i=1,2,\cdots,|P|-1$.
\vskip 0.20cm
The advantage of partially using adelic language is that	
the values of $p_P^g$ may be written down precisely.
Indeed, using Langlands decompositon $g=n\cdot m\cdot a(g)\cdot k$ 
with $n\in N_P(\mathbb A), m\in M_P(\mathbb A)^1, a\in A_+$ and 
$k\in K:=\prod_pSL(\mathcal O_{{\mathbb Q}_p})\times SO(r)$. Write 
$$a=a(g)=\mathrm{diag}\Big(a_1I_{d_1},a_2I_{d_2},\cdots,a_{|P|}I_{d_{|P|}}
\Big)$$ where $r=d_1+d_2+\cdots+d_{|P|}$ is the partition corresponding to $P$.
Then it is a standard fact that
$$\mathrm{deg}\Big(\Lambda_i^{g,P}\Big)=-\log\Big(\prod_{j=1}^i a_j^{d_j}\Big)
=-\sum_{j=1}^i d_j\log a_j,\qquad i=1,\cdots,|P|.$$

Set now ${\bf 1}(p_P^*>_Pp)$ to be the characteristic function 
of the subset
of $g$'s such that $p_P^g>_Pp$. Then by a certain calculation, we obtain the following
\vskip 0.30cm
\noindent
{\bf Micro Bridge.} (Lafforgue, Weng) {\it For a fixed convex normalized polygon 
$p:[0,r]\to\mathbb R$, and $g\in SL_r(\mathbb A)$, with respect to
 any parabolic subgroup $P$, we have}
 $$\hat\tau_P\Big(-H_0(g)-T(p)\Big)=\bold 1\Big(p_P^g>_Pp\Big).$$

\subsection{Analytic Truncations and Stability}
\vskip 0.30cm
With the micro bridge above, now we are ready to state the following
analytic characterization of stability.
\vskip 0.30cm
\noindent
{\bf Global Bridge.} (Lafforgue, Weng) {\it For a fixed normalized convex polygon $p:[0,r]\to 
{\mathbb R}$, let $T(p)\in\frak a_0$ be the associated vector defined by
$$\Big(p(1), p(2)-p(1),
\cdots, p(i)-p(i-1),\cdots,p(r-1)-p(r-2),-p(r-1)\Big).$$
If $T(p)$ is sufficiently positive,
then} $$\bold 1(\overline p^g\leq p)=\Big(\Lambda^{T(p)}\bold 1\Big)(g).$$
\vfill\eject
\centerline
{\Large\bf Chapter VII. Non-Abelian $L$-Functions}
\section{High Rank Zetas and Eisenstein Series}
\subsection{Epstein Zeta Functions and High Rank Zetas}

Recall that the rank $r$ non-abelian zeta function 
$\xi_{{\mathbb	Q},r}(s)$ of ${\mathbb	Q}$ is given by
$$\xi_{{\mathbb	 Q},r}(s)=
\int_{{\mathcal	 M}_{{\mathbb  Q},r}}\left(e^{h^0({\mathbb  Q},\Lambda)}-1
\right)
\cdot \big(e^{-s}\big)^{\mathrm{deg}(\Lambda)}\, d\mu(\Lambda),\qquad 
\mathrm{Re}(s)>1,$$
 with $e^{h^0({\mathbb	Q},\Lambda)}:=\sum_{x\in \Lambda}
\exp\big(-\pi|x|^2\big)$ and $\mathrm{deg}(\Lambda)=-\log\,
\mathrm{Vol}\big(\mathbb R^r/\Lambda\big)$.

Decompose according to their volumes to get
${\mathcal  M}_{{\mathbb  Q},r}=\cup_{T>0}{\mathcal  M}_{{\mathbb  Q},r}[T]$.
Using the natural morphism
${\mathcal  M}_{{\mathbb  Q},r}[T]\to {\mathcal	 M}_{{\mathbb  Q},r}[1],
\ \Lambda\mapsto T^{1\over r}\cdot\Lambda,$
we obtain $$\begin{aligned} \xi_{{\mathbb  Q},r}(s)=&
\int_{\cup_{T>0}{\mathcal  M}_{{\mathbb	 Q},r}[T]}\left( e^{h^0({\mathbb  Q},
\Lambda)}-1
\right)\cdot \big(e^{-s}\big)^{\mathrm{deg}(\Lambda)}\, d\mu(\Lambda)\\
=&\int_0^\infty T^s{{dT}\over T}
\int_{{\mathcal	 M}_{{\mathbb  Q},r}[1]}\left( 
e^{h^0({\mathbb	 Q},T^{1\over r}\cdot \Lambda)}-1\right)
\cdot d\mu(\Lambda).\end{aligned}$$
But $$h^0({\mathbb  Q},T^{1\over r}\cdot \Lambda)=
\log\left(\sum_{x\in \Lambda}
\exp\big(-\pi|x|^2\cdot T^{2\over r}\big)\right)$$
and  $$\int_0^\infty e^{-AT^B}T^s{{dT}\over T}={1\over B}\cdot A^{-{s\over B}}
\cdot\Gamma({s\over B}),\qquad	B\not=0,$$ we have
$\xi_{{\mathbb	Q},r}(s)={r\over 2}\cdot\pi^{-{r\over 2}\, s}
\Gamma({r\over 2}\, s)\cdot\int_{{\mathcal  M}_{{\mathbb  Q},r}[1]}
\left(\sum_{x\in\Lambda\backslash\{0\}}|x|^{-rs}\right)\cdot d\mu_1(\Lambda).$
Accordingly, introduce the completed Epstein zeta function for
$\Lambda$ by 
$$\hat E(\Lambda;s):=\pi^{-s}\Gamma(s)\cdot \sum_{x\in \Lambda\backslash \{0\}}
|x|^{-2s}.$$

\noindent
{\bf Proposition.} (Weng) ({\bf Eisenstein Series and High Rank Zetas}) 
$$\xi_{{\mathbb	 Q},r}(s)={r\over 2}\int_{{\mathcal  M}_{{\mathbb  Q},r}[1]}
\hat E(\Lambda,{r\over 2}s)\,d\mu_1(\Lambda).$$

\subsection{Rankin-Selberg Method: An Example with $SL_2$}
Consider the action of $\mathrm{SL}(2,{\mathbb	Z})$ on the upper half plane 
${\mathcal  H}$.
Then a standard	 \lq fundamental domain\rq$\ $ is given by 
$D=\{z=x+iy\in {\mathcal  H}:|x|\leq {1\over 2},y>0,x^2+y^2\geq 1\}.$ 
Recall also the completed standard Eisenstein series 
$$\hat E(z;s):=\pi^{-s}\Gamma(s)\cdot\sum_{(m,n)\in {\mathbb  Z}^2
\backslash \{(0,0)\}}{{y^{s}}\over {|mz+n|^{2s}}}.$$
Naturally, we are led to  the integral
$\int_D\hat E(z,s){{dx\,dy}\over {y^2}}.$
However, this integration diverges. Indeed,
near the only cusp $y=\infty$, $\hat E(z,s)$ has the  Fourier expansion
$$\hat E(z;s)=\sum_{n=-\infty}^\infty a_n(y,s)e^{2\pi i nx}$$ with 
$$a_n(y,s)=\begin{cases} \xi(2s)y^s+\xi(2-2s)y^{1-s},&	\text{if $n=0$;}\\
2|n|^{s-{1\over 2}}\sigma_{1-2s}(|n|){\sqrt y}K_{s-{1\over 2}}(2\pi|n|y),& 
\text{if $n\not=0$,} \end{cases}$$
where $\xi(s)$ is the completed Riemann zeta function,
$\sigma_s(n):=\sum_{d|n}d^s$, and 
$K_s(y):={1\over 2}\int_0^\infty e^{-y(t+{1\over t})/2}t^s{{dt}\over t}$ 
is the K-Bessel function. Moreover, 
$$|K_s(y)|\leq e^{-y/2}K_{\mathrm{Re}(s)}(2),\ \mathrm{if}\  y>4,\qquad
\mathrm{and}\qquad K_s=K_{-s}.$$ 
So $a_{n\not=0}(y,s)$ decay exponentially, and the problematic 
term comes from $a_0(y,s)$, which is of slow growth.

Therefore, to make the original integration meaningful, 
we need to cut-off the slow growth part. Recall from the discussions
in previous three chapters, 
we have two different ways to do so: 
one is geometric and hence rather direct and simple; while the other is 
analytic, and hence rather technical and traditional, dated back to 
Rankin-Selberg.
\vskip 0.30cm
\noindent
(a) {\bf Geometric Truncation}

Draw a horizontal line $y=T\geq 1$ and set
 $$D_T=\{z=x+iy\in D:y\leq T\},\qquad D^T=\{z=x+iy\in D:y\geq T\}.$$ 
Then $D=D_T\cup D^T$.
Introduce a well-defined integration
$$I^{\mathrm{Geo}}_T(s):=\int_{D_T}\hat E(z,s)\,{{dx\,dy}\over {y^2}}.$$

\noindent
(b) {\bf Analytic Truncation}

Define a truncated Eisenstein series $\hat E_T(z;s)$ by
$$\hat E_T(z;s):=\begin{cases} \hat E(z;s),&\text{if $y\leq T$;}\\
\hat E(z,s)-a_0(y;s),&\text{if $y>T$.}\end{cases}$$
Introduce a well-defined integration
$$I_T^{\mathrm{Ana}}(s):=\int_D\hat E_T(z;s)\,{{dx\,dy}\over {y^2}}.$$

With this, from the Rankin-Selberg method, 
one checks that we have the following:
\vskip 0.30cm
\noindent
{\bf Proposition.} ([W2,3,5]) ({\bf Analytic Truncation=Geometric Truncation in Rank Two}) 
{\it $$I_T^{\mathrm{Geo}}(s)=\xi(2s){{T^{s-1}}\over {s-1}}-\xi(2s-1){{T^{-s}}
\over {s}}=I_T^{\mathrm{Ana}}(s).$$}

Each of the above two integrations has its own merit: for the geometric one, 
we keep the Eisenstein series unchanged, while for the analytic one, we keep
the original fundamental domain of ${\mathcal  H}$ under 
$\mathrm{SL}(2,{\mathbb	 Z})$
as it is.

Note that the nice point about the fundamental domain is that it admits a
modular interpretation. Thus it would be very idealistic if we could at the
same time keep the Eisenstein series unchanged, while offer some integration
domains which appear naturally in certain moduli problems. Guided by this,
in the follows, we will introduce non-abelian $L$-functions using integrations of 
Eisenstein series over generalized moduli spaces.
\vskip 0.30cm
\noindent
(c) {\bf  Arithmetic Truncation}

Now we explain why above discussion and Rankin-Selberg method have 
anything to do with our non-abelian zeta functions. For this, we introduce 
yet another truncation, the algebraic, or better arithmetic, one.

So back to the moduli space  of rank 2 lattices of volume 1 over 
${\mathbb  Q}$.
Then classical reduction theory gives a	 natural map from this moduli
space to the fundamental domain $D$ of $\mathrm{SL}(2,\mathbb Z)$ on $\mathcal H$: 
For any lattice $\Lambda$, fix ${\bold x}_1\in \Lambda$ such that its length	
gives the first Minkowski minimum $\lambda_1$ of $\Lambda$ ([Min]).
Then via rotation, we may assume that ${\bold x}_1=(\lambda_1,0)$. 
Further, from the reduction theory 
${1\over {\lambda_1}}\Lambda$ may be viewed as the lattice of the volume
$\lambda_1^{-2}=y_0$ which is generated by $(1,0)$ and $\omega=x_0+iy_0\in D$. 
That is to say, 
the points in $D_T$ are in one-to-one corresponding to
the  rank two lattices of volume one whose
first Minkowski minimum $\lambda_1^{-2}\leq T$, i.e,
$\lambda_1\geq T^{-{1\over 2}}$. Set 
${\mathcal  M}_{{\mathbb  Q},2}^{\leq {1\over 2}\log T}[1]$  be the moduli 
space of
rank 2 lattices $\Lambda$ of volume 1 over ${\mathbb  Q}$ whose sublattices 
$\Lambda_1$ of rank 1 have degrees 
$\leq {1\over 2}\log T$. As a direct consequence, we have the following
\vskip 0.30cm
\noindent
{\bf Proposition.} (Geometric Truncation = Algebraic Truncation) {\it 
There is a natural one-to-one, onto morphism 
$${\mathcal  M}_{{\mathbb  Q},2}^{\leq {1\over 2}\log T}[1]\simeq D_T.$$
In particular, $${\mathcal  M}_{{\mathbb  Q},2}^{\leq 0}[1]=
{\mathcal  M}_{{\mathbb	 Q},2}[1]\simeq D_1.$$}

Consequently,  we have the following
\vskip 0.30cm
\noindent
{\bf Example in Rank 2}. {\it	
$\xi_{{\mathbb	Q},2}(s)={{\xi(2s)}\over {s-1}}-{{\xi(2s-1)}\over {s}}.$}

\section{Non-Abelian $L$-Functions: Definitions}
\subsection{Automorphic Forms and Eisenstein Series}

To faciliate our ensuing discussion, we make the following preparations.
Here,  as usual, instead of parabolic subgroups $P$, 
we adopt their Levi subgroups $M$ as running symbols. 
For details, see e.g., [MW] and [W-1].

Fix a connected reduction group $G$ defined over $F$, denote by $Z_G$ its 
center. Fix a minimal parabolic subgroup $P_0$ of $G$. Then $P_0=M_0N_0$, 
where as usual we fix once and for all the Levi $M_0$ and  the unipotent
 radical $N_0$. Recall that a parabolic subgroup $P$ is $G$ is called standard if 
$P\supset P_0$. For such groups, write $P=MN$ with $M_0\subset M$ the 
standard Levi and $N$ the unipotent radical. Denote by $\mathrm {Rat}(M)$ 
the group of rational characters of $M$, i.e, the morphism 
$M\to {\mathbb G}_m$ 
where ${\mathbb G}_m$ denotes the multiplicative group. Set 
$\frak a_M^*:=\mathrm {Rat}(M)\otimes_{\mathbb Z}{\mathbb C},\ \frak a_M
:=\mathrm{Hom}_{\mathbb Z}(\mathrm {Rat}(M),{\mathbb C}),$ and 
$\mathrm {Re}\frak 
a_M^*:=\mathrm {Rat}(M)\otimes_{\mathbb Z}{\mathbb R},\ 
\mathrm{Re}\frak a_M
:=\mathrm{Hom}_{\mathbb Z}(\mathrm {Rat}(M),{\mathbb R}).$ For any $\chi\in 
\mathrm {Rat}(M)$, we obtain a (real) character $|\chi|:M({\mathbb A})\to 
{\mathbb R}^*$ defined by $m=(m_v)\mapsto m^{|\chi|}:=\prod_{v\in S}
|m_v|_v^{\chi_v}$ with $|\cdot|_v$ the $v$-absolute values. Set then 
$M({\mathbb A})^1:=\cap_{\chi\in \mathrm {Rat}(M)}\mathrm{Ker}|\chi|$, 
which is 
a normal subgroup of $M({\mathbb A})$. Set $X_M$ to be the group of complex 
characters which are trivial on $M({\mathbb A})^1$. Denote by 
$H_M:=\log_M:M({\mathbb A})\to \frak a_M$ the map such that 
$\forall\chi\in \mathrm {Rat}(M)\subset \frak a_M^*,\langle\chi,
\log_M(m)\rangle:=\log(m^{|\chi|})$. Clearly, $M({\mathbb A})^1
=\mathrm{Ker}(\log_M);\ \log_M(M({\mathbb A})/M({\mathbb A})^1)\simeq 
\mathrm{Re}\frak a_M.$ Hence in particular there is a natural isomorphism 
$\kappa:\frak a_M^*\simeq X_M.$
Set $\mathrm{Re}X_M:=\kappa (\mathrm{Re}\frak a_M^*),\ \mathrm{Im}X_M:=
\kappa (i\cdot \mathrm{Re}\frak a_M^*).$ Moreover define our working space 
$X_M^G$ to be the subgroup of $X_M$ consisting of complex characters of 
$M({\mathbb A})/M({\mathbb A})^1$ which are trivial on $Z_{G({\mathbb A})}$.

Fix a maximal compact subgroup $K$ such that for all standard 
parabolic
subgroups $P=MN$ as above, $P({\mathbb A})\cap K=M({\mathbb A})
\cap K\cdot
U({\mathbb A})\cap{K}.$ Hence we get the Langlands decomposition 
$G({\mathbb A})=M({\mathbb A})\cdot N({\mathbb A})\cdot {K}$. 
Denote by $m_P:G({\mathbb A})\to M({\mathbb A})/M({\mathbb A})^1$ the map 
$g=m\cdot n\cdot k\mapsto M({\mathbb A})^1\cdot m$ where $g\in G({\mathbb A}), 
m\in M({\mathbb A}), n\in N({\mathbb A})$ and $k\in {K}$.
\vskip 0.20cm
Fix Haar measures on $M_0({\mathbb A}), N_0({\mathbb A}), {K}$ 
respectively such that	the induced measure on $N_0(F)$ is the counting 
measure and the volumes 
of $N(F)\backslash N_0({\mathbb A})$ and $K$ are all 1.

Such measures then also induce Haar measures via $\log_M$ to $\frak a_{M_0}, 
\frak a_{M_0}^*$, etc. Furthermore, if we denote by $\rho_0$ the half of the 
sum of the positive roots of  the maximal split torus $T_0$ of the central 
$Z_{M_0}$
of $M_0$, then $f\mapsto \int_{M_0({\mathbb A})\cdot N_0({\mathbb A})\cdot 
{K}}f(mnk)\,dk\,dn\,m^{-2\rho_0}dm$ defined for continuous functions 
with compact supports on $G({\mathbb A})$  defines a Haar measure $dg$ on 
$G({\mathbb A})$. This in turn gives measures on $M({\mathbb A}), 
N({\mathbb A})$ 
and hence on $\frak a_{M}, \frak a_{M}^*$, $P({\mathbb A})$, etc, for all 
parabolic subgroups $P$. In particular, the following 
compactibility condition
$$\begin{aligned}\int_{M_0({\mathbb A})\cdot N_0({\mathbb A})\cdot {K}}&
f(mnk)\,dk\,dn\,m^{-2\rho_0}dm\\
=&\int_{M({\mathbb A})\cdot N({\mathbb A})\cdot {K}}f(mnk)\,dk\,
dn\,m^{-2\rho_P}dm\end{aligned}$$ holds 
for all continuous functions $f$ with compact supports 
on $G({\mathbb A})$, where $\rho_P$ denotes the half of the sum of the 
positive roots of  the maximal split torus $T_P$ of the central $Z_{M}$ of 
$M$. For later use, 
denote also by $\Delta_P$ the set of positive roots determined by $(P,T_P)$ 
and $\Delta_0=\Delta_{P_0}$.
\vskip 0.20cm
Fix an isomorphism $T_0\simeq {\mathbb G}_m^R$. Embed ${\mathbb R}_+^*$ by 
the map 
$t\mapsto (1;t)$. Then we obtain a natural injection $({\mathbb R}_+^*)^R
\hookrightarrow T_0({\mathbb A})$ which splits. Denote by 
$A_{M_0({\mathbb A})}$ 
the unique connected subgroup of $T_0({\mathbb A})$ which projects onto 
$({\mathbb R}_+^*)^R$. More generally, for a standard parabolic subgroup 
$P=MN$, set $A_{M({\mathbb A})}:=
A_{M_0({\mathbb A})}\cap Z_{M({\mathbb A})}$ where as used above $Z_*$ denotes 
the center of the group $*$. Clearly, $M({\mathbb A})=A_{M({\mathbb A})}\cdot 
M({\mathbb A})^1$. For later use, set also 
$A_{M({\mathbb A})}^G:=\{a\in A_{M({\mathbb A})}:\log_Ga=0\}.$ 
Then $A_{M({\mathbb A})}=A_{G({\mathbb A})}\oplus
A_{M({\mathbb A})}^G.$

Note that ${K}$, $M(F)\backslash M({\mathbb A})^1$ and $N(F)\backslash 
N({\mathbb A})$ are all compact, thus with the Langlands decomposition 
$G({\mathbb A})=N({\mathbb A})M({\mathbb A}){K}$ in mind, the 
reduction theory 
for $G(F)\backslash G({\mathbb A})$ or more generally 
$P(F)\backslash G({\mathbb A})$ is reduced to that for 
$A_{M({\mathbb A})}$ since $Z_G(F)\cap Z_{G({\mathbb A})}\backslash 
Z_{G({\mathbb A})}\cap G({\mathbb A})^1$ is compact as well. As such 
for $t_0\in M_0({\mathbb A})$ set
$A_{M_0({\mathbb A})}(t_0):=\{a\in A_{M_0({\mathbb A})}:a^\alpha>t_0^\alpha
\forall\alpha\in\Delta_0\}.$ Then, for a fixed compact subset 
$\omega\subset P_0({\mathbb A})$, we have the corresponding Siegel set 
$S(\omega;t_0):=\{p\cdot a\cdot k:p\in \omega, a\in A_{M_0({\mathbb A})}
(t_0),k\in {K}\}.$ In particular, for big enough $\omega$ and small 
enough $t_0$, i.e, $t_0^\alpha$ is very close to 0 for all 
$\alpha\in\Delta_0$, the classical reduction theory may be restated as 
$G({\mathbb A})=G(F)\cdot S(\omega;t_0)$. More generally set 
$A_{M_0({\mathbb A})}^P(t_0):=\{a\in A_{M_0({\mathbb A})}:a^\alpha>t_0^\alpha
\forall\alpha\in\Delta_0^P\},$ and $S^P(\omega;t_0):=\{p\cdot a\cdot k:
p\in \omega, a\in A_{M_0({\mathbb A})}^P(t_0),k\in {K}\}.$	
Then similarly as above
for big enough $\omega$ and small enough $t_0$, $G({\mathbb A})=P(F)\cdot 
S^P(\omega;t_0)$. (Here 
$\Delta_0^P$ denotes the set of positive roots for $(P_0\cap M,T_0)$.)
\vskip 0.20cm
Fix an embedding $i_G:G\hookrightarrow SL_n$ sending $g$ to $(g_{ij})$. 
Introducing a hight function on $G({\mathbb A})$ by setting 
$\|g\|:=\prod_{v\in S}\mathrm{sup}\{|g_{ij}|_v:\forall i,j\}$. It is 
well-known 
that up to $O(1)$, hight functions are unique. This implies that the following 
growth conditions do not depend on the height function we choose.

A function $f:G({\mathbb A})\to {\mathbb C}$ is said to have {\it moderate 
growth} if 
 there exist 
$c,r\in {\mathbb R}$ such that $|f(g)|\leq c\cdot \|g\|^r$ for all
 $g\in G({\mathbb A})$. Similarly, for a standard parabolic subgroup $P=MN$, 
a function $f:N({\mathbb A})M(F)\backslash G({\mathbb A})\to {\mathbb C}$ 
is said 
to have moderate growth if there exist $c,r\in {\mathbb R},\lambda\in 
\mathrm{Re}X_{M_0}$ such that for any $a\in A_{M({\mathbb A})},k\in 
{K}, 
m\in M({\mathbb A})^1\cap S^P(\omega;t_0)$,
$|f(amk)|\leq c\cdot \|a\|^r\cdot m_{P_0}(m)^\lambda.$

Also a function $f:G({\mathbb A})\to {\mathbb C}$ is said to be {\it smooth} 
if for any 
$g=g_f\cdot g_\infty \in G({\mathbb A}_f)\times G({\mathbb A}_\infty)$, 
there exist 
open neighborhoods $V_*$ of $g_*$ in $G({\mathbb A})$ and a 
$C^\infty$-function 
$f':V_\infty\to {\mathbb C}$ such that $f(g_f'\cdot g_\infty')=
f'(g_\infty')$ for 
all $g_f'\in V_f$ and $g_\infty'\in V_\infty$.

By contrast, a function $f: S(\omega;t_0)\to {\mathbb C}$ is said to be 
{\it rapidly decreasing} if there exists $r>0$ and for all 
$\lambda\in \mathrm{Re}X_{M_0}$ there exists $c>0$ such that for 
$a\in A_{M({\mathbb A})}, g\in G({\mathbb A})^1\cap 
S(\omega;t_0)$, $|\phi(ag)|\leq c\cdot\|a\|\cdot m_{P_0}(g)^\lambda$. And a 
function $f:G(F)\backslash G({\mathbb A})\to {\mathbb C}$ is said to be 
rapidly decreasing if $f|_{S(\omega;t_0)}$ is so.
\vskip 0.20cm
By definition,	a function 
$\phi:N({\mathbb A})M(F)\backslash G({\mathbb A})\to {\mathbb C}$ is 
called {\it automorphic} if

\noindent
(i) $\phi$ has moderate growth; 

\noindent
(ii) $\phi$ is smooth;

\noindent
(iii) $\phi$ is ${K}$-finite, i.e, the ${\mathbb C}$-span of all 
$\phi(k_1\cdot *\cdot k_2)$ parametrized by $(k_1,k_2)\in {K}\times 
{K}$ is finite dimensional; and 

\noindent
(iv) $\phi$ is $\frak z$-finite, i.e, the ${\mathbb C}$-span of all 
$\delta(X)\phi$ parametrized by all $X\in \frak z$ is finite dimensional. Here 
$\frak z$ denotes the center of the universal enveloping algebra 
$\frak u:=\frak U(\mathrm{Lie}G({\mathbb A}_\infty))$ of the Lie algebra of	
$G({\mathbb A}_\infty)$ and $\delta(X)$
denotes the derivative of $\phi$ along $X$.

For automorphic function $\phi$, set $\phi_k:M(F)\backslash M({\mathbb A})\to 
{\mathbb C}$ by $m\mapsto m^{-\rho_P}\phi(mk)$ for all $k\in {K}$. 
Then one checks that $\phi_k$ is an automorphic form in the usual sense. Set
 $A(N({\mathbb A})M(F)\backslash G({\mathbb A}))$ be the space of automorphic 
forms on $N({\mathbb A})M(F)\backslash G({\mathbb A})$.
\vskip 0.20cm
For a measurable locally $L^1$-function $f:N(F)\backslash G({\mathbb A})\to 
{\mathbb C}$, define its {\it constant term} along with the standard parabolic 
subgroup 
$P=NM$ to be the function $f_P:N({\mathbb A})\backslash G({\mathbb A})\to
 {\mathbb C}$ 
given by
$g\to\int_{N(F)\backslash G({\mathbb A})}f(ng)dn.$ 
By definition, an automorphic form $\phi\in A(N({\mathbb A})M(F)\backslash 
G({\mathbb A}))$ is 
called {\it cuspidal} if for any standard parabolci subgroup $P'$ properly 
contained in $P$, $\phi_{P'}\equiv 0$. Denote by
 $A_0(N({\mathbb A})M(F)\backslash G({\mathbb A}))$  the space of cusp 
forms on $N({\mathbb A})M(F)\backslash G({\mathbb A})$. Obviously, 
all cusp forms are rapidly decreasing. Hence, there is a natural pairing
$$\langle\cdot,\cdot\rangle:A_0(N({\mathbb A})M(F)\backslash 
G({\mathbb A}))\times 
A(N({\mathbb A})M(F)\backslash G({\mathbb A}))\to {\mathbb C}$$ defined by
$$\langle\psi,\phi\rangle:=\int_{Z_{M({\mathbb A})}N({\mathbb A})
M(F)\backslash 
G({\mathbb A})}\psi(g)\bar\phi(g)\,dg.$$

Moreover, for a (complex) character 
$\xi:Z_{M({\mathbb A})}\to {\mathbb C}^*$, set
$$\begin{aligned}A(N({\mathbb A})M(F)\backslash G({\mathbb A}))_\xi&:=
\Big\{\phi\in A(N({\mathbb A})M(F)
\backslash G({\mathbb A})):\\
&\phi(zg)=z^{\rho_P}\cdot\xi(z)\cdot\phi(g),
\forall z\in Z_{M({\mathbb A})}, g\in G({\mathbb A})\Big\},\end{aligned}$$ 
and $A_0(N({\mathbb A})M(F)\backslash G({\mathbb A}))_\xi$ its subspace consisting	
of cusp forms.

Set now $$A_{(0)}(N({\mathbb A})M(F)\backslash G({\mathbb A}))_Z:=\sum_{\xi\in 
\mathrm{Hom}(Z_{M({\mathbb A})},{\mathbb C}^*)}A_{(0)}
(N({\mathbb A})M(F)\backslash G({\mathbb A}))_\xi.$$ Then the natural morphism 
$$\begin{aligned}{\mathbb C}
[\mathrm{Re}\frak a_M]\otimes 
A_{(0)}(N({\mathbb A})M(F)\backslash G({\mathbb A}))_Z
&\to& A_{(0)}(N({\mathbb A})M(F)\backslash G({\mathbb A}))\\
(Q,\phi)&\mapsto& \big(g\mapsto Q(\log_M(m_P(g))\big)\cdot \phi(g)\end{aligned}$$ is an 
isomorphism.

Let $\Pi_0(M({\mathbb A}))_\xi$ be isomorphism classes of irreducible 
representations of $M({\mathbb A})$ occuring in the space 
$A_0(M(F)\backslash M({\mathbb A}))_\xi$, and
$$\Pi_0(M({\mathbb A}):=\cup_{\xi\in \mathrm{Hom}(Z_{M({\mathbb A})},
{\mathbb C}^*)}\Pi_0(M({\mathbb A}))_\xi.$$ (In fact, we should 
use $M({\mathbb A}_f)\times (M({\mathbb A})\cap {K},
\mathrm{Lie}(M({\mathbb A}_\infty))\otimes_{\mathbb R}{\mathbb C}))$ 
instead of 
$M({\mathbb A})$.) For any 
$\pi\in \Pi_0(M({\mathbb A}))_\xi$, set $A_0(M(F)\backslash 
M({\mathbb A})_\pi$ to 
be the isotypic component of type $\pi$ of $A_0(M(F)\backslash 
M({\mathbb A})_\xi$, i.e, the set of cusp forms of $M({\mathbb A})$ 
generating a 
semi-simple isotypic $M({\mathbb A}_f)\times (M({\mathbb A})\cap {K},
\mathrm{Lie}(M({\mathbb A}_\infty))\otimes_{\mathbb R}{\mathbb C}))$-module 
of type $\pi$.
Set $$\begin{aligned}A_0(N({\mathbb A})M(F)\backslash G({\mathbb A}))_\pi
:=&\Big\{\phi\in 
A_0(N({\mathbb A})M(F)\backslash G({\mathbb A})):\\
&\phi_k\in A_0(M(F)\backslash 
M({\mathbb A}))_\pi,\forall k\in {K}\Big\}.\end{aligned}$$ It is quite clear that
$$A_0(N({\mathbb A})M(F)\backslash G({\mathbb A}))_\xi=\oplus_{\pi\in 
\Pi_0(M({\mathbb A}))_\xi} A_0(N({\mathbb A})M(F)\backslash 
G({\mathbb A}))_\pi.$$

More generally, let $V\subset A(M(F)\backslash M({\mathbb A}))$ be an 
irreducible 
$M({\mathbb A}_f)\times (M({\mathbb A})\cap {K},\mathrm{Lie}
(M({\mathbb A}_\infty))
\otimes_{\mathbb R}{\mathbb C}))$-module with $\pi_0$ the induced 
representation of 
$M({\mathbb A}_f)\times (M({\mathbb A})\cap {K},
\mathrm{Lie}(M({\mathbb A}_\infty))
\otimes_{\mathbb R}{\mathbb C}))$. Then we call $\pi_0$ an automorphic 
representation of $M({\mathbb A})$. Denote by $A(M(F)\backslash 
M({\mathbb A})_{\pi_0}$	 the isotypic 
subquotient module of type $\pi_0$ of $A(M(F)\backslash M({\mathbb A})$. 
One checks that
$$\begin{aligned}V\otimes \mathrm{Hom}_{M({\mathbb A}_f)\times (M({\mathbb A})\cap 
{K},
\mathrm{Lie}(M({\mathbb A}_\infty))\otimes_{\mathbb R}{\mathbb C}))}(V,&A(M(F)
\backslash M({\mathbb A})))\\
&\simeq A(M(F)\backslash M({\mathbb A}))_{\pi_0}.\end{aligned}$$
Set $$\begin{aligned}A(N({\mathbb A})M(F)\backslash G({\mathbb A}))_{\pi_0}
:=&\Big\{\phi\in 
A(N({\mathbb A})M(F)\backslash G({\mathbb A})):\\
&\phi_k\in A(M(F)\backslash 
M({\mathbb A}))_{\pi_0},\forall k\in {K}\Big\}.\end{aligned}$$
Moreover if 
$A(M(F)\backslash M({\mathbb A}))_{\pi_0}\subset
A_0(M(F)\backslash M({\mathbb A}))$, we call $\pi_0$ cuspidal.

Automorphic representations $\pi$ and $\pi_0$ of $M({\mathbb A})$ are said 
to be equivalent if $\pi\simeq \pi_0\otimes\lambda$ for some 
$\lambda\in X_M^G$. This, in practice, means that $A(M(F)\backslash 
M({\mathbb A}))_\pi
=\lambda\cdot A(M(F)\backslash M({\mathbb A}))_{\pi_0}.$  Consequently, 
$$A(N({\mathbb A})M(F)\backslash G({\mathbb A}))_\pi=(\lambda\circ m_P)\cdot 
A(N({\mathbb A})M(F)\backslash G({\mathbb A}))_{\pi_0}.$$ Denote by 
$\frak P:=[\pi_0]$ the equivalence class of $\pi_0$. Then $\frak P$ is an 
$X_M^G$-principal homogeneous space, hence admits a natural complex structure. 
Usually we call $(M,\frak P)$ a cuspidal datum of $G$ if $\pi_0$ is cuspidal. 
Also for $\pi\in\frak P$ set $\mathrm{Re}\pi:=
\mathrm{Re}\chi_\pi=|\chi_\pi|\in 
\mathrm{Re}X_M$, where $\chi_\pi$ is the central character of $\pi$, and 
$\mathrm{Im}\pi:=\pi\otimes(-\mathrm{Re}\pi)$.

For $\phi\in A(N({\mathbb A})M(F)\backslash G({\mathbb A}))_\pi$ with 
$\pi$ an irreducible automorphic representation of $M({\mathbb A})$, 
define the 
associated {\it Eisenstein series} $E(\phi,\pi):G(F)\backslash G({\mathbb A})\to 
{\mathbb C}$ by $$E(\phi,\pi)(g):=\sum_{\delta\in P(F)\backslash G(F)}\phi
(\delta g).$$ Then there is an open cone ${\mathcal C}\subset 
\mathrm{Re}X_M^G$ such that if $\mathrm{Re}\pi\in {\mathcal C}$, 
$E(\lambda\cdot \phi,
\pi\otimes\lambda)(g)$ converges uniformly for $g$ in a compact subset of 
$G({\mathbb A})$ and $\lambda$ in an open neighborhood of 0 in $X_M^G$. For 
example, if $\frak P=[\pi]$ is cuspidal, we may even take ${\mathcal C}$ to be 
the cone $\{\lambda\in \mathrm{Re}X_M^G:\langle\lambda-\rho_P,
\alpha^\vee\rangle>0,\forall\alpha\in \Delta_P^G\}$.
As a direct consequence, then $E(\phi,\pi)\in A(G(F)\backslash 
G({\mathbb A}))$ is an automorphic form.

\subsection{Non-Abelian L-Functions}

Being automorphic forms, Eisenstein series are of moderate growth. 
Consequently, they are not integrable over $G(F)\backslash G({\mathbb A})^1$ in general.
On the other hand, Eisenstein series are also smooth and hence integrable over
compact subsets of $G(F)\backslash G({\mathbb A})^1$. So it is very natural
for us to search for compact domains which are intrinsically defined.

As such, let us now return to the group $G=GL_r$. 
Then, we obtain compact moduli spaces 
$${\mathcal M}_{F,r}^{\leq p}[\Delta_F^{r\over 2}]:=
\Big\{g\in GL_r(F)\backslash GL_r({\mathbb  A}):
\mathrm{deg}\, g=0,\bar p^g\leq p\Big\}$$ 
for a fixed convex polygon $p:[0,r]\to{\mathbb	R}$. For example, 
${\mathcal M}_{\mathbb Q,r}^{\leq 0}[1]=\mathcal M_{\mathbb Q,r}[1]$,
(the adelic inverse image of) the moduli space of rank $r$
semi-stable $\mathbb Z$-lattices of volume 1.

More generally, for the standard parabolic subgroup $P$ of $GL_r$, we 
introduce the moduli spaces $${\mathcal M}_{F,r}^{P;\leq p}
[\Delta_F^{r\over 2}]
:=\Big\{g\in P(F)\backslash GL_r({\mathbb  A}):\mathrm{deg}\, g=0, 
\bar p_P^g\leq p, \bar p_P^g\geq -p\Big\}.$$ One checks that these 
moduli spaces ${\mathcal M}_{F,r}^{P;\leq p}[\Delta_F^{r\over 2}]$ 
are all compact.

As usual, we fix the minimal parabolic subgroup 
$P_0$ corresponding to the partition $(1,\cdots,1)$ with $M_0$ 
consisting of diagonal matrices. Then	
$P=P_I=N_IM_I$ corresponds to a certain partition 
$I=(r_1,\cdots,r_{|P|})$ of $r$ with $M_I$ the standard Levi and $N_I$ the 
unipotent radical.

Now for a fixed irreducible automorphic representation $\pi$ of 
$M_I({\mathbb  A})$, choose $$\begin{aligned}\phi\in A&
(N_I({\mathbb  A})M_I(F)\backslash 
G({\mathbb  A}))_\pi\cap 
L^2(N_I({\mathbb  A})M_I(F)\backslash G({\mathbb  A}))\\
:=&A^2(N_I({\mathbb  A})M_I(F)
\backslash G({\mathbb  A}))_\pi,\end{aligned}$$ with 
$L^2(N_I({\mathbb  A})M_I(F)\backslash	
G({\mathbb  A}))$ the space of $L^2$ functions on the space 
$Z_{G({\mathbb	A})}N_I({\mathbb  A})M_I(F)\backslash G({\mathbb  A})$. 
Denote the associated Eisenstein series by 
$E(\phi,\pi)\in A(G(F)\backslash G({\mathbb  A}))$.
\vskip 0.20cm
\noindent
{\bf Definition.} (Weng) {\it The rank $r$ non-abelian $L$-function 
$L_{F,r}^{\leq p}(\phi,\pi)$
associated to the $L^2$-automorphic form 
$\phi\in A^2(N_I({\mathbb  A})M_I(F)\backslash 
G({\mathbb  A}))_\pi$ for the number field $F$
is defined by the following integration}
$$L_{F,r}^{\leq p}(\phi,\pi):=
\int_{{\mathcal M}_{F,r}^{\leq p}[\Delta_F^{r\over 2}]}E(\phi,\pi)(g)\,dg,
\qquad \text {Re}\,\pi\in {\mathcal C}.$$

More generally, for any standard parabolic subgroup $P_J=N_JM_J\supset P_I$ 
(so that the partition $J$ is a refinement of $I$), we obtain a
relative Eisenstein series $$E_I^J(\phi,\pi)(g):=\sum_{\delta\in P_I(F)
\backslash P_J(F)}\phi(\delta g),\qquad\forall g\in P_J(F)\backslash 
G({\mathbb  A}).$$ There is  an open cone ${\mathcal C}_I^J$
in $\text {Re}X_{M_I}^{P_J}$ s.t. if $\text {Re}\pi\in 
{\mathcal C}_I^J$, then 
$E_I^J(\phi,\pi)\in A(P_J(F)\backslash G({\mathbb  A})),$ 
where $X_{M_I}^{P_J}$ is defined similarly as $X_M^G$ with $G$ replaced by 
$P_J$. As such, we are able to define 
the associated non-abelian $L$-function by
$$L_{F,r}^{P_J;\leq p}(\phi,\pi):=
\int_{{\mathcal M}_{F,r}^{P_J;\leq p}[\Delta_F^{r\over 2}]}
E_I^J(\phi,\pi)(g)\,dg, \qquad \text {Re}\pi\in {\mathcal C}_I^J.$$ 

\noindent
{\it Remark.} Here when defining non-abelian $L$-functions we assume that 
$\phi$ comes from a single irreducible automorphic representations. But this 
restriction is rather artifical and can be removed easily:
such a restriction only serves the purpose of giving the 
constructions and results in a very neat form. 

\section{Basic Properties of Non-Abelian $L$-Functions}
\subsection{Meromorphic Extension and Functional Equations}

With the same notation as above, set $\frak P=[\pi]$. For $w\in W$ the Weyl 
group of $G=GL_r$, fix once and for all representative $w\in G(F)$ of $w$. Set 
$M':=wMw^{-1}$ and denote the associated parabolic subgroup by $P'=N'M'$. 
$W$ acts naturally on the automorphic representations, from which we obtain 
an equivalence
classes $w{\frak P}$ of automorphic representations of $M'({\mathbb  A})$. As 
usual, define the associated {\it intertwining operator} $M(w,\pi)$ by 
$$(M(w,\pi)\phi)(g):=\int_{N'(F)\cap wN(F)w^{-1}\backslash N'({\mathbb	A})}
\phi(w^{-1}n'g)dn',\qquad
\forall g\in G({\mathbb	 A}).$$ One checks that if $\langle \text {Re}\pi,
\alpha^\vee\rangle\gg 0,\forall \alpha\in\Delta_P^G$,

\noindent
(i) for a fixed $\phi$, $M(w,\pi)\phi$ depends only on the double coset 
$M'(F)wM(F)$. So $M(w,\pi)\phi$ is well-defined for $w\in W$;

\noindent
(ii) the above integral converges absolutely and uniformly for $g$ varying 
in a compact subset of $G({\mathbb  A})$;

\noindent
(iii) $M(w,\pi)\phi\in A(N'({\mathbb  A})M'(F)\backslash 
G({\mathbb  A}))_{w\pi}$;
and if $\phi$ is $L^2$, which from now on we always assume, so is	
$M(w,\pi)\phi$.
\vskip 0.30cm
\noindent
{\bf Basic Facts of Non-Abelian $L$-Functions.} (Langlands, Weng)

\noindent
$\bullet$ ({\bf Meromorphic Continuation}) {\it $L_{F,r}^{\leq p}(\phi,\pi)$ for 
$\mathrm{Re}\pi\in {\mathcal C}$ is well-defined and admits a unique 
meromorphic continuation to the whole space $\frak P$;}

\noindent
$\bullet$ ({\bf Functional Equation}) {\it As meromorphic functions on $\frak P$,}
$$L_{F,r}^{\leq p}(\phi,\pi)=L_{F,r}^{\leq p}(M(w,\pi)\phi,w\pi),
\qquad \forall w\in W.$$
This is a direct consequence of the fundamental results of Langlands 
on Eisenstein series and spectrum decompositions and explains why only 
$L^2$-automorphic forms are used in the definition of non-abelian $L$s. 
(See e.g, [Ar1], [La1], [MW] and/or [W2,5]).

\subsection{Holomorphicity and Singularities}

Let $\pi\in \frak P$ and $\alpha\in\Delta_M^G$. Define the function 
$h:\frak P\to {\mathbb	C}$ by $\pi\otimes\lambda\mapsto \langle \lambda,
\alpha^\vee\rangle,\forall\lambda\in X_M^G\simeq \frak a_M^G$. Here as usual, 
$\alpha^\vee$ denotes the coroot associated to $\alpha$. Set 
$H:=\{\pi'\in\frak P:h(\pi')=0\}$ and call it a root hyperplane. Clearly the 
function $h$ is determined by $H$, hence we also denote	 $h$ by $h_H$. 
Note also that root hyperplanes depend on the base point $\pi$ we choose.

Let $D$ be a set of root hyperplanes. Then 

\noindent
(i) the singularities of a meromorphic function $f$ on $\frak P$ is said to 
be supported by $D$ if for all $\pi\in\frak P$, there exist $n_\pi:D\to 
{\mathbb  Z}_{\geq 0}$ zero almost everywhere such that $\pi'\mapsto 
\big(\Pi_{H\in D}h_H(\pi')^{n_\pi(H)}\big)\cdot f(\pi')$ is holomorphic at
$\pi'$;

\noindent
(ii) the singularities of $f$ are said to be without multiplicity at $\pi$ if 
$n_\pi\in\{0,1\}$;

\noindent
(iii) $D$ is said to be locally finite,	 if for any compact subset 
$C\subset\frak P$, $\{H\in D:H\cap C\not=\emptyset\}$ is finite.
\vskip 0.30cm
\noindent
{\bf Basic Facts of Non-Abelian $L$-Functions.} (Langlands, Weng) 

\noindent
$\bullet$ ({\bf Holomorphicity}) (i) {\it When $\mathrm{Re}\pi\in {\mathcal C}$, 
$L_{F,r}^{\leq p}(\phi,\pi)$ is holomorphic;}

\noindent
(ii) {\it $L_{F,r}^{\leq p}(\phi,\pi)$ is holomorphic at $\pi$ where 
$\mathrm{Re}\pi=0$;}

\noindent
$\bullet$ ({\bf Singularities}) {\it Assume further that $\phi$ is a cusp form. Then}

\noindent
(i) {\it There is a locally finite set of root hyperplanes 
$D$ such that the singularities of $L_{F,r}^{\leq p}(\phi,\pi)$ are supported by $D$;}

\noindent
(ii) {\it Singularities of $L_{F,r}^{\leq p}(\phi,\pi)$ are without 
multiplicities at $\pi$ if $\langle \mathrm{Re}\pi,\alpha^\vee\rangle\geq 0,
\forall \alpha\in\Delta_M^G$;}

\noindent
(iii) {\it There are only finitely many of singular hyperplanes of 
$L_{F,r}^{\leq p}(\phi,\pi)$ which intersect $\{\pi\in\frak P:
\langle\mathrm{Re}\pi,\alpha^\vee\rangle\geq 0,\forall\alpha\in\Delta_M\}$.}

\noindent
As above, this is a direct consequence of the fundamental results of 
Langlands on Eisenstein series and spectrum decompositions. (See e.g, [Ar1], 
[La1], [MW] and/or [W2,5]).
\vfill
\eject
\centerline{\Large\bf Chapter VIII. Symmetries and the Riemann Hypothesis}
\section{Abelian Parts of High Rank Zetas}
\subsection{Analytic Studies of High Rank Zetas}

Associated to a number field $F$ is the genuine high rank zeta function
$\xi_{F,r}(s)$ for every fixed $r\in\mathbb Z_{>0}$.
Being natural generalizations of (completed)
Dedekind zeta functions, these functions satisfy canonical properties for zetas as well.
Namely, they admit meromorphic continuations to the whole complex $s$-plane,
satisfy the functional equation $\xi_{F,r}(1-s)=\xi_{F,r}(s)$ and have only two singularities,
all simple poles, at $s=0,\,1$. Moreover, we expect that the Riemann Hypothesis holds for all zetas $\xi_{F,r}(s)$, namely,
all zeros of $\xi_{F,r}(s)$ lie on the central line $\mathrm {Re}(s)=\frac{1}{2}$.

Recall that $\xi_{F,r}(s)$ is defined by
\begin{equation*}
\xi_{F,r}(s):=
\Big(|\Delta_F|\Big)^{\frac{rs}{2}} \int_{\mathcal M_{F,r}}
\Big(e^{h^0(F,\Lambda)}-1\Big)
\cdot\big(e^{-s}\big)^{\mathrm{deg}(\Lambda)}\,d\mu(\Lambda),
\ \mathrm{Re}(s)>1
\end{equation*}
where $\Delta_F$ denotes the discriminant of $F$,
$\mathcal M_{F,r}$ the moduli space of semi-stable $\mathcal O_F$-lattices of rank $r$
(here $\mathcal O_F$ denotes the ring of integers),
$h^0(F,\Lambda)$ and $\mathrm{deg}(\Lambda)$ denote the 0-th geo-arithmetic cohomology
and the Arakelov degree of the lattice $\Lambda$, respectively,
and $d\mu(\Lambda)$ a certain Tamagawa type measure on $\mathcal M_{F,r}$.
Defined using high rank lattices,
these zetas then are expected to be naturally related with non-abelian aspects of number fields.

On the other hand, algebraic groups associated to $\mathcal O_F$-lattices
are general linear group $GL$ and special linear group $SL$.
A natural question then is whether principal lattices associated to other reductive groups $G$
and their associated zeta functions can be introduced and studied.

While arithmetic approach using stability seems to be complicated, analytic one using analytic truncation is ready to be exposed.
To start with, let us go back to high rank zetas.
For simplicity, take $F$ to be the field $\mathbb Q$ of rationals. Then, via
a Mellin transform, high rank zeta $\xi_{\mathbb Q,r}(s)$ can be written
as $$\xi_{\mathbb Q,r}(s)=\int_{\mathcal M_{\mathbb Q,r}[1]}
\widehat E(\Lambda,s)\,d\mu(\Lambda),\quad \mathrm{Re}(s)>1,$$ where 
$\mathcal M_{\mathbb Q,r}[1]$ denotes the moduli space of	
$\mathbb Z$-lattices of rank $r$ and volume 1 and $\widehat E(\Lambda,s)$ the completed
 Epstein zeta functions associated to $\Lambda$.
Recall that the moduli space $\mathcal M_{\mathbb Q,r}[1]$ may be viewed as a compact subset in $SL(r,\mathbb Z)\backslash SL(r,\mathbb R)/SO(r)$ and 
Epstein zeta functions may be written as the relative Eisenstein series
$E^{SL(r)/P_{r-1,1}}(\bold 1;s;g)$ associated to the constant function $\bold 1$ on 
the maximal parabolic subgroup $P_{r-1,1}$ corresponding to the partition $r=(r-1)+1$ of $SL(r)$, we have
$$\begin{aligned}\frac{2}{r}\cdot\xi_{\mathbb Q,r}(\frac{2}{r}\cdot s)=&\int_{\mathcal M_{\mathbb Q,r}[1]\subset SL(r,\mathbb Z)\backslash SL(r,\mathbb R)/SO(r)}\widehat E(\Lambda,s)\,d\mu(g)\\
=&
\int_{SL(r,\mathbb Z)\backslash SL(r,\mathbb R)/SO(r)}{\bold 1}_{\mathcal M_{\mathbb Q,r}[1]}(g)
\cdot\widehat E(\bold 1;s;g)\,d\mu(g)\end{aligned}$$ where 
${\bold 1}_{\mathcal M_{\mathbb Q,r}[1]}(g)$ denotes the characteristic function of the compact 
subset $\mathcal M_{\mathbb Q,r}[1]$.

Recall also that, in parallel, to remedy the divergence of  integration
$$\int_{SL(r,\mathbb Z)\backslash SL(r,\mathbb R)/SO(r)}\widehat E(\bold 1;s;g)\,d\mu(g),$$ in theories of automorphic forms and trace formula, Rankin, Selberg and Arthur introduced an analytic truncation
for smooth functions $\phi(g)$ over $SL(r,\mathbb Z)\backslash SL(r,\mathbb R)/SO(r)$. Simply put, Arthur's analytic truncation
is a device to get rapidly decreasing functions from slowly increasing functions by cutting off slow growth parts near all types of cusps uniformly. Being truncations near cusps, a rather large, or better, sufficiently regular, new parameter $T$ must be introduced. In particular, when applying to Eisenstein series
$\widehat E(\bold 1;s;g)$ and to $\bold 1$ on $SL(r,\mathbb R)$, we get the truncated function 
$\Lambda^T\widehat E(\bold 1;s;g)$ and $(\Lambda^T\bold 1)(g)$, respectively. Consequently,
by using basic properties on Arthur's truncation,
we obtain the following well-defined integrations
$$
\begin{aligned}\int_{SL(r,\mathbb Z)\backslash SL(r,\mathbb R)/SO(r)}&
\Lambda^T\widehat E(\bold 1;s;g)\,d\mu(g)\\
=&\int_{SL(r,\mathbb Z)\backslash SL(r,\mathbb R)/SO(r)}(\Lambda^T\bold 1)(g)\cdot 
\widehat E(\bold 1;s;g)\,d\mu(g)\\
=&\int_{\frak F(T)\subset SL(r,\mathbb Z)\backslash SL(r,\mathbb R)/SO(r)}\widehat E(\bold 1;s;g)\,d\mu(g)\end{aligned}
$$ where $\frak F(T)$ is the compact subset in (a fundamental domain for the quotient space)
$SL(r,\mathbb Z)\backslash SL(r,\mathbb R)/SO(r)$ whose characteristic function is given by $(\Lambda^T\bold 1)(g)$.

\subsection{Advanced Rankin-Selberg and Zagier Methods}
As such, we find an analytic way to understand our high rank zetas, provided that the above analytic discussion for sufficiently positive parameter $T$ can be
further strengthened so as to work for smaller $T$, in particular, for $T=0$, as well. 
In general, it is very difficult. Fortunately, as recalled in the previous two chapters, 
in the case of $SL$, 
this can be achieved based on an intrinsic
geo-arithmetic result, called the Micro-Global Bridge, an analogue of the following basic principle in Geometric Invariant Theory for unstability: A point is not stable, then there is a parabolic subgroup which destroys the stability. Consequently,
we have $$\frac{2}{r}\cdot\xi_{\mathbb Q,r}(\frac{2}{r}\cdot s)=\Big(\int_{G(\mathbb Z)\backslash G(\mathbb R)/K}\Lambda^T\widehat E(\bold 1;s;g)\,d\mu(g)\Big)\Big|_{T=0}.$$
This then leads to evaluation of the
special Eisenstein periods $$\int_{G(\mathbb Z)\backslash G(\mathbb R)/K}\Lambda^T\widehat E(\bold 1;s;g)\,d\mu(g),$$
and more generally the evaluation of {\it Eisenstein periods}
$$\int_{G(\mathbb Z)\backslash G(\mathbb R)/K}\Lambda^T E(\phi;\lambda;g)\,d\mu(g),$$ where 
$K$ a certain maximal compact subgroup of a reductive group $G$, $\phi$ is a $P$-level automorphic forms with $P$	 parabolic,  and $E(\phi;\lambda;g)$ the relative Eisenstein 
series from $P$ to $G$ associated to a $P$-level $L^2$ form $\phi$.

Unfortunately, in general, it is quite difficult to find a close formula for Eisenstein periods. 
But, when $\phi$ is cuspidal, then the corresponding Eisenstein period can be calculated, 
thanks to the work of [JLR] and [W4], 
an advanced version of Rankin-Selberg \& Zagier method.

\subsection{Discovery of Maximal Parabolics: SL, Sp and $G_2$}

Back to high rank zeta functions, the bad news is that this powerful calculation cannot be applied directly, since in the specific Eisenstein series, i.e., the classical Epstein zeta, used, the function $\bold 1$, corresponding to $\phi$ in general picture, on the maximal parabolic $P_{r-1,1}$  is only $L^2$, far from being cuspidal. To overcome this technical difficulty, 
we, partially also motivated by our earlier work on the so-called abelian part of 
high rank zeta functions ([W2,4]) and Venkov's trace formula for $SL(3)$ ([Ve]), introduce Eisenstein series $E^{G/B}(\bold 1;\lambda;g)$
associated to the constant function  $\bold 1$ on $B=P_{1,1,\dots,1}$, the Borel, 
into our study, since 

\noindent
1) being over the Borel, the constant function $\bold 1$ is cuspidal. So the associated Eisenstein period $\omega_{\mathbb Q}^{G;T}(\lambda)$ can be evaluated following [JLR]/[W4]; and

\noindent
2)  $E(\bold 1;s;g)$ used in high rank zetas can be realized as residues of 
$E^{G/B}(\bold 1;\lambda;g)$ along with suitable singular hyper-planes, a result 
essentially due to Siegel and Langlands, but carried out by Diehl ([D]). 

\noindent
In particular, for 1), we now know that $$\omega_{\mathbb Q}^{G;T}(\lambda)=\sum_{w\in W}\Bigg(\frac{e^{\langle w\lambda-\rho,T\rangle}}{\prod_{\alpha\in\Delta_0}\langle w\lambda-\rho,\alpha^\vee\rangle}\cdot\prod_{\alpha>0,w\alpha<0}\frac{\xi_{\mathbb Q}(\langle\lambda,\alpha^\vee\rangle)}
{\xi_{\mathbb Q}(\langle\lambda,\alpha^\vee\rangle+1)}\Bigg).$$ Here $W$ denotes the associated Weyl
group, $\Delta_0$ the collection of simple roots, $\rho:=\frac{1}{2}\sum_{\alpha >0}\alpha$, 
and $\alpha^\vee$ the co-root associated to $\alpha$.

With all this, it is clear that to get genuine
zetas associated to reductive groups $G$, it may be more
economical to use  the period $\omega_{\mathbb Q}^G(\lambda)$ defined by
$$\omega_{F}^G(\lambda):=\sum_{w\in W}
\Bigg(\frac{1}{\prod_{\alpha\in\Delta_0}\langle w\lambda-\rho,\alpha^\vee\rangle}
\cdot\prod_{\alpha>0,w\alpha<0}\frac{\xi_{F}(\langle\lambda,\alpha^\vee\rangle)}
{\xi_{F}(\langle\lambda,\alpha^\vee\rangle+1)}\Bigg)$$
which make sense for all reductive groups $G$ defined over $F$. Here as usual, $\xi_F(s)$ denotes the completed Dedekind zeta function
of $F$.

\bigskip
Back to the field of rationals.
The period $\omega_{\mathbb Q}^G(\lambda)$ of $G$ over $\mathbb Q$
is of $\mathrm{rank}(G)$ variables. To get a single variable zeta out from it, we need to 
take residues along with $\mathrm{rank}(G)-1$ (linearly independent) singular hyper-planes. So
proper choices for singular spaces should be made.	
This is done for $SL$ and $Sp$ in [W7], thanks to Diehl's paper ([D]). 
(In fact, Diehl dealt with $Sp$ only. But due to the fact that positive definite matrices are 
naturally associated to $\mathbb Z$-lattices and Siegel upper spaces, $SL$ can be also treated 
successfully with a bit extra care.) Simply put, for each $G=SL(r)$ (or $=Sp(2n)$), within the framework of classical Eisenstein series, 
there exists {\it only one} choice of  $\mathrm{rank}(G)-1$ singular hyper-planes 
$H_1=0,\,H_2=0,\dots,H_{\mathrm{rank}(G)-1}=0$. Moreover, after taking residues along with them, 
that is,  $$\mathrm{Res}_{H_1=0,\,H_2=0,\dots,H_{\mathrm{rank}(G)-1}=0}\,\omega_{\mathbb Q}^G(\lambda),$$ 
with suitable normalizations, we can get a new zeta $\xi_{G;\mathbb Q}(s)$ for $G$. 

At this point, the role played in new zetas $\xi_{G;\mathbb Q}(s)$ by maximal parabolic subgroups has not yet emerged. It is only after the study done for $G_2$ that we understand such a key role. Nevertheless, what we do observe from these discussions on $SL$ and $Sp$ is the follows: all singular hyper-planes are taken from only a single term appeared in the period 
$\omega_{\mathbb Q}^G(\lambda)$. More precisely, the term corresponding to $w=\mathrm{Id}$,
 the  Weyl element Identity. In other words, singular hyper-planes are taken from the denominator of 
the expression $$\frac{1}{\prod_{\alpha\in\Delta_0}\langle \lambda-\rho,\alpha^\vee\rangle}.$$ 
(Totally, there are $\mathrm{rank}(G)$ factors, among which we have carefully chosen $\mathrm{rank}(G)-1$ for $G=SL,\, Sp$.)
In particular, for the exceptional $G_2$, being a rank two group and hence an obvious  choice for our next test, this reads as
$$\frac{1}{\langle \lambda-\rho,\alpha_{\mathrm{short}}^\vee\rangle\cdot \langle \lambda-\rho,\alpha_{\mathrm{long}}^\vee\rangle}$$ where $\alpha_{\mathrm{short}}, \,\alpha_{\mathrm{long}}$ denote
the short and long roots of $G_2$ respectively. So two possibilities,

a) $\mathrm{Res}_{\langle \lambda-\rho,\alpha_{\mathrm{short}}^\vee\rangle=0}\,\omega_{\mathbb Q}^{G_2}(\lambda)$, leading to
$\xi_{\mathbb Q}^{G_2/P_{\mathrm{long}}}(s)$ after suitable normalization; and

b) $\mathrm{Res}_{\langle \lambda-\rho,\alpha_{\mathrm{long}}^\vee\rangle=0}\,\omega_{\mathbb Q}^{G_2}(\lambda)$, leading to
$\xi_{\mathbb Q}^{G_2/P_{\mathrm{short}}}(s)$ after suitable normalization.

\noindent
With this, by the fact that there exists a natural one-to-one and onto correspondence between collection of conjugation classes of 
maximal parabolic groups and simple roots, we are able to detect in [W7] the crucial role played
by maximal parabolic subgroups and hence are able to offer the proper definition for
the genuine zetas associated to pairs of reductive groups and their maximal parabolic subgroups.

\section{Abelian Zetas for $(G,P)$}
\subsection{Definition}

Motivated by the above discussion, we can introduce a genuine abelian zeta function for pairs $(G,P)$ defined over  number fields,
consisting of reductive groups $G$ and their maximal reductive groups. As the details is explained in [W7] collected in this volume, we here only sketch key features of such zetas.

Thus let $G$ be a reductive group and $P$ a maximal parabolic subgroup of $G$ both defined over $\mathbb Q$. Denote by $\Delta_0$ the collection of simple roots. For any root $\alpha$ denotes by $\alpha^\vee$ the corresponding co-root and $\rho:=\frac{1}{2}\sum_{\alpha>0}\alpha$.
Denote by $W$ the associated Weyl group. The for any $\lambda$ in 
a  suitable positive chamber of the root space, define the abelian zeta function associated to 
$(G,P)$ over $\mathbb Q$ by
$$\xi_{\mathbb Q}^{G/P}(s):=\mathrm{Norm}\Big[
\mathrm{Res}_{\langle\lambda-\rho,\alpha^\vee\rangle=0,\alpha\in\Delta_0\backslash\{\alpha_P\}}
\omega_{\mathbb Q}^G(\lambda)\Big)\Big]
$$ where as above, $$\omega_{\mathbb Q}^G(\lambda):=
\sum_{w\in W}\frac{1}{\prod_{\alpha\in\Delta_0}\langle w\lambda-\rho,\alpha^\vee\rangle}
\cdot\prod_{\alpha>0, w\alpha<0}\frac{\xi_{\mathbb Q}(\langle\lambda,\alpha^\vee\rangle)}
{\xi_{\mathbb Q}(\langle\lambda,\alpha^\vee\rangle+1)},$$  $\alpha_P$ denotes the unique simple root corresponding to the maximal parabolic subgroup $P$, 
$s:=\langle\lambda-\rho,\alpha_P^\vee\rangle$, and Norm means a certain normalization, 
the details of which may be found in [W7].

\subsection{Conjectural FE and the RH}

As such, then easily, $\xi_{\mathbb Q}^{G/P}(s)$ is a well-defined meromorphic function on the whole complex $s$-plane. And strikingly, the structures of all this zetas can be summarized by 
the following
\vskip 0.20cm
\noindent
{\bf Main Conjecture.} (i) ({\bf Functional Equation}) $$\xi_{\mathbb Q}^{G/P}(1-s)
=\xi_{\mathbb Q}^{G/P}(s);$$

\noindent
(ii) ({\bf The Riemann Hypothesis}) 
\vskip 0.20cm
\centerline{$\xi_{\mathbb Q}^{G/P}(s)=0$ implies that
$\mathrm{Re}(s)=\frac{1}{2}.$}
\vskip 0.20cm
\noindent
{\it Remarks.} (i) Funational equation is checked in [W7]
for 10 examples listed in the appendix there, namely for the groups 
$SL(2,3,4,5), Sp(4)$ and $G_2$;	 More generally, in April 2008, Henry Kim 
 in a joint effort with the author
obtained a proof of the functional equation 
for $\xi^{SL(n)/P_{n-1,1}}_{\mathbb Q}(s)$  ([KW2]);
Independently, in June, 2009, Yasushi Komori ([Ko]) 
found an elegant proof of the functional equation 
for all zetas $\xi^{G/P}_{\mathbb Q}(s)$:

\noindent
{\bf Functional Equation.}
 {\it For zeta functions $\xi_{\mathbb Q}^{G/P}(s)$, we have
 $$\xi_{\mathbb Q}^{G/P}(1-s)
=\xi_{\mathbb Q}^{G/P}(s).$$}

\noindent
(ii) Based on symmetries, the Riemann Hypothesis for the above 10 examples is solved
partially by J. Lagarias-M. Suzuki, Suzuki, and fully by H. Ki. Ki's method is expected 
to have more applications. For details, please go to ([LS], [Su1,2], [SW], [Ki1,2]).

\section{Abelian Parts of High Rank Zetas}

In a certain sense, $\xi_{\mathbb Q}^{SL(r)/P_{r-1,1}}(s)$ may be viewed as
an abelian part of the 
high rank zeta $\xi_{{{\mathbb Q}},r}(s)$, since it is naturally related to
the so-called constant terms of the Eisenstein series $E^{SL/B}({\bold 1};\lambda;g)$.
Formally, starting from Eisenstein series $E^{G/B}({\bold 1};\lambda;g)$,
we can get high rank zeta functions by first taking the residues along suitable singular
hyperplanes, then taking integration over moduli spaces of semi-stable lattices. 
That is to say,
$\xi_{{{\mathbb Q}},r}(s)$ corresponds to $({\rm Res} \to \int)$-ordered construction. In this sense,
the zeta function $\xi_{{\rm SL}(r),{{\mathbb Q}}}(s)$
corresponds to $(\int \to {\rm Res})$-ordered construction.
 
Since there is no needs to  take residues,  for $SL(2)$,
we have $\xi_{{{\mathbb Q}},2}(s)=\xi_{\mathbb Q}^{{\rm SL}(2)/P_{1,1}}(s)$.
However, in general, there is a discrepancy between $\xi_{{{\mathbb Q}},r}(s)$ and $\xi^{{\rm SL}(r)/P_{r-1,1}}_{\mathbb Q}(s)$,
because of the obstruction for the exchanging of $\int$ and ${\rm Res}$. 
\vskip 0.20cm
\noindent
{\it Remarks.} (i) Non-abelain zetas were essentially introduced around 2000. Contrary to
the publishing order, the zetas for number fields was first introduced, and it was for the purpose to get some concrete feelings that we started our examples with function fields;

\noindent
(ii) There are a few flaws in our works
on the zeta associated to SL(3) in the final chapter of [W2]. More precisely, what we have done there is the abelian zeta $\xi^{SL(3)/P_{2,1}}_{\mathbb Q}(s)$, instead of the non-abelian rank 3 zeta
$\xi_{\mathbb Q,3}(s)$; Moreover, there are sign mistakes in the formula for $\xi^{SL(3)/P_{2,1}}_{\mathbb Q}(s)$. The right one should be
\begin{equation}
\begin{aligned}\xi^{SL(3)/P_{2,1}}_{\mathbb Q}(s)=&\xi_{\mathbb Q}(2)\cdot\frac{1}{3s-3}\cdot\xi_{\mathbb Q}(3s)\\
&-\xi_{\mathbb Q}(2)\cdot\frac{1}{3s}\cdot\xi_{\mathbb Q}(3s-2)\\
&-\frac{1}{3}\cdot\frac{1}{3s-3}\cdot\xi_{\mathbb Q}(3s-1)\\
&+\frac{1}{3}\cdot\frac{1}{3s}\cdot\xi_{\mathbb Q}(3s-1)\\
&+\frac{1}{2}\cdot\frac{1}{3s-1}\cdot\xi_{\mathbb Q}(3s-2)\\
&-\frac{1}{2}\cdot\frac{1}{3s-2}\cdot\xi_{\mathbb Q}(3s)\end{aligned}
\end{equation} 

\noindent
(iii) Combinatorial techniques used by Arthur for reduction theory and 
analytic truncations are discussed in details in our 
preprint (arXiv:Math/\linebreak 0505016). 
But we remind the reader that $\tau$, 
the characteristic function in \S13.4, does not work well for analytic truncations.

\vfill
\eject
\centerline{\Large\bf Part C. General CFT and Stability}
\vskip 0.80cm
In this part, we will propose a general CFT for $p$-adic number fields using stability of what we call filtered $(\varphi,N;\omega)$-modules, built on Fontaine's theory of $p$-adic Galois representations.
The key points are
\vskip 0.20cm
\noindent
1) (Fontaine$\|$Berger) $p$-adic monodromy theorem 
for $p$-adic representations which claims that a de Rham representation is a potentially 
semi-stable representation; 

\noindent
2)  (Fontaine$\|$Fontaine, Colmez-Fontaine) characterization of semi-stable representations in terms of weakly admissible filtered $(\varphi,N)$-modules;

\noindent
3) a notion of $\omega$-structures measuring (higher) ramifications of de Rham representations;

\noindent
4) a conjectural Micro Reciprocity Law, characterizing de Rham representations in terms of semi-stable filtered $(\varphi,N;\omega)$-modules of slope zero.

\vskip 0.7cm
\centerline{\Large\bf Chapter IX. $l$-adic Representations for $p$-adic Fields}
\section{Finite Monodromy and Nilpotency}
\subsection{\bf Absolute Galois Group and Its pro-$l$ Structures}

Let $K$ be a $p$-adic number field, i.e., a finite extension of $\mathbb Q_p$.
Denote by $k$ its residue field. Fix an algebraic closure $\overline K$ of $K$. 
Let $G_K:=\mathrm{Gal}(\overline K/K)$ be the absolute Galois group of $K$ with $I_K$ its inertial 
subgroup and $P_K$ its wild ramification group.
Then from the theory of local fields, we have the following structural 
exact sequences
$$1\to I_K\to G_K\to G_k\to 1\qquad\mathrm{and}\qquad 1\to P_K\to I_K\to\prod_{l(\not=p)}\mathbb Z_l(1)\to 1.$$

 With its application to $l$-adic representation in mind, let us fix a prime $l\not=p.$ 
To avoid the pro-$l$ part systematically, define $P_{K,l}$ to be the inverse image of 
$\prod_{l'(\not=p,l)}\mathbb Z_{l'}(1)$. Accordingly,  we have an induced exact sequence 
$$1\to P_K\to P_{K,l}\to \prod_{l'\not=p,l}\mathbb Z_{l'}(1)\to 1.$$ By contrast, the pro-$l$ part can be read from the exact sequence $$1\to \mathbb Z_l(1)\to G_{K,l}\to G_k\to 1,$$
where the group $G_{K,l}$ is defined via the exact sequence $$1\to P_{K,l}\to G_K\to G_{K,l}\to 1.$$  Consequently, $g\in G_k$ acts naturally on $\gamma\in P_{K,l}$ via
$$\gamma\mapsto g\gamma g^{-1}.$$

We are ready to state one of the most intrinsic relations for Galois groups of local
fields:
\vskip 0.20cm
\noindent
{\bf Tame Relation.} (Iwasawa) {\it For any $\gamma\in \mathbb Z_l(1)$ and $\mathrm{Fr}_k\in G_k$ the absolute arithmetic Frobenius, a topological generator,  we have $$\mathrm{Fr}_k\cdot \gamma\cdot \mathrm{Fr}_k^{-1}=\gamma^q$$ where $q:=\#k$.}

\subsection{\bf Finite Monodromy}

\noindent
We say that a representation $\rho:G_K\to \mathrm{Aut}_{\mathbb Q_l}(V)$ is {\it  a $l$-adic
representation of} $G_K$ if $V/\mathbb Q_l$ is finite dimensional and $\rho$ is continuous.	
The following is the basic result on the structure of $l$-adic Galois representations:
\vskip 0.20cm
\noindent
{\bf Finite Monodromy.} (Grothendieck) {\it If $\rho:G_K\to \mathrm{Aut}_{\mathbb Q_l}(V)$ is a $l$-adic representation, then
$\rho(P_{K,l})$ is finite.}

\noindent
{\it Sketch of a proof.} Since it is a profinite group, the Galois group $G_K$ is compact. Consequently, there exists a maximal $G_K$-stable $\mathbb Z_l$-lattice $\Lambda$ in $V$ 
such that $\rho$ admits an integral form $$\rho_{\mathbb Z_l}:G_K\to\mathrm{Aut}_{\mathbb Z_l}(\Lambda).$$

As such, for any $n\in\mathbb N$, define a subgroup $N_n$  of $\mathrm{Aut}_{\mathbb Z_l}(\Lambda)$ to be the kernel of mod $l^n$ map
$$1\to N_n\to \mathrm{Aut}_{\mathbb Z_l}(\Lambda)\to
\mathrm{Aut}_{\mathbb Z_l}(\Lambda/l^n\Lambda)\to 1.$$
Clearly, $N_1/N_n$ is a finite group of order equal to a power of $l$
and hence $\displaystyle{N_1=\lim_{\leftarrow n}N_n}$ is a pro-$l$ group.
 
On the other hand, by definition, $P_{K,l}$ is a projective limit of finite groups whose orders are prime to $l$, thus $\rho_{\mathbb Z_l}(P_{K,l})\cap N_1=\{1\}$. Consequently, 
$\rho(P_{K,l})=\rho_{\mathbb Z_l}(P_{K,l})$ is naturally embedded in $\mathrm{Aut}_{\mathbb Z_l}(\Lambda/l\Lambda)$ which is a finite group.

\subsection{\bf Unipotency}
Based on finite monodromy property, we further have the following
\vskip 0.20cm
\noindent
{\bf Monodromy Theorem.} (Grothendieck)
{\it Let $\rho:G_K\to \mathrm{Aut}_{\mathbb Q_l}(V)$ be a $l$-adic representation. Then there exists a finite Galois extension $L/K$ such that
for the induced representation $\rho|_{G_L}:G_L(\subset G_K)\to 
\mathrm{Aut}_{\mathbb Q_l}(V)$, the inertial subgroup $I_L(\subset G_L)$ acts unipotently.} 

\noindent
{\it Sketch of a proof.} This is a direct consequence of the Tame relation. Indeed, by the finite monodromy result in the previous subsection,
replacing $K$ by a finite Galois extension, we may assume that $P_{K,l}$
acts on $V$ trivially. Consequently, since $G_K/ P_{K,l}=G_{K,l},$  the representation $\rho$ factors through $G_{K,l}$:
$$\rho:G_K\twoheadrightarrow G_{K,l}\buildrel\bar\rho\over\to \mathrm{Aut}_{\mathbb Q_l}(V).$$

Recall now that we have the following structural exact sequence
$$1\to\mathbb Z_l(1)\to G_{K,l}\to G_k\to 1$$
and the tame relation, recalled above, implies that 
for any $t\in\mathbb Z_l(1), n\in\mathbb N,$
$$\mathrm{Fr}_k^n \cdot t\cdot \mathrm{Fr}_k^{-n}=t^{nq},$$
with $\mathrm{Fr}_k$ the absolute Frobenius of $k$ and $q=\#k$.
Consequently, if $\lambda$ is an eigenvalues of $\bar \rho(t)=\rho(t)$, then so is $\lambda^n$. This implies that all such $\lambda$'s are roots of unity.
Namely, all elements of $\mathbb Z_l(1)\subset G_{K,l}$ act unipotently.
But $\mathbb Z_l(1)$ is rank one, so if we choose $t_0$ as a topological generator, then the topological closure $\overline{\langle t_0\rangle}$ of the subgroup generated by $t_0$
acts unipotently on $V$. Since $\overline{\langle t_0\rangle}$ is clearly
 an open subgroup of $\mathbb Z_l(1)$, so
the whole $\mathbb Z_l(1)$ acts on $V$ unipotently.
With this, to complete the proof, it suffices to note that the induced action of inertia subgroup $I_K$ factors through $\mathbb Z_l(1)$. From the exact sequences
$$0\to P_K\to P_{K,l}\to \prod_{l'\not=p,l}\mathbb Z_{l'}(1)\to 0\ \mathrm{and}\ 0\to P_{K,l}\to G_K\to G_{K,l}\to 0,$$
we conclude that the induced action on $I_K$ factors through $\mathbb Z_l(1)$ via the natural projection map $$I_K\twoheadrightarrow I_K/P_K\simeq\mathbb Z_l(1)\times\prod_{l'\not=p,l}\mathbb Z_{l'}(1)\twoheadrightarrow\mathbb Z_l(1),$$
and hence is unipotent.
\vskip 0.20cm
\noindent
{\bf Example.} If $V/\mathbb Q_l$ is  one dimensional, from the Monodromy Theorem above, there exists a finite Galois extension $L/K$ such that the induced action of $I_L$ on $V$ is unipotent.
That means that the image of $I_L$ is a finite group. As such, replacing $L$ with a further extension, we may assume that $I_L$ acts trivially on $V$. Particularly, this works for the Tate module $\mathbb Z_l(1)$.
\vskip 0.20cm
\noindent
{\bf Definition.} Let $\rho:G_K\to\mathrm{Aut}(V)$ be a $l$-adic representation. Then $\rho$ is called

\noindent
{\bf 1.}a) {\it unramified} if $I_K$ acts on $V$ trivially;

\noindent
{\bf 1.}b) {\it potentially unramified} if there exists a finite Galois extension $L/K\subset \overline K/K$ such that the induced action of $I_L$ on $V$ is trivial;
\vskip 0.20cm
\noindent
{\bf 2.}a) {\it semi-stable}   if $I_K$ acts on $V$ unipotently;

\noindent
{\bf 2.}b) {\it potentially semi-stable} if there exists a finite Galois extension 
$L/K\subset \overline K/K$ such that the induced action of $I_L$ on $V$ is unipotent.
\vskip 0.20cm
In terms of this language, then Grothendieck's Monodromy Theorem	
claims that all $l$-adic Galois representation of a $p$-adic number field, $l\not=p$, is potentially semi-stable.

\vfill
\eject
\centerline{\Large\bf Chapter X. Primary Theory of $p$-adic Representations}
\vskip 0.80cm
In this chapter, we expose some elementary structures of $p$-adic Galois representations following [FO].
\section{Preliminary Structures of Absolute Galois Groups}
\subsection{Galois Theory: A $p$-adic Consideration}
Let $K$ be a $p$-adic number field with $k$ its residue field. Fix an
algebraic closure $\overline K$. $\overline K$ is not complete with respect to the 
natural extension of the $p$-adic valuation of $K$. Denote the corresponding completion of 
$\overline K$ by $\mathbb C=\mathbb C_p$.

Denote by $G_K:=\mathrm{Gal}(\overline K/K)$ the absolute Galois group of $K$. Then, from $p$-adic theory point of view,	$G_K$  can be naturally decomposed into two parts, namely arithmetic one corresponding to the cyclotomic extensions by $p^n$-th roots of unity, and the geometric one, corresponding to the so-called field of norms. 

More precisely, let $K_n:=K(\mu_{p^n})$ where $\mu_{p^n}$ denotes the collection of 
$p^n$-th roots of unity in $\overline K$ and set $K_\infty:=\cup_n K_n.$ Denote 
the corresponding Galois groups by $H_K:=\mathrm{Gal}(\overline K/K_\infty)$ and
$\Gamma_K:=\mathrm{Gal}(K_\infty/K).$ Clearly, $G_K/H_K\simeq\Gamma_K.$

\subsection{Arithmetic Structure: Cyclotomic Character}

Denote by $K_0:=\mathrm{Fr}\,W(k)$ the fractional field of the ring of Witt vectors 
with coefficients in $k$. Then it is known that $K_0$ is the maximal unramified extension of $\mathbb Q_p$ contained in $K$ and $\Gamma_{K_0}$ is canonically isomorphic to $\mathbb Z_p^*$ via the cyclotomic character $\chi_{\mathrm{cyc},p}=\chi_{\mathrm{cyc}}.$ Clearly, $\Gamma_K$ may be viewed as an open subgroup of $\Gamma_{K_0}$ via $\chi_{\mathrm{cyc}}$.

The natural exponential map gives a $\mathbb Z_p$-module structure on $\mathbb Z_p^*$. One 
can easily checks that it is of rank one and its torsion part is given by $$(\mathbb Z_p^*)_{\mathrm{tor}}
=\begin{cases} \mathbb F_p^*,&	p\not=2\\
\mathbb Z/2\mathbb Z, & p=2.\end{cases}$$
Consequently, if we denote by $\Delta_K$ the torsion subgroup of $\Gamma_K$, then $K_\infty^{\Delta_K}=(K_{0,\infty})^{\Delta_{K_0}}\cdot K/K$ is a $\mathbb Z_p$-extension with the same residue field $k$ of $K$.

For later use, denote by $k'$ the residue field of $K_\infty$.
From the discussion above, we see that it may happen that
$k'$ is different from	$k$.

\subsection{Geometric Structure: Fields of Norms}
\subsubsection{Definition}
With $\Gamma_K$ understood, let us turn our attention to $H_K$ part. This then leads to
the theory of fields of norms due to Witenberger. Roughly speaking, this theory
says that the arithmetically defined Galois group $H_K:=\mathrm{Gal}(\overline K/K_\infty)$ 
of fields of characteristic zero admits a natural geometric interpretation in terms of Galois 
group of localizations of function fields over finite fields, due to the fact that the natural 
norm map $N_{K_n/K_{n-1}}$ is quite related with the $p$-th power map.

More precisely, motivated by a work of Tate, for  fields $K_n$, consider  norm maps 
$N_{K_n/K_{n-1}}$. Clearly, $\big\{(K_n,N_{K_n/K_{n-1}})\big\}_{n\in \mathbb N}$ forms a projective system. Let $\displaystyle{\mathcal N_K:=\lim_{\leftarrow_n}K_n}$ be the corresponding limit. That is,

\noindent
(i) as a set, $$\mathcal N_K=\Big\{(x^{(0)},x^{(1)}, \dots, x^{(n)},\dots):x^{(n)}\in K_n,\ N_{K_n/K_{n-1}} \big(x^{(n)}\big)= x^{(n-1)}\Big\};$$
(ii) for the ring structure, the addition and multiplication on $\mathcal N_K$ are given by
$$\begin{aligned}\big(x+y\big)^{(n)}:=&\lim_{m\to\infty} N_{K_{n+m}/K_{n}} \Big(x^{(n+m)}+ y^{(n+m)}\Big)\\ 
\big(x\cdot y\big)^{(n)}:=& x^{(n)}\cdot y^{(n)}\end{aligned}$$ for $x=\big(x^{(n)}\big),\ y=\big(x^{(n)}\big)\in\mathcal N_K.$
\vskip 0.20cm
Much more holds:

\noindent
{\bf Theorem.} (Wintenberger) {\it $\mathcal N_K$ is a field, the so-called field of norms of $K_\infty/K$, such that its separable closure $\mathcal N_K^s$ is given by
$$\bigcup_{L/K:\mathrm{finite\,Galois}}\mathcal N_L,$$ and $G_{\mathcal N_K}:=\mathrm{Gal}(\mathcal N_K^s/\mathcal N_K)$ is isomorphic to $H_K$. 

In particular,

\noindent
(i) for every finite Galois extension $L/K$ in $\overline K/K$, $\mathcal N_L/\mathcal N_K$ is a finite Galois extension with $$\mathrm{Gal}\Big(\mathcal N_L/\mathcal N_K\Big)\simeq\mathrm{Gal}\Big(L_\infty/K_\infty\Big);$$

\noindent
(ii) for every finite Galois extension $\mathcal N_{*}/\mathcal N_K$, there
exists	a finite Galois extension $L/K$ such that $\mathcal N_L=\mathcal N_{*}.$}
\vskip 0.20cm
\subsubsection{Geometric Interpretation}

To  give a geometric interpretation of $\mathcal N_K$, let us start with $\mathcal N_{K_0}.$ 
If we set $E_{K_0}:=k\big((\pi_{K_0})\big)$
for a certain indeterminant $\pi_{K_0}$ over $k$, then $$\mathcal N_{K_0}\simeq E_{K_0}=k\big((\pi_{K_0})\big).$$
And more generally, for a certain indeterminant $\pi_K$ over $k'$, $$\mathcal N_{K}\simeq E_{K}=k'((\pi_K)).$$

To be more precise, this is realized via the following consideration.
First, by ramification theory, we see that the norm map $N_{K_n/K_{n-1}}$ is not far away from being the $p$-th
power map. Accordingly, it is natural to introduce the ring
$$\widetilde{\mathbb E^+}:=\lim_{\begin{subarray}{1}\ 
\longleftarrow\\
{x\mapsto x^p}\end{subarray}}\mathcal O_{\mathbb C}:=
\Big\{\big(x^{(0)},x^{(1)},\dots \big): x^{(n)}\in
\mathcal O_{\mathbb C},\ 
\Big(x^{(n+1)}\Big)^p=x^{(n)}\Big\}$$
where $\mathcal O_{\mathbb C}$ denotes the ring of integers of $\mathbb C$. Define the ring structure on $\widetilde{\mathbb E^+}$ by
$$\big(x+y\big)^{(n)}:=\lim_{m\to\infty}  \Big(x^{(n+m)}+ y^{(n+m)}\Big)^{p^m}\quad \&\quad 
\big(x\cdot y\big)^{(n)}:= x^{(n)}\cdot y^{(n)}$$ for $x=\big(x^{(n)}\big),\ y=\big(x^{(n)}\big)\in\widetilde{\mathbb E^+}.$

One can easily check that $\widetilde{\mathbb E^+}$ is {\it perfect}. It is also of characteristic $p$. Indeed, there is a bijection
$$\lim_{\begin{subarray}{1}\ \longleftarrow\\
{x\mapsto x^p}\end{subarray}}\mathcal O_{\mathbb C}\simeq \lim_{\begin{subarray}{1}\ \longleftarrow\\
{x\mapsto x^p}\end{subarray}}\mathcal O_{\mathbb C}/p\mathcal O_{\mathbb C}.$$ This implies that
$$\widetilde{\mathbb E^+}\simeq \lim_{\begin{subarray}{1}\ \longleftarrow\\
{x\mapsto x^p}\end{subarray}}
\mathcal O_{\overline K}/p\mathcal O_{\overline K},$$ since $\mathcal O_{\mathbb C}/p\mathcal O_{\mathbb C}\simeq \mathcal O_{\overline K}/p\mathcal O_{\overline K}$.

Moreover, if we set $\varepsilon=(\varepsilon^{(n)})\in \widetilde{\mathbb E^+}$ with $\varepsilon^{(0)}=1, \varepsilon^{(1)}\not=1$ defined by primitive $p^n$-th roots of unity, and set
$$\widetilde{\mathbb E}=\widetilde{\mathbb E^+}
\big[(\varepsilon-1)^{-1}\big].$$ Then this is the completion of the algebraic (yet non-separable) closure of $\mathbb F_p\big((\varepsilon-1)\big)$. 

By definition, there is a natural action of $H_K$ on $\widetilde{\mathbb E}$.
With the interpretation of $\widetilde{\mathbb E^+}\simeq \lim_{\leftarrow_{x\mapsto x^p}}
\mathcal O_{\overline K}/p\mathcal O_{\overline K}$ in terms of $\mathcal O_{\overline K}$ (not the one from the definition in terms of the completion
$\mathcal O_{\mathbb C}$), this action can be read clearly as follows:

We have a natural injective morphism $$\begin{matrix}\mathcal N_K&\to&\widetilde{\mathbb E}&\\
\big(x^{(n)}\big)&\mapsto&\Big(y^{(n)}&:=\lim_{m\to\infty}(x^{(n+m)})^{p^m}\Big)\end{matrix}$$ 
Moreover, one checks that

\noindent
(i) $\displaystyle{\mathcal N_{K_0}\simeq k\big((\pi)\big)}$ with $\pi=\varepsilon-1$;

\noindent
(ii)
$\mathbb E_K=\Big(\widetilde{\mathbb E}\Big)^{H_K}$	
coincides with the image of $\mathcal N_K$;

\noindent
(iii) $H_{L/K}:=H_K/H_L=\mathrm{Gal}\Big(L_\infty/K_\infty\Big)\simeq 
\mathrm{Gal}\Big(\mathcal N_L/\mathcal N_K\Big)\simeq \mathrm{Gal}\Big(\mathbb E_L/\mathbb E_K\Big)$. 

\section{Galois Representations: Characteristic $p$-theory}

In this section we concentrate on  Galois representations of fields of characteristic 
$p$, motivated by the geometric interpretation of $H_K$.

\subsection{$\mathbb F_p$-Representations}

Assume that $E$ is a field of characteristic $p>0$. Fix a separable
closure $E^s$ and let $G_E:=\mathrm{Gal}(E^s/E)$ be the 
corresponding absolute Galois group. Denote by $\sigma:\lambda
\mapsto\lambda^p$  the absolute Frobenius of $E$. Let $V$ be a 
mod $p$ representation of $G_E$ of dimension $d$, i.e., a $\mathbb F_p$-vector space $V$
of  dimension $d$
 equipped with a linear and continuous 
action of $G_E$. 

Since $G_E$ acts naturally on $E^s$, it makes sense to talk about the 
$E^s$-representation $E^s\otimes_{\mathbb F_p}V$ equipped with $G_E$. The 
advantage of taking this extension of scalars is that, by Hilbert Theorem 
90, one checks that if we set 
$\mathbb D(V):=\Big(E^s\otimes_{\mathbb F_p}V\Big)^{G_E}$, then 

(i) $\mathbb D(V)$ is a	 $E$-vector space of dimension $d$; and

(ii) the natural map $$\alpha_V:E^s\otimes_E\mathbb D(V)\to E^s\otimes_{\mathbb F_p}V$$ is an isomorphim of $G_E$-modules. Here, as usual, on the left hand side, the action concentrates on the coefficients $E^s$, while on the right, it is given by the diagonal action.

Moreover, since the absolute Frobenius $\sigma$ commutes with the action of $G_E$, via the natural definition $\varphi:\lambda\otimes v\mapsto\sigma(\lambda)\otimes v$, we obtain a Frobenius on $E^s\otimes_{\mathbb F_p}V$ such that if $x\in \mathbb D(V)$ then so is $\varphi(x)$. Consequently, we obtain a natural Frobenius $\varphi:\mathbb D(V)\to\mathbb D(V)$.

\subsection{Etale $\varphi$-modules}

Motivated by the above discussion, we call a finite dimensional $E$-vector space $M$ equipped with a $\sigma$-semi-linear map $\varphi:M\to M$ a
$\varphi$-{\it module} over $E$.

We call a $\varphi$-module {\it etale} if  $M=E\cdot\varphi(M)$.
\vskip 0.20cm
\noindent
{\bf Proposition.} (See e.g.,[FO]) {\it If $V$ is a 
$\mathbb F_p$-representation of $G_E$ of dimension $d$, then 
$\mathbb D(V):=\Big(E^s\otimes V\Big)^{G_E}$ is an etale 
$\varphi$-module of dimension $d$ over $E$. Moreover, as $G_E$-modules, 
we have an isomorphism} $$\alpha_V:E^s\otimes_E\mathbb D(V)
\to E^s\otimes_{\mathbb F_p}V.$$

\subsection{Characteristic $p$ Representation and Etale $\varphi$-Module}

Denote by $\mathrm{\bf Rep}_{\mathbb F_p}(G_E)$ the category of all mod $p$ representations of 
$G_E$ and $\mathcal M_{\varphi}^{\mathrm{et}}(E)$ the category of etale $\varphi$-modules over 
$E$ with morphisms being $E$-linear maps which commute with $\varphi$. Then from the paragraph
 above
we have a natural functor $$\mathbb D_E:\mathrm{\bf Rep}_{\mathbb F_p}(G_E)\to \mathcal M_{\varphi}^{\mathrm{et}}(E).$$

\noindent
{\bf Proposition.} (Fontaine) {\it The
natural functor $$\begin{matrix}\mathbb D_E:&\mathrm{\bf Rep}_{\mathbb F_p}(G_E)&\to& \mathcal M_{\varphi}^{\mathrm{et}}(E)&\\
&V&\mapsto&\mathbb D_E(V)&:=\Big(E^s\otimes_{\mathbb F_p}V\Big)^{G_E}\end{matrix}$$
gives an equivalence of categories and its quasi-inverse is given by}
$$\begin{matrix}\mathbb V_E:&\mathcal M_{\varphi}^{\mathrm{et}}(E)&\to& \mathrm{\bf Rep}_{\mathbb F_p}(G_E)&\\
&M&\mapsto&\mathbb V_E(M)&:=\Big(E^s\otimes_EM\Big)^{\varphi=1}.\end{matrix}$$

\section{Lifting to Characteristic Zero}

As our final aim is to study $p$-adic representations of Galois groups 
of local fields, it is natural to see how the discussions above
on $\mathbb F_p$-representations, a characteristic $p$-theory, 
can be lifted to $p$-adic representations, a characteristic zero theory.
We present the ralated materials following [FO] (and [Ber2]).

\subsection{\bf Witt Vectors and Teichm\"uller Lift} 

Let us start with a preparation on the coefficients, particularly, the theory of Witt vectors.

So let $R$ be a perfect ring of characteristic $p$. We want to	
construct a ring $W(R)$, the so-called {\it ring of Witt vectors with coefficients in} $R$, such that $p$ is not nilpotent 
and $W(R)$ is separated and complete for the topology defined by $p^nW(R)$. The main result on Witt rings
is  that {\it such a ring $W(R)$ does exists, unique up to isomorphism,
and has $R$ as its residual ring}. Consequently, if $\sigma:R\to S$ is a morphism, then $\sigma$ 
lifts to a morphism $W(\sigma)=:\sigma:W(R)\to W(S)$. Particularly, all Witt ring admits a lift 
of Frobenius $\sigma$!
\vskip 0.20cm
\noindent
{\bf Examples:}

\noindent 
(i) $W(\mathbb F_p)=\mathbb Z_p$; 

\noindent
(ii) If $k$ is a finite field, then $W(k)$ is the ring of integers of the unique unramified extension of $\mathbb Q_p$ whose residue field is $k$. Consequently,

\noindent
(iii) $W(\overline{\mathbb F_p})=\mathcal O_{\widehat{\mathbb Q_p^{\mathrm{un}}}}$ is the ring of integers of the $p$-adic completion of the maximal unramified extension $\mathbb Q_p^{\mathrm{un}}$ of $\mathbb Q_p$.
\vskip 0.20cm
For $x=x_0\in R$, since $R$ is perfect, it makes sense to talk about $x^{p^{-n}}$ in $R$ for all $n$. (This is in fact the key condition for a field to be perfect.) Up to $W(R)$, choose then an element $\widetilde{x_n}\in W(R)$ such that its residue class coincides with $x^{p^{-n}}$. Then the sequence $\{\widetilde{x_n}\}_{n\geq 0}$ converges in $W(R)$, say, to an element $[x]$. This $[x]$ is known to depend only on $x$, not on the choices of 
$\widetilde{x_n}$. As such, we obtain a multiplicative map, the so-called {\it Teichm\"uller lift}:
$$\begin{matrix}[\cdot]:&R&\to& W(R)\\
&x&\mapsto&[x].\end{matrix}$$ Clearly, 

\noindent
(i) the Teichm\"uller lift is a special section to the natural reduction map;

\noindent
(ii) every element $x\in W(R)$ can be written uniquely as $x=\sum_{n=0}^\infty p^n[x_n]$ wth $x_n\in R$.
Moreover, 

\noindent
(iii)  there exist universal homogeneous polynomials 
\newline $S_n,P_n\in\mathbb Z[X_i^{p^{-n}},Y_i^{p^{-n}}:i=0,1,\dots,n]$ of degree 1 (where $\mathrm{deg}X_i:=1=:\mathrm{deg}Y_i$)
such that for all $x,y\in W(R)$, we have $$\begin{aligned}x+y=&\sum_{n=0}^{\infty}p^n\Big[S_n(x_0,y_0,\dots,x_n,y_n)\Big]\\ 
 xy=&\sum_{n=0}^{\infty}p^n\Big[P_n(x_0,y_0,\dots,x_n,y_n)\Big].\end{aligned}\eqno(*)$$ For instance, $$\begin{aligned}
S_0(X_0,Y_0):=&X_0+Y_0;\\
S_1(X_0,Y_0,X_1,Y_1):=&X_1+Y_1+p^{-1}\Big((X_0^{1/p}+Y_0^{1/p})^p-X_0-Y_0\Big)\end{aligned}$$
Indeed, with the help of the polynomials $S$ and $P$, we can construct $W(R)$ by setting 

\noindent
(a) as a set, $W(R):=\prod_{n=0}^\infty R$, and 

\noindent
(b) for the ring structure, set the addition and the multiplication 
according to the above relations $(*)$.
\vskip 0.20cm
Furthermore, the concept of Witt ring can be extended to the case when $R$ is {\it not} perfect. In this later case, we call the result ring a {\it Cohen ring} $C(R)$. Cohen rings are not really unique, but still  they are  of characteristic zero with residual ring $C(R)/pC(R)=R.$ For example, $C\Big(\mathbb F_p\big[[X]\big]\Big)=\mathbb Z_p\big[[X]\big]$.

\subsection{\bf $p$-adic  Representations of Fields of Characteristic 0}

\subsubsection{Lift of base fields}

Let $\mathbb E_K\subset \widetilde {\mathbb E}$ be the field isomorphic to the field of norms $\mathcal N_K$ introduced before. 
It is of characteristic $p$ and may not be perfect. Denote its associated Cohen ring $\mathcal C(\mathbb E_K)$ by $\mathcal O_{\mathcal E_K}$ and write $\mathcal E_K$ the associated fraction field which is of characteristic 0. Denote by $\varphi:\mathcal E_K\to \mathcal E_K$ a lift of the Frobenius
$\sigma:\mathbb E_K\to \mathbb E_K$. 
Consequently, $$\mathcal O_{\mathcal E_K}=\lim_{\leftarrow_n}\mathcal O_{\mathcal E_K}/p^n\mathcal O_{\mathcal E_K},\ \ \mathcal O_{\mathcal E_K}/p\mathcal O_{\mathcal E_K}=\mathbb E_K\ \ \mathrm{and}\ \ {\mathcal E_K}=\mathcal O_{\mathcal E_K}[\frac{1}{p}].$$

Let $\mathcal F$ be a finite extension of $\mathcal E_K$ and $\mathcal O_{\mathcal F}$ be the ring of integers. We say that $\mathcal F/\mathcal E_K$ is {\it unramified} if 

\noindent
(i) $p$ is a generator of the maximal ideal of $\mathcal O_{\mathcal F}$; and

\noindent
(ii) $F=\mathcal O_{\mathcal F}/p\mathcal O_{\mathcal F}$ is a separable extension of $\mathbb E_K$.

For any finite separable extension $F$ of $\mathbb E_K$, the inclusion $\mathbb E_K\hookrightarrow F$ induces a local homomorphism $\mathcal C(\mathbb E_K)\to\mathcal C(F)$ through which we may identify $\mathcal C(\mathbb E_K)$ and a subring of $\mathcal C(F)$ and $\mathrm{Fr}\,\mathcal C(F)$ as a field extension of $\mathrm{Fr}\,\mathcal C(\mathbb E_K)$, which in particular is unramified. 
Much more is correct: By the field of norms, all finite unramified extensions of $\mathcal E_K$ 
are obtained in this way.  If we let 
$\displaystyle{\mathcal E^{\mathrm{ur}}:=\lim_{\to F\in S}\mathcal E_F}$ and let
  $\widehat{\mathcal E^{\mathrm{ur}}}$ be the $p$-adic completion 
 of $\mathcal E^{\mathrm{ur}}$ with $\mathcal O_{\widehat{\mathcal E^{\mathrm{ur}}}}$  its ring of integers, then	 
 $\mathcal O_{\widehat{\mathcal E^{\mathrm{ur}}}}$ is a local ring 
 and $$\mathcal O_{\widehat{\mathcal E^{\mathrm{ur}}}}=
 \lim_{\longleftarrow}\mathcal O_{\mathcal E^{\mathrm{ur}}}/
 p^n\mathcal O_{\mathcal E^{\mathrm{ur}}}.$$
Clearly, all are equipped with Frobenious $\varphi$ which commute with the natural action of $H_K$. Moreover, one checks directly the following holds:

\noindent
(i) $\Big(\widehat{\mathcal E^{\mathrm{ur}}}\Big)^{H_K}=\mathcal E_K,\ \Big(\mathcal O_{\widehat{\mathcal E^{\mathrm{ur}}}}\Big)^{H_K}=\mathcal O_{\mathcal E_K};$

\noindent
(ii) $\Big(\widehat{\mathcal E^{\mathrm{ur}}}\Big)^{\varphi=1}=\mathbb Q_p,\ \Big(\mathcal O_{\widehat{\mathcal E^{\mathrm{ur}}}}\Big)^{\varphi=1}=\mathbb Z_p$.

\vskip 0.20cm
\subsubsection{$p$-adic Representations}
For simplicity, write $\mathcal E$ for $\mathcal E_K$.
We say that a $\varphi$-{\it module} $M$ over $\mathcal E$ is a finite dimensional $\mathcal E$-vector space equipped with a $\sigma$-semi-linear morphism $\varphi:M\to M$; and a $\varphi$-module is called {\it etale} if $M=\mathcal E\cdot\varphi(M)$. 
One can easily check that {\it for a $p$-adic representation $V$ of $H_K$, 
$$\mathbb D(V):=\Big(\widehat{\mathcal E^{\mathrm{ur}}}\otimes_{\mathbb Q_p}V\Big)^{H_K}$$ is an 
etale $\varphi$-module over $\mathcal E$ such that
the natural map $$\widehat{\mathcal E^{\mathrm{ur}}}
\otimes_{\mathcal E}\mathbb D(V)\to 
\widehat{\mathcal E^{\mathrm{ur}}}\otimes_{\mathbb Q_p}V$$ 
is a $H_K$-equivariant isomorphism}. 

\subsection{$p$-adic Representations and Etale $(\varphi,\Gamma)$-Modules}

Let $V$ be a $\mathbb Q_p$-representation of $G_K$, set $$\mathbb D(V):=\Big(\widehat{\mathcal E^{\mathrm{ur}}}\otimes_{\mathbb Q_p}V\Big)^{H_K},$$
then  $\mathbb D(V)$ admits  natural $\Gamma_K$-actions. 
We say that $D$ is a $(\varphi,\Gamma)$-{\it module} over $\mathcal O_{\mathcal E}$ (resp. over 
$\mathcal E$) if it a $\varphi$-module over 
$\mathcal O_{\mathcal E}$ (resp. over $\mathcal E$) together with a $\sigma$-semi-linear action of $\Gamma_K$ commuting with $\varphi$. Moreover, $D$ is called {\it etale} if it is an etale $\varphi$-module and the action of $\Gamma_K$ is continuous.

Denote by	
$\mathrm{\bf Rep}_{\mathbb Q_p}(G_K)$ the category of  $p$-adic representations of $G_K$ and	
$\mathcal M_{\varphi,\Gamma}^{\mathrm{et}}({\mathcal E})$ the category of etale $(\varphi,\Gamma)$-modules over  $\mathcal E$. Then we have the following

\noindent
{\bf Corollary.} (Fontaine) {\it The
natural functor $$\begin{matrix}\mathbb D:&\mathrm{\bf Rep}_{\mathbb Q_p}(G_K)&\to& \mathcal M_{\varphi,\Gamma}^{\mathrm{et}}(\mathcal E)&\\
&V&\mapsto&\mathbb D(V)&:=\Big(\widehat{\mathcal E^{\mathrm{ur}}}\otimes_{\mathbb Q_p}V\Big)^{H_K}\end{matrix}$$
gives an equivalence of categories and its quasi-inverse is given by
$$\begin{matrix}\mathbb V:&\mathcal M_{\varphi,\Gamma}^{\mathrm{et}}(\mathcal E)&\to& \mathrm{\bf Rep}_{\mathbb Q_p}(G_K)&\\
&M&\mapsto&\mathbb V(M)&:=\Big(\widehat{\mathcal E^{\mathrm{ur}}}\otimes_{\mathcal E}M\Big)^{\varphi=1}.\end{matrix}$$}
\vfill
\eject
\centerline{\Large\bf Chapter XI. $p$-adic Hodge and Properties of Periods}
\vskip 0.80cm
To expose basic structures of $p$-adic Galois representations, 
we shift our attentions to the so-called $p$-adic
Hodge theory, based on the following reason:
etale cohomology not only  offers natural examples of Galois representations, 
but provides all the fine structures which play key roles in the theory of $p$-adic Galois representations.

\section{Hodge Theory over $\mathbb C$}

Let $X$ be a projective smooth variety over a field $E$ of characteristic zero. Then we have 
the associated complex of sheaf of differential forms $$\Omega_{X/E}^*:\mathcal O_{X/E}\to \Omega_{X/E}^1\to \Omega_{X/E}^2\to\cdots.$$
By definition, the de Rham cohomology groups $H_{\mathrm{dR}}^m(X/E)$ are the hyper-cohomology groups $\mathbb H^m(\Omega_{X/E}^*)$ for all $m$.

On the other hand, for any embedding $E\hookrightarrow \mathbb C$,
since $X(\mathbb C)$ is a compact complex manifold, the singular cohomology $H^m(X(\mathbb C),\mathbb Q)$, being the dual of $H_m(X(\mathbb C))$, is a finite dimensional $\mathbb Q$-vector space. The comparison theorem in the classical Hodge theory then says that there exists a canonical isomorphism
$$\mathbb C\otimes_{\mathbb Q}H^m(X(\mathbb C),\mathbb Q)\simeq\mathbb C\otimes_EH_{\mathrm{dR}}^m(X/E).$$
Thus without loss of generality, we may assume that $E$ is simply $\mathbb C$.

For a complex smooth projective variety $X$, denote by $A^n(X)$, resp. by $A^{p,q}(X)$, the space of $C^\infty$ $n$-forms, resp. $C^\infty$ $(p,q)$-forms. Clearly, $A^n(X)=\bigoplus_{p+q=n}A^{p,q}(X)$. With respect to the total differential operator $d:A^n(X)\to A^{n+1}(X)$, we have the cohomology groups
$$H^{p,q}(X):=\Big\{\phi\in A^{p,q}(X):d\phi=0\Big\}\Big/dA^{n-1}(X)\cap A^{p,q}(X).$$
Then the Hodge decomposition theorem in the classical Hodge theory claims that there exists a canonical isomorphism
$$H^n_{\mathrm{dR}}(X,\mathbb C)=\bigoplus_{p+q=n}H^{p,q}(X).$$ Furthermore, there is a decreasing filtration on $A^n(X)$ defined by $$\mathrm{Fil}^pA^n(X):=A^{n,0}(X)\oplus
A^{n-1,1}(X)\oplus\cdots\oplus A^{p,n-p}(X)$$ and the induced decreasing filtration of $H_{\mathrm{dR}}^n(X)$ defined by
$$\mathrm{Fil}^pH_{\mathrm{dR}}^n(X):=H^{n,0}(X)\oplus
H^{n-1,1}(X)\oplus\cdots\oplus H^{p,n-p}(X).$$
Clearly, $$\begin{aligned}\mathrm{Fil}^pH_{\mathrm{dR}}^n(X)=&\Big\{\phi\in \mathrm{Fil}^pA^n(X):d\phi=0\Big\}\Big/dA^{n-1}(X)\cap \mathrm{Fil}^pA^n(X),\\
H^{p,q}(X)=&\overline{H^{q,p}(X)},\\
H^{p,q}(X)=&\mathrm{Fil}^pH_{\mathrm{dR}}^n(X)\cap\overline{\mathrm{Fil}^qH_{\mathrm{dR}}^n(X)}.\end{aligned}$$

\section{Admissible Galois Representations}

Before we go to the essentials of $p$-adic Hodge theory, let us make a further preparation. 

Let $G$ be a  topological group and $B$ a topological commutative ring equipped with a 
continuous $G$ action. Then by a $B$-{\it representation} $V$ of $G$, we mean a free $B$-module
 $V$ of finite rank $d$ together with a semi-linear and continuous action of $G$. Such a representation is said to be {\it trivial} if there exists a basis of $V$ consisting of only elements of $V^G$, the invariants of $V$ with respect to the action of $G$.

Assume that $E:=B^G$ is a field and let $F$ be a closed subfield of $E$. 
Then $B$ is called $(F,G)$-{\it regular} if

\noindent
(1) $B$ is a domain;

\noindent
(2) $B^G=\mathrm{Fr}\,B^G$, where the action of $G$ on $B$ extends naturally on its fraction field;

\noindent
(3) all elements $$\Big\{b\in B-\{0\}:\ \forall g\in G,\ \exists
\lambda(g)\in F\ \mathrm{s.t.}\ g(b)=\lambda(g)\cdot b\,\Big\}$$
are invertible in $B$.
\vskip 0.20cm
Let  $V$ be a $F$-representation of $G$. Set then $\mathbb D_B(V):=(B\otimes_FV)^G$. Accordingly, we have a natural $B$-linear and $G$-equivariant morphism $$\begin{matrix}\alpha_V:&B\otimes_E\mathbb D_B(V)&\to&B\otimes_FV\\
&\lambda\otimes x&\mapsto&\lambda x.\end{matrix}$$
We say that
$V$ is $B$-{\it admissible} if $B\otimes_FV$ is a trivial $B$-representation of $G$. 
\vskip 0.20cm
\noindent
{\bf Lemma.} (See e.g., [FO]) {\it Assume $B$ is $(F,G)$-regular and let $V$ be a
$F$-representation of $G$. Then

\noindent
(1) The map $\alpha_V$ is injective and $$\mathrm{dim}_E\mathbb D_B(V)\leq\mathrm{dim}_FV;$$

\noindent
(2) The following things are equivalent:

(i) $V$ is $B$-admissible;

(ii) $\mathrm{dim}_E\mathbb D_B(V)=\mathrm{dim}_FV$;

(iii) $\alpha_V$ is an isomorphism.} 

\section{Basic Properties of Various Periods}

With the above discussion and the $p$-adic Hodge structures (to be stated below) in mind, we then can summarize
the essential properties of various $p$-adic periods rings. Our treatment follows [Tsu2]. 

\subsection{Hodge-Tate Periods}

Define the {\it ring of Hodge-Tate periods} to be the graded ring $$\mathbb B_{\mathrm{HT}}:=\bigoplus_{i\in\mathbb Z}\mathbb B_{\mathrm{HT}}^i$$
where, 

\noindent
(i$)_{\mathrm{HT}}$ the $i$-th piece is given by $\mathbb B_{\mathrm{HT}}^i :=\mathbb C(i)$; and 

\noindent
(ii$)_{\mathrm{HT}}$ the ring structure is 
given by the natural multiplication $$\mathbb C(i)\otimes_{\mathbb C}\mathbb C(j)
\to\mathbb C(i+j).$$

\subsection{de Rham Periods}

Fix a $p$-adic number field $K$. 
Denote by $\mathbb B_{\mathrm{dR}}$ the ring of de Rham periods.

\noindent
{\bf Basic Properties of $\mathbb B_{\mathrm{dR}}$:}

\noindent
(i$)_{\mathrm{dR}}$ $\mathbb B_{\mathrm{dR}}$ is a complete discrete valuation field with $\mathbb C_p$ its residue field; 

\noindent
(ii$)_{\mathrm{dR}}$ $\mathbb B_{\mathrm{dR}}$ admits a natural decreasing filtration $$\mathrm{Fil}_{\mathrm{HT}}^i\mathbb B_{\mathrm{dR}}:=\Big\{x\in \mathbb B_{\mathrm{dR}}:v(x)\geq i\Big\}$$ (reflecting  the structure of Hodge filtration). Here we have normalized the valuation so that $v(\mathbb B_{\mathrm{dR}}^*)=\mathbb Z;$ 

\noindent
(iii$)_{\mathrm{dR}}$ $\mathbb B_{\mathrm{dR}}$
admits a natural $G_K$ action which not only preserves the above filtration, but is compatible with the natural induced projection $\mathrm{Fil}^0\mathbb B_{\mathrm{dR}}\to\mathbb C;$

\noindent
(iv$)_{\mathrm{dR}}$ $\mathbb B_{\mathrm{dR}}$ satisfies the following additional fine structures/properties:
\vskip 0.20cm
\noindent
(1)$_{\mathrm{dR}}$ There is a natural $G_K$-equivariant embedding 

\centerline{$P_0:=K_0^{\mathrm{ur}}\otimes_{K_0}\overline{K}\hookrightarrow
\mathcal O_{\mathbb B_{\mathrm{dR}}}=:\mathbb B_{\mathrm{dR}}^+$}

\noindent 
such that its composition with the residue map $\mathbb B_{\mathrm{dR}}^+\twoheadrightarrow\mathbb C$ coincides with the natural embedding
$K_0^{\mathrm{ur}}\otimes_{K_0}\overline{K}\hookrightarrow\mathbb C$;

\noindent
(2)$_{\mathrm{dR}}$ There is a natural $G_K$-equivariant injection $\mathbb Q_p(i)\hookrightarrow \mathrm{Fil}^i\mathbb B_{\mathrm{dR}}$ such that  one (and hence  all) $a\in\mathbb Q_p(1), a\not=0$, maps into a prime element of $\mathbb B_{\mathrm{dR}}$.
In particular, 

(2.1)$_{\mathrm{dR}}$ there are natural $G_K$-equivariant injections 
 $\mathbb Q_p(i)\hookrightarrow \mathrm{Fil}^i\mathbb B_{\mathrm{dR}}$; 

(2.2)$_{\mathrm{dR}}$ there are natural $G_K$-equivariant  isomorphisms

\centerline{$\mathbb C(i)\simeq \mathrm{Gr}_{\mathrm{HT}}^i\mathbb B_{\mathrm{dR}}:=\mathrm{Fil}_{\mathrm{HT}}^i\mathbb B_{\mathrm{dR}}/\mathrm{Fil}_{\mathrm{HT}}^{i+1}\mathbb B_{\mathrm{dR}};$}

\noindent
(3)$_{\mathrm{dR}}$ $\mathbb B_{\mathrm{dR}}^{G_K}=K$.

It appears that $\mathbb B_{\mathrm{dR}}$ depends on $K$. For this, we have

\noindent
(v)$_{\mathrm{dR}}$ If $L/K$ is a finite Galois extension contained in $\overline K/K$, then $$\Big(\mathbb B_{\mathrm{dR}}(L),G_L\Big)\simeq \Big(\mathbb B_{\mathrm{dR}}(K), G_L(\subset G_K)\Big).$$ That is to say, 
$\mathbb B_{\mathrm{dR}}(L)$ together with its Galois action $G_L$ coincides with $\mathbb B_{\mathrm{dR}}(K)$ associated to $K$ together with the induced action of $G_L$ as the restriction from $G_K$ to its subgroup $G_L$.

\subsection{Crystalline Periods}
Denote by $\mathbb B_{\mathrm{crys}}$ the ring of crystalline periods.

\noindent
{\bf Basic Properties of $\mathbb B_{\mathrm{crys}}$:}

\noindent
(i$)_{\mathrm{crys}}$ $\mathbb B_{\mathrm{crys}}$ is a $G_K$-stable subring 
of $\mathbb B_{\mathrm{dR}}$ such that the induced decreasing filtration $\mathrm{Fil}^i\mathbb B_{\mathrm{crys}}:=\mathbb B_{\mathrm{crys}}\cap \mathrm{Fil}^i\mathbb B_{\mathrm{dR}}$ has the same graded pieces $\mathbb C(i)$;

\noindent
(ii$)_{\mathrm{crys}}$ $\mathbb B_{\mathrm{crys}}$ satisfies the following additional structures/properties:

\noindent
(1)$_{\mathrm{crys}}$ There is a natural $\sigma$-semi ($P_0$-)linear
action of $G_K$
and a $G_K$-equivariant	 injective morphism $\varphi:\mathbb B_{\mathrm{crys}}\to\mathbb B_{\mathrm{crys}}$, the so-called Frobenius, such that the following holds

(1.1)$_{\mathrm{crys}}$ For  $t\in\mathbb Q_p(1)\subset \mathbb B_{\mathrm{crys}}$, $\varphi(t)=pt$;

(1.2)$_{\mathrm{crys}}$ $\mathrm{Fil}^0\mathbb B_{\mathrm{crys}}\cap \mathbb B_{\mathrm{crys}}^{\varphi=1}=\mathbb Q_p;$

(1.3)$_{\mathrm{crys}}$ $\forall x\in\mathbb Q_p(i)$, 
$\varphi(x)=p^ix$ and $\mathrm{Fil}^i\mathbb B_{\mathrm{crys}}\cap \mathbb B_{\mathrm{crys}}^{\varphi=p^i}=\mathbb Q_p(i);$

\noindent
(2)$_{\mathrm{crys}}$ The natural map $K\otimes_{K_0}\mathbb B_{\mathrm{crys}}\to \mathbb B_{\mathrm{dR}}$ is injective;

\noindent
(3)$_{\mathrm{crys}}$ $\mathbb B_{\mathrm{crys}}^{G_K}=K_0$;

\noindent
(4)$_{\mathrm{crys}}$ All one dimensional $G_K$-stable $\mathbb Q_p$-vector subspaces of $\mathbb B_{\mathrm{crys}}$ are contained in $P_0\cdot\mathbb Q_p(i), i\in\mathbb Z$.

Similarly, as for $\mathbb B_{\mathrm{dR}}$, we have

\noindent
(iii)$_{\mathrm{crys}}$ If $L/K$ is a finite Galois extension contained in $\overline K/K$, then $$\Big(\mathbb B_{\mathrm{crys}}(L),G_L\Big)\simeq \Big(\mathbb B_{\mathrm{crys}}(K), G_L(\subset G_K)\Big).$$

\subsection{Semi-Stable Periods}

Denote by $\mathbb B_{\mathrm{st}}$ the ring of semi-stable periods.

\noindent
{\bf Basic Properties of  $\mathbb B_{\mathrm{st}}$:}

\noindent
(i)$_{\mathrm{st}}$
$\mathbb B_{\mathrm{st}}$ may be understood as a $G_K$-stable subring of $\mathbb B_{\mathrm{dR}}$. However,
different from $\mathbb B_{\mathrm{crys}}$, such an embedding of $\mathbb B_{\mathrm{st}}$ in $\mathbb B_{\mathrm{crys}}$ depends on the choices of prime element $\pi$ of $K$.

\noindent
(ii)$_{\mathrm{st}}$
$\mathbb B_{\mathrm{st}}$ satisfies the following additional structures/properties:

\noindent
(1)$_{\mathrm{st}}$ Corresponding to a systematic choice of $p^n$-th root of $\pi$ in 
$\overline K$: $s=(s_n)_{n\in\mathbb N}, s_0=\pi, s_{n+1}^p=s_n$, there is a natural element 
$u_s\in \mathbb B_{\mathrm{st}}$ such that

(1.1)$_{\mathrm{st}}$ $\mathbb B_{\mathrm{st}}=\mathbb B_{\mathrm{crys}}\big[u_s\big]$;

(1.2)$_{\mathrm{st}}$	$\forall g\in G_K$, $g(u_s)=u_{g(s)}$, where
 $g(s)=\big(g(s_n)\big)_{n\in \mathbb N}$;

(1.3)$_{\mathrm{st}}$ If $s'=(s_n')$ is another choice, then $u_{s'}=u_s+t$, where

\hskip 1.0cm $(s_n's_n^{-1})_{n\in \mathbb N}=:t\in\mathbb Q_p(1)\subset \mathbb B_{\mathrm{crys}}$;

\noindent
(2)$_{\mathrm{st}}$ $\mathbb B_{\mathrm{st}}$ admits a natural $G_K$-equivariant Frobenius $\varphi(u_s)=p\cdot u_s$ extending the Frobenius $\varphi$ on $\mathbb B_{\mathrm{crys}}$;

\noindent
(3)$_{\mathrm{st}}$ $\mathbb B_{\mathrm{st}}$ admits a natural monodromy operator $N:\mathbb B_{\mathrm{st}}\to \mathbb B_{\mathrm{st}}$ satisfying

(3.0)$_{\mathrm{st}}$ $N$ is a $\mathbb B_{\mathrm{crys}}$-derivation and $N(u_s)=1$; 

(3.1)$_{\mathrm{st}}$ $N$ is $G_K$-equivariant;

(3.2)$_{\mathrm{st}}$ $N\varphi=p\varphi N$;

(3.3)$_{\mathrm{st}}$ $\mathbb B_{\mathrm{st}}^{N=0}=\mathbb B_{\mathrm{crys}};$ and

(3.4)$_{\mathrm{st}}$ $\mathrm{Fil}^0\mathbb B_{\mathrm{dR}}\cap \mathbb B_{\mathrm{st}}^{N=0,\varphi=1}=\mathbb Q_p;$

\noindent
(4)$_{\mathrm{st}}$ The natural map $K\otimes_{K_0}\mathbb B_{\mathrm{st}}\to \mathbb B_{\mathrm{dR}}$ is injective; and

\noindent
(5)$_{\mathrm{st}}$ $\mathbb B_{\mathrm{st}}^{G_K}=K_0$;

\noindent
(6)$_{\mathrm{crys}}$ All one dimensional $G_K$-stable $\mathbb Q_p$-vector subspaces of $\mathbb B_{\mathrm{st}}$ are contained in $P_0\cdot\mathbb Q_p(i), i\in\mathbb Z$.

Similarly,

\noindent
(iii)$_{\mathrm{crys}}$ If $L/K$ is a finite Galois extension contained in $\overline K/K$, then $$\Big(\mathbb B_{\mathrm{st}}(L),G_L, e(L/K)^{-1}N\Big)\simeq \Big(\mathbb B_{\mathrm{st}}(K), G_L(\subset G_K),N\Big).$$ Here 
$e(L/K)$ denotes the ramification index of the extension $L/K$.

\section{Hodge-Tate, de Rham, Semi-Stable and Crystalline Reps}

\subsection{Definition}
Let $V$ be a $p$-adic representation of $G_K$, and let $$\mathbb D_{\bullet}(V):=\Big(\mathbb B_{\bullet}\otimes_{\mathbb Q_p}V\Big)^{G_K}$$ where $\bullet$ is the running symbol for
HT, dR, st, crys, and $G_K$ acts on $\mathbb B_{\bullet}\otimes_{\mathbb Q_p}V$ via diagonal action of $G_K$. Clearly, from the natural structure of the ring of periods, 
there is an induced structures on $\mathbb D_{\bullet}(V)$. In particular,
since $$\mathbb C^{G_K}=\mathbb B_{\mathrm{HT}}^{G_K}= \mathbb B_{\mathrm{dR}}^{G_K}=K, \quad\mathrm{and}\quad\mathbb B_{\mathrm{st}}^{G_K}= \mathbb B_{\mathrm{crys}}^{G_K}=K_0,$$

\noindent
(i) $\mathbb D_{\mathrm{HT}}(V), \mathbb D_{\mathrm{dR}}(V)$ are $K$-vector spaces; and	 

\noindent
(ii) $\mathbb D_{\mathrm{st}}(V), \mathbb D_{\mathrm{crys}}(V)$ are $K_0$-vector spaces. 

One checks easily that $\mathbb B_{\bullet}$ is $\Big(\mathbb B_{\bullet}^{G_K},G_K\Big)$-regular.
Accordingly, following Fontaine,  we call a $p$-adic Galois representation $V$ of 
$G_K$ a $\bullet$-{\it representation}, where $\bullet$=Hodge-Tate, de Rham, 
semi-stable, crystalline, if
$V$ is $\mathbb B_{\bullet}$-admissible, that is to say, if 
$$\mathrm{dim}_{\mathbb B_{\bullet}^{G_K}}\mathbb D_{\bullet}(V)=\mathrm{dim}_{\mathbb Q_p}(V).$$

\subsection{Basic Structures of $\mathbb D_{\bullet}(V)$}

Induced from Fontaine's	 rings of various periods, there are natural structures on the space
$\mathbb D_{\bullet}(V)$ associated to a $p$-adic Galois representation $V$ of $G_K$. 

\noindent
$\bullet$ {\bf Hodge-Tate}: The graded structure on $\mathbb B_{\mathrm{HT}}$ induces  a natural graded structure on $K$-vector space $\mathbb D_{\mathrm{HT}}(V)$. More precisely,
$$\mathbb D_{\mathrm{HT}}(V)=\bigoplus_{i\in\mathbb Z}
 \mathbb D_{\mathrm{HT}}^i(V)\quad\mathrm{where}\quad \mathbb D_{\mathrm{HT}}^i(V):=\Big(\mathbb C(i)\otimes_{\mathbb Q_p}V\Big)^{G_K}.$$
 
\noindent
$\bullet$ {\bf de Rham}:  The decreasing filtration structure on $\mathbb B_{\mathrm{dR}}$ induces  a natural decreasing filtration of $K$-vector subspaces on $\mathbb D_{\mathrm{dR}}(V)$. More precisely,$$\mathrm{Fil}_{\mathrm{HT}}\mathbb D_{\mathrm{dR}}^i(V):=\Big(\mathrm{Fil}^i_{\mathrm{HT}}\mathbb B_{\mathrm{dR}}\otimes_{\mathbb Q_p}V\Big)^{G_K}.$$ This filtration
 is {\it exhaustive} and {\it separated}, that is, we have 
 $$\bigcup_{i\in\mathbb Z}\mathrm{Fil}_{\mathrm{HT}}^i
 \mathbb D_{\mathrm{dR}}(V)=\mathbb D_{\mathrm{dR}}(V)
 \quad\mathrm{and}\quad 
 \bigcap_{i\in\mathbb Z}\mathrm{Fil}_{\mathrm{HT}}^i
 \mathbb D_{\mathrm{dR}}(V)=0.$$ Moreover, 
by (2.2)$_{\mathrm{dR}}$, we have the following natural injection 
of $K$-vector spaces
$$\mathrm{Gr}_{\mathrm{HT}}\mathbb D_{\mathrm{dR}}(V)
:=\bigoplus_{i\in\mathbb Z}\mathrm{Fil}_{\mathrm{HT}}^i\mathbb D_{\mathrm{dR}}(V)\big/\mathrm{Fil}_{\mathrm{HT}}^{i+1}\mathbb D_{\mathrm{dR}}(V)
\hookrightarrow\mathbb D_{\mathrm{HT}}(V).\eqno(*)$$
 
\noindent
$\bullet$ {\bf Semi-Stable}: By (4)$_{\mathrm{st}}$, we have a non-canonical embedding of 
$K\otimes_{K_0}\mathbb B_{\mathrm{st}}\hookrightarrow \mathbb B_{\mathrm{dR}},$ and hence 
a natural inclusion $$K\otimes_{K_0}\mathbb D_{\mathrm{st}}(V)\hookrightarrow\mathbb D_{\mathrm{dR}}(V).\eqno(**)$$ Consequently,
there is a natural decreasing filtration by $K$-vector subspaces on 
$K\otimes_{K_0}\mathbb D_{\mathrm{st}}(V)$. Moreover, from the Frobenius structure 
$\varphi$ and monodromy operator $N$ on $\mathbb B_{\mathrm{st}}$, we get a natural 
Frobeinus structure $\varphi: \mathbb D_{\mathrm{st}}(V)\to \mathbb D_{\mathrm{st}}(V)$ 
and a monodromy operator $N:\mathbb D_{\mathrm{st}}(V)\to \mathbb D_{\mathrm{st}}(V)$ 
which are all $K_0$-linear and satisfy the relation
$$N\varphi=p\cdot\varphi N.$$

\noindent
$\bullet$ {\bf Crystalline}: By (2)$_{\mathrm{crys}}$, we have a canonical embedding 
$K\otimes_{K_0}\mathbb B_{\mathrm{crys}}\hookrightarrow \mathbb B_{\mathrm{dR}},$ 
and hence a natural inclusion $$K\otimes_{K_0}\mathbb D_{\mathrm{crys}}(V)
\hookrightarrow\mathbb D_{\mathrm{dR}}(V).\eqno(*_3)$$ Consequently,
there is a natural decreasing filtration by $K$-vector subspaces on 
$K\otimes_{K_0}\mathbb D_{\mathrm{crys}}(V)$. Moreover, from the Frobenius 
structure $\varphi$  on $\mathbb B_{\mathrm{crys}}$, we get a 
natural Frobeinus $\varphi: \mathbb D_{\mathrm{crys}}(V)\to \mathbb D_{\mathrm{crys}}(V)$ 
which is  $K_0$-linear. 

Finally, by (3.3)$_{\mathrm{st}}$, we have $\mathbb B_{\mathrm{st}}^{N=0}=\mathbb B_{\mathrm{crys}}$ and hence
$$\mathbb D_{\mathrm{st}}(V)^{N=0}=\mathbb D_{\mathrm{crys}}(V).\eqno(*_4)$$

\subsection{Relations among Various $p$-adic Representations}

 Let $V$ be a $p$-adic representation of $G_K$. Then from Lemma in \S28, 
 $(*,**,*_3,*_4)$, and the fact that the natural $\mathbb C$-linear morphism 
$$\bigoplus_{i\in\mathbb Z}\mathbb C(-i)\otimes_K\Big(\mathbb C(i)
\otimes_{K}V\Big)^{G_K}\to \mathbb C\otimes_{\mathbb Q_p}V$$ is an injection, we obtain the following inequalities:
$$\begin{aligned}\mathrm{dim}_{K_0}&\mathbb D_{\mathrm{crys}}(V)\leq\mathrm{dim}_{K_0}
\mathbb D_{\mathrm{st}}(V)\\
&\leq
\mathrm{dim}_{K}\mathbb D_{\mathrm{dR}}(V)\leq \mathrm{dim}_{K}\mathbb D_{\mathrm{HT}}(V)\\
&\leq\mathrm{dim}_{\mathbb Q_p}V.\end{aligned}$$
Consequently,

\noindent
(i) $\mathbb D_{\bullet}(V)$ are all finite dimensional $\mathbb B_\bullet^{G_K}$-vector spaces;

\noindent
(ii) $\varphi$, whenever makes sense, is an isomorphism; and most importantly, 

\noindent
(iii) there are simple implications that 
$$\mathrm{crystalline}\ \Rightarrow\ \mathrm{semi\, stable}\ \Rightarrow\ \mathrm{de\, Rham}\Rightarrow\mathrm{Hodge\ Tate}.$$
\eject
\noindent
{\bf Proposition.} (Fontaine) {\it Let $V$ be a $p$-adic representation of $G_K$. Use $\bullet$ as the runing symbol for HT, dR, st, crys. Then

\noindent
(1) the natural $\mathbb B_{\bullet}$-linear map $$
\mathbb B_{\bullet}\otimes_{\mathbb B_{\bullet}^{G_K}}\mathbb D_{\bullet}(V)\hookrightarrow\mathbb B_{\bullet}\otimes_{\mathbb Q_p}V$$
is  a $G_K$-equivariant morphism which preserves the grads, where

\noindent
(i)  $G_K$ acts on the left hand side via the action on $\mathbb B_{\bullet}$ and on the right hand side via the diagonal one; and 

\noindent
(ii) the graded structures are given  on the left hand side by	$$\sum_{i=i_0+i_1}\mathrm{Fil}^{i_0}\mathbb B_{\bullet}\otimes_K\mathrm{Fil}^{i_0}\mathbb D_{\bullet}^{i_1}(V)$$	 and on the right hand	by $\mathrm{Fil}^{i}\mathbb B_{\bullet}\otimes_{\mathbb Q_p}V$;

\noindent
(2) If $V$ is $\bullet$-admissible, then the $\mathbb B_{\bullet}$-linear map $$
\mathbb B_{\bullet}\otimes_{\mathbb B_{\bullet}^{G_K}}\mathbb D_{\bullet}(V)\to\mathbb B_{\bullet}\otimes_{\mathbb Q_p}V$$
 is an isomorphism; Moreover, 

\noindent
(3)
(i) If $V$ is of Hodge-Tate, then by considering the degree zero parts, 
we get a natural isomorphism, the so-called Hodge-Tate decomposition,
$$\oplus_{i\in\mathbb Z}\mathbb C(-i)\otimes_{\mathbb Q_p}
\mathbb D_{\mathrm{HT}}^i(V)\simeq \mathbb C\otimes_{\mathbb Q_P}V;$$

\noindent
(ii) If $V$ is semi-stable, then  the natural $\mathbb B_{\mathrm{st}}$-linear map $$
\mathbb B_{\mathrm{st}}\otimes_{K_0}\mathbb D_{\mathrm{st}}(V)
\simeq\mathbb B_{\mathrm{st}}\otimes_{\mathbb Q_p}V$$ commutes with $\varphi$ and $N$;

\noindent
(iii) If $V$ is crystalline, then the natural $\mathbb B_{\mathrm{crys}}$-linear map $$
\mathbb B_{\mathrm{crys}}\otimes_{K_0}\mathbb D_{\mathrm{crys}}(V)\simeq\mathbb B_{\mathrm{crys}}\otimes_{\mathbb Q_p}V$$
 commutes with $\varphi$.}
	
\subsection{Examples}

\noindent
(1) {\it Tate Twist}: $\mathbb Q_p(i)$ given by cyclotomic characters 
$\chi_{\mathrm{cyclo}}^i, i\in\mathbb Z$. All are crystalline. Indeed, $D=\mathbb D_{\mathrm{crys}}(\mathbb Q_p(i))=K_0\cdot e$
with $e=t^{-i}\otimes t^i$ and $$\varphi(e)=p^{-i}e,\ \ \mathrm{Fil}_{\mathrm{HT}}^{-i}D=D,\ \ \mathrm{Fil}_{\mathrm{HT}}^{-i+1}D=0.
$$

\noindent
(2) {\it Unramified Representations}: A unramified $p$-adic Galois representation, i.e., where
the inertial group $I_K$ acts trivially, is crystalline. Moreover, 
a crystalline representation is unramified if and only if its associated Hodge-Tate 
filtration satisfies
$\mathrm{Fil}^0_{\mathrm{HT}}\mathbb D_{\mathrm{dR}}(V)=\mathbb D_{\mathrm{dR}}(V)$ and $\mathrm{Fil}^1_{\mathrm{HT}}\mathbb D_{\mathrm{dR}}(V)=0$.
\vskip 0.20cm
\noindent
(3) {\it Semi-Stable Representations}:
All Tate modules $T_p(E)$ for Tate curves $E$ are semi-stable representations.
\vskip 0.20cm
\noindent
(4) {\it de Rham and Hodge-Tate Representations}: 

\noindent
(i) Extension of $\mathbb Q_p$ by $\mathbb Q_p(1)$ is de Rham, but

\noindent
(ii) Non-trivial  extension of $\mathbb Q_p(1)$ by $\mathbb Q_p$ is Hodge-Tate but not de Rham.
\vskip 0.20cm
\noindent
(5) {\it One Dimensional Galois Representations}: In this case, there are following equivalences

\noindent
(i) Hodge-Tate $\Leftrightarrow$ de Rham

\noindent
$\Leftrightarrow$ There is an open subgroup $I_L$ of $I_K$ and an integer $i$ such that the induced action of $I_L$ on $V(-i)$ is trivial;

\noindent
(ii) Semi-stable $\Leftrightarrow$ Crystalline

\noindent
$\Leftrightarrow$  the induced action of $I_K$ on $V(-i)$ is trivial.
\vskip 0.20cm
\noindent
(6) {\it Not Even Hodge-Tate}: $V$ is a two dimensional $\mathbb Q_p$-vector space 
equipped with an action of $G_K$ given by $\begin{pmatrix}1&\log_p\chi(g)\\
0&1\end{pmatrix}.$ By \S45, the Sen operator $\Theta_V=\begin{pmatrix}0&1\\ 0&0\end{pmatrix},$
it is not of Hodge-Tate.

\section{$p$-adic Hodge Theory}
Fix a $p$-adic number field $K$ with absolute Galois group $G_K$.
\vskip 0.20cm
\noindent
{\bf Therem.} ({\bf $p$-adic Hodge Theory}) {\it  Let $X$ be a $n$-dimensional proper regular variety defined over $K$. Denote by $\Big(V:=H_{\mathrm{et}}^m(X_{\overline K},\mathbb Q_p),\rho\Big)$ the induced representation of $G_K$, where $H_{\mathrm{et}}^m(X_{\overline K},\mathbb Q_p)$ denotes the $m$-th $p$-adic etale cohomology group of $X$. Then the following conjectures hold:}
\vskip 0.20cm
\noindent
$\bullet$ {\bf Hodge-Tate} {\it (i) The Galois representation $\Big(H_{\mathrm{et}}^m(X_{\overline K},\mathbb Q_p),\rho_K\Big)$ is of Hodge-Tate type; and 

\noindent
(ii) There is a natural graded preserved isomorphism $$\mathbb D_{\mathrm{HT}}\Big(H_{\mathrm{et}}^m(X_{\overline K},\mathbb Q_p)\Big)\simeq\oplus_{i\in\mathbb Z}H^{m-i}(X,\Omega_{X/K}^i),$$ and hence the following $G_K$-equivariant Hodge-Tate decomposition} 
$$
\mathbb C\otimes_{\mathbb Q_p}H_{\mathrm{et}}^m(X_{\overline K},\mathbb Q_p)
\simeq\oplus_{i=0}^m\mathbb C(-i)\otimes_KH^{m-i}(X,\Omega_{X/K}^i);$$
$\bullet$ {\bf de Rham} {\it (i) The Galois representation	
$\Big(H_{\mathrm{et}}^m(X_{\overline K},\mathbb Q_p),\rho_K\Big)$ is of de Rham type. Moreover, 

\noindent
(ii) $\mathbb D_{\mathrm{dR}}(V)$ together with its associated Hodge 
filtration is isomorphic to the de Rham cohomology $H_{\mathrm{dR}}^m(X_K/K)$ equipped with the Hodge filtration;}
\vskip 0.20cm
\noindent
$\bullet$ {\bf Semi-Stable} {\it (i) If $X$ has	 a semi-stable reduction 
$(Y,D)$, then the Galois representation	 $\Big(H_{\mathrm{et}}^m(X_{\overline K},\mathbb Q_p),\rho_K\Big)$ is semi-stable. Moreover,

\noindent
(ii) The associated filtered $(\varphi,N)$-module $\mathbb D_{\mathrm{st}}(V)$
is canonically isomorphic to the following filtered $(\varphi,N)$-module on the	 log crystalline cohomology $H^m_{\log}((Y,D)/K_0)$: Choose a semi-stable model $\frak X\to\mathcal O_K$ of $X/K$ so that we obtain 
a log geometric structure $(Y,D)$ on the special fiber. Then induced from the log crystalline cohomology 
of the special fiber, there is a natural weakly admissible filtered 
$(\varphi,N)$-module structure on $\Big(H^m_{\log}((Y,D)/K_0),H_{\mathrm {dR}}^m(X_K/K)\Big)$;}
\vskip 0.20cm
\noindent
$\bullet$ {\bf Crystalline} {\it (i) If $X$ has	 a good reduction, 
the Galois representation  $\Big(H_{\mathrm{et}}^m(X_{\overline K},\mathbb Q_p),\rho_K\Big)$ is crystalline. Moreover,

\noindent
(ii) The filtered $\varphi$-module $\mathbb D_{\mathrm{crys}}(V)$
is canonically isomorphic to the following filtered $\varphi$-module on the  crystalline	
cohomology $H^m_{\mathrm{crys}}(Y/K_0)$: Choose a proper regular model $\frak X\to\mathcal O_K$ of $X/K$. Then induced from	
crystalline cohomology of the special fiber, there is a natural weakly admissible filtered 
$\varphi$-module structure on  $\Big(H^m_{\mathrm{crys}}(Y/K_0),
H_{\mathrm {dR}}^m(X_K/K)\Big)$.}
\vskip 0.20cm
The Hodge-Tate conjecture, due to mainly Tate, is a $p$-adic analogus of the standard Hodge theory for projective complex manifolds. This conjectures was solved by Tate for abelian varieties with good reduction, by Raynaud for all abelian varieties, by Bloch-Kato ([BK1])
for varieties with good reduction and finally by Faltings ([Fa1]) in general.

The de Rham conjecture and the crystalline conjecture are due to Fontaine ([Fon3]) and are 
solved by Fontaine-Messing ([FMe]) when $K=K_0$, $\mathrm{dim}X\leq p-1$ and $X$ has good 
reduction, and by Faltings ([Fa2]) in gereral. The above filtered $\varphi$-module structure on 
the de Rham cohomology is due to Berthelot-Ogus ([Ber1,2], [BO1,2]), and the independence 
issue for the filtered $\varphi$-structure on the de Rham cohomology on the model used is 
established by Gillet-Messing ([GM]).

The semi-stable conjecture is due to Fontaine and U. Jannsen ([Fo6]), 
solved by Fontaine  for abelian varieties ([Fo6]), by Kato	
when $\mathrm{dim}X\leq(p-1)/2$ ([K]), by 
Tsuji ([Tsu1]), Niziol ([Ni1,2]) and Faltings ([Fa4]) independently in general. The above 
filtered $(\varphi,N)$-structure on the de Rham cohomology is due to Hyodo-Kato ([HK])
and the independence of the model chosen can be established via de Jong's alternation theory ([dJ]).
\vfill
\eject
\centerline{\Large\bf Chapter XII. Fontaine's Rings of Periods}
\vskip 0.80cm
In this chapter, for completeness, we explain the essentials of various rings of periods following Fontaine (e.g. [FO]).
\section{The Ring of de Rham Periods $\mathbb B_{\mathrm{dR}}$}
To have a reasonable theory of $p$-adic Galois representations, the standard
$p$-adic cyclotomic character should be involved in a natural way. 
Accordingly, to construct the ring of good period $\mathbb B_{\bullet}$, 
we need to find an element $t\in \mathbb B_{\bullet}$ which is a period for the cyclotomic 
character. That is to say, there should be an element 
$t\in \mathbb B_{\bullet}$ such that
$$g(t)=\chi_{\mathrm{cycl}}(g)\cdot t\qquad\mathrm{for\ all}\ g\in G_K.$$ 

As a starting point, one may naively try $\mathbb C$. However it does not work 
since $$\Big\{x\in\mathbb C:g(x)=\chi_{\mathrm{cycl}}(g)\cdot x,\forall g\in G_K\Big\}=\{0\}.$$
Thus we need to enlarge it. This then leads to the Tate module $\mathbb Z(1)$ 
and hence the ring of Hodge-Tate periods
$$\mathbb B_{\mathrm{HT}}:=\oplus_{i\in\mathbb Z}\mathbb C(i),$$ which in a certain sense
is the simplest ring of periods.

With the simplest one found, it is then very natural for us to seek a sort of
$\lq$universal' one. 
With the theory of field of norms,  we are led	
to the Cohen ring $\widetilde{\mathbb A}^+:=W(\widetilde{\mathbb E}^+)$ associated to 
$\widetilde{\mathbb E}^+$, or better, to its fractional field 
$\widetilde{\mathbb B}^+:=\widetilde{\mathbb A}^+[\frac{1}{p}]$. While this basically works, an essential modification should be made.

To be more precise, let $\varepsilon=(\varepsilon^{(n)})\in \widetilde{\mathbb E}^+$ with $\varepsilon^{(0)}=1$, 
$\varepsilon^{(1)}\not=1$. Assume that 
$$t:=\log[\varepsilon]=-\sum_{n=1}^\infty\frac{(1-[\varepsilon])^n}{n}$$ makes sense. That is to say, assume that the infinite power series above converges. 
Then, formally, we have for all $g\in G_{K_0}$,
$$\begin{aligned}g(t)=&g\big(\log\,[\varepsilon]\big)=\log\Big(\big[g(\varepsilon^{(0)},\varepsilon^{(1)},\dots)\big]\Big)\\
=&\log\big([\varepsilon^{\chi_{\mathrm{cycl}}(g)}]\big)= \chi_{\mathrm{cycl}}(g)\cdot t.\end{aligned}$$
In other words, whenever it makes sense, $t=\log\,[\varepsilon]$ is a cyclotomic period. 
Thus, we need to create a ring within which the above series defining $\log\,[\varepsilon]$ converges.

For the infinite series defining $\log[\varepsilon]$ 
to converge, it suffices to make $1-[\varepsilon]$ small. However, in $\widetilde{\mathbb E}^+$, 
we have $$v_{\mathbb E}(\varepsilon-1)=\lim_{n\to\infty}v_p(\varepsilon^{(n)}-1)^{p^n}=\frac{p}{p-1}.$$ 
In other words, within $\widetilde{\mathbb E}^+$,
$\varepsilon-1$ is not really very small. To overcome this difficulty, following Fontaine, 
we go as follows:

From the natural isomorphism $\displaystyle{\widetilde{\mathbb E}\simeq\lim_{\begin{subarray}{1}\ \longleftarrow\\
{x\mapsto x^p}\end{subarray}}\mathcal O_{\mathbb C}/p\mathcal O_{\mathbb C}}$, we obtain an induced homomorphism $$\begin{matrix}\theta:&\widetilde{\mathbb E}&\to&\mathcal O_{\mathbb C}/p\mathcal O_{\mathbb C}\\
&(x^{(n)})&\mapsto& x^{(0)}.\end{matrix}$$ 
Lift this construction to the characteristic zero world. Since 
$$\widetilde{\mathbb B}^+:=\widetilde{\mathbb A}^+\Big[\frac{1}{p}\Big]:=
\Big\{\sum_{k\gg-\infty}p^k[x_k]: x_k\in\widetilde{\mathbb E}^+\Big\}$$ where 
$[x]\in \widetilde{\mathbb A}^+$ denotes the Teichm\"uller lift of $x\in \widetilde{\mathbb E}^+$,
we obtain a natural morphism, a lift of $\theta$, $$\begin{matrix}\theta:&\widetilde{\mathbb B}^+&\to&\mathbb C\\
&\sum p^k[x_k]&\mapsto&\sum p^kx_k^{(0)}.\end{matrix}$$ (Here we have used 
the isomorphism $$\widetilde{\mathbb E}^+\simeq \lim_{\begin{subarray}{1}\ \longleftarrow\\
{x\mapsto x^p}\end{subarray}}\mathcal O_{\mathbb C}/p\mathcal O_{\mathbb C}
\simeq \lim_{\begin{subarray}{1}\ \longleftarrow\\
{x\mapsto x^p}\end{subarray}}\mathcal O_{\mathbb C},$$ namely, a shift from $\mathcal O_{\mathbb C}/p\mathcal O_{\mathbb C}$
a characteristic $p$ one to  $\mathcal O_{\mathbb C}$ a characteristic zero world,
 so that elements $x$ take the forms
 $x=(x^{(n)})$ with $x^{(n)}\in\mathcal O_{\mathbb C}$.)

Recall that $\varepsilon=(\varepsilon^{(n)})\in \widetilde{\mathbb E}^+$ with 
$\varepsilon^{(0)}=1$, $\varepsilon^{(0)}\not=1$. Set 
$$\varepsilon_1:=\varepsilon^p=(\varepsilon^{(1)},\varepsilon^{(2)},\dots)
\in \widetilde{\mathbb E}^+\qquad\mathrm{and}\qquad
\omega:=\frac{[\varepsilon]-1}{[\varepsilon_1]-1}.$$ Then $\theta(\omega)=1+\varepsilon^{(1)}
+\cdots+(\varepsilon^{(1)})^{p-1}=0.$ In other words, $\langle\omega\rangle\subset \mathrm{Ker}(\theta).$ 
\vskip 0.20cm
\noindent
{\bf Lemma.} (Fontaine) {\it $\mathrm{Ker}(\theta)=\langle\omega\rangle.$}

\noindent
{\it Proof.} Obviously, $\mathrm{Ker}(\theta)$ is an ideal of $\widetilde{\mathbb E}^+$ whose 
elements satisfying $v_{\mathbb E}(x)\geq 1$. Note that $\omega\in \mathrm{Ker}(\theta)$ 
with its modulo $p$ reduction $\overline \omega$ satisfies $v_{\mathbb E}(\overline \omega)=1$. 
Thus the natural injection map $\langle \omega\rangle\to \mathrm{Ker}(\theta)$ is	
surjective modulo $p$. Since both sides are complete for the $p$-adic topology, 
this has to be an isomorphism.
\vskip 0.20cm
Note that in particular $\theta([\varepsilon]-1)=0$ i.e., 
$[\varepsilon]-1\in\mathrm{Ker}(\theta)=\langle\omega\rangle$. Thus in order to make 
$[\varepsilon]-1$ small, it suffices to introduce the $\mathrm{Ker}(\theta)$-adic, 
or the same $\omega$-adic, topology. 
Accordingly, let $$\mathbb B_{\mathrm{dR}}^+:=\lim_{\leftarrow_n}\widetilde{\mathbb B}^+/(\mathrm{Ker}\theta)^n,$$ namely, define
$\mathbb B_{\mathrm{dR}}^+$ to be the ring obtained by completing $\widetilde{\mathbb B}^+$ 
with respect to the $\mathrm{Ker}(\theta)$-adic topology. 

Clearly, $t=\log([\varepsilon])\in \mathbb B_{\mathrm{dR}}^+.$ Indeed, we have the follows.
\vskip 0.20cm
\noindent
{\bf Lemma.} (Fontaine) {\it (1) $\mathbb B_{\mathrm{dR}}:=\mathbb B_{\mathrm{dR}}^+[\frac{1}{t}]$ is a field;

\noindent
(2) There is a natural filtration $\mathrm{Fil}_{HT}^i\mathbb B_{\mathrm{dR}}
=t^i\cdot \mathbb B_{\mathrm{dR}}^+$ such that
$$\mathrm{Gr}_{\mathrm{HT}}\mathbb B_{\mathrm{dR}}\simeq\oplus_{i\in\mathbb Z}\mathbb C(i);$$

\noindent
(3) There is a natural $G_K$ action on $\mathbb B_{\mathrm{dR}}$ with $\mathbb B_{\mathrm{dR}}^{G_K}=K.$}

\section{The Ring of Crystalline Periods $\mathbb B_{\mathrm{crys}}$}

The point here is to create a subring  $\mathbb B_{\mathrm{crys}}$ of 
$\mathbb B_{\mathrm{dR}}$ which contains the cyclotomic period $t$ and is equipped with a natural Frobenius structure.  Its construction is essentially based on the following two relations:

\noindent
(1) $\displaystyle{\varphi(t)=\log([\varepsilon^p])=\log([\varepsilon]^p)
=p\log([\varepsilon])=p\cdot t,}$ and

\noindent
(2) $\displaystyle{\varphi\big(\mathrm{Ker}(\theta)+p\cdot W(\widetilde{\mathbb E^+})\big)\subset 
\mathrm{Ker}(\theta)+p\cdot W(\widetilde{\mathbb E^+}).}$
\vskip 0.20cm
Indeed, in order to have $t\in \mathbb B_{\mathrm{crys}}$, we need to analyze the terms $\frac{([\varepsilon]-1)^n}{n}$ appeared in the defining series of $t=\log[\varepsilon]$. Note that
$$\frac{([\varepsilon]-1)^n}{n}=(n-1)!([\varepsilon_1]-1)^n\frac{\omega^n}{n!}.$$
Since,	(i) $p$-adically, $(n-1)!\to 0$, and (ii) both 
$[\varepsilon_1]-1$ and $\varphi([\varepsilon_1]-1)$ are in 
$W(\widetilde{\mathbb E}^+)$,
we need to understand how all $\varphi\Big(\frac{\omega^n}{n!}\Big)$ behave.

For this, recall that on $W(\widetilde{\mathbb E}^+)$, we have a Frobenius map $$\varphi:(a_0,a_1,\dots,a_n,\dots)\mapsto (a_0^p,a_1^p,\dots,a_n^p,\dots).$$ So, for all $b\in 
W(\widetilde{\mathbb E}^+)$, $\varphi(b)\equiv b^p\ \mathrm{mod}\,p.$ In particular,
$$\varphi(\omega)=\omega^p+p\eta=p\Big(\eta+(p-1)!\frac{\omega^p}{p!}\Big)$$ for a certain 
$\eta\in W(\widetilde{\mathbb E}^+)$. Consequently, 
$$\varphi\Big(\frac{\omega^m}{m!}\Big)=\frac{p^m}{m!}\cdot \Big(\eta+(p-1)!
 \frac{\omega^p}{p!}\Big)^m$$ which are contained in $W(\widetilde{\mathbb E}^+)\Big[\frac{\omega^p}{p!}\Big].$

All this then leads to the following constructions:
\vskip 0.20cm
\noindent
(1) Starting from 
$\widetilde{\mathbb A^+}=W(\widetilde{\mathbb E}^+)$, we introduce the ring $\mathbb A_{\mathrm{crys}}^0$ by adding all elements $\frac{a^m}{m!}$ for $a\in \mathrm{Ker}(\theta)$, the so-called {\it divided power envelope  of} $\widetilde{\mathbb A^+}=W(\widetilde{\mathbb E}^+)$
{\it with respect to} $\mathrm{Ker}(\theta)$;

\noindent
(2) To make $(n-1)!$ small, we need to use $p$-adic topology and hence to obtain  the ring
$$\mathbb A_{\mathrm{crys}}:=
\lim_{\leftarrow_n}\mathbb A_{\mathrm{crys}}^0/p^n\mathbb A_{\mathrm{crys}}^0=\Big\{\sum_{n=0}^\infty a_n\frac{\omega^n}{n!}:
a_n\to 0\ p\mathrm{-adically\ in\ } W(\widetilde{\mathbb E}^+)\Big\};$$

\noindent
(3) By inverting $p$, we get $$\mathbb B_{\mathrm{crys}}^+:=\mathbb A_{\mathrm{crys}}[\frac{1}{p}]=\Big\{\sum_{n=0}^\infty a_n\frac{\omega^n}{n!}:
a_n\to 0\ p\mathrm{-adically\ in}\ W(\widetilde{\mathbb E}^+)\Big[\frac{1}{p}\Big]\Big\}.$$
Clearly, $\mathbb B_{\mathrm{crys}}^+$ contains $t$ and is naturally contained in $\mathbb B_{\mathrm{dR}}^+$;
(Indeed, we have
$$\mathbb B_{\mathrm{crys}}^+=\Big\{\sum_{n=0}^\infty a_n\frac{\omega^n}{n!}
\in\mathbb B_{\mathrm{dR}}: a_n\to 0\ \ \mathrm{in}\ \ W(\widetilde{\mathbb E}^+)
\Big[\frac{1}{p}\Big]\Big\}.)$$

\noindent
(4) Finally, define the ring of crystalline periods by $\mathbb B_{\mathrm{crys}}:=\mathbb B_{\mathrm{crys}}^+[\frac{1}{t}]$ with the extension of Frobenius via $\varphi(\frac{1}{t}):=\frac{1}{pt}$.
\vskip 0.20cm
\noindent
{\it Remark.} The domain $\mathbb B_{\mathrm{crys}}$ is	 
not a field. For example,  $\omega-p$ is in $\mathbb B_{\mathrm{crys}}\backslash \mathbb B_{\mathrm{crys}}^*$.

\section{The Ring of Semi-Stable Periods $\mathbb B_{\mathrm{st}}$}

Since for semi-stable periods,
$\mathbb B_{\mathrm{st}}^{N=0}=\mathbb B_{\mathrm{crys}}$, a natural way to construct	
$\mathbb B_{\mathrm{st}}$ is to enlarge
$\mathbb B_{\mathrm{crys}}$.
For this purpose, motivated by analysis, 
we may simply try 
to find a transcendental element $T$ over $\mathbb B_{\mathrm{crys}}$, or better, 
over its fraction field $\mathbb C_{\mathrm{crys}}:=\mathrm{Fr}\,\mathbb B_{\mathrm{crys}}$, such that

\noindent
(1) $\varphi(T)=pT$;

\noindent
(2) $N(T)=1$, which implies $N\big(\sum a_nT^n\big)=\sum na_nT^{n-1}$ for all $a_n\in \mathbb C_{\mathrm{crys}}$; and

\noindent
(3) There is a natural action of $G_K$ on $T$ which commutes
with the operators $\varphi$ and $N$. That is to say, for all $g\in G_K$, 
$$g\big(\varphi(T)\big)=\varphi\big(g(T)\big)\qquad\mathrm{and}\qquad g\big(N(T)\big)=N\big(g(T)\big).$$

This, by (1) and (2), shows that, for all $g\in G_K$, 
$$\varphi\big(g(T)\big)=p\cdot g(T)\qquad\mathrm{and}\qquad N\big(g(T)\big)=1.$$
Consequently, if such a $T$ exists, $g(T)$ should satisfy an additive relation 
$$g(T)=T+\eta(g)$$ for a certain $\eta(g)\in \mathbb B_{\mathrm{crys}}$ such that
$\varphi(\eta(g))=p\cdot \eta(g)$. A good choice of
$\eta(g)$ is $\chi_{\mathrm{cycl}}(g)\cdot t$. This then leads to finding an element $T\in\mathbb B_{\mathrm{dR}}$ such that

\noindent
(i) $T$ is transcendental over $\mathbb B_{\mathrm{crys}}$;

\noindent
(ii) $\varphi(T)=pT$; and

\noindent
(iii) $g(T)=T+\chi_{\mathrm{cycl}}(g)\cdot t$ for all $g\in G_K$.

From our experience, a natural way to obtain 
transcendental element is via logarithmic map. 
Thus, by applying the exponential map, we must
 find an element $\varpi\in \widetilde{\mathbb E}^+$ 
satisfying the multiplicative relation $$g(\varpi)=\varpi\cdot\varepsilon^{\chi_{\mathrm{cycl}}(g)}.$$ But this is relatively easy 
since the element $\varpi:=(\varpi^{(n)})\in \widetilde{\mathbb E}^+$ with $\varpi^{(0)}=p$ does the job. Indeed, $\theta\Big(\frac{[\varpi]}{p}-1\Big)=\frac{p}{p}-1=0$. Thus
$$\log[\varpi]:=\log\Big(\frac{[\varpi]}{p}\Big)=\sum_{i=0}^\infty(-1)^{n+1}\frac{\Big(\frac{[\varpi]}{p}-1\Big)^n}{n}=-\sum_{n=0}^\infty\frac{\omega^n}{np^n},$$ 
which is clearly convergent in $\mathbb B_{\mathrm{dR}}$.
As a by-product, this also offers us a (non-canonical) embedding
 $$\mathbb B_{\mathrm{st}}:=\mathbb B_{\mathrm{crys}}\big[\log[\varpi]\big] \hookrightarrow\mathbb B_{\mathrm{dR}}.$$

\vfill
\eject
\centerline{\Large\bf Chapter XIII. Micro Reciprocity Laws and General CFT}
\vskip 0.80cm
\section{\bf Filtered $(\varphi,N)$-Modules and Semi-Stable Reps}

\subsection{Definition}

Let $\rho:G_K\to\mathrm{GL}(V)$ be a $p$-adic Galois representation.
Following Fontaine, define the associated  spaces of periods by
$$\mathbb D_{\mathrm{HT}}(V):=\Big(\mathbb B_{\mathrm{HT}}\otimes_{{\mathbb Q}_p}V\Big)^{G_K},\qquad
\mathbb D_{\mathrm{dR}}(V):=\Big(\mathbb B_{\mathrm{dR}}\otimes_{{\mathbb Q}_p}V\Big)^{G_K},$$
$$\ \,\mathbb D_{\mathrm{st}}(V):=\Big(\mathbb B_{\mathrm{st}}\otimes_{{\mathbb Q}_p}V\Big)^{G_K},\qquad
\mathbb D_{\mathrm{crys}}(V):=\Big(\mathbb B_{\mathrm{crys}}\otimes_{{\mathbb Q}_p}V\Big)^{G_K}.$$

Then by the properties of corresponding rings of periods $\mathbb B_\bullet$,
we know that $\mathbb D_{\mathrm{HT}}(V)$ (resp. $\mathbb D_{\mathrm{dR}}(V)$, 
resp. $\mathbb D_{\mathrm{st}}(V)$, resp. $\mathbb D_{\mathrm{crys}}(V)$) is finite 
dimensional $K$ (resp. $K$, resp. $K_0$, resp. $K_0$)-vector space. Moreover, it is known that

\noindent
(1) the following inequalities hold:
$$\begin{aligned}\mathrm{dim}_{K_0}&\mathbb D_{\mathrm{crys}}(V)\leq \mathrm{dim}_{K_0}\mathbb D_{\mathrm{st}}(V)\\
\leq& \mathrm{dim}_{K}\mathbb D_{\mathrm{dR}}(V)\leq \mathrm{dim}_{K}\mathbb D_{\mathrm{HT}}(V)\\
\leq& \mathrm{dim}_{\mathbb Q_p}(V);\end{aligned}$$

\noindent
(2) Refined structures on the rings of periods $\mathbb B_\bullet$, where $\bullet=$HT, dR, st, crys, naturally induce additional structures on $\mathbb D_{\bullet}(V)$ as well. More precisely,

\noindent
$\bullet$ {\bf Hodge-Tate} On $D:=\mathbb D_{\mathrm{HT}}(V)$, there is a natural filtration structure, the {\it Hodge-Tate filtration}, given by $$\mathrm{Fil}^i_{\mathrm{HT}}\mathbb D_{\mathrm{HT}}(V):=\Big(\mathrm{Fil}^i_{\mathrm{HT}}\mathbb B_{\mathrm{HT}}\otimes_{{\mathbb Q}_p}V\Big)^{G_K}.$$
This is a decreasing filtration by $K$-vector subspaces on $\mathbb D_{\mathrm{HT}}(V).$ Define then the associated graded piece by
$$\mathrm{Gr}_{\mathrm{HT}}^i(D):=\mathrm{Fil}^i_{\mathrm{HT}}(D)/\mathrm{Fil}^{i+1}_{\mathrm{HT}}(D),$$ and its associated
{\it Hodge-Tate 
slope} by $$\mu_{\mathrm{HT}}(D):=\frac{1}{\mathrm{dim}_KD}\cdot\sum_{i\in\mathbb Z}i\cdot\mathrm{dim}_K\mathrm{Gr}_{\mathrm{HT}}^i(D).$$

\noindent
$\bullet$ {\bf de Rham} On $D:=\mathbb D_{\mathrm{dR}}(V)$, there is a natural filtration structure given by $$\mathrm{Fil}^i_{\mathrm{HT}}\mathbb D_{\mathrm{dR}}(V):=\Big(\mathrm{Fil}^i_{\mathrm{HT}}\mathbb B_{\mathrm{dR}}\otimes_{{\mathbb Q}_p}V\Big)^{G_K}.$$
This is a decreasing filtration by $K$-vector subspaces on $\mathbb D_{\mathrm{dR}}(V).$ Define then the associated graded piece by
$$\mathrm{Gr}_{\mathrm{dR}}^i(D):=\mathrm{Fil}^i_{\mathrm{HT}}(D)/\mathrm{Fil}^{i+1}_{\mathrm{HT}}(D),$$ and its associated
slope by $$\mu_{\mathrm{dR}}(D):=\frac{1}{\mathrm{dim}_KD}\cdot\sum_{i\in\mathbb Z}i\cdot \mathrm{dim}_K\mathrm{Gr}_{\mathrm{dR}}^i(D).$$
Since  $$\mathrm{Gr}_{\mathrm{dR}}^i\mathbb B_{\mathrm{dR}}
\simeq\mathbb C(i)\simeq \mathrm{Gr}_{\mathrm{HT}}^i\mathbb B_{\mathrm{HT}},$$	quite often, we will call the above filtration and its associated slope the {\it Hodge-Tate filtration} and the {\it Hodge-Tate slope} respectively;
\vskip 0.20cm
\noindent
$\bullet$ {\bf Crystalline} Being naturally embedded in $\mathbb B_{\mathrm{dR}}$,  $K\otimes_{K_0}\mathbb B_{\mathrm{crys}}$ admits a natural filtration. Consequently, this induces a natural Hodge-Tate filtration on $D:=K\otimes_{K_0}\mathbb D_{\mathrm{crys}}(V)$ by $K$-vector subspaces. 

Denote by $D_0:=\mathbb D_{\mathrm{crys}}(V)$. Since $\mathbb B_{\mathrm{crys}}$ admits a 
natural Frobenius $\varphi$, we obtain a natural $\varphi$-module structure on $D_0$ induced from 
$\varphi\otimes\mathrm{Id}_V$ on $\mathbb B_{\mathrm{crys}}\otimes_{\mathbb Q_p}V$. Thus 
working over $\overline K_0^{\mathrm{ur}}$, or better, via the natural residue map, working 
over $\overline k$, the associated algebraic closure of the residue field $k$ of $K$, according to Dieudonne, there is a natural decomposition
$$\overline D_0=\oplus_{l\in\mathbb Q} \overline{D_{0,l}}.$$ Here
$\overline{D_{0,l=\frac{s}{r}}}$ denotes the $p^s$ eigen-space of $\varphi^r$
with $l=\frac{s}{r}$ the reduced expression of $l\in\mathbb Q$ in terms of
quotient of integers $r,s$, i.e., $r,s\in\mathbb Z$, $(r,s)=1$ and $r> 0$.
Introduce then the associated $\mathbb Q$-indexed filtration by $K_0$ vector subspaces, called the {\it Dieudonne filtration} associated to the $\varphi$-module $D_0$,
by $$\mathrm{Fil}_{\mathrm{Dieu}}^{l_0}D_0:=\oplus_{l\geq l_0}\overline{D_{0,l}}.
$$ Accordingly, define the associated graded $K_0$-vector space by 
$$\mathrm{Gr}_{\mathrm{Dieu}}^{l_0}(D):=
\mathrm{Fil}^{l_0}_{\mathrm{Dieu}}(D_0)\big/
\cup_{l<l_0}\mathrm{Fil}^{l}_{\mathrm{Dieu}}(D_0),$$ and the {\it 
Dieudonne slope} by $$\mu_{\mathrm{Dieu}}(D_0):=\frac{1}
{\mathrm{dim}_{K_0}D_0}\cdot\sum_{l\in\mathbb Q}l\cdot 
\mathrm{dim}_{K_0}\mathrm{Gr}_{\mathrm{Dieu}}^l(D_0).$$

\noindent
$\bullet$ {\bf Semi-Stablility} Unlike $\mathbb B_{\mathrm{crys}}$, there is no natural
embedding of $K\otimes_{K_0}\mathbb B_{\mathrm{st}}$
in $\mathbb B_{\mathrm{dR}}$. But still we can embed $K\otimes_{K_0}\mathbb B_{\mathrm{st}}$ in $\mathbb B_{\mathrm{dR}}$. Fix such an embedding. Then we obtain a filtration by 
$K$-vector subspaces on $D:=K\otimes_{K_0} \mathbb D_{\mathrm{st}}(V)$. One can easily 
check that
this filtration does not depend on the choice used above, thus it is well-justified to call such a filtration the Hodge-Tate filtration on $D$.

Similarly, the Frobenius structure on $\mathbb B_{\mathrm{st}}$ induces a natrual $\varphi$-module structure on the finite dimensional $K_0$-vector space
$D_0:=\mathbb D_{\mathrm{st}}(V)$, or better,  on $\overline{D_0}/\overline{K_0^{\mathrm{ur}}}$. Accordingly, we can introduce the Dieudonne filtration on $D_0$ and hence its associated Dieudonne slope $\mu_{\mathrm{Dieu}}(V)$.
\vskip 0.20cm
Moreover, the natural monodromy operator $N$ on $\mathbb B_{\mathrm{st}}$ 
 introduces a nilpotent monodromy operator $N$ on $D_0$ via $N\otimes\mathrm{Id}_V$ on 
$\mathbb B_{\mathrm{st}}\otimes_{\mathbb Q_p}V$. Motivated by this, we say that $\mathbb D=
(D_0,D)$ is a {\it filtered $(\varphi,N)$-module} if it consists of a
finite dimensional $K_0$-vector space $D_0$ and a finite dimensional $K$-vector space 
$D$, equipped with a exhaustive and separated filtration by $K$-vector subspaces on $D$, 
a $\varphi$-module structure on $D_0$, and a monodromy operator $N$ satisfying the following compatibility conditions:

\noindent
(i) $D\simeq K\otimes_{K_0}D_0$;

\noindent
(ii) $N\circ \varphi=p\varphi\circ N$.

Set 
$$\mu_{\mathrm{HT}}(\mathbb D):=\mu_{\mathrm{HT}}(D),\qquad
\mu_{\mathrm{Dieu}}(\mathbb D):=\mu_{\mathrm{Dieu}}(D_0).$$
It is known that $\mu_{\mathrm{Dieu}}(\mathbb D)$ is equal to the Newton slope
 $\mu_{\mathrm{N}}(\mathbb D)$ of $\mathbb D$. Here, 
 $\mu_{\mathrm{N}}(\mathbb D)$ is defined as follows:

\noindent
(a) If $D_0$ is of dimension 1 over $K_0$, say, $D_0=K_0\cdot d$. Then, we can see that we have
$N=0$ and there exists a non-zero $\lambda\in K_0$ such that $\varphi(d)=\lambda\cdot d$. 
Consequently, we have $$\mu_{\mathrm{N}}(\mathbb D)=v_{K_0}(\lambda);$$
(b) In general, we have $$\mu_{\mathrm{N}}(\mathbb D)=
\mu_{\mathrm{N}}(\mathrm{det}\,\mathbb D),$$ where
$\mathrm{det}\,\mathbb D$ denotes the determinant of $\mathbb D$ obtained by taking the maximal exterior products of $D_0$ and $D$.

Tautologically, we have also the notion of {\it saturated filtered $(\varphi,N)$-submodules}.

\subsection{Weak Admissibility and Semi-Stablility}
Clearly, if $V$ is a $p$-adic semi-stable representation of $G_K$, then
$\mathbb D(V):=\Big(\mathbb D_{\mathrm{st}}(V),\mathbb D_{\mathrm{dR}}(V)\Big)$ admits a natural	filtered $(\varphi,N)$-module structure, since in this case
$$\mathbb D_{\mathrm{dR}}(V)=K\otimes_{K_0}\mathbb D_{\mathrm{st}}(V).$$
Hence it makes sense to talk about the corresponding Hodge-Tate slopes and Newton slopes. Along with this line, an important discovery of Fontaine is the following basic:
\vskip 0.20cm
\noindent
{\bf Theorem.} (Fontaine) {\it Let $\rho:G_K\to\mathrm{GL}(V)$ be a semi-stable $p$-adic representation of $G_K$ and set $\mathbb D:=\big(D_0,D\big)$ with
$$D_0:=\mathbb D_{\mathrm{st}}(V)\quad{and}\quad
D:=\mathbb D_{\mathrm{dR}}(V).$$
Then 

\noindent
(i) $\mu_{\mathrm{HT}}(\mathbb D)=\mu_{\mathrm{N}}(\mathbb D);$ and

\noindent
(ii) $\mu_{\mathrm{HT}}(\mathbb D')\leq\mu_{\mathrm{N}}(\mathbb D')$ for any saturated filtered $(\varphi,N)$-submodule $\mathbb D'=(D_0',D')$ of $\mathbb D=(D_0,D).$} 
\vskip 0.20cm
If a filtered $(\varphi,N)$-module $(D_0,D)$ satisfies the above two 
conditions (i) and (ii), following Fontaine, we call it a 
{\it weakly admissible filtered $(\varphi,N)$-module}. 
So the above result then simply says that for a semi-stable representation $V$, its 
associated periods $\mathbb D:=\Big(\mathbb D_{\mathrm{st}}(V), 
\mathbb D_{\mathrm{dR}}(V)\Big)$ is weakly admissible.
More surprisingly,   the converse holds correctly. That is to say, we also have the following
\vskip 0.20cm
\noindent
{\bf Theorem.} (Fontaine$\|$Colmez-Fontaine) {\it If $(D_0,D)$ is a weakly admissible 
filtered $(\varphi,N)$-module. Then there exists a semi-stable	
representation $V$ of $G_K$ such that}
$$D=\mathbb D_{\mathrm{dR}}(V)\quad\mathrm{and}\quad D_0=\mathbb D_{\mathrm{st}}(V).$$

\noindent
{\it Remark.}	(A$\|$B), for contributors, means
that the assertion is on one hand conjectured by A and on the other   proved by B.

\section{Monodromy Theorem for $p$-adic Galois Representations}

We have already explained two of fundamental results on $p$-adic Galois 
representations, namely, Theorems in \S31 and \S35. Here we introduce another one, 
the so-called Monodromy Theorem for $p$-adic Galois Representations. 

To explain this, let us recall that a $p$-adic Galois representation $\rho:G_K\to\mathrm{GL}(V)$ is called {\it potentially semi-stable}, if there exists a finite Galois extension $L/K$ such that the induced Galois representation
$\rho|_{G_L}:G_L(\hookrightarrow G_K)\to\mathrm{GL}(V)$ is a semi-stable representation. One can 
easily check that every potentially semi-stable representation is de Rham. As a $p$-adic 
analogue of the Monodromy Theorem for $l$-adic Galois Representations, we have the 
following fundamental thing:
\vskip 0.20cm
\noindent
{\bf Monodromy Theorem for $p$-adic Galois Reps.}\,(Fontaine$\|$Berger)
{\it All de Rham representations are potentially semi-stable representations.}
\vskip 0.20cm
Started with Sen's theory for $\mathbb B_{\mathrm{dR}}$ of Fontaine,
bridged by over-convergence of $p$-adic representations due to (Cherbonnier$\|$Cherbonnier-Colmez),
Berger's proof is based on the so-called $p$-adic monodromy theorem (for $p$-adic differentials equations) of (Crew, Tsuzuki$\|$Crew, Tsuzuki, Andre, Kedelaya, Menkhout).
 For more details, please refer to the final chapter.

\section{Semi-Stability of Filtered $(\varphi,N;\omega)$-Modules}

\subsection{Weak Admissibility = Stability and of Slope Zero}

With the geometric picture in mind, particularly the works of Weil, Grothendieck, Mumford, Narasimhan-Seshadri and Seshadri, we then notice that 
weakly admissible condition for filtered $(\varphi,N)$-module $\mathbb D=(D_0,D)$ is an arithmetic analogue of the condition on semi-stable bundles of slope zero. Indeed, if we set 
$\mu_{\mathrm{total}}(\mathbb D):=\mu_{\mathrm{HT}}(D)-\mu_{\mathrm{Dieu}}(D_0),$
then the first condition of weak admissibility, namely,

\noindent
(i) $\mu_{\mathrm{HT}}(D)=\mu_{\mathrm{Dieu}}(D_0)$

\noindent
is {\it equivalent} to the slope zero condition 

\noindent
(i)$'$	$\mu_{\mathrm{total}}(\mathbb D)=0$; 

\noindent
and the second condition

\noindent
(ii) $\mu_{\mathrm{HT}}(D')\leq\mu_{\mathrm{Dieu}}(D_0')$ for any saturated filtered $(\varphi,N)$-submodule $(D_0',D')$ of $(D_0,D),$

\noindent
is {\it equivalent} to the semi-stability condition

\noindent
(ii)$'$ $\mu_{\mathrm{total}}(\mathbb D')\leq \mu_{\mathrm{total}}(\mathbb D)=0$
for all saturated filtered $(\varphi,N)$-submodule $\mathbb D'$ of $\mathbb D$.

\noindent
Put in this way, the above correspondence between semi-stable Galois representations and weakly admissible filtered $(\varphi,N)$-modules
may be understood as an arithmetic analogue of the Narasimhan-Seshadri 
correspondence between irreducible unitary representations and stable bundles of degree zero 
over compact Riemann surfaces. 

Accordingly, in order to establish a general class field theory for $p$-adic 
number fields,	we need 
to introduce some new structures to tackle ramifications. Recall that in 
algebraic geometry, as explained in Part A, there are two parallel theories for this purpose, 
namely, the $\pi$-bundle one on the covering space using 
Galois groups; and the parabolic bundle one on the base space 
using parabolic structures. Hence, in our current arithmetic setting, we would like to develop corresponding theories.

The $\pi$-bundle analogue is easy, based on Monodromy theorem for $p$-adic Galois 
Representations. In fact, we have the following orbifold version:
\vskip 0.20cm
\noindent
{\bf Theorem.} (Fontaine$\|$Fontaine, Colmez-Fontaine, Berger)
{\it There exists a natural one-to-one and onto correspondence 
 \vskip 0.30cm
 \centerline  {\Big\{de Rham Galois representations of $G_K$\Big\}}
 \vskip 0.20cm
\centerline{$\Updownarrow$}
 \vskip 0.20cm
\centerline {\Big\{semi-stable filtered $(\varphi,N;G_{L/K})$ of slope zero: $\exists$\ L/K finite Galois\Big\}.}}

\subsection{Ramifications}

In geometry, parabolic structures take care of ramifications. Recall that if 
$M^0\hookrightarrow M$ is a punctured Riemann surface, then around 
the punctures $P_i\in M\backslash M^0, i=1,2,\dots,N,$ the associated 
monodromy groups generated by parabolic elements $S_i$ are 
isomorphic to $\mathbb Z$, an abelian group. Thus for  a unitary 
representation $\rho:\pi_1(M^0;*)\to \mathrm{GL}(V)$, the images of $
\rho(S_i)$ are given by	 diagonal matrices with diagonal entries $
\exp(2\pi \sqrt{-1}\alpha_{i;k})$, that is to say, they are determined 
by  unitary characters $\exp(2\pi\sqrt{-1}\alpha),\,\alpha\in\mathbb Q$. 
As such, to see the corresponding ramifications, one usually
choose a certain cyclic covering with ramifications around $P_i$'s 
such that the orbifold semi-stable bundles can be characterized by semi-stable 
parabolic bundles on $(M^0,M)$.

However, in arithmetic side, the picture is much more complicated since there is no simple 
way to make each step abelian. By contrast, the good news is that there is a 
well-established theory 
in number theory to measure ramifications, namely, the theory of high ramification groups.

Let then $G_K^{(r)}$ be the upper-indexed high ramification groups of $G_K$, parametrized by 
non-negative reals $r\in\mathbb R_{\geq 0}$.  (See e.g., [Se3].)
Denote then by $V^{(r)}:=V^{G_K^{(r)}}$ the invariant subspace of $V$ under  $G_K^{(r)}$, 
and $K^{(r)}:={\overline K}^{G_K^{(r)}}.$
For a $p$-adic Galois representation $V$, define the associated $r$-th 
graded piece by $$\mathrm{Gr}^{(r)}V:=\bigcap_{s\geq r}V^{(s)}
\big/\bigcup_{s<r}V^{(s)},$$ and its
{\it Swan conductor} by
$$\mathrm{Sw}(\rho):=\sum_{r\in
\mathbb R_{\geq 0}}r\cdot \mathrm{dim}_{{\mathbb Q}_p}\mathrm{Gr}^{(r)}V.$$
{\bf Proposition.} {\it Let $\rho:G_K\to\mathrm{GL}(V)$ be a de Rham representation.

\noindent
(i) ({\bf Hasse-Arf Lemma}) All jumps of $\mathrm{Gr}^{(r)}V$ are rational;

\noindent
(ii) ({\bf Artin, Fontaine}) There exists a Swan representation 
$\rho_{\mathrm{Sw}}:G_K\to\mathrm{GL}(V_{\mathrm{Sw}})$ such that
$$\langle\rho_{\mathrm{Sw}},\rho\rangle=\mathrm{Sw}(\rho).$$
In particular, $\mathrm{Sw}(\rho)\in\mathbb Z_{\geq 0}.$}

\subsection{$\omega$-structures}

Recall that in geometry ([MY]),	 parabolic structures, taking care of ramifications, can also be
characterized via an $\mathbb R$-index filtration
$$E_t:=\Big(p_*\Big(W\otimes\mathcal O_Y\big(-[\#\Gamma\cdot t]D\big)\Big)\Big)^\Gamma,$$ and its associated parabolic degree is measured by 
$$\sum_i\alpha_i\cdot\mathrm{dim}_{\mathbb C}\mathrm{Gr}^iV.$$
Moreover, it is known that the filtration $E_t$ is 

\noindent
(i) left continuous; 

\noindent
(ii) has jumps only at $t=\alpha_i-\alpha_{i-1}\in\mathbb Q$; and

\noindent
(iii) with parabolic degree in $\mathbb Z_{\geq 0}$.
\vskip 0.20cm
Even we have not yet checked with geometers whether their ramification filtration constructions are motivated by the arithmetic one related to
the filtration of upper indexed high ramification groups, the similarities
between both constructions are quite apparent. 
Indeed, it is well-known that,
for the filtrations on Galois groups $G_K$ and on representations $V$
induced from that of high ramification groups $G_K^{(r)}$,

\noindent
(i)  by definition, $G_K^{(r)}$ and hence $V^{(r)}$ are left continuous;

\noindent
(ii)  from the Hasse-Arf Lemma, all jumps of $G_K^{(r)}$ and hence of $V^{(r)}$ are rational; and

\noindent
(iii) according to essentially a result of Artin, the Artin/Swan conductors
are non-negative integers.
\vskip 0.20cm
Motivated by this, for a finite dimensional $K$-vector space $D$,  a
$\omega$-{\it filtration} $\mathrm{Fil}_\omega^rD$ is defined to be a 
$\mathbb R_{\geq 0}$-indexed {\it increasing} and exhausive filtration 
by finite dimensional $K$-vector subspaces on $D$ satisfying the following properties:

\noindent
(i) ({\bf Continuity}) it is left continuous;

\noindent
(ii) ({\bf Hasse-Arf's Rationality}) it has all jumps at rationals; 

Define then the associated $r$-th graded piece by
$$\mathrm{Gr}^{(r)}_\omega D:=\bigcap_{s\geq r}\mathrm{Fil}_\omega^{(s)}D
\big/\bigcup_{s<r}\mathrm{Fil}_\omega^{(s)}D,$$ and its
$\omega$-{\it slope} by
$$\mu_\omega(D):=\frac{1}{\mathrm{dim}_KD}\cdot
\sum_{r\in \mathbb R_{\geq 0}}r\cdot \mathrm{dim}_{K}\mathrm{Gr}_\omega^{(r)}D.$$

\noindent
(iii) ({\bf Artin's Integrality}) The $\omega$-degree $$\mathrm{deg}_\omega(D):=
\sum_{r\in\mathbb R_{\geq 0}}r\cdot\mathrm{dim}_K\mathrm{Gr}_\omega^{(r)}D
=\mathrm{dim}_KD\cdot\mu_\omega(D)$$ is a non-negative integer.

\subsection{Semi-Stability of Filtered $(\varphi,N;\omega)$-Modules}

By the monodromy theorem of $p$-adic Galois representations, 
for a de Rham representation $V$ of $G_K$, there exists a finite Galois extension $L/K$
such that $V$, as a representation of $G_L$, is semi-stable. As such, 
then, over the extension field $L$,
the weakly admissible filtered $(\varphi,N)$-structure on 
$\Big(\mathbb D_{\mathrm{st},L}(V),\mathbb D_{\mathrm{dR},L}(V)\Big)$ 
is equipped with a compatible Galois action of $G_{L/K}$. By contrast, 
motivated by the non-abelian class field theory for Riemann surfaces, 
we expect that the $\omega$-structures would play a similar role 
in our approach to a general CFT in arithmetic 
as that of parabolic structures in geometry. Accordingly, we make the following
\vskip 0.20cm
\noindent
{\bf Definition.} {\it (i) A filtered $(\varphi,N;\omega)$-module 
${\bf D}:=\Big(D_0,D;\mathrm{Fil}_\omega^rD\Big)$ is 
 a filtered $(\varphi,N)$-module $(D_0,D)$ equipped with a compatible $\omega$-structure on $D$;

\noindent 
(ii) Tautologically, we	 have the notion of a  saturated filtered $(\varphi,N;\omega)$-submodule ${\bf D}':=\Big(D_0',D';\mathrm{Fil}_\omega^rD'\Big)$ of ${\bf D}=\Big(D_0,D;\mathrm{Fil}_\omega^rD\Big)$;

\noindent 
(iii) Define the total slope of a filtered $(\varphi,N;\omega)$-module ${\bf D}:=\Big(D_0,D;\mathrm{Fil}_\omega^rD\Big)$ by
$$\mu_{\mathrm{total}}({\bf D}):=\mu_{\mathrm{HT}}(D)-\mu_{\mathrm{Dieu}}(D_0)-\mu_{\omega}(D);$$
(iv) A filtered $(\varphi,N;\omega)$-module ${\bf D}=\Big(D_0,D;\mathrm{Fil}_\omega^rD\Big)$ is called  semi-stable and of slope zero if

\noindent
(a) ({\bf Slope 0}) it is of total slope  zero,	 i.e.,
$$\mu_{\mathrm{total}}({\bf D})=0;$$

\noindent
(b) ({\bf Semi-Stability}) For every saturated filtered 
$(\varphi,N;\omega)$-module ${\bf D}'$ of ${\bf D}$, we have} 
$$\mu_{\mathrm{total}}({\bf D}')\leq \mu_{\mathrm{total}}({\bf D}).$$

\section{General CFT for $p$-adic Number Fields}

\subsection{Micro Reciprocity Law}

With all these preparations, we are now ready to make the following:
\vskip 0.20cm
\noindent
{\bf Conjectural Micro Reciprocity Law.} {\it There exists a canonical one-to-one  correspondence}
\vskip 0.20cm
\centerline{\Big\{de Rham representations of $G_K$\Big\}}
\vskip 0.20cm
\centerline{$\Updownarrow$}
\vskip 0.20cm
\centerline{\Big\{semi-stable filtered $(\varphi,N;\omega)$-modules of slope zero over $K$\Big\}.}

\subsection{General CFT for $p$-adic Number Fields}
Denote	the category of semi-stable filtered $(\varphi,N;\omega)$-modules of slope zero over 
$K$  by $\mathrm{FM}^{\mathrm{ss};0}_K(\varphi,N;\omega)$. Assuming the MRL, i.e., 
the micro reciprocity law, then 
we can easily show that, with respect to natural structures, 
$\mathrm{FM}_K^{\mathrm{ss};0}(\varphi,N;\omega)$ becomes a Tannakian category. 
Denote by $\mathbb F$ the natural fiber functor	 to the category of finite 
$K$-vector spaces. Then, from the standard Tannakian category theory, we obtain the following
\vskip 0.20cm
\noindent
{\bf General CFT for $p$-adic Number Fields}
\vskip 0.20cm
\noindent
$\bullet$ {\bf Existence Theorem} {\it There exists a canonical one-to-one correspondence}
\vskip 0.25cm
\centerline{\Big\{Finitely Generated Sub-Tannakian Categories $\Big(\Sigma,\mathbb F|_{\Sigma}\Big)$\Big\}}
$$\Updownarrow\Pi$$
\centerline{\Big\{Finite Galois Extensions L/K\Big\};}

\noindent
{\it Moreover,

\noindent
$\bullet$ {\bf Reciprocity Law} The canonical correspondence above 
induces a natural isomorphism}
$$\mathrm{Aut}^{\otimes}\Big(\Sigma,\mathbb F|_{\Sigma}\Big)\simeq\mathrm{Gal}\Big(\Pi\Big(\Sigma,\mathbb F|_{\Sigma}\Big)\Big).$$

In fact much refined result holds: By using $\omega$-filtration, 
for all $r\in\mathbb R_{\geq 0}$, we may form a
sub-Tannakian category $(\Sigma^{(r)},\mathbb F|_{\Sigma^{(r)}})$
of $(\Sigma,\mathbb F|_{\Sigma})$, consisting of objects admitting trivial 
$\mathrm{Fil}_\omega^{r'}$ for all $r'\geq r$. 
\vskip 0.20cm
\noindent
$\bullet$ {\bf Refined Reciprocity Law} {\it The natural correspondence $\Pi$ induces, 
for all $r\in\mathbb R_{\geq 0}$, a canonical isomorphism}
$$\mathrm{Aut}^{\otimes}\Big(\Sigma^{(r)},\mathbb F|_{\Sigma^{(r)}}\Big) \simeq \mathrm{Gal}\Big(\Pi\Big(\Sigma,\mathbb F|_{\Sigma}\Big)\Big)\Big/\mathrm{Gal}^{(r)}\Big(\Pi\Big(\Sigma,\mathbb F|_{\Sigma}\Big)\Big).$$

\vfill
\eject
\centerline{\Large\bf Chapter XIV. GIT Stability, Moduli and Invariants}
\vskip 0.80cm
\section{Moduli Spaces}
Let ${\bf D}:=\Big(D_0,D;\mathrm{Fil}_\omega^r(D)\Big)$ be a 
filtered $(\varphi,N;\omega)$-module of rank $d$ over $K$. Then
$D_0$ is a  $d$-dimensional $K_0$-vector 
space equipped with a $(\varphi,N)$-module structure, which induces
a $K_0$-vector subspace filtration of 
$D_0$, namely, the $\mathbb Q$-indexed Dieudonne filtration
$\big\{\mathrm{Fil}_{\mathrm{Dieu}}^l(D_0)\big\}_{l\in\mathbb Q},$ 
$D=K\otimes_{K_0}D_0$, and there are two $K$-vector subspace filtrations of $D$, 
namely, the decreasing Hodge-Tate filtration 
$\big\{\mathrm{Fil}_{\mathrm{HT}}^i(D)\big\}_{i\in\mathbb Z}$, and 
the increasing $\omega$-filtration 
$\big\{\mathrm{Fil}_{\omega}^r(D)\big\}_{r\in\mathbb R_{\geq 0}}$ 
which is compatible with $\varphi$ and $N.$ 

Let $P(\kappa_{\mathrm{Dieu}})$ and $P(\kappa_{\mathrm{HT}})$ 
be the corresponding parabolic subgroups of 
$\mathrm{GL}(D_0)$ and of $\mathrm{GL}(D)$. Define the character 
$L_{\kappa_{\mathrm{HT}}}$ of $P(\kappa_{\mathrm{HT}})$ by
$$L_{\kappa_{\mathrm{HT}}}:=\bigotimes_{i\in \mathbb Z}
\Big(\mathrm{det}\,\mathrm{Gr}_{\mathrm{HT}}^i(D)
\Big)^{\otimes -i}.$$
Similarly, define the (rational)
character $L_{\kappa_{\mathrm{Dieu}}}$ of 
$P(\kappa_{\mathrm{Dieu}})$ by
$$L_{\kappa_{\mathrm{Dieu}}}:=\bigotimes_{l\in\mathbb Q}
\Big(\mathrm{det}\,
\mathrm{Gr}_{\mathrm{Dieu}}^l(D_0)\Big)^{\otimes -l}.$$
(Unlike $L_{\kappa_{\mathrm{HT}}}$, which is an element of the group 
$X^*(P_{\kappa_{\mathrm{HT}}})$ of characters of 
$P_{\kappa_{\mathrm{HT}}}$, being rationally indexed, 
$L_{\kappa_{\mathrm{Dieu}}}$ is in general not an element of $X^*(P_{\kappa_{\mathrm{Dieu}}})$, but a rational character, i.e., it belongs to $X^*(P_{\kappa_{\mathrm{Dieu}}})\otimes\mathbb Q.$)

Moreover, since all jumps of an $\omega$-structure are rationals,
it makes sense to define the associated parabolic subgroup 
$P(\kappa_{\omega})$ and a (rational) character 
$L_{\kappa_{\omega}}$ of $P(\kappa_{\omega})$ by 
$$L_{\kappa_{\omega}}:=\bigotimes_{r\in \mathbb R_{\geq 0}}
\Big(\mathrm{det}\,\mathrm{Gr}_{\omega}^r(D)\Big)^{\otimes -r}.$$

As usual, identify $L_{\kappa_{\mathrm{HT}}}$
with an element of 
$\mathrm{Pic}^{\mathrm{GL}(D)}\Big(\mathrm{Flag}
(\kappa_{\mathrm{HT}})\Big)$, where 
$\mathrm{Flag}(\kappa_{\mathrm{HT}})$ denotes the partial 
flag variety consisting of all filtrations of $D$ with the same graded piece 
dimensions $\mathrm{dim}_{K}\mathrm{Gr}_{\mathrm{HT}}^k(D).$ 
(We have identified $\mathrm{Flag}(\kappa_{\mathrm{HT}})$ with $\mathrm{GL}(D)
\big/P_{\kappa_{\mathrm{HT}}}$.) Similarly, we get an element
$L_{\kappa_{\omega}}$
of $\mathrm{Pic}^{\mathrm{GL}(D)}\Big(\mathrm{Flag}(\kappa_{\omega})\Big)\otimes\mathbb Q$, with
$\mathrm{Flag}(\kappa_{\omega})$
the partial flag variety consisting 
of all filtrations of $D$ with the same $\mathrm{dim}_{K}\mathrm{Gr}_{\omega}^r(D).$ 
Thus, it makes sense to talk about the rational line bundle $\Big(L_{\kappa_{\mathrm{HT}}}\boxtimes L_{\kappa_{\omega}}\Big)\otimes L_{\kappa_{\mathrm{Dieu}}}$ on the product variety $\mathrm{Flag}(\kappa_{\mathrm{HT}})\times
\mathrm{Flag}(\kappa_{\omega})$. Moreover,
define $J=J_K$ be an algebraic group  whose $\mathbb Q_p$-rational points consist of 
automorphisms of the filtered $(\varphi,N;\omega)$-module ${\bf	 D}$ over $K$. 
We infer the following Proposition essentially from the works of
Langton, Mehta-Seshadri, Rapoport-Zink, and particularly, Totaro.
\vskip 0.20cm
\noindent
{\bf Proposition.} ([Lan], [MS], [To]) {\it Assume  $k$ is algebrically closed.	 
Then $\Big(D_0,D;\mathrm{Fil}_\omega^r(D)\Big)$ 
is semi-stable of slope zero if and only if the 
corresponding point $$\Big(\mathrm{Fil}_{\mathrm{HT}}^i(D),
\mathrm{Fil}_{\omega}^r(D)\Big)\in \mathrm{Flag}(\kappa_{\mathrm{HT}})\times
\mathrm{Flag}(\kappa_{\omega})$$
is semi-stable with respect to all one-parameter subgroups 
$\mathbb G_m\to J$ defined over $\mathbb Q_p$ and the 
rational $J$-line bundle
$$\Big(L_{\kappa_{\mathrm{HT}}}\boxtimes L_{\kappa_{\omega}}\Big)\otimes L_{\kappa_{\mathrm{Dieu}}}$$ on $\mathrm{Flag}(\kappa_{\mathrm{HT}})\times
\mathrm{Flag}(\kappa_{\omega})$.} 
\vskip 0.20cm
As a direct consequence, following Mumford's  Geometric Invariant Theory ([M]),
we then obtain the moduli space $\frak M_{K;d,0}^{\varphi,N;\omega}$
 of rank $d$ semi-stable filtered $(\varphi,N;\omega)$-modules of slope zero over $K$.
 In particular, when there is no $\omega$-structure involved, 
 we denote the corresponding moduli space simply by
 $\frak M_{K;d,0}^{\varphi,N}$.
\vskip 0.20cm
\noindent
{\it Remark.} The notion of semi-stable filtered $(\varphi,N;\omega)$-modules of slope $s$ and the associated moduli space $\frak M_{K;r,s}^{\varphi,N;\omega}$ for arbitrary $s$ can also be introduced similarly.	 We leave the details to the reader.

\section{Polarizations and Galois Cohomology}

With moduli spaces of semi-stable filtered $(\varphi,N;\omega)$-modules built, next we want to introduce various invariants (using these spaces). Recall that in (algebraic) geometry for semi-stable vector bundles, this process is divided into two: First we construct natural polarizations via the so-called Mumford-Grothendieck determinant
line bundles of cohomologies; then we study the cohomologies of these polarizations.

Moduli spaces of semi-stable filtered
$(\varphi,N;\omega)$-modules, being projective, admit
natural geometrized polarizations as well. However, such  geometric
polarizations, in general, are quite hard to be used arithmetically, due to the fact that it is difficult to reinterpret them in terms of arithmetic structures involved. To overcome this difficulty, we here want to use Galois cohomologies of  $p$-adic representations, motivated by the $(\frak g,K)$-modules interpretations of cohomology of
(certain types of) vector bundles over homogeneous spaces. 

On the other hand, as said, such polarizations, or better,
determinant line bundles, if exist, should be understood as 
arithmetic analogues of Grothendieck-Mumford determinant line 
bundles constructed using cohomologies of vector bundles. 
Accordingly, if we were seeking a perfect theory, we should first develop
an analogue of sheaf cohomology for filtered 
$(\varphi,N;\omega)$-modules. We will discuss this elsewhere, but merely
point out here the follows:

\noindent
(i) a good cohomology theory in the simplest abelian case 
of $r=1$ is already very interesting since it would naturally
lead to a true arithmetic analogue of the theory of 
Picard varieties, an understanding of which is expected to play 
a key role in our intersectional approach to the Riemann Hypothesis 
proposed in our Program paper [W2];

\noindent
(ii) the yet to be developed cohomology theory would help us to build	
up $p$-adic $L$-functions algebrically.	 This algebraically defined 
$L$-function for filtered $(\varphi, N;\omega)$-modules
then should be compared to $p$-adic $L$-functions for Galois representations
defined using Galois cohomology ([PR]).	 We expect that these two different types of $L$'s 
correspond to each other in a canonical way and further 
can be globalized within the framework of the thin theory
of adelic Galois representations proposed in the introduction.

\section{Iwasawa Cohomology and Dual Exp Map}

In this section, we recall some basic facts about Iwasawa cohomology needed 
in defining $p$-adic $L$-functions following [Col1, Col3].

\subsection{Galois Cohomology}

Let $M$ be a $\mathbb Z_p$-representation of $G_K$. As usual, for any $n\in\mathbb N$, denote by $C_c^n(G_K,M)$ the collections of continuous maps $G_K^n\to M$, called $n$-cochains of $G_K$ with coefficients in $M$. (Thus $C_c^0(G_K,M)$ is simply $M$.)
Define the boundary map $d_n:C_c^n(G_K,M)\to C_c^{n+1}(G_K,M)$ by
$$\begin{aligned}(d_0a)(g)&:=g(a)-a;\\
(d_1f)(g_1g_2)&:=g_1(f(g_2))-f(g_1g_2)+f(g_1);\\
\dots&\dots\end{aligned}$$
$$\begin{aligned}(d_nf)(g_1,g_2,\dots,g_{n+1})&:=g_1(f(g_2,g_3,\dots,g_{n+1})\\
+\sum_{i=1}^n(-1)^if(&g_1,g_2,\dots, g_{i-1}, g_ig_{i+1},g_{i+2},\dots, g_n, g_{n+1})\\
&+(-1)^{n+1}f(g_1,g_2,\dots,g_n).\end{aligned}$$
One can easily check that $\Big(C_c^*(G_K,M),d_*\Big)$ forms 
a complex  of abelian groups. Set $Z_c^n(G_K,M):=\mathrm{Ker}\,d_n$ 
be the collections of
$n$-th cocycles, and $B^n(G_K,M):=\mathrm{Im}\,d_{n-1}$ the collections of $n$-th coboundaries. 
Then, the $n$-th
 {\it Galois cohomology} of $M$ is defined by	
$$H_c^n(G_K,M):=H^n(C_c^*(G_K,M),d_*):=Z_c^n(G_K,M)\big/B_c^n(G_K,M).$$
For examples, 
$H^0(G_K,M)=M^{G_K}$,	
$$Z_c^1(G_K,M)=\Big\{f:G_K\to M: f\ \mathrm{continuous},\ f(g_1g_2)=g_1f(g_2)+f(g_1)\Big\}$$ and
$\displaystyle{B_c^1(G_K,M):=\Big\{f_m:g\mapsto gm-m:\ \exists\ m\in M\Big\}}.$

As usual, for	a $p$-adic representation $V$ of $G_K$, choose a maximal $G_K$-stable $\mathbb Z_p$-lattice $M$,	and set 
$$H^n(G_K,V):=H^{n}(G_K,M)\otimes\mathbb Q_p.$$

\noindent
{\bf Proposition.} (See e.g., [Hi]) {\it Let $V$ be a $p$-adic representation of $G_K$. Then

\noindent
(i) $H^{n\geq 3}(G_K,V)=\{0\}$;

\noindent
(ii) $\displaystyle{H^2(G_K,V)=H^0(G_K,V^\vee(1))^\vee}$ and $\displaystyle{H^1(G_K,V)=H^1(G_K,V^\vee(1))^\vee}$;

\noindent
(iii) $\sum_{n=0}^2(-1)^n\,\mathrm{dim}_{\mathbb Q_p}H^n(G_K,V)=-[K:\mathbb Q_p]\cdot\mathrm{dim}_{\mathbb Q_p}V.$}

\subsection{$(\varphi,\Gamma)$-Modules and Galois Cohomology}

We already knew	 that the category of etale $(\varphi,\Gamma)$-
modules is equivalent to that of $p$-adic Galois representations. Thus, in 
principle, it is possible to compute Galois cohomologies in terms of 
$(\varphi,\Gamma_K)$-modules.

Let $K$ be a $p$-adic number field. As usual, denote by
$K_n:=K(\mu_{p^n}),$ $n\geq 1$ and $K_\infty=\cup_{n\geq 1} K_n$ with 
$\mu_{p^n}$ the $p^n$-th roots of unity. 
Set $\Gamma_n:=\mathrm{Gal}(K_\infty/K_n)$. For simplicity, in the sequel,  assume that $\Gamma_n$ is free and hence of  rank 1 over $\mathbb Z_p$. 

Let $V$ be a $p$-adic representation of $G_K$. For a fixed generator $\gamma\in\Gamma_K$, introduce the complex $C_{\varphi,\gamma}(K,V)$ via:
$$0\to\mathbb D(V)\buildrel{(\varphi-1,\gamma-1)}\over 
\longrightarrow \mathbb D(V)\oplus\mathbb D(V)\buildrel{(\gamma-1)\mathrm{pr}_1-(\varphi-1)\mathrm{pr}_2}\over\longrightarrow \mathbb D(V)\to 0.$$

\noindent
{\bf Lemma.} (Herr) {\it Let $V$ be a $p$-adic representation of $G_K$.
Then the cohomology of the complex $C_{\varphi,\gamma}(K,V)$ is naturally isomorphic to the 
Galois cohomology of $V$.}

\subsection{Iwasawa Cohomology $H_{\mathrm{Iw}}^i(K,V)$}

Choose a system of generators $\gamma_n$ of $\Gamma_n$
such that $\gamma_n=\gamma_1^{p^{n-1}}$. Then,	
$\mathbb Z_p[[\Gamma_K]]$, the so-called the Iwasawa algebra, may be realized as the topological ring $\mathbb Z_p[[T]]$ with the $(p,T)$-adic topology  ($T\,\leftrightarrow\, \gamma-1$), and $$\mathbb Z_p\big[[\Gamma_K]\big]\big/(\gamma_n-1)\simeq \mathbb Z_p\Big[\mathrm{Gal}(K_n/K)\Big].$$ Moreover, via the quitient map $G_K\to \Gamma_K$, we obtain a natural $G_K$ action on $\mathbb Z_p[[\Gamma_K]]$ and hence a $G_K$-action on $\mathbb Z_p[\mathrm{Gal}(K_n/K)]$.

Recall that for a $\mathbb Z_p[G_K]$-module $M$, using Shapiro's lemma, see e.g., [Hi], we have  canonical isomorphisms
 $$H^i\big(G_{K_n},M\big)\simeq H^i\Big(G_K,\mathbb Z_p[\mathrm{Gal}(K_n/K)]\otimes M\Big),$$
which then make	 the corestriction maps
$H^i(G_{K_{n+1}},M)\to H^i(G_{K_n},M)$ a projective system. Consequently,  associated to a 
$\mathbb Z_p$-representation $M$ of $G_K$, we obtain the well-defined 
{\it Iwasawa cohomology groups}
 $$H_{\mathrm{Iw}}^i(K,M):=\lim_{\longleftarrow}
H^i(G_{K_n},M).$$ 
Moreover, for a $p$-adic representation $V$ of $G_K$,  define its 
associated Iwasawa cohomology by
$$H_{\mathrm{Iw}}^i(K,V):=H_{\mathrm{Iw}}^i(K,\Lambda)\otimes_{\mathbb Z_p}\mathbb Q_p,$$ where $\Lambda$ is a (maximal) $G_K$-stable $\mathbb Z_p$-lattice of $V$.

\subsection{Two Descriptions of $H_{\mathrm{Iw}}^i(K,V)$}

There are various ways to describe Iwasawa cohomologies.
For example, we have the following:
 \vskip 0.20cm
 \noindent
{\bf Proposition.}  $\displaystyle{H_{\mathrm{Iw}}^i(K,V)=H^i\Big(G_K,\mathbb Z_p[[\Gamma_K]]\otimes V\Big).}$

Consequently, Iwasawa cohomologies admit natural 
$Z_p[[\Gamma_K]]$-module structures. Quite often we also call $H_{\mathrm{Iw}}^i(K,V)$ Iwasawa modules associated to $V$. Moreover, recall that there is a natural bijection
$$\begin{matrix}\mathbb Z_p[[\Gamma_K]]\otimes V&\simeq&\mathcal D_0(\Gamma_K,V)\\
\gamma\otimes v&\mapsto& \delta_\gamma\otimes v,\end{matrix}$$
where $\mathcal D_0(\Gamma_K,V)$ denotes
the set of $p$-adic measures from $\Gamma_K$ to $V$, and 
$\delta_\gamma$ denotes the Dirac measure at $\gamma$.	
Therefore, we can interpret elements of	 $H^1_{\mathrm{Iw}}(K,V)$ 
in terms of  $p$-adic measures.
In particular, if $\eta:\Gamma_K\to\mathbb Q_p^*$ is a continuous character, 
then, for any $n\geq 1$, we obtain a natural map
$$\begin{matrix}H^1_{\mathrm{Iw}}(K,V)&\to& H^1(G_K,V\otimes\eta)\\
\mu&\mapsto&\int_{\Gamma_{K_n}}\eta\,\mu.\end{matrix}$$

We can also interpret Iwasawa modules in terms of 
$(\varphi,\Gamma)$-modules. Denote by $\psi$ 
a left inverse of the Frobenius $\varphi$.
If $V$ is a $\mathbb Z_p$-representation of $G_K$, 
then there exists a unique operator $\psi:
\mathbb D(V)\to\mathbb D(V)$ such that 
$\psi(\varphi(a)x)=a\psi(x)$ and $\psi(a\varphi(x))=\psi(a)x$ for
$a\in A_K, x\in\mathbb D(V)$ and $\psi$ commutes 
with the action of $\Gamma_K$. Similarly, if $D
$ is an etale $(\varphi,\Gamma)$-module over $A_K$ or $B_K$, 
there exists a unique operator $
\psi:D\to D$ as above. In particular, for any $x\in D$,
$x$ can be written as $x=\sum_{i=0}^{p^n-1}[\varepsilon]^i\varphi^n(x_i)$ 
where $x_i:=\psi^n([\varepsilon]^{-i}x)$.
\vskip 0.20cm
\noindent
{\bf Lamma.} (See e.g., [Col3]) 
{\it (1) If $D$ is an etale $\varphi$-module over $A_K$ (resp. over $B_K$), then

(i) $D^{\psi=1}$ is compact (resp. locall compact);

(ii) $D/(\psi-1)$ is finitely generated over $\mathbb Z_p$ (resp. over $\mathbb Q_p)$.

\noindent
(2) Let $V$ be a $p$-adic representation of $G_K$. Let $C_{\psi,\gamma}$ be the complex
$$0\to\mathbb D(V)\buildrel{(\psi-1,\gamma-1)}\over\longrightarrow
\mathbb D(V)\oplus \mathbb D(V)\buildrel{(\gamma-1)\mathrm{pr}_1-(\psi-1)\mathrm{pr}_2}\over\longrightarrow\mathbb D(V)\to 0.$$
Then we have a commutative diagram of between complexes $C_{\varphi,\gamma}$ and $C_{\psi,\gamma}:$
$$\begin{matrix}0\to&\mathbb D(V)&\buildrel{(\varphi-1,\gamma-1)}\over\longrightarrow&
\mathbb D(V)\oplus \mathbb D(V)&\buildrel{(\gamma-1)\mathrm{pr}_1-(\varphi-1)\mathrm{pr}_2}\over\longrightarrow&\mathbb D(V)&\to 0\\
&\mathrm{Id}\downarrow&&-\psi\oplus \mathrm{Id} \downarrow&&\downarrow -\psi&\\
0\to&\mathbb D(V)&\buildrel{(\psi-1,\gamma-1)}\over\longrightarrow&
\mathbb D(V)\oplus \mathbb D(V)&\buildrel{(\gamma-1)\mathrm{pr}_1-(\psi-1)\mathrm{pr}_2}\over\longrightarrow&\mathbb D(V)&\to 0,\end{matrix}$$
which induces an isomorphism on cohomologies.}
\vskip 0.20cm
\noindent
{\bf Corollary.} (See e.g., [Col3]) {\it If $V$ is a $\mathbb Z_p/\mathbb Q_p$-representation of $G_K$, then $C_{\psi,\gamma}(K,V)$ computes the Galois cohomology of $V$. More precisely,
 
\noindent
(i) $H^0(G_K,V)=\mathbb D(V)^{\psi=1,\gamma=1};$
 
 \noindent
 (ii) $H^2(G_K,V)=\mathbb D(V)/(\psi-1,\gamma-1);$ and
 
\noindent
(iii) there exists a short exact sequence
$$0\to\mathbb D(V)/(\gamma-1)\to H^1(G_K,V)\to\Big(\mathbb D(V)/(\psi-1)\Big)^{\gamma=1}\to 0.$$
Consequently, $H_{\mathrm{Iw}}^i(K,V)=0$ if $i\not=1,\,2$,
and
there are canonical isomorphisms}
$$\mathrm{Exp}^*:H_{\mathrm{Iw}}^1(K,V)=\mathbb D(V)^{\psi=1},\quad\mathrm{and}\quad H_{\mathrm{Iw}}^2(K,V)=\mathbb D(V)\big/(\psi-1).$$

\subsection{Dual Exponential Maps}

From now on, assume that $V$ is de Rham. Then we have the following natural isomorphisms
$$\mathbb B_{\mathrm{dR}}\otimes_{\mathbb Q_p}V\simeq 
\mathbb B_{\mathrm{dR}}\otimes_K\mathbb D_{\mathrm{dR}}(V),\qquad \mathbb D_{\mathrm{dR}}(V)
=H^0(G_K,\mathbb B_{\mathrm{dR}}\otimes_{\mathbb Q_p}V)$$ and
$$H^1(G_K,\mathbb B_{\mathrm{dR}}\otimes_{\mathbb Q_p}V)=H^1(G_K,\mathbb B_{\mathrm{dR}})\otimes_K\mathbb D_{\mathrm{dR}}(V).$$
Recall also  that 

\noindent
(i) for all $k\not=0$, $H^i(G_K,\mathbb C_p(k))=0$ for all $i$;

\noindent
(ii) for all $i\geq 2$, $H^i(G_K,\mathbb C_p)=0$, 
$H^0(G_K,\mathbb C_p)=K$,  and

\noindent
(iii) $H^1(G_K,\mathbb C_p)$ is a one-dimensional $K$-vector space generated by $\log\,\chi\in H^1(G_K,\mathbb Q_p).$

\noindent
Consequently, the cup product $x\mapsto x\cup\log\,\chi$ gives
isomorphisms
$$H^0(G_K,\mathbb C_p)\simeq H^1(G_K,\mathbb C_p)\quad\mathrm{and}\quad
\mathbb D_{\mathrm{dR}}(V)\simeq H^1(G_K,\mathbb B_{\mathrm{dR}}\otimes V).$$
All this then leads to the so-called {\it Bloch-Kato dual exponential map} ([BK2])
for a de Rham representation $V$ of $G_K$, i.e., the composition
$$\mathrm{exp}^*:H^1(G_K,V)\to H^1(G_K,B_{\mathrm{dR}}\otimes V)\simeq \mathbb D_{\mathrm{dR}}(V).$$
Consequently, for any $\mu\in H^1_{\mathrm{Iw}}(K,V)$,
for any $k\in\mathbb Z$, we obtain a natural element
$$\mathrm{exp}^*\Big(\int_{\Gamma_{K_n}}\chi^k\,\mu\Big)\in t^{-k}K_n\otimes_K\mathbb D_{\mathrm{dR}}(V),$$
which is zero when $k\gg 0$.

Moreover, from the overconvergent theory ([CC]),
there exists  $n(V)$ such that, for all $n\geq n(V)$, the natural map $\varphi^{-n}$ sends $\mathbb D(V)^{\psi=1}$ into 
$$\varphi^{-n}\Big(\mathbb D(V)^{\psi=1}\Big)\subset K_n((t))\otimes_K\mathbb D_{\mathrm{dR}}(V).$$

\noindent
{\bf Exp versus exp.} {\it Let $V$ be a de Rham representation of $G_K$, and 
$\mu\in H^1_{\mathrm{Iw}}(K,V)$. Then,
for all $n\geq n(V)$,
$$p^{-n}\cdot\varphi^{-n}\Big(\mathrm{Exp}^*(\mu)\Big)
=\sum_{k\in\mathbb Z}\exp^*\Big(\int_{\Gamma_{K_n}}\chi^k\,\mu\Big).$$}

That is to say, when $V$ is de Rham, the isomorphism
$$\mathrm{Exp}^*:H^1_{\mathrm{Iw}}(K,V)\simeq \mathbb D(V)^{\psi=1}$$ and the Bloch-Kato dual exponential map
admit much more refined arithmetic
structures. This is particularly so when the representation is semi-stable.
In fact, following Perrin-Riou ([PR]),	it is known that they are related to
theory of $p$-adic $L$-functions.
We leave the details to the literatures. Instead,
to end this discussion of polarizations, let us simply point out that
the associated determinants, or better, exterior products, are very important invariants 
and hence should be investigated from a more board point of view.

\vfill
\eject

\centerline{\Large\bf Chapter XV. Two Approaches to Conjectural MRL}
\vskip 0.80cm

\section{Algebraic and Geometric Methods}

There are two different approaches to establish the conjectural Micro Reciprocity Law. 
Namely, algebraic one and geometric one. 

Let us start with algebraic approach. Here, we want to establish a correspondence between 
filtered $\big(\varphi,N;G\big)$-modules $M$ and filtered $\big(\varphi,N;\omega\big)$-modules $D$.
Obviously, this is an arithmetic analogue  of Seshadri's correspondence
between $\pi$-bundles and parabolic bundles over Riemann surfaces. 
Therefore, we expect further that
our correspondence satisfies the following two compatibility conditions:

\noindent
(i) it induces a natural correspondence between saturated subobjects $M'$
and $D'$ of $M$ and $D$; and

\noindent
(ii) it scales the slopes by a constant multiple of $\#G$. Namely,
$$\mu_{\mathrm{total}}(M')=\#G\cdot \mu_{\mathrm{total}}(D').$$

Assume the existence of such a correspondence. Then,
as a direct consequence of the compatibility conditions, 
semi-stable  filtered $\big(\varphi,N;G\big)$-modules $M$ of slope zero 
correspond naturally to semi-stable filtered $\big(\varphi,N;\omega\big)$-modules  $D$
of slope zero. Indeed, if $M$ is a semi-stable filtered $\big(\varphi,N;G\big)$-module 
of slope zero, then using the correspondence, we obtain a 
filtered $\big(\varphi,N;\omega\big)$-module $D$. Clearly, by (ii), we conclude that
the slope of $D$ is zero. Furthermore, $D$ is semi-stable as well:
Let $D'$ be a saturated submodule of $D$. Then, via the induced
correspondence (i) for saturated submodules, there exists a saturated submodule $M'$ of $M$ 
such that the slope of $M'$ is a positive multiple of the slope of $D'$. On the other hand, 
since $M$ is semi-stable, the slope of $M'$ is at most zero. Consequently,
the slope of $D'$ is at most zero too. So $D$ is semi-stable of slope zero. We are done.
Conversely, if $D$ is a semi-stable filtered $\big(\varphi,N;\omega\big)$-module 
of slope zero, then the corresponding
filtered $\big(\varphi,N;G\big)$-module $M$ can be similarly proved 
to be semi-stable of slope zero.

In this way, via the MRL with limited ramifications and the Monodromy Theorem for
$p$-adic Galois Representations, we are able to establish the conjectural MRL.
\vskip 0.20cm
With algebraic approach roughly discussed, let us say also a few wrods on
the geometric approach here. Simply put, the main point we want to
establish there is a direct correspondence between $p$-adic representations with
finite monodromy around marks of fundamental groups of curves defined over finite fields of characteristic $p$ and what we call semi-stable rigid parabolic $F$-bundles in what should be called logarithmic rigid analytic geometry.

\section{MRL with Limited Ramifications}

Before we give more details on our algebraic approach, for completion, let us in this section 
recall some of the key ingredients in establishing the natural connection between semi-stable 
Galois representations and weakly admissible filtered $(\varphi,N)$-modules.

\subsection{Logarithmic Map}

We start with a description of refined structures of $\mathbb B_{\mathrm{crys}}$.

Set $\mathcal O_{\mathbb C}^{**}:=\big\{x\in\mathcal O_{\mathbb C}:\|x-1\|<1\big\}$ be a subgroup of
units of $\mathcal O_{\mathbb C}$. Clearly, 

\noindent
(i) if $x\in \mathcal O_{\mathbb C}^{**}$, then $x^{p^r}\to 1$ as $r\to+\infty$; and

\noindent
(ii) for all $r\in\mathbb Z_{\geq 0}$, the map $x\mapsto x^{p^r}$ induces a surjective
morphism from $\mathcal O_{\mathbb C}^{**}$ into itself with kernel $\mu_{p^r}(\mathbb C).$
Consequently, any element in $\mathcal O_{\mathbb C}^{**}$ has exactly $p^r$ numbers of
$p^r$-th roots in $\mathcal O_{\mathbb C}^{**}$.

Let $$\begin{aligned}U^*:=&\Big\{(x^{(n)})\in\widetilde{\mathbb E^+}:x^{(0)}\in 
\mathcal O_{\mathbb C}^{**}\Big\},\\
U_1^*:=&\Big\{(x^{(n)})\in\widetilde{\mathbb E^+}:x^{(0)}\in 1+2p\mathcal O_{\mathbb C}\Big\}.
\end{aligned}$$
From above, one can easily check that 

\noindent
(iii) the multiplicative group $U_1^*$, resp. $U^*$, admits a natural 
$\mathbb Z_p$-module structure,	 resp. $\mathbb Q_p$-vector space structure, such that
$$U^*\simeq \mathbb Q_P\otimes_{\mathbb Z_p}U_1^*;$$

\noindent
(iv) if $x\in U_1^*$,  then $[x]-1\in\mathrm{Ker}\,\theta+p\cdot W(\widetilde{\mathbb E^+})$. 

Consequently,
the series $$\log[x]:=-\sum_{n=1}^\infty (-1)^n\frac{([x]-1)^n}{n}$$ converges in 
$\mathbb A_{\mathrm{crys}}$. Hence, we get a logarithmic map 
$\log[\ ]:U_1^*\to \mathbb A_{\mathrm{crys}}$ which can also be extended to a 
logarithmic map $\log[\ ]:U^*\to \mathbb B_{\mathrm{crys}}^+$. Denote its image by $U$. 
Clearly, $\varphi([x])=(x^p)$ and $\varphi(\log[x])=p\cdot\log[x]$ for all $x\in U^*$.

\subsection{Basic Structures of $\mathbb B_{\mathrm{crys}}^{\varphi=1}$}

As usual, let 
$$\mathbb B_{\mathrm{crys}}^{\varphi=1}:=\Big\{x\in \mathbb B_{\mathrm{crys}};\varphi(x)=x\Big\}.$$
Also fix an element $v\in U(-1)-\mathbb Q_p$. Then, we have the following

\noindent
{\bf Theorem.} ([CF]) {\it (i.a) $\mathrm{Fil}^0\mathbb B_{\mathrm{crys}}^{\varphi=1}=\mathbb Q_p$;

\noindent
(i.b) $\mathrm{Fil}^i\mathbb B_{\mathrm{crys}}^{\varphi=1}=0$ for all $i>0$;

\noindent
(i.c) $\mathrm{Fil}^{-1}\mathbb B_{\mathrm{crys}}^{\varphi=1}=U(-1);$

\noindent
(i.d)  All elements $b\in\mathrm{Fil}^{-i}\mathbb B_{\mathrm{crys}}^{\varphi=1}$, $i\geq 1$, 
can be written 
in the form $$b=b_0+b_1v+\cdots+b_{r-1}v^{r-1}$$ where $b_0,b_1,\dots,b_{r-1}\in U(-1);$
\vskip 0.15cm
\noindent
(ii.a) For all $r\geq 1$, there is an exact sequence
$$0\to\mathbb Q_p\to\mathrm{Fil}^{-r}\mathbb B_{\mathrm{crys}}^{\varphi=1}\to \Big(\mathrm{Fil}^{-r}\mathbb B_{\mathrm{dR}}/
\mathbb B_{\mathrm{dR}}^+\Big)\to 0;$$

\noindent
(ii.b) There is an exact sequence}
$$0\to\mathbb Q_p\to\mathbb B_{\mathrm{crys}}^{\varphi=1}\to \mathbb B_{\mathrm{dR}}/
\mathbb B_{\mathrm{dR}}^+\to 0.$$

\subsection{Rank One Structures}

Let $V$ be a $p$-adic Galois  representation, then we can form a filtered 
$(\varphi,N)$-module via
$$\mathbb D_{\mathrm {st}}(V):=\Big(\mathbb B_{\mathrm{st}}\otimes_{\mathbb Q_p}V\Big)^{G_K},\quad \mathbb D_{\mathrm {dR}}(V):=\Big(\mathbb B_{\mathrm{dR}}\otimes_{\mathbb Q_p}V\Big)^{G_K}.$$ 
Following Fontaine, if $V$ is semi-stable, then $\big(\mathbb D_{\mathrm {st}}(V),
\mathbb D_{\mathrm {dR}}(V)\big)$ is weakably admissible.
Conversely, for a filtered $(\varphi,N)$-module ${\bf D}=(D_0,D)$, we can introduce a Galois representation
via the functor $$\mathbb V_{\mathrm{st}}({\bf D}):=\Big\{v\in\mathbb B_{\mathrm{st}}\otimes D:
\varphi v=v,\ Nv=0\ \&1\otimes v\in\mathrm{Fil}_{\mathrm{HT}}^0\Big(\mathbb B_{\mathrm{dR}}
\otimes_KD\Big)\Big\}.$$ \
Moreover, following Colmez-Fontaine,
if $(D_0,D)$ is weakly admissible, then $\mathbb V_{\mathrm{st}}({\bf D})$ is semi-stable.
 
While for general ranks, the proof of this equivalence between semi-stable representations 
and weakly admissible filtered $(\varphi,N)$-module is a bit twisted, the rank one case is rather transparent, thanks to the structural result above on $\mathbb B_{\mathrm{crys}}^{\varphi=1}$.
As the statement, together with its proof, is a good place to understand the essentials involved, we decide to include full details.
\vskip 0.20cm
\noindent
{\bf Proposition.} ([CF]) {\it Let ${\bf D}=(D_0,D)$ be a filtered $(\varphi,N)$-module of 
dimension 1 over $K$.

\noindent
(i) If $t_H(D)<t_N(D_0)$, $\mathbb V_{\mathrm {st}}({\bf D})=\{0\}$;

\noindent
(ii) If	 $t_H(D)=t_N(D_0)$,  $\mathrm{dim}_{\mathbb Q_p}\bf V_{\mathrm {st}}({\bf D})=1$.
If $\mathbb V_{\mathrm {st}}({\bf D})$ is generated by $\alpha\cdot\bold x$, 
$\alpha$ is an invertible  element of $\mathbb B_{\mathrm {st}}$;

\noindent
(iii) If  $t_H(D)>t_N(D_0)$,  $\mathbb V_{\mathrm {st}}({\bf D})$ is infinite 
dimensional over $\mathbb Q_p$.}
\vskip 0.20cm
\noindent
{\it Proof.}  The core is really the 
structural result on ${\bf B}_{\mathrm {crys}}^{\varphi=1}$ stated in the previous subsection.
(In fact, only (i) and (ii) will be used.) 

\noindent
Step One: {\it Twisted by $\mathbb Q(-m)$ to make Hodge-Tate weight zero}.
Since $\dim_{K_0}D_0=1$ and $N$ 
is nilpotent, we have $D_0=K_0{\bold x}$ with $N\bold x=0$. 
Let $\varphi(\bold x)=a\cdot \bold x=p^ma_0\cdot 
\bold x$ with $m=v_p(a)=t_N(D)$ and $a_0\in K_0$ satisfying $v_p(a_0)=0$. Then there 
exists an element $\alpha_0\in W(\overline k)$ satisfying
$\varphi(\alpha_0)=a_0\alpha_0$. Set accordingly $\alpha=\alpha_0^{-1}\cdot t^{-m}.$
Clearly, $\alpha$ is an invertible element in $\mathbb B_{\mathrm {crys}}$. 

\noindent
Step Two: {\it Deduced to Crystalline Periods}. 
If 
$\beta\bold x\in \mathbb V_{\mathrm {st}}(D)$ with $\beta\not=0$, then

a) $0=N(\beta \bold x)=N(\beta)\bold x+\beta\cdot N(\bold x)=N(\beta)\bold x$. Hence $N(\beta)=0$;

b) $\beta\in \mathrm{Fil}^{-t_H(D)}\mathbb B_{\mathrm {st}}$ by definition; And

c) $\beta\bold x=\varphi(\beta)\varphi(\bold x)=\varphi(\beta)\cdot a\bold x$.
So $\varphi(\beta)=a^{-1}\cdot\beta.$

\noindent
Therefore, 
$$\mathbb V_{\mathrm {st}}(D)=\Big\{\beta\bold x\,\big|\,\beta\in \mathrm{Fil}^{-t_H(D)}
\mathbb B_{\mathrm {st}}, 
N(\beta)=0, \varphi(\beta)=a^{-1}\beta\Big\}.$$
Set then $\beta=y\alpha\in\mathbb B_{\mathrm {st}}$, (since 
$\alpha\in \mathbb B_{\mathrm {crys}}$ is invertible, 
this is possible,) and we have 
$$\varphi(\beta)=\varphi(y)\varphi(\alpha)=
\varphi(y)\cdot\varphi(\alpha_0^{-1})\cdot\varphi(t^{-m})=\varphi(y)
\cdot\varphi(\alpha_0)^{-1}p^{-m}t^{-m}$$ since $\varphi(t^{-1})=(pt)^{-1}.$ 
On the other hand, $$\begin{aligned}
\varphi(\beta)=&a^{-1}\beta=a^{-1}y\alpha=a^{-1}y\cdot 
\alpha_0^{-1}t^{-m}\\
=&y\cdot p^{-m}a_0^{-1}\alpha_0^{-1}\cdot t^{-m}=y\cdot 
\varphi(\alpha_0)^{-1}\cdot p^{-m}t^{-m}.\end{aligned}$$ 
Consequently, $\varphi(y)=y$.
Therefore, 
$$\begin{aligned}\mathbb V_{\mathrm {st}}(D)=&\Big\{y\cdot\alpha\bold x\,\big|\, 
y\in \mathrm{Fil}^{t_N(D)-t_H(D)}\mathbb B_{\mathrm {st}}, N(y)=0, \varphi(y)=1\Big\}\\
=&\Big\{y\cdot\alpha\bold x\,\big|\, y\in \mathrm{Fil}^{t_N(D)-t_H(D)}
\mathbb B_{\mathrm {crys}}, \varphi(y)=1\Big\}\\
=&\Big\{y\cdot\alpha\bold x\,\big|\, y\in \mathrm{Fil}^{t_N(D)-t_H(D)}
\mathbb B_{\mathrm {st}}^{\varphi=1}\Big\}.\end{aligned}$$
This then completes the proof of the Proposition.

\section{Filtration of Invariant Lattices}
Now let come back to our algebaric approach to the conjectural MRL.

Let then ${\bf D}_L:=(D_0,D)$ be a filtered $(\varphi,N;G_{L/K})$-module.
So $D_0$ is defined over $L_0$ and $D$ is over $L$. 
By the compactness of the Galois groups, there exists a lattice version of $(D_0,D)$ which we denote by $(\Lambda_0,\Lambda)$. In particular, $\Lambda_0$ is an $\mathcal O_{L_0}$-lattice with a group action $G_{L_0/K_0}$. 
Consider then the finite covering map $$\pi_0:\mathrm{Spec}\,\mathcal O_{L_0}\to \mathrm{Spec}\,\mathcal O_{K_0}.$$ We identify $\Lambda_0$ with its associated coherent sheaf on $\mathrm{Spec}\,\mathcal O_{L_0}$.
Set $$\Lambda_{0,K}:=\Big((\pi_0)_*\Lambda_0\Big)^{\mathrm{Gal}(L_0/K_0)}.$$ Clearly, there is a natural $(\varphi,N)$-structure on $\Lambda_{0,K}$.	

Moreover, for the natural covering map $$\pi:\mathrm{Spec}\,\mathcal O_{L}\to \mathrm{Spec}\,\mathcal O_{K},$$ view $\Lambda$ as a coherent sheaf on $\mathrm{Spec}\,\mathcal O_{L}$ and form the coherent sheaf $\mathcal O_L\Big(-[\mathrm{deg}(\pi)\cdot t]\frak m_{L}\Big)$, where 
$t\in\mathbb R_{\geq 0}$ and $\frak m_{L}$ denotes the maximal idea of $\mathcal O_{L}$. Consequently, it makes sense to talk about
$$\Lambda_{K}(t):=\bigg(\pi_*\Big(\Lambda\otimes\mathcal O_L\Big(-\big[\mathrm{deg}(\pi)\cdot t\big]\frak m_{L}\Big)\Big)\bigg)^{\mathrm{Gal}(L/K)}.$$ Or equivalently, in pure algebaric language,
$$\Lambda_{K}(t):=\Big(\Lambda\otimes\frak m_L^{\big[t\cdot \#G_{L/K}\big]}
\Big)^{\mathrm{Gal}(L/K)}.$$

Even we can read ramification information involved from	 
this decreasing filtration consisting of invariant $\mathcal O_K$-lattices, unfortunately, we have not yet been able to obtain its relation with $\omega$-structure wanted. 
\vskip 0.20cm
On the other hand, to go  back from 
filtered $(\varphi,N;\omega)$-modules to filtered
$(\varphi,N;G_{L/K})$-modules,	a solution to the inverse Galois problem 
for $p$-adic number fields is needed. (Alternatively, as pointed by Hida, we can first use 
an independent geometric approach to be explained below to establish the conjectural MRL and hence
the general CFT for $p$-adic number fields and then turn back as an application of our CFT to solve 
the inverse Galois problem for $p$-adic number fields.)

\section{Tate-Sen Theory and Its Generalizations}

From now on, we explain what is involved in our second approach to the conjectural Micro 
Reciprocity Law. As said, this approach is an arithmetic-geometrical one, 
with the main aim to characterize {\it $p$-adic representations} of fundamental groups 
with finite monodromy around marks of algebraic curves defined 
over finite fields of characteristic $p$ in terms of what we call semi-stable 
parabolic rigid $F$-bundles on the logarithmic rigid analytic spaces associated to logarithmic
formal schemes whose special fibers are the original marked curves.
For this purpose, also for the completeness, we start with some preparations.

\subsection{Sen's Method}

Consider then the natural action of $G_K$ on $\widehat{\overline K}=\mathbb C$. For a closed subgroup $H$ of $G_K$, clearly, $\overline K^H\subset\mathbb C^H$, which implies in particular that $\widehat{\overline K^H}\subset \mathbb C^H$. In fact much strong result holds:
\vskip 0.20cm
\noindent
{\bf Ax-Sen-Tate Theorem.} {\it For every closed subgroup $H$ of $G_K$, we have $\widehat{\overline K^H}= \mathbb C^H$. In particular,
$\widehat {K_\infty}=\mathbb C^{H_K}$.}
 
With this, to understand the action of $G_K$ on $\mathbb C$, we are led to the study of the 
residual action of $\Gamma_K$ on $\widehat {K_\infty}$. By using the so-called Tate-Sen 
decompletion process, this can be reduced to the study of 
the action of $\Gamma_K$ on $K_\infty$, which is known to be given by the cyclotomic character.
\vskip 0.20cm
Motivated by this, following Sen, for a general $\mathbb C$-representation of $G_K$, we first concentrate on its $H_K$-invariant part, which offers a natural $\widehat {K_\infty}$-representation of $\Gamma_K$; then by
the decomposition technique just mentioned, 
we are led to a $K_\infty$-representation of $\Gamma_K$.	
$\Gamma_K$ is a rather simple $p$-adic Lie group, namely,
abelian of rank 1 over $\mathbb Z_p$. This final residual $K_\infty$-representation 
of $\Gamma_K$  can be described via its infinitesimal action of $\mathrm{Lie}\,\Gamma_K$, which
in turn is controlled by a single differential operator (modulo a certain finite extension):
\vskip 0.20cm
\noindent
{\bf Theorem.} (Sen) {\it (1) $H^1(H_K,GL_d(\mathbb C))=1$;

\noindent
(2) The natural map $H^1(\Gamma_K,GL_d(K_\infty))\rightarrow 
H^1(\Gamma_K,GL_d(\widehat{K_\infty}))$
induced by the natural inclusion $K_\infty\hookrightarrow\widehat{K_\infty}$ is a bijection;

\noindent
(3) Denote by $\mathbb D_{\mathrm{Sen}}(V)$ the union of all $K_\infty$-vector subspaces 
of $(\mathbb C\otimes_{\mathbb Q_p}V)^{H_K}$
which are $\Gamma_K$-stable and	 finite dimensional (over $K_\infty$). 
Then for $\gamma\in \Gamma_K$ close enough to 1,
the series operator on $\mathbb D_{\mathrm{Sen}}(V)$ defined by
$$\Theta:=-\frac{1}{\log_p\chi_{\mathrm{cyc}}(\gamma)}\cdot\sum_{n\geq 1}\frac{(1-\gamma)^n}{n}$$ converges and is independent of the choice of $\gamma$.}
\vskip 0.20cm
Consequently, for a $\mathbb C$-representation $V$ of $G_K$ of dimension $d$, 
we have the following associated structures:

\noindent
(1) The $H_K$-invariants  $(\mathbb C\otimes_{\mathbb Q_p}V)^{H_K}$ is a 
 $\widehat{K_\infty}$-vector space of dimension $d$;

\noindent
(2) $\mathbb D_{\mathrm{Sen}}(V)$ is a	$K_\infty$-vector space of dimension $d$;

Therefore, the natural map $$\widehat{K_\infty}\otimes_{K_\infty}
\mathbb D_{\mathrm{Sen}}(V)\to (\mathbb C\otimes_{\mathbb Q_p}V)^{H_K}$$ is an isomorphism and we have a natural residual action
of $\Gamma_K$ on $\mathbb D_{\mathrm{Sen}}(V)$.

\noindent
(3) The action of $\mathrm{Lie}(\Gamma_K)$ on $\mathbb D_{\mathrm{Sen}}(V)$ is given by the operator $\Theta:=\frac{\log(\gamma)}{\log_p\chi_{\mathrm{cyc}}(\gamma)}$ (where $\gamma\in\Gamma_K$
is chosen to be close enough to 1) defined as above, which is  $K_\infty$-linear.

Due to the fact that $\Theta$ is defined only for $\gamma$ close enough to 1,
the Lie action is only defined for a certain open subgroup of $\Gamma_K$. 
This is why in literature quite often we have to shift our discussion from $K$-level to $K_n$-
level for a certain $n$.

\subsection{Sen's Theory for $\mathbb B_{\mathrm{dR}}$}

The above result of Sen is based on the so-called Sen-Tate method. This method has been
generalized by Colmez to a much more general context. (See e.g., [Col3], [FO].)
This then leads Fontaine to obtain Sen's theory
for $\mathbb B_{\mathrm{dR}}$ and Cherbonnier-Colmez to the theory of overconvergence, both of 
which play  key roles in Berger's solution to Fontaine's Monodromy Conjecture for $p$-adic 
Galois representations.

\vskip 0.20cm
\noindent
{\bf Theorem.} (Fontaine)
{\it Let $V$ be a $p$-adic representation of $G_K$ of dimension $d$. Then we have 
the following associated structures:

\noindent
(i) There is a maximal element $\mathbb D_{\mathrm{Fon}}^+(V)$ in the set of 
finitely generated  $\Gamma_K$-stable $\mathbb K_\infty[[t]]$-submodules of $(\mathbb B_{\mathrm{dR}}^+\otimes_{\mathbb Q_p}V)^{H_K}$;

\noindent
(ii) The $\mathbb K_\infty[[t]]$-submodule $\mathbb D_{\mathrm{Fon}}^+(V)$ is a free 
$\mathbb K_\infty[[t]]$ of rank $d$ equipped with a natural residual $\Gamma_K$-action whose infinitesimal action via $\mathrm{Lie}(\Gamma_K)$ is given by a differential 
operator $\nabla_V$;

\noindent
(iii) $V$ is de Rham if and only if $\nabla_V$ has a full set of solutions in 
$\mathbb D_{\mathrm{Fon}}^+(V)$;

\noindent
(iv) Natural residue map $\theta:\mathbb B_{\mathrm{dR}}^+\to \mathbb C$ when applying to $(\mathbb D_{\mathrm{Fon}}^+(V),\nabla_V)$ gives rise naturally to $(\mathbb D_{\mathrm{Sen}}(V),\Theta_V).$}

\subsection{Overconvergency}

By the work of Fontaine, for a $p$-adic representation $V$ of $G_K$, we can associate it to
an etale $(\varphi,\Gamma)$-module $\mathbb D(V):=\big(\mathbb B\otimes_{\mathbb Q_p}V\big)^{H_K}.$
While useful, this etale $(\varphi,\Gamma)$-module $\mathbb D(V)$ is only a 
first approximation to the Galois representation $V$ since $\mathbb B$ is too rough. 
Thus, 
certain refined structures should be introduced. This leads to the theory of overconvergence. 

Let $\mathbb B^{\dag,r}$ be the subring of $\mathbb B$ defined by
$$\begin{aligned}\mathbb B^{\dag,r}:=\Big\{x\in\mathbb B:&x=\sum_{k\gg -\infty}p^k[x_k],\\
&x_k\in\widetilde{\mathbb E}, \lim_{k\to\infty}\Big(k+\frac{p-1}{p}\cdot\frac{1}{r}\cdot v_E(x_k)\Big)=+\infty\Big\}.\end{aligned}$$
One checks that $$\begin{aligned}\mathbb B_K^{\dag,r}:=&(\mathbb B^{\dag,r})^{H_K}
=\Big\{\sum_{k=-\infty}^\infty a_k\pi_K^k:\ a_k\in K_\infty\cap F^{\mathrm{ur}},\\
&\sum_{k=-\infty}^\infty a_k X^k\mathrm{convergent\ and\ bounded\ on}\ p^{-1/e_Kr}\leq|X|<1\Big
\}\end{aligned}$$ where $e_K$ denotes the ramification index of $K_\infty/K_{0,\infty}$.

We say that a $p$-adic representation $V$ of $G_K$ is {\it overconvergent} 
if, for some $r\gg 0$, $\mathbb D(V):=\big(\mathbb B\otimes_{\mathbb Q_p}V\big)^{H_K}$ 
has a basis consisting of elements of 
$\mathbb D^{\dag,r}(V):=(\mathbb B^{\dag,r}\otimes_{\mathbb Q_p}V)^{H_K}.$
In other words, there exists a basis of $\mathbb D(V)$ whose corresponding matrix $\mathrm{Mat}(\varphi)$ for the Frobenius $\varphi$ belongs to $M(d,\mathbb B^{\dag,r})$ for some $r\gg 0$.
\vskip 0.20cm
\noindent
{\bf Theorem.} (Cherbonnier$\|$Cherbonnier-Clomez) {\it Every $p$-adic representation of $G_K$ 
is overconvergent.}

\section{$p$-adic Monodromy Theorem}

Now we are ready to recall Berger's proof of Monodromy Theorem for $p$-adic Representations.

Let $V$ be a $p$-adic Galois representation of $G_K$. Following Fontaine, we obtain
an etale $(\varphi,\Gamma)$-module $\mathbb D(V)$.	
This, together with the overconvergence of $\mathbb D(V)$, naturally gives raise to 
the question whether the differential operator
$\nabla(:=\log(\gamma)/\log_p(\chi(\gamma))$ for $\gamma\in \Gamma_K$  close enough to 1, reflecting the Lie action of  $\Gamma_K$,)
 makes sense on the overconvergent subspace 
$\mathbb D^\dag(V):=\big(\mathbb B^\dag\otimes_{\mathbb Q_p}V\big)^{H_K}$. Thus, 
we need to check how $\nabla$ acts on the periods 
$\mathbb B_K^\dag:=\cup_{r\gg 0}\mathbb B_{K}^{\dag,r}$. Unfortunately,
$\mathbb B_K^\dag$ is not $\nabla$-closed: Easily one finds that 
$$\nabla\big(f(\pi)\big)=\log(1+\pi)\cdot(1+\pi)\cdot df/d\pi.$$ 
In particular, with the appearence of the factor $\log(1+\pi)$, 
boundness condition for the elements involved in the definition of
$\mathbb B_K^\dag$ becomes clearly too restricted and hence should be removed. 
To remedy this, 
we make the following extension of periods (from $\mathbb B_K^\dag$) to
$$\begin{aligned}
\mathbb B_{\mathrm{rig},K}^{\dag,r}:=\Big\{f(\pi_K)=&\sum_{k=-\infty}^\infty a_k\pi_K^k:\, a_k\in \mathrm{Fr}\,W\big(k_{K_\infty}\big)\\
 &\&\,
f(X)\ \mathrm{convergent\ on}\ p^{-1/e_Kr}\leq |X|<1\Big\}\end{aligned}$$ 
 to include $\log(1+\pi)$. In fact, much more has been
 achieved, namely, we now have a natural geometric interpretation for the periods:
 The union $\displaystyle{\cup_{r\gg 0}\mathbb B_{\mathrm{rig},K}^{\dag,r}=:\mathbb B_{\mathrm{rig},K}^{\dag}}$ 
 is exactly the so-called {\it Robba ring} used in the theory of
$p$-adic differential equations. Consequently, $\mathbb B_K^\dag$ is the subring
of $\mathbb B_{\mathrm{rig},K}^{\dag}$ consisting of those functions which are bounded; and 
 $\nabla$ naturally acts
on the periods $$\mathbb D_{\mathrm{rig}}^\dag(V):=
\mathbb B_{\mathrm{rig},K}^\dag\otimes_{\mathbb B_K^\dag}\mathbb D^\dag(V).$$

For general $p$-adic representations $V$, the differential operators $\nabla$ do not behave nicely. 
However, for de Rham representations, the situation changes dramatically: 
\vskip 0.20cm
\noindent
{\bf Theorem.} (Berger) {\it  Let $V$ be a  $p$-adic Galois representation of dimension 
$d$. Then 

\noindent
(i) $V$ is de Rham
if and only if there exists a free $\mathbb B_{\mathrm{rig},K}^\dag$-submodule 
$\mathbb N_{\mathrm{BW}}(V)$ of rank $d$ of $\mathbb D_{\mathrm{rig}}^\dag(V)\big[\frac{1}{t}\big]$ which is stable under the differential operator 
$\partial_V:=\frac{1}{\log(1+\pi)}\cdot\nabla_V$
and  the Frobenius operator $\varphi$ such that 
$\varphi^*\mathbb N_{\mathrm{BW}}(V)=\mathbb N_{\mathrm{BW}}(V)$;

\noindent
(ii) $V$ is semi-stable if and only if $\big(\mathbb B_{\mathrm{log},K}^\dag[\frac{1}{t}]
\otimes_{\mathbb B_K^\dag}\mathbb D^\dag(V)\big)^{\Gamma_K}$ is a  $K_0$-vector space of dimension $d$, where, as usual,	$\mathbb B_{\mathrm{log},K}^\dag:=\mathbb B_{\mathrm{log},K}^\dag[\log\pi]$; and

\noindent
(iii) $V$ is crystalline if and only if $\big(\mathbb B_{\mathrm{rig},K}^\dag[\frac{1}{t}]\otimes_{\mathbb B_K^\dag}\mathbb D^\dag(V)\big)^{\Gamma_K}$ is a  $K_0$-vector space of dimension $d$.}

In fact, (ii) and (iii) may be obtained by using Sen's method, that is, 
the so-called regularization and decompletion processes.
\vskip 0.20cm
\noindent
{\bf Examples}. (1) When $V$ is crystalline, we have  $\mathbb N_{\mathrm{BW}}(V)=\mathbb B_{\mathrm{rig},K}^\dag\otimes_F\mathbb D_{\mathrm{crys}}(V)$, a result essentially due to Wach [Wa1,2];

\noindent
(2)  When $V$ is semi-stable, we have $\mathbb N_{\mathrm{BW}}(V)=\mathbb B_{\mathrm{rig},K}^\dag\otimes_F\mathbb D_{\mathrm{st}}(V)$.
\vskip 0.20cm
Berger-Wach modules $\mathbb N_{\mathrm{BW}}(V)$ above are examples of the so-called
$p$-adic differential equation with Frobenius structure. In this language, 
Berger's theorem claims that $V$ is de Rham if and only if there exists a
$p$-adic differential equation $\mathbb N_{\mathrm{BW}}(V)\subset \mathbb D_{\mathrm{rig}}^\dag(V)\big[\frac{1}{t}\big]$ with Frobenius structure.
\vskip 0.20cm
\noindent
{\it Remark.} We say that a {\it $p$-adic differential equation} is a free  module $M$
of finite rank over the Robba ring $\mathbb B_{\mathrm{rig},K}^\dag$ equipped with a connection $\partial_M:M\to M$; $M$ is equpied with a {\it Frobenius structure} if there is a semi-linear Frobenius 
$\varphi_M:M\to M$ which commutes with $\partial_M$; and $M$ is called 
{\it quasi-unipotent} if there exists a finite extension $L/K$ such that $\partial_M$ has a full set 
of horizontal solutions in $\mathbb B_{\mathrm{rig},L}^\dag\big[\log(\pi)\big]\otimes_{\mathbb B_{\mathrm{rig},K}^\dag}M$. 

With all this, then we are in a position to recall the following 
 fundamental result on $p$-adic differential equations.
\vskip 0.20cm
\noindent
{\bf $p$-adic Monodromy Thm.}\,(Crew,Tsuzuki$\|$Andre,Kedlaya,Mebkhout) 
{\it Every $p$-adic differential equation with a Frobenius structure is quasi-unipotent.}
\vskip 0.20cm
Consequently, if $V$ is a de Rham representation of dimension $d$, 
then following Berger, we obtain a $p$-adic differential equation
$\mathbb N_{\mathrm{BW}}(V)$ equipped with a Frobenius structure. Thus 
there exists a finite extension $L/K$ such that 
$\big(\mathbb B_{\mathrm{rig},L}^\dag\big[\log(\pi)\big]\otimes_{\mathbb B_{\mathrm{rig},K}^\dag}\mathbb N_{\mathrm{BW}}(V)\big)^{G_K}$ is a $K_0$-vector space of dimension $d$. Therefore, 
by Theorem (ii), $V$  is a semi-stable representation of $G_L$. 
In other words, $V$ itself is a potentially semi-stable representation of $G_K$.
This is nothing but the statement of Fontaine$\|$Berger's Monodromy Theorem for $p$-adic Galois Representations.

\section{Infinitesimal, Local and Global}

In this section, we briefly recall how micro arithmetic objects of Galois representations are naturally related with global geometric objects of the so-called overconvergent $F$-isocrystals.

\subsection{\bf From Arithmetic to Geometry}

The shift from arithmetic to geometry, as said, is carried out via 
Fontaine-Winterberger's	 fields of norms.

Let $K$ be a $p$-adic number field
with $\overline K$ a fixed algebraic closure and $K_\infty=\cup_n K_n$ with $K_n:=K(\mu_{p^n})$ the cyclotomic extension of $K$ by adding $p^n$-th root of unity. Denote by $k$ 
its residue field, and $K_0:=\mathrm{Fr}\,W(k)$ the maximal unramified extension of $\mathbb Q_p$ contained in $K$.
Set $\varepsilon:=(\varepsilon^{(n)})$ with $\varepsilon^{(n)}\in\mu_{p^n}$
satisfying $\varepsilon^{(1)}\not=1,\  (\varepsilon^{(n+1)})^p=\varepsilon^{(n)}$, and introduce
the base field $E_{K_0}:=k_K((\varepsilon-1))$. Then, from the theory of fields of norms,
associated to $K$, there exists a finite extension $E_K$ of $E_{K_0}$ in a fixed separable closure $E_{K_0}^{\mathrm{sep}}$
such that we have a canonical isomorphism
$$H_K:=\mathrm{Gal}\Big(\overline K/K_\infty\Big)\simeq\mathrm{Gal}\Big(E_K^{\mathrm{sep}}/E_K\Big),$$ where
$E_K^{\mathrm{sep}}:=\bigcup_{L/K:\mathrm{finite\,Galois}} E_L$ is	
a separable closure of $E_K$. In this way,  the arithmetically defined Galois group
$H_K$ for $p$-adic field $K_\infty$ is 
tranformed into the geometrically defined Galois group 
$\mathrm{Gal}\Big(E_K^{\mathrm{sep}}/E_K\Big)$ 
for the field $E_K$ of power series defined over finite field.

\subsection{From Infinitesimal to Global}

Let $\rho:G_K\to GL(V)$ be a $p$-adic representation of $G_K$. Then, following Fontaine, 
we obtain
an etale $(\varphi,\Gamma)$-module $\mathbb D(V)$. Moreover, by a result of Cherbonnier-Colmez [CC], 
$\mathbb D(V)$ is an overconvergent representation. 
Note that now $\Gamma_K$, being the Galois group of $K_\infty/K$, is abelian and may be viewed 
as an open subgroup of $\mathbb Z_p^*$ via cyclotomic character. This,
following Sen, leads naturally to a certain connection. In this way, 
we are able to transform our initial arithmetic 
objects of Galois representations into the corresponding structures in
geometry, namely, that of $p$-adic differential equations with Frobenius structure, 
following Berger [B1]. 
However, despite of this successful transformation, we now face a new challenge -- 
In general, the $p$-adic differential equations obtained have singularities. 
This finally leads to the category of 
de Rham representations: thanks to the works of Fontaine and Berger, for de Rham representations, there areonly removable singularities.

On the other hand, contrary to this  infinitesimal theory, 
thanks to the works of	Levelt and Katz ([Le], [Ka2]), we are led to a corresponding global theory, 
the framework of which was first built up by Crew ([Cre]) based on
Berthelot's overconvergent isocrystals ([B2], [BO1,2], [O]). For more details, see the discussion below.
Simply put, the up-shot is the follows: 
If $X^0\hookrightarrow X$ is a marked regular algebraic curve defined over $\mathbb F_q$, then, 
Crew (for rank one) ([Cre]) and Tsuzuki (in general) ([Ts1]) show that 
there exists a canonical one-to-one  correspondence between $p$-adic representations of 
$\pi_1(X^0,*)$ with finite monodromy along $Z=X\backslash X^0$ and the so-called unit-root 
$F$-isocrystals on $X^0$ overconvergent around $Z$.
This result is an arithmetic-geometric analogue of the result of Weil recalled in Part A on correspondence between complex representations of fundamental groups and flat bundles over compact Riemann surfaces, at least when $Z$ is trivial.

Conversely, to go from global overconvergent isocrystals to micro $p$-adic Galois representations, 
aiming at establishing the conjectural MRL relating de Rham representations to semi-stable
filtered $(\varphi,N;\omega)$-modules, additional works should be done. 
We suggest the reader to go to the papers [Ber3], [Tsu2] and [Mar].

\section{Convergent $F$-isocrystals and Rigid Stable $F$-Bundles}

Recall that the $p$-adic Monodromy Theorem 
is built up on Crew and Tsuzuki's works on overconvergent unit-root $F$-isocrystals.
To understand it, in this section, we make some preparations following [Cre].
Along with this same line, we also offer a notion
called semi-stable rigid $F$-bundles of slope zero in rigid analytic geometry, 
which is the key to our algebraic characterization of $p$-adic representations 
of fundamental groups of complete, regular,  geometrically irreducible 
curves defined over finite fields.

\subsection{Rigid Analytic Spaces}

Let $R$ be a complete DVR of characteristic zero with perfect residue field $k$ of
 characteristic $p$ and fraction field $K$. Let $\frak X/R$ be a flat $p$-adic formal $R$-scheme 
 with closed fiber $X=\frak X\otimes k$ and generic fiber the rigid analytic space 
 $\frak X^{\mathrm{an}}/K$. Following Raynaud, the points of $\frak X^{\mathrm {an}}$ then are naturally 
 in bijection corresponding to the set of closed subschemes of $\frak X$ which are integral, 
 finite and flat over $R$. Therefore, we have 
the so-called {\it specialization map} $\mathrm{sp}:\frak X^{\mathrm{an}}\to X$ sending a point of $\frak X^{\mathrm{an}}$, viewed as a subscheme $\frak Z\subset\frak X$, to its support $\frak Z\otimes k$, which is a closed point of $X$. Define, for any subscheme of $X$ (or of $\frak X$), its tube $]Z[_{\frak X}:=]Z[:=\mathrm{sp}^{-1}(Z).$ One can easily	check that 

\noindent
(i) if $Z\subset X$ is open then $]Z[_{\frak X}=\frak Z^{\mathrm{an}}$ where $\frak Z$ is a flat lifting of $Z$ over $R$. In particular, $]X[=\frak X^{\mathrm{an}}$;

\noindent
(ii) if $Z\subset X$ is closed, say, defined by $f_1,\dots, f_n\in\Gamma(\mathcal O_{\frak X})$, then $$]Z[_{\frak X}=\big\{x\in \frak X^{\mathrm{an}}:|f_i(x)|<1\ \forall i\big\}.$$

\subsection{Convergent $F$-Isocrystals}
Let $X/k$ be a separated $k$-scheme of finite type, 
and $X\hookrightarrow \frak Y$ a closed immersion into a flat $p$-adic formal $R$-scheme that is 
formally smooth in a neighborhood of $X$. Then the diagonal embedding
gives us two natural projections
$p_1,p_2:]X[_{\frak Y\times\frak Y}\to]X[_{\frak Y}$. Following Berthelot ([B2]), 
a {\it convergent isocrystal} on $(X/K,\frak Y)$ is a locally free 
sheaf $\mathcal E$ of $\mathcal O_{]X[_{\frak Y}}$-modules endowed with an isomorphism 
$$p_1^*\mathcal E\simeq p^*_2\mathcal E\eqno(*)$$ restricting to the identity on the image of the diagonal and satisfying the usual compatibility conditions (for more involved copies). A morphism of convergent isocrystals on $(X/K,\frak Y)$ is just a morphism of locally free sheaves compatible with $(*)$.
\vskip 0.20cm
\noindent
{\bf Theorem.} (Berthelot) {\it The category of convergent isocrystals on 
\newline $(X/K,\frak Y)$ is 

\noindent
(i) independent, up to	canonical equivalence, of the choice of $X\hookrightarrow \frak Y$;

\noindent
(ii) functorial in $X/K$; and

\noindent
(iii) of local nature on $X$.}

Consequently, since every separated $X/k$ of finite type always admits such embeddings locally on $X$, we obtain the category of convergent isocrystals on a general $X/k$ by glueing.

\subsection{Integrable and Convergent Connections}

Let $\mathcal E$ be a locally free $\mathcal O_{]X[}$-sheaf.
Then an integrable connection $\nabla:\mathcal E\to\mathcal E\otimes\Omega^1_{]X[}$ on $\mathcal E$ may be obtained via an isomorphismn $$q_1^*\mathcal E\to q_2^*\mathcal E,\qquad q_{1,2}:\Delta_1
\to]X[$$ where $\Delta_1$ is the {\it first} infinitesimal neighborhood $\Delta_1$ of the diagonal 
$]X[\subset ]X[\times]X[$, satisfying the usual cocycle conditions (above). Motivated by 
this, an integrable connection $\nabla$ on $\mathcal E$ is called {\it convergent} if 
the associated isomorphism above can be extended to $(*)$, i.e., from the first 
infinitesimal neighborhood to all levels of infinitesmial neighborhoods.

\subsection{Frobenius Structure}
Now assume that $k\supset\mathbb F_q$ and let $F=F_q$ be a fixed power of the absolute Frobenius 
of $k$. Choose once and for all	 a homomorphism $\sigma:K\to K$ extending the $p$-adic lifting 
of $F_q$ on $W(k)$ and fixing a uniformizer $\pi$ of $R$. Then by the functorial property of 
categories of convergent isocrystals, the pair $(F_q,\sigma)$ gives rise to a semi-linear 
functor $F_\sigma^*$.  An {\it F-isocrystal} on $X/K$ is defined to be a
convergent isocrystal $\mathcal E$ equipped with an isomorphism 
$$\Phi:F_\sigma^*\mathcal E\to \mathcal E.$$ We can see that
if $\nabla$ is the integral connection with
$\nabla(s)=:\sum_is_i\otimes\eta_i, \ \eta_i\in\Omega_{]X[}^1$, then $$\nabla\big(\Phi(s)\big)
=\sum_i\Phi(s_i)\otimes\sigma^*\eta_i.$$

\subsection{Unit-Root $F$-Isocrystals}

In the case when $X=\mathrm{Spec}(k)$, an $F$-isocrystal on $X/K$ is simply a finite-dimensional $K$-vector space
endowed with a $\sigma$-linear automorphism $\Phi:\sigma^*V\simeq V$. Since we assume that $\sigma(\pi)=\pi$, following Dieudonne (see e.g., [Man], [Dem]), there is a natural decomposition of $F$-isocrystals
$V=\oplus_lV_l$ indexed by a finite set of $l\in\mathbb Q$, where, if $l=a/b$ with $a,b\in\mathbb Z, (a,b)=1$ and $b>0$, $V_l\otimes_{W(k)}W(\bar k)$ is simply a $\pi^a$-eigenspace of $\Phi^b$.
We call the number
$\frac{l}{\mathrm{dim}V}$ the {\it Dieudonne slope} of $\Phi$ in $V$. If all slopes are the same,	 
$V$ is {\it pure}; moreover, $V$ is called a {\it unit-root isocrystal}, 
if it is pure of slope zero.

More generally, if $(\mathcal E,\Phi)$ is an $F$-isocrystal on $X/R$, then for any point $x\to X$ with values in a perfect field, there exists a formal covering 
$\mathrm{Spf}(R')\to\mathrm{Spf}(R)$ for
$\mathrm{Spf}W\big(k(x)\big)\to \mathrm{Spf}\big(W(k)\big)$. Denote by $\sigma':R'\to R'$ a compatible lift of $F$. Then the pull-back of $(\mathcal E,\Phi)$ to $x/R'$ is an $F$-isocrystal on $x/R'$, the so-called {\it fiber} of $(\mathcal E,\Phi)$ at $x$. We say that 
an $F$-isocrystal
$(\mathcal E,\Phi)$ is called {\it unit-root} if all its fibers are.
\vskip 0.20cm
\noindent
{\bf Theorem.} (Crew) {\it Let $X/k$ be a smooth $k$-scheme and suppose that $\mathbb F_q\subset k$. Then there exists a natural equivalence of categories
$$\mathbb G:\mathbb R\mathrm{ep}_K\big(\pi_1(X)\big)\simeq \mathrm{Isoc}^{F;\mathrm{ur}}(X/K)$$
where $\mathbb R\mathrm{ep}_K\big(\pi_1(X)\big)$ denotes the category of $K$-representations of the fundamental group $\pi_1(X)$ of $X$, and $\mathrm{Isoc}^{F;\mathrm{ur}}(X/K)$ denotes the category of unit-root $F$-isocrystals on $X/K$.}

This result is based on Katz's work on the correspondence between $R$-representations of $\pi_1(X)$ and the so-called unit-root $F$-lattices on $\frak X/R$ ([Ka1]).
Here, as usual, by an $F$-{\it lattice} on $\frak X/(R,\phi)$, we mean
a locally free $R\otimes\mathcal O_{\frak X}$-modules $\mathbb E$ equipped with a map $\Phi:\phi^*\mathbb E\to\mathbb E$ such that $\Phi\otimes\mathbb Q$ is an isomorphism ($\phi:\frak X\to\frak X$ a lifting of the absolute Frobenius of $X$).

The key to Crew's proof is the following Langton type result:
\vskip 0.20cm
\noindent
{\bf Lemma.} ([Cre]) {\it Let $X/k$ be a smooth affine $k$-scheme and $(\mathcal E,\Phi)$ be 
a unit-root $F$-isocrystal on $X/K$. Then there is a unit-root $F$-lattice $(\mathbb E,\Pi)$ on 
$\frak X/R$ such that}
$(\mathcal E,\Phi)=(\mathbb E,\Pi)^{\mathrm{an}}.$

\subsection{\bf Stability of Rigid $F$-Bundles}

The above result of Crew may be viewed as an arithmetic analogue of Weil's result on the 
correspondence between representations of fundamental groups and flat bundles over compact Riemann 
surfaces. However now the context is changed to curves defined over finite fields of characteristic 
$p$, the representations are $p$-adic, and, \
accordingly the flat bundles are replaced by unit-root $F$-isocrystals. In fact, the 
arithmetic result is a bit more refined: since the associated fundamental group is pro-finite, 
the actural analogue in geometry is better to be
understood as the one for unitary representations and unitary flat bundles.

With this picture in mind, it is then very naturally to ask whether an arithmetic structure
in parallel with Narasimhan-Seshadri correspondence between unitary representations and semi-stable bundles of  slope zero can be established in the current setting. This is our next topic.

With the same notationa as above, assume in addition that $X$ is completed. Then it makes 
sense to talk about locally free $F$-sheaves $\mathcal E$ of $\mathcal O_{]X[}$-modules. If 
$X=\mathrm{Spec}(k)$, then $\mathcal E$ is nothing but a finite-dimensional $K$-vector space $V$ endowed with a $\sigma$ automorphism $\Phi:\sigma^*V\simeq V$. Similarly, 
 we have its associated Dieudonne slope. Consequently, for general $X$, if 
$\mathcal E$ is a locally free $F$-sheaves $\mathcal E$ of $\mathcal O_{]X[}$-modules, 
then we can talk above its fibers at points of $X$ with values in a perfect field. We say that 
a locally free $F$-sheaf $\mathcal E$ of $\mathcal O_{]X[}$-modules is of {\it slope} 
$s\in\mathbb Q$, denoted by $\mu(\mathcal E)=s$, if all its fibers have slope 
$s$; and $\mathcal E$ is called {\it semi-stable} if
for all saturated $F$-submodules $\mathcal E'$, we have all slopes of the fibers of 
$\mathcal E'$ is at most $\mu(\mathcal E)$. As usual, if the slopes satisfy the strict inequalities, then we call $\mathcal E$ {\it stable}.
For simplicity, we call such locally free objects semi-stable (resp. stable) rigid 
$F$-bundles on $X/K$ of slope $s$.
\vskip 0.20cm
\noindent
{\bf Conjectural MRL in Rigid Analytic Geometry.} {\it Let $X$ be a regular projective curve 
defined over $k$. There is a natural one-to-one 
correspondence between absolutely irreducible $K$-representations of $\pi_1(X)$ and stable rigid 
$F$-bundles on $X/K$ of slope zero.}
\vskip 0.20cm
\noindent
{\it Remark.} It is better to rename the above as a Working Hypothesis: 
There are certain points here which have not yet been completed understood
due to lack of time.  (For example, in terms of intersection, the so-called Hodge polygon is better than Newton polygon adopted here. ...) See however [Ked1,2].

\section{Overconvergent $F$-Isocrystals, Log Geometry and Stability}
\subsection{Overconvergent Isocrystals}
Suppose that $j:X\hookrightarrow \bar X$ is an open immersion, $\bar X\hookrightarrow \frak Y$ is a closed immersion with $\frak Y/R$ smooth in a neighborhood of $X$ and let $Z:=\bar X- X$. 
If $Z$ is locally defined by $f_1,\dots, f_n\in\Gamma(\mathcal O_{\frak Y})$, set, for 
$\lambda<1$, $$Z_\lambda:=\big\{x\in]\bar X[:|f_i(x)|<\lambda\ \forall i\big\},\qquad X_\lambda:=
]\bar X[ - Z_\lambda,$$
 and let $j_\lambda:X_\lambda\hookrightarrow ]\bar X[$ be the natural inclusion. 
 It is well-known that
the pro-object $\{X_\lambda\}_{\lambda\to 1}$ does not depend on the choice of $f_i$. 
So, for any coherent sheaf $\mathcal E$ on $]\bar X[$, it makes sense to talk about 
$j^\dag\mathcal E:=\lim_{\to}(j_\lambda)_*j_\lambda^*\mathcal E.$ For example, the sheaf 
$j^\dag\mathcal O_{]\bar X[}\subset\mathcal O_{]X[}$ is the ring of germs of functions on $]X[$ 
extending into the tube $]Z[$. Denote by $p_1^*, p_2^*$ the two functors from the category of 
$j^\dag\mathcal O_{]\bar X[_{\frak Y}}$-modules to the category of 
$j^\dag\mathcal O_{]\bar X[_{\frak Y\times \frak Y}}$-modules.	An {\it overconvergent isocrystal} 
$\mathcal E$ on $(X/K,\frak Y,Z)$ is defined to be 
a locally free sheaf of $j^\dag\mathcal O_{]\bar X[_{\frak Y}}$-module $\mathcal E$
endowed with an isomorphism $p_1^*\mathcal E\simeq p_2^*\mathcal E$ satisfying the 
standard cocyle conditions.
\vskip 0.20cm
\noindent
{\bf Theorem.} (Berthelot) {\it	 The category of overconvergent isocrystals on $(X/K,\frak Y,Z)$ is 

\noindent
(i) independent of $\frak Y$, up to canonical equivalence;

\noindent
(ii) of local nature on $\bar X$; and

\noindent
(iii) functorial in the pair $(X\subset\bar X)$.}

Consequently, we define a category of overconvergent isocrystals on $(X/K,Z)$ for any 
$X\subset\bar X$ with $\bar X/k$ separated of finite type by glueing. In fact, 
much stronger result holds:
\vskip 0.20cm
\noindent
{\bf Theorem.} (Berthelot) {\it If $X/k$ is separated and of finite type and $X\subset \bar X$ is a compactification of $X$, then the category of overconvergent isocrystals on $(X/K,\bar X)$ (i) 
depends, up to canonical equivalence, on $X/K$ only; (ii) is of local nature on $X$; and (iii) is functorial in $X/K$.}

Due to this, we often call it the category of overconvergent isocrystals on $X/K$ simply.

Similarly,  an {\it overconvergent $F$-isocrystal} on $X/K$ is defined to be 
an overconvergent isocrystal $\mathcal E$ equipped with an isomorphism $\Phi:F_\sigma^*\mathcal E\simeq\mathcal E.$
Denote	by $\mathrm{OIsoc}^{F;\mathrm{ur}}(X/K)$ the category of unit-root overconvergent
$F$-isocrystals on $X/K$.

\subsection{$p$-adic Reps with Finite Local Monodromy}

From now on assume that $X/k$ is a regular geometrically connected curve with regular compatification $\bar X$. Let $Z:=\bar X-X$. We say that
 a $p$-adic representation $\rho:\pi_1(X)\to GL(V)$ is	
 having {\it finite (local) monodromy around $Z$} if for each $x\in Z$, the image under $\rho$ of the inertia group at $x$ is finite. Denote
by $\mathbb R\mathrm{ep}_K\big(\pi_1(X)\big)^{\mathrm{fin}}$ the associated Tannakian category.
\vskip 0.20cm
\noindent
{\bf Theorem.} (Crew$\|$Crew for rank one, Tsuzuki in general) {\it The restriction of the Crew equivalence $\mathbb G$
induces a natural equivalence}
$$\mathbb G^\dag:\mathbb R\mathrm{ep}_K\big(\pi_1(X)\big)^{\mathrm{fin}}\to 
\mathrm{OIsoc}^{F;\mathrm{ur}}(X/K).$$

More generally, instead of unit-root condition, there is a notion of quasi-unipotency. In this language, then the $p$-adic Monodromy Theorem is nothing but the following
\vskip 0.20cm
\noindent
{\bf $p$-adic Monodromy Theorem.} (Crew, Tsuzuki$\|$Crew, Tsuzuki, Andre, Kedlaya, Mebkhout) {\it Every overconvergent $F$-isocrystal is quasi-unipotent.}

In addition, quasi-unipotent overconvergent $F$-isocrystal has been beautifully 
classified by Matsuda ([Mat]). Simply put, we now have the following structural
\vskip 0.20cm
\noindent
{\bf Theorem.} (Crew, Tsuzuki, MA(C)K, Matsuda) {\it Every overconvergent $F$-isocrystal is  Matsudian, i.e., admits a natural decomposition to the so-called
Matsuda blocks defined by tensor products of etale and unipotent objects.}

\noindent
In a certain sense, while unit-root objects are coming from representations of fundamental groups, quasi-unipotent objects are related with representations of central extension of fundamental groups.
Finally, we would like to recall that overconvergent isocrystals have been used by Shiho to define	
crystalline fundamental groups for high dimensional varieties ([Sh1,2]).

\subsection{Logarithmic Rigid Analytic Geometry}

The above result of Crew \& Tsuzuki is built up from the open part $X$ of $\bar X$, a kind of arithmetic analogue of local constant systems over $\mathbb C$. As we have already seen, 
in Part A,
to have a complete theory, it is even better if such a theory can be studied over the whole 
$\bar X$: After all, for representation side, $\mathbb R\mathrm{ep}_K\big(\pi_1(X)\big)^{\mathrm{fin}}$ is nothing but $\mathbb R\mathrm{ep}_K\big(\pi_1(X)\big)^{Z}$, that is, $p$-adic representations of $\pi_1(X)$
with finite local monodromy around every mark $P\in Z$. For doing so, we propose two different approaches, namely, analytic one and algebraic one.

Let us start with the analytic approach. As said, the analytic condition of 
unit-root $F$-isocrystals on $X$ overconvergent around $Z$ is defined over (infinitesmal neighborhood of) $X$. We need to extend it to the total space $\bar X$. As usual,
this can be done if we are willing to pay the price, i.e., allowing singularities along the boundary. Certainly, in general term, singularities are very hard to deal with. However, with our experience over ${\bf C}$, particularly, the work of Deligne on local constant systems ([De1]),
for the case at hands, fortunately, we can expect that singularities involved are very mild -- 
There are only logarithmic singularities appeared. This leads to the notion of logarithmic convergent $F$-isocrystals $\mathcal E$ over $(\bar X, Z)$: Simply put, it is 
an overconvergent $F$-isocrystal that can be extended and hence
realized as a locally free sheaf of $\mathcal O_{]\bar X[}$-module $\mathcal E$, endowed with an integral connection $\nabla$ with logarithmic singularities along $Z$
$$\nabla:\mathcal E\to\mathcal E\otimes\Omega_{]\bar X[}^1(\log Z),$$
not only defined over the first infinitesimal neighborhood but over all levels of infinitesimal neighborhoods.

Let us next turn to algebraic approach. With the notion of semi-stable rigid $F$-bundles 
introduced previously, it is not too difficult to introduce the notion of what should be
 called semi-stable parabolic rigid $F$-bundles.

Even we understand that additional work has to be
done here using what should be called logarithmic formal, rigid analytic geometry, 
but with current level of understanding of mathematics involved, we decide to leave the details
to the ambitious reader. Nevertheless, we 
would like to single out the following
\vskip 0.20cm
\noindent
{\bf Correspondence I.} {\it There is a natural one-to-one  correspondence between 
unit-root $F$-isocrystals on $X$ overconvergent around $Z:=\bar X- X$ and what should be called 
unit-root logarithmic overconvergent $F$-isocrystals on $(X,Z)/K$.}

\noindent
{\bf Correspondence II.} {\it There is a natural one-to-one  correspondence between
unit-root $F$-isocrystals on $X$ overconvergent around $Z:=\bar X-X$ and what should be called poly-semi-stable parabolic rigid $F$-bundles of slope zero
on $\big(\frak X^{\mathrm{an}},\frak Z^{\mathrm{an}}\big)$. Here $(\frak X,\frak Z)$ denotes a logarithmic formal scheme associated to $(X,Z)$.}

Moreover, by comparing the theory to be developed here with that for $\pi$-bundles
of algebraic geometry for Riemann surfaces recalled in Part A, for a fixed finite Galois
covering $\pi: Y\to X$ ramified at $Z$, branched at $W:=\pi^{-1}(Z)$,
 it is also natural 
for us to expect the following
\vskip 0.20cm
\noindent
{\bf Correspondence III.} {\it There is a natural one-to-one  correspondence between
orbifold rigid $F$-bundles on $\big(\frak Y^{\mathrm{an}},\frak W^{\mathrm{an}}\big)$
and   rigid parabolic $F$-bundles 
on $\big(\frak X^{\mathrm{an}},\frak Z^{\mathrm{an}}\big)$
  satisfying the following compatibility conditions:

\noindent
(i) it induces a natural correspondences among saturated sub-objects;

\noindent
(ii) it scales the slopes by a constant multiple $\mathrm{deg}(\pi)$.}
\eject
\vskip 0.20cm
Assuming all this, then we can obtain the following
\vskip 0.20cm
\noindent
{\bf Micro Reciprocity Law in Log Rigid Analytic  Geometry.} {\it There is a natural
one-to-one and onto correspondence}

\vskip 0.25cm
\centerline{\Big\{irreducible $p$-adic representations of $\pi_1(X,*)\qquad\quad$}
\centerline{$\qquad\quad$ with finite monodromy along $Z:=\bar X\backslash X$\Big\}}
$$\Updownarrow$$
\centerline{\Big\{stable parabolic rigid $F$-bundles of slope 0 on $\big(\frak X^{\mathrm{an}},\frak Z^{\mathrm{an}}\big)$\Big\}.}

\vskip 10.0cm

Institute for Fundamental Research,

The $L$-Academy

 \& 
 
Graduate School of Mathematics,

Kyushu University, 

Fukuoka 819-0395, 

Japan

E-mail: weng@math.kyushu-u.ac.jp

\end{document}